\documentclass[english,11pt]{article}
\usepackage[T1]{fontenc}
\usepackage[utf8]{inputenc}
\usepackage[pdftex]{graphicx}
\usepackage{a4wide}
\usepackage[english]{babel}
\usepackage{amsthm,amsmath,amssymb,bbm,bm,dsfont,natbib,verbatim,xr,natbib,color,float}
\usepackage[draft]{hyperref}
\usepackage{pdfpages}

\numberwithin{equation}{section}
\newtheorem{thm}{Theorem}[section]
\newtheorem{prop}{Proposition}[section]
\newtheorem{lem}{Lemma}[section]
\newtheorem{hyp}{Assumption}

\newcommand{\convD}{\stackrel{d}{\longrightarrow}}

\newcommand{\convAS}{\stackrel{\text{as}}{\longrightarrow}}
\newcommand{\convASstar}{\stackrel{\text{as}*}{\longrightarrow}}
\newcommand{\E}{\mathbb E}
\renewcommand{\P}{\mathbb P}
\newcommand{\I}{\mathbb I}
\newcommand{\G}{\mathbb G}
\newcommand{\V}{\mathbb V}
\newcommand{\R}{\mathbb R}
\newcommand{\N}{\mathbb N}

\newcommand{\eps}{\varepsilon}

\newcommand{\norm}[1]{\left\Vert#1\right\Vert}
\renewcommand{\j}{\bm{j}}

\newcommand{\un}{\boldsymbol{1}}
\newcommand{\unb}{\bm{1}'}
\newcommand{\deux}{\boldsymbol{2}_j}

\newcommand{\Cov}{\mathbb{C}ov}
\newcommand{\e}{\bm{e}}
\newcommand{\n}{\bm{n}}

\renewcommand{\i}{\bm{i}}


\date{}

\begin{document}

\title{Empirical Process Results for Exchangeable Arrays\thanks{We are grateful to anonymous referees and an associate editor  for their thoughtful comments that improved the paper. We would also like to thank St\'{e}phane Bonhomme, Bryan Graham, Isabelle M\'ejean, Pedro Sant' Anna and participants at various seminars and conferences for their remarks.}}
\author{Laurent Davezies\thanks{CREST-ENSAE, laurent.davezies@ensae.fr}
\and Xavier D'Haultf\oe uille
\thanks{CREST-ENSAE. xavier.dhaultfoeuille@ensae.fr}
\and Yannick Guyonvarch
\thanks{CREST-ENSAE. yannick.guyonvarch@ensae.fr}}
\maketitle

\begin{abstract}
Exchangeable arrays are natural tools to model common forms of dependence between units of a sample. Jointly exchangeable arrays are well suited to dyadic data, where observed random variables are indexed by two units from the same population. Examples include trade flows between countries or relationships in a network. Separately exchangeable arrays are well suited to multiway clustering, where units sharing the same cluster (e.g. geographical areas or sectors of activity when considering individual wages) may be dependent in an unrestricted way. We prove uniform laws of large numbers and central limit theorems for such exchangeable arrays. We obtain these results  under the same moment restrictions and conditions on the class of functions as those typically assumed with i.i.d. data. We also show the convergence of bootstrap processes adapted to such arrays.

\textbf{Keywords:} exchangeable arrays, empirical processes, bootstrap.
\end{abstract}

\section{Introduction}

Taking into account dependence between observations is crucial for making correct inference. For instance, different observations may face common shocks, tending to correlate them positively and thus leading to overly optimistic inference when ignored \citep{bertrand2004}. Such common shocks may arise if the data are polyadic (e.g., dyadic), namely they involve interactions between several units of a given population. An example is international trade, where each observation corresponds to a pair of countries, one exporting and the other importing. We can then expect that two such pairs may be dependent whenever they share at least one country, because of that country's specificities in terms of international trade. Common shocks may also correspond to aggregate fluctuations that affect all units sharing some characteristics. For instance, wages of two individuals may be correlated either because they live in the same geographical area, or because they work in the same sector. We refer to multiway clustering when there are several dimensions along which units may be correlated.

\medskip
\cite{holland1976local}, \cite{fafchamps2007formation} derived variance formulas for linear regressions with dyadic data, while \cite{cameron2011}  propose similar formulas for  multiway clustering. The Stata command \texttt{ivreg2} and the R package \texttt{multiwaycov} are now used routinely to report standard errors accounting for multiway clustering. However, theory has lagged behind this practice. \cite{Tabor2019} shows the asymptotic validity of inference based on \citeauthor{holland1976local}'s suggestion for dyadic data, but for OLS estimators only. \cite{graham2018} and \cite{graham2019kernel} study respectively parametric regressions and density estimation with dyadic data. Regarding multiway clustering, the only papers we are aware of are the recent works of \cite{menzel2017} and \cite{mackinnon2017}. Again, they focus on linear parameters.\footnote{On the other hand and interestingly,  \cite{menzel2017} studies inference both with and without asymptotically normality. He also shows that refinements in asymptotic approximations are possible using the wild bootstrap.} 

\medskip
In this paper, we establish uniform laws of large numbers (LLN) and central limit theorems (CLT) for such type of data. Uniform LLNs and CLTs are key in showing consistency and asymptotic normality of nonlinear estimators under weak regularity conditions. As such, they have been studied extensively with i.i.d. but also dependent data. We refer to, e.g., \cite{vanderVaartWellner1996} and \cite{GineNickl2015} for overviews with i.i.d. data,  and \cite{dehling2002empirical} for the case of time series \citep[see also, e.g.,][for recent results on sampling designs]{bertail2017, han2019complex}.  Noteworthy, we obtain these uniform LLNs and CLTs under the same moment restrictions and conditions on the class of functions as those usually considered with i.i.d. data. Thus, statistical results deducted from the uniform LLNs and CLTs with i.i.d. data  directly extend to the exchangeable arrays we consider. As a proof of concept, we consider Z-estimators and smooth functionals of the empirical cumulative distribution function (cdf).

\medskip
We also study consistency of a direct generalization of  the standard bootstrap for i.i.d. data to polyadic data. A related bootstrap scheme for multiway clustering is the so-called pigeonhole bootstrap, suggested by \cite{mccullagh2000resampling} and studied by \cite{owen2007pigeonhole}, but for which no uniform result has been established so far. For both, we establish weak convergence of the corresponding process. These results imply the validity of the corresponding bootstrap schemes in a wide range of setting, including Z-estimators and smooth functionals of the empirical cdf.

\medskip
To prove these results, we first argue that polyadic data correspond to dissociated, jointly exchangeable arrays. Similarly, multiway clustering corresponds to dissociated separately exchangeable arrays. We then rely extensively on the so-called Aldous-Hoover-Kallenberg representation \citep{Hoover1979, Aldous1981,kallenberg1989} for such arrays. This representation allows us in particular to prove a symmetrization lemma, which is very useful to derive the uniform LLNs and CLTs. This lemma generalizes a similar result for i.i.d. data, but also for U-processes \citep[see, e.g.][Theorem 3.5.3]{delapena1999}. Note that simple LLNs  and CLTs have been already proved, or are direct consequences of known results on dissociated, jointly exchangeable arrays. For LLNs, we refer to \cite{Eaglesonweber78} and Lemma 7.35 in \cite{kallenberg05}. For CLTs, see \cite{silverman1976limit}. But to our knowledge, no abstract uniform LLNs and CLTs have been proved so far for such arrays. 

\medskip
Finally, we illustrate our results with two applications to international trade. In the first, we test whether international trade remains stable from one year to another, using a Kolmogorov-Smirnov test. Given the dependence structure over pairs of countries and through  time, the asymptotic distribution of the test under the null is complicated, making the  bootstrap attractive. We show that neglecting the dependence between dyads leads to important overrejection of the null hypothesis. Next, we  estimate the so-called gravity equation, a very popular model for explaining trade between countries. Since \cite{silva2006log}, this equation has often been estimated with Poisson pseudo maximum likelihood, an estimator for which our results apply. Again, much fewer explanatory variables are significant at usual levels when accounting for dependence between pairs of countries than when considering such pairs  to be i.i.d. observations \citep[as in][]{silva2006log}. 

\medskip
The paper is organized as follows. Section \ref{sec:gen_results} describes the set-up and gives our main results for jointly exchangeable arrays. In addition to uniform LLNs and CLTs, we prove weak convergence of our bootstrap scheme. We also show results for Z-estimators and smooth functionals of the empirical cdf. Section \ref{sec:extensions} considers a few extensions. In particular, we study separately exchangeable arrays. An important difference for such arrays is that the multiple dimensions, corresponding to different sources of clustering, may not grow at the same rate. We show that our results still hold in this case. We also study ``degenerate'' cases (in the same sense as with  U-processes) and consider another bootstrap scheme. The two applications to international trade are developed in Section \ref{sec:application}. The appendix presents three key lemmas. In the supplementary material, we present additional extensions. In particular, we generalize our main results to cases where the number of observations for each $k$-tuple (e.g., the number of matches between two sport players) varies. We also display Monte Carlo simulations and all the proofs of our results.

\section{The set up and main results}
\label{sec:gen_results}

\subsection{Set up}
\label{sub:setup}

Before formally defining our data generating process, we introduce some notation. For any $A\subset \R$ and $B\subset \R^k$ for some $k\geq 2$, we let $A^+=A \cap (0,\infty)$ and
$$\overline{B}=\left\{b=(b_1,...b_k)\in B: \; \forall (i,j)\in \{1,...,k\}^2, i\neq j, b_i \neq b_j\right\}.$$
We then let $\mathbb{I}_k=\overline{\N^{+k}}$ denote the set of $k$-tuples of $\mathbb{N}^+$ without repetition. Similarly, for any $n\in \N^+$, we let $\I_{n,k}=\overline{\{1,...,n\}^k}$. For any $\i=(i_1,...,i_k)$  and $\j=(j_1,...,j_k)$ in $\N^k$, we let $\i \odot \j = (i_1 j_1,...,i_k j_k)$. With a slight abuse of notation, we also let, for any $\i=(i_1,...,i_k)\in \N^k$, $\{\i\}$ denote the set of distinct elements of $(i_1,...i_k)$. For any $r\in\{1,...,k\}$, we let
 $$\mathcal{E}_r=\left\{(e_1,...,e_k) \in \{0,1\}^k: \sum_{j=1}^k e_j = r\right\}.$$
Finally, for any $A\subset \N^+$, we let $\mathfrak{S}(A)$ denote the set of permutations on $A$. For any $\i=(i_1,...,i_k)\in \N^{+k}$ and $\pi\in \mathfrak{S}(\N^+)$, we let $\pi(\i)=(\pi(i_1),...,\pi(i_k))$.

\medskip
We are interested in polyadic data, that is to say random variables $Y_{\i}$ (whose support is denoted by $\mathcal{Y}$)  indexed by $\i \in \I_k$. Dyadic data, which are the most common case, correspond to $k=2$. For instance, when considering trade data, $Y_{i_1, i_2}$ corresponds to export flows from country $i_1$ to country $i_2$. In network data, $Y_{i_1, i_2}$ could be a dummy for whether there is a link from $i_1$ to $i_2$. In directed networks, $Y_{i_1,i_2}\neq Y_{i_2,i_1}$, while $Y_{i_1,i_2}=Y_{i_2,i_1}$ in undirected networks. Similarly,  $Y_{i_1, i_2, i_3}$ could capture whether $(i_1,i_2,i_3)$ forms a triad or not \citep[see, e.g.][for a motivation on triad counts]{wasserman1994social}. $Y_{\i}$ could also correspond to data subject to multiway clustering. Then  $i_1$,..., $i_k$ are the indexes corresponding to the different dimensions of clustering, for instance geographical areas and sectors of activity. In such cases, however, adaptations of our set-up are needed, and we postpone this discussion to Section \ref{sub:unbalanced} below.

\medskip
We assume that the random variables are generated according to a jointly exchangeable and dissociated array, defined formally as follows:

\begin{hyp}\label{as:dgp}
For any $\pi\in \mathfrak{S}(\mathbb{N}^+)$, $(Y_{\i})_{\i \in \I_k}\overset{d}{=}(Y_{\pi(\i)})_{\i \in \I_k}$. Moreover, for any $A,B$ disjoint subsets of $\mathbb{N}^+$ with $\min(|A|,|B|)\geq k$, $(Y_{\i})_{\i \in \overline{A^k}}$ is independent of $(Y_{\i})_{\i \in \overline{B^k}}$.
\end{hyp}

The first part imposes that the labelling conveys no information: the joint distribution of the data remains identical under any possible permutation of the labels. The second part states that the array is dissociated: the variables are independent if they share no unit in common. For instance, $Y_{(i_1,i_2)}$ must be independent of $Y_{(j_1,j_2)}$ if $\{i_1,i_2\}\cap \{j_1,j_2\}=\emptyset$. On the other hand, Assumption \ref{as:dgp} does not impose independence otherwise. This is important in many applications. In the international trade example, $Y_{i_1, i_2}$ and $Y_{i_1,i_3}$ are likely to be dependent because if $i_1$ is open to international trade, it tends to export more than the average to any other country. It may also import more from other countries, meaning that $Y_{i_1, i_2}$ and $Y_{i_3,i_1}$ could also be dependent.

\medskip
Lemma \ref{lem:repr} below is very helpful to better understand the dependence structure imposed by joint exchangeability and dissociation. It may be seen as an extension of de Finetti's theorem to arrays satisfying such restrictions. It is also key in establishing our asymptotic results below.

\begin{lem}\label{lem:repr}
	Assumption \ref{as:dgp} holds if and only if there exist i.i.d. variables $(U_{J})_{J\subset \N^+, 1\leq |J|\leq k}$ and a measurable function $\tau$ such that almost surely,\footnote{In this formula, the $(U_{\{\i \odot \e\}^+})_{\bm{e}\in \cup_{r=1}^k \mathcal{E}_r}$ appear according to a precise ordering, which we let nonetheless implicit as it bears no importance hereafter.}
	\begin{equation}
Y_{\i}=\tau\left((U_{\{\i \odot \e\}^+})_{\bm{e}\in \cup_{r=1}^k \mathcal{E}_r}\right) \quad \forall \i\in\I_k.		
		\label{eq:AHK}
	\end{equation}
\end{lem}

This result is due to \cite{kallenberg1989} but a weaker version, where the equality only holds in distribution, is known as Aldous-Hoover representation \citep{Aldous1981,Hoover1979}. Accordingly, we refer to \eqref{eq:AHK} as the AHK representation hereafter. To illustrate it, let us consider dyadic data ($k=2$). Then, according to Lemma \ref{lem:repr}, we have, for every $i_1<i_2$,
\begin{equation}
Y_{i_1,i_2}=\tau(U_{i_1},U_{i_2},U_{\{i_1,i_2\}}).	
	\label{eq:ex_AH}
\end{equation}
Thus, in the example of trade flows, the volume of exports  from $i_1$ to $i_2$ depends on factors specific to $i_1$ and $i_2$, such as  their own GDP, but also on factors relating both, such as the distance between the two countries. \eqref{eq:ex_AH} has been also used by \cite{bickel2009nonparametric} and \cite{bickel2011method} to model network formation (in which case $Y_{i_1,i_2}=1$ if there is a link between $i_1$ and $i_2$, 0 otherwise). Note also the link between \eqref{eq:ex_AH} and U-statistics: $Y_{i_1,i_2}$ would correspond to such a statistic if $\tau$ did not depend on its third argument.

\medskip
Under Assumption \ref{as:dgp}, the $(Y_{\i})_{\i\in \I_k}$ have a common marginal probability distribution, which we denote by $P$. We are interested in estimating and making inference on features of this distribution, such as its expectation or a quantile, based on observing the first $n$ units only, namely the sample $(Y_{\i})_{\i \in \I_{n,k}}$, with $n\geq k$.

\subsection{Uniform laws of large numbers and central limit theorems}
\label{sub:emp_process}

Let $\mathcal{F}$ denote a class of real-valued functions admitting a first moment with respect to the distribution $P$ and let $Pf$ denote the corresponding moment $\E\left[f(Y_{\un})\right]$ (with $\un$ the $k-$tuple $(1,...,k)$). To avoid measurability issues and the use of outer expectations subsequently, we maintain the following assumption:

\begin{hyp}\label{as:measurability}
There exists a countable subclass $\mathcal{G}\subset \mathcal{F}$ such that elements of $\mathcal{F}$ are pointwise limits of sequences of elements of $\mathcal{G}$.
\end{hyp}

Assumption \ref{as:measurability} is not necessary but often imposed \citep[see, e.g.][]{cherno2014,Kato2019}. We refer to \citeauthor{Kosorok2006} (\citeyear{Kosorok2006}, pp.137-140) for further discussion.

\medskip
In this section, we study the empirical measure $\mathbb{P}_{n}$ and the empirical process $\mathbb{G}_{n}$ defined on $\mathcal{F}$ by
$$\mathbb{P}_{n}f=\frac{(n-k)!}{n!}\sum_{\i \in \I_{n,k}}f(Y_{\i}),$$
$$\mathbb{G}_{n}f=\sqrt{n}\left(\mathbb{P}_{n}f - Pf\right).$$
Let $\ell^\infty(\mathcal{F})$ denote the set of bounded functions on $\mathcal{F}$. We prove below that under restrictions on $\mathcal{F}$, $\mathbb{P}_{n}f$ converges almost surely to $Pf$ uniformly over $f\in \mathcal{F}$, while $\mathbb{G}_{n}$ converges weakly in $\ell^\infty(\mathcal{F})$ to a Gaussian process. We  refer to, e.g., \cite{vanderVaartWellner1996} for a formal definition of weak convergence of empirical processes. These results, stronger than  pointwise convergence of $\mathbb{P}_{n}f$ and $\mathbb{G}_{n}f$, are key in establishing the consistency and asymptotic normality of, e.g., smooth functionals of the empirical cdf or Z- and M-estimators. We consider briefly applications in Section \ref{sub:applis} below, and refer to Part 3 of \cite{vanderVaartWellner1996} for a more comprehensive review of statistical applications of empirical process results.

\medskip
We use the rate $\sqrt{n}$ to normalize $\P_nf -Pf$, though we have $n!/(n-k)!$ different random variables. In general, we cannot expect a better rate of convergence. To see this, let $(X_i)_{i\in\N^+}$ be i.i.d. random variables and let $Y_{\i}=\sum_{j\in \{\i\}} X_j$. Then $(Y_{\i})_{\i\in \I_k}$ satisfies Assumption \ref{as:dgp}, and  $\P_nf$ boils down to an average over $n$ i.i.d. terms only. In some cases, however, for instance if the $(Y_{\i})_{\i\in \I_k}$ are i.i.d., the convergence rate is faster than $\sqrt{n}$.\footnote{\label{foot:different_rates} As with U-statistics, we expect different rates depending on the degree of ``degeneracy''.} Theorem \ref{thm:unif} below remains valid in such cases, but the limit Gaussian process is then degenerate. We come back in more details to such cases in Section \ref{sub:degenerate} below.

\medskip
Let us now introduce the restrictions on $\mathcal{F}$ that we use to obtain uniform laws. We require additional notation for that purpose. For any $\eta>0$ and any seminorm $||\cdot||$ on a space containing $\mathcal{F}$, $N(\eta,\mathcal{F},||\cdot||)$ denotes the minimal number of $||\cdot||$-closed balls of radius $\eta$ with centers in $\mathcal{F}$ needed to cover $\mathcal{F}$. $N_{[\;]}(\eta,\mathcal{F},||\cdot||)$ denotes the minimal number of $\eta$-brackets needed to cover $\mathcal{F}$, where an $\eta$-bracket for $f\in\mathcal{F}$ is a pair of functions $(\ell, u)$ such that $\ell \leq f \leq u$ and $||u-\ell||<\eta$. The seminorms we consider hereafter are $\|f\|_{\mu,r}=(\int |f|^rd\mu)^{1/r}$ for any $r\geq 1$ and probability measure or cdf $\mu$. Hereafter, an envelope of $\mathcal{F}$ is a measurable function $F$ satisfying $F(u)\geq\sup_{f\in \mathcal{F}}|f(u)|$. Finally, we let $\mathcal{Q}$ denote the set of probability measures with finite support on $\mathcal{Y}$.

\begin{hyp}\label{as:entropy}
The class $\mathcal{F}$ either:
\begin{enumerate}
	\item[(i)] admits an envelope $F$ with $PF<\infty$ and $\forall \eta>0$, $$\sup_{Q\in\mathcal{Q}}N\left(\eta||F||_{Q,1},\mathcal{F},||\cdot||_{Q,1}\right)<\infty;$$
	\item[(ii)] or satisfies $N_{[\;]}\left(\eta,\mathcal{F},||\cdot||_{L_1(P)}\right)<\infty$ for all $\eta>0$.
\end{enumerate}
\end{hyp}

\begin{hyp}\label{as:vc}
The class $\mathcal{F}$ either:
\begin{enumerate}
	\item[(i)] admits an envelope $F$ with $PF^2<\infty$ and
\begin{equation*}
    \int_0^{\infty}\sup_{Q\in\mathcal{Q}}\sqrt{\log N\left(\eta||F||_{Q,2},\mathcal{F},||\cdot||_{Q,2}\right)}d\eta<\infty;
\end{equation*}
\item[(ii)] or satisfies $\int_0^{\infty}{\sqrt{\log N_{[\;]}\left(\eta,\mathcal{F},||\cdot||_{L_2(P)}\right)}d\eta}<\infty$.
\end{enumerate}
\end{hyp}

Assumptions \ref{as:entropy} and \ref{as:vc} are exactly the same as the conditions often imposed with i.i.d. data to show uniform LLNs and CLTs \citep[see, e.g., Theorems 19.4, 19.5, 19.13 and 19.14 in ][]{vanderVaart2000}.\footnote{In \cite{vanderVaart2000}, the supremum in Assumptions~\ref{as:entropy} and \ref{as:vc} is taken over the set of probability measures $Q$ with finite support on $\mathcal{Y}$ and such that $||F||_{Q,2}>0$. This additional restriction is simply due to a different convention in constructing covering numbers, as \citeauthor{vanderVaart2000}  considers open balls while we use closed balls, following, e.g., \cite{Kato2019}.} In particular, Assumption \ref{as:vc}-(i) (resp. (ii))  imposes a condition on what is usually referred to as the uniform (resp. bracketing) entropy integral, see, e.g., \cite{vanderVaartWellner1996}. Finiteness of the uniform entropy integral is satisfied by any VC-type class of functions \citep[see][for a definition]{cherno2014}, or by the convex hull of such classes under some restrictions. The bracketing entropy integral is finite for instance for classes of monotone or H\"older continuous functions \citep[see, e.g.][]{vanderVaartWellner1996}.

The following theorem establishes uniform LLNs and CLTs under these two conditions. We  denote by $\unb$ the $k-$tuple $(1,k+1,...,2k-1)$.

\begin{thm}\label{thm:unif}
	Suppose that Assumptions \ref{as:dgp}-\ref{as:measurability} hold. Then:
\begin{enumerate}
	\item If Assumption \ref{as:entropy} holds, $\sup_{f\in \mathcal{F}}\left|\mathbb{P}_{n}f-Pf\right|$ tends to 0 a.s. and in $L^1$.
	\item If Assumption \ref{as:vc} holds, the process $\mathbb{G}_n$ converges weakly in $\ell^\infty(\mathcal{F})$ to a centered Gaussian process $\mathbb{G}$ on $\mathcal{F}$ as $n$ tends to infinity. Moreover, the covariance kernel $K$ of $\mathbb{G}$ satisfies:
$$K(f_1,f_2) =  \frac{1}{(k-1)!^2} \sum_{(\pi,\pi')\in \mathfrak{S}(\{\un\})\times \mathfrak{S}(\{\unb\})}\Cov\left(f_1(Y_{\pi(\un)}),f_2(Y_{\pi'(\unb)})\right).$$
\end{enumerate}	
\end{thm}

The proof is in Section \ref{sub:proof_thm_unif} of the supplement. When Assumption \ref{as:entropy}-(ii) holds, Part 1 can be proved by essentially combining Theorem 3 in \cite{Eaglesonweber78} and Lemma 7.35 in \cite{kallenberg05}. Part 2 was also proved for a finite $\mathcal{F}$ by \cite{silverman1976limit}. But the weak convergence result under the bracketing entropy condition, and the uniform laws under the uniform entropy conditions, do not follow from such results. To prove the former, we adapt a maximal inequality in \citeauthor{GineNickl2015} (2015, see their Lemma 3.5.12) to our context. To this end, we show that Hoeffding's bound on U-statistic \citep[][Section 5.a]{Hoeffding1963} still applies to our context.

\medskip
To prove the results under the uniform entropy conditions, the key ingredient, as with i.i.d. data, is a symmetrization lemma stated in Appendix \ref{sec:proof_lem_sym} below and proved in the supplement. Its proof relies extensively on Lemma \ref{lem:repr} and a decoupling inequality that may be of independent interest (see Lemma \ref{lem:coup}). The latter result generalizes a similar inequality for U-processes \citep[see][]{de1992decoupling}. In the proofs of both lemmas, we follow similar strategies as with U-processes, with two complications.  First, even with $k=2$, $Y_{\i}$ does not only depend on $U_{i_1}$ and $U_{i_2}$, but also on $U_{\{i_1,i_2\}}$. Second, when $k\geq 3$, dependence between observations arises not only because of single-unit terms such as $U_{i_1}$ or $U_{i_2}$, but also because of multiple-unit terms such as $U_{\{i_1,i_2\}}$. 

\medskip
As in the i.i.d. case, Assumption \ref{as:entropy} is actually stronger than necessary to obtain the uniform law of large numbers. The following proposition gives an exact characterization, where, for simplicity, we restrict to $k=2$. It is similar to the characterization for i.i.d. data \citep[see, e.g. Theorem 3.7.4 in][]{GineNickl2015} or for U-processes \citep[see Theorem 5.2.2 in][]{delapena1999}. Let us introduce the following norms:
\begin{align*}
	\norm{f}_{1,1} & = \frac{1}{n}\sum_{i_1=1}^n\left|\frac{1}{n-1}\sum_{i_2\neq i_1}f(Y_{i_1,i_2}) + f(Y_{i_2,i_1})\right|, \\
\norm{f}_{1,2} & =\frac{1}{n(n-1)}\sum_{1\leq i_1 < i_2\leq n}\left|\E\left[f(Y_{i_1,i_2})+f(Y_{i_2,i_1})\mid U_{\{i_1,i_2\}}\right]\right|.
\end{align*}

\begin{prop}\label{prop:charact_LLN}
	Suppose that Assumptions \ref{as:dgp}-\ref{as:measurability} hold and $\mathcal{F}$ admits an envelop $F$ with $PF<\infty$. Then  $\sup_{f\in \mathcal{F}}\left|\mathbb{P}_{n}f-Pf\right| \convAS 0$ if and only if both
$\log N(\varepsilon,\mathcal{F},||\cdot||_{1,2})/n^2$ and $\log N(\varepsilon,\mathcal{F},||\cdot||_{1,1})/n$ tend to $0$ in outer probability.\footnote{For a definition of convergence in outer probability or outer almost-sure convergence considered below, see e.g. Chapter 1.9 in \cite{vanderVaartWellner1996}.}
\end{prop}

Proposition \ref{prop:charact_LLN} emphasizes the two aspects of dissociated, exchangeable arrays. The first is i.i.d. variations, through the random entropy term related to $||\cdot||_{1,2}$, which only involves $(U_{\{i_1,i_2\}})_{\i\in\I_{n,2}}$. The second is U-statistic like variations, through the random entropy term related to $||\cdot||_{1,1}$: up to negligible terms, $||f||_{1,1}$ only depends on $(U_{i_1})_{1\leq i_1\leq n}$.
Key in establishing the necessity of these two conditions is a weak converse of the symmetrization lemma for $k=2$, see Equation \eqref{eq:desym} in the supplement. 


\subsection{Convergence of the bootstrap process}
\label{sub:boot}

We now study the properties of the following bootstrap sampling scheme, which extends the pigeonhole bootstrap \citep{mccullagh2000resampling,owen2007pigeonhole} to jointly separable arrays:
\begin{enumerate}
	\item $n$ units are sampled independently in $\{1,...,n\}$ with replacement and equal probability. $W_{i}$ denotes the number of times unit $i$ is sampled.
	\item the $k-$tuple $\i=(i_1,...,i_k) \in \I_{n,k}$ is then selected $W_{\i}=\prod_{j=1}^kW_{i_j}$ times in the bootstrap sample.
\end{enumerate}
Then we consider $\mathbb{P}^{\ast}_{n}$ and $\mathbb{G}_{n}^{\ast}$, defined on $\mathcal{F}$ by
$$\mathbb{P}_{n}^{\ast}f=\frac{(n-k)!}{n!}\sum_{\i \in \I_{n,k}} W_{\i} f(Y_{\i}),$$
$$\mathbb{G}^{\ast}_{n}f=\sqrt{n} \left(\mathbb{P}_{n}^{\ast}f - \mathbb{P}_{n}f \right).$$
Asymptotic validity of the bootstrap amounts to showing that conditional on the data $(Y_{\i})_{\i \in \I_k}$, $\mathbb{G}_n^{\ast}$ converges weakly to the process $\mathbb{G}$ defined in Theorem \ref{thm:unif}.\footnote{For the sake of brevity, we focus afterwards on convergence results under the sole uniform entropy condition (Assumption \ref{as:vc}-(i)).} 
As discussed in, e.g., \citeauthor{vanderVaartWellner1996} (1996, Chapter 3.6),  the outer almost-sure conditional weak convergence boils down to proving
\begin{equation}
\sup_{h\in \text{BL}_1} \left|\E\left(h(\mathbb{G}_n^{\ast})\big|(Y_{\i})_{\i \in \I_k}\right)-\E\left(h(\mathbb{G})\right)\right|\convASstar 0,	
\label{eq:conv_boot_proc}
\end{equation}
where $\text{BL}_1$ is the set of bounded and Lipschitz functions from $\ell^{\infty}(\mathcal{F})$ to $[0,1]$ and ``$\convASstar$'' denotes outer almost-sure convergence.

\begin{thm}\label{thm:unifboot}
	If Assumptions \ref{as:dgp}-\ref{as:measurability} and  \ref{as:vc}-(i) hold, the process $\mathbb{G}^{\ast}_n$ converges weakly in $\ell^\infty(\mathcal{F})$ to $\mathbb{G}$, conditional on $(Y_{\i})_{\i \in \I_k}$ and outer almost surely.
\end{thm}

This theorem ensures the asymptotic validity of the bootstrap above not only for sample means, but also for smooth functionals of the empirical cdf and nonlinear estimators, as we shall see below. The proof of Theorem \ref{thm:unifboot}, in Section \ref{ssec:boot_gen_case}  of the supplement,  follows the same lines as that of Theorem \ref{thm:unif}, though some of the corresponding steps are more involved, as often with the bootstrap. In particular, to prove pointwise convergence, we use arguments in Lindeberg's proof of the CLT for triangular arrays, Theorem \ref{thm:unif}.1 and Urysohn's subsequence principle, combined with Prohorov's theorem.

\medskip
Note that in contrast with the standard bootstrap for i.i.d. data,
$$\E\left(\mathbb{P}^{\ast}_{n}(f)\big| (Y_{\i})_{\i \in \I_k}\right)= \frac{1}{n^k}\sum_{\i \in \I_{n,k}}f(Y_{\i})\neq \P_nf.$$
However, the difference between $\P_n$ and $\P'_n$, the empirical measure with weights $1/n^k$, becomes negligible as $n\rightarrow\infty$. Accordingly, we also show in the proof of Theorem \ref{thm:unifboot} the almost-sure conditional convergence of $\sqrt{n} \left(\mathbb{P}_{n}^{\ast}f - \mathbb{P}'_{n}f \right)$, in addition to that of $\G^*_n$.

\subsection{Application to nonlinear estimators}
\label{sub:applis}

Theorem \ref{thm:unif} ensures consistency and asymptotic normality of a large class of estimators. In turn, Theorem \ref{thm:unifboot} shows that using the bootstrap for such estimators is asymptotically valid. To illustrate these points, we consider here two popular classes of estimators, namely Z-estimators and smooth functionals of the empirical cdf.  Similar results could be obtained for, e.g., M-estimators \citep[see, e.g.][]{cheng2010} or generalized method of moments estimators \citep[see, e.g.][]{hansen1982large}.

\medskip
Let us first consider Z-estimators. Let $\Theta$ denote a normed space, endowed with the norm $\|\cdot\|_{\Theta}$ and let $(\psi_{\theta,h})_{(\theta,h)\in \Theta\times \mathcal{H}}$ denote a class of real, measurable functions.  Let $\Psi(\theta)(h)=P\psi_{\theta,h}$, $\Psi_n(\theta)(h)=\P_n\psi_{\theta,h}$ and $\Psi^*_n(\theta)(h)=\P^*_n\psi_{\theta,h}$. We let, for any real function $g$ on $\mathcal{H}$, $\|g\|_{\mathcal{H}}=\sup_{h\in \mathcal{H}} |g(h)|$. The parameter of interest $\theta_0$, which satisfies $\Psi(\theta_0)=0$, is estimated by $\widehat{\theta}=\arg\min_{\theta\in\Theta} \|\Psi_n(\theta)\|_{\mathcal{H}}$. We also define $\widehat{\theta}^*=\arg\min_{\theta\in\Theta} \|\Psi^*_n(\theta)\|_{\mathcal{H}}$ as the bootstrap counterpart of $\widehat{\theta}$. The following theorem  extends Theorem 13.4 in \cite{Kosorok2006} to jointly exchangeable and dissociated arrays. For related results on Z-estimators in the i.i.d. case, see Section 3.2 in  \cite{vanderVaartWellner1996} and \cite{wellner1996}.

\begin{thm}\label{thm:Z_est}
	Suppose that Assumption \ref{as:dgp} holds and:
	\begin{enumerate}
		\item $\|\Psi(\theta_m)\|_{\mathcal{H}}\rightarrow 0$ implies $\|\theta_m-\theta_0\|_{\Theta} \rightarrow 0$ for every $(\theta_m)_{m\in \N}$ in $\Theta$;
		\item The class $\{\psi_{\theta, h}: (\theta,h)\in \Theta\times \mathcal{H}\}$ satisfies Assumptions \ref{as:measurability}-\ref{as:entropy}, with the envelope function $F$ satisfying $PF<\infty$;
		\item There exists $\delta>0$ such that the class $\{\psi_{\theta, h}: \|\theta-\theta_0\|_{\Theta}<\delta, h\in\mathcal{H}\}$  satisfies Assumptions \ref{as:measurability} and \ref{as:vc}, with an envelope function $F_\delta$ satisfying $PF^2_\delta<\infty$;
		\item $\lim_{\theta\rightarrow \theta_0} \sup_{h\in \mathcal{H}} P\left(\psi_{\theta,h}-\psi_{\theta_0,h}\right)^2=0$;
		\item $\|\Psi_n(\widehat{\theta})\|_{\mathcal{H}}=o_p(n^{-1/2})$ and $P\left(\|\sqrt{n}\Psi^*_n(\widehat{\theta}^*)\|_{\mathcal{H}}>\eta | (Y_{\i})_{\i\in\I_k}\right)=o_p(1)$ for every $\eta>0$;
		\item $\theta\mapsto \Psi(\theta)$ is Fr\'echet-differentiable at $\theta_0$, with continuously invertible derivative $\dot{\Psi}_{\theta_0}$.
	\end{enumerate}
	Then $\sqrt{n}(\widehat{\theta}-\theta_0)$ converges in distribution to a centered Gaussian process $\G$. Moreover, conditional on $(Y_{\i})_{\i\in\I_k}$ and almost surely, $\sqrt{n}(\widehat{\theta}^*-\widehat{\theta})$ converges in distribution to $\G$.
\end{thm}

Next, we consider smooth functionals of $F_Y$, the cdf of $Y_{\i}$. Suppose that $\mathcal{Y}\subset \R^p$ for some $p\in \N^+$ and $\theta_0=g(F_Y)$, where $g$ is Hadamard differentiable \citep[for a definition, see, e.g.,][Section 3.9.1]{vanderVaartWellner1996}. We estimate $\theta_0$ with $\widehat{\theta}=g(\widehat{F_Y})$, where $\widehat{F_Y}$ denotes the empirical cdf of  $(Y_{\i})_{\i \in\I_{n,k}}$. Finally, we let $\widehat{\theta}^*$ denote the bootstrap counterpart of $\widehat{\theta}$.

\begin{thm}\label{thm:smooth}
Suppose that $g$ is Hadamard differentiable at $F_Y$ tangentially to a set $\mathbb{D}_0$, with derivative equal to $g'_{F_Y}$. Suppose also that Assumption \ref{as:dgp} holds. Then:
\begin{enumerate}
	\item $\sqrt{n}(\widehat{F_Y}-F_Y)$ converges weakly, as a process indexed by $y$, to a Gaussian process $\mathbb{G}$ with kernel $K$ satisfying	
\begin{align*}
  K(y_1,y_2) =&  \frac{1}{(k-1)!^2} \sum_{(\pi,\pi')\in \mathfrak{S}(\{\un\})\times \mathfrak{S}(\{\unb\})} \Cov\big(\mathds{1}_{\{Y_{\pi(\un)}\leq y_1\}}, \mathds{1}_{\{Y_{\pi'(\unb)}\leq y_2\}}\big).
\end{align*}
	\item If $\mathbb{G} \in \mathbb{D}_0$ with probability one,
$$\sqrt{n} \left(\widehat{\theta}-\theta_0\right) \convD \mathcal{N}(0,\V(g'_{F_Y}(\mathbb{G}))).$$
Moreover, conditional on $(Y_{\i})_{\i\in\I_k}$ and almost surely, $\sqrt{n}(\widehat{\theta}^*-\widehat{\theta})$ converges in distribution to the same limit.
\end{enumerate}
\end{thm}

In practice, $\mathbb{D}_0$ often corresponds to the set of functions that are continuous everywhere or at a certain point $y_0$. This is the case for instance with $g:F_Y\mapsto F_Y^{-1}(\tau)$ for $\tau \in (0,1)$. In such cases, one can show that $\mathbb{G} \in \mathbb{D}_0$ under the same condition as for i.i.d. data, namely that $F_Y$ is continuous everywhere or at the point $F_Y^{-1}(\tau)$.

\section{Extensions}\label{sec:extensions}

We now consider several extensions to our main results. First, we study the asymptotic behavior of the properly normalized empirical process in degenerate cases where $K(f,f)=0$. Second, we establish additional results on the bootstrap. Third, we study separately, rather than jointly, separable arrays. Other extensions to arrays with multiple observations per $k$-tuple and arrays where $Y_{\i}$ is defined even if there are identical indices in $\i$ are considered in the supplement. We also develop therein a test that the data are in fact i.i.d.

\subsection{Degenerate cases} 
\label{sub:degenerate}

We consider here situations where $K(f,f)=0$ for all $f\in\mathcal{F}$, focusing for simplicity on $k=2$.\footnote{If $K(f,f)=0$ for only some $f\in\mathcal{F}$, we focus on $\mathcal{F}'=\{f\in \mathcal{F}: K(f,f)=0\}$.} Such a degeneracy appears for instance if the variables in the array are actually i.i.d., in which case $\sqrt{n}\G_n$ converges to a Gaussian process with covariance kernel $K(f_1,f_2)=\Cov(f_1(Y_{1,2}),f_2(Y_{1,2}))$. As another example \citep[see][]{menzel2017,bretagnolle1983}, suppose that  $Y_{i_1,i_2}=X_{i_1}X_{i_2}$, with $(X_i)_{i\in\N^+}$ i.i.d. variables with $\E(X_1)=0$, $\V(X_1)=1$. Let also  $\mathcal{F}=\{f_\lambda(x)=\lambda x, \lambda\in I\}$ for a compact $I\subset\R$. Then one can easily see that $\sqrt{n}\G_n$ converges weakly in $\ell^\infty(\mathcal{F})$ to $\G(f_\lambda)=\lambda (Z^2 - 1) $, with $Z$ a standard normal variable.

\medskip
More generally and as with U-processes \citep[see, e.g.][]{ArconesGine1993}, when $K(f,f)=0$, the rate of convergence of $\P_nf-Pf$ is $n^{-1}$ rather than $n^{-1/2}$ and the asymptotic distribution may not be normal. For any $(i_1,i_2)\in\I_2$, let $Y_{i_1,i_2}=\tau(U_{i_1},U_{i_2}, U_{\{i_1,i_2\}})$ be the Aldous-Hoover-Kallenberg representation where, without loss of generality, the variables in $\tau(\cdot,\cdot,\cdot)$ are assumed to be uniform on $[0,1]$. Let  $\psi_m(u)=\left(1+\mathds{1}_{\{m\geq 2\}}\right)^{1/2}\cos\left(m\pi u\right)$ for $m$ even and $\psi_m(u)=\sqrt{2}\sin((m+1)\pi u)$ for $m$ odd. Then $(\psi_m)_{m\in\N}$ forms an orthonormal basis of $L^2[0,1]$. For all $\bm{m}\in\N^3$ and any $f\in\mathcal{F}$, we define  $\mu_{\bm{m}}(f)$ by
$$\mu_{\bm{m}}(f)=\E\left[\left[f(Y_{1,2})-\E\left(f(Y_{1,2})\right)\right]\psi_{m_1}(U_{1})\psi_{m_2}(U_{2}) \psi_{m_3}(U_{\{1,2\}})\right].$$
Let $(Z_m)_{m \in \N^+}$, $(Z_{m_1,m_2})_{(m_1,m_2)\in\N\times\N^+}$ and $(Z_{\{m_1,m_2\},m_3})_{(m_1,m_2,m_3)\in\N^2\times\N^+:m_1<m_2}$   denote independent standard normal variables. We then define the process $\G^d$ on $\mathcal{F}$ by
\begin{align*}\G^d(f)=& 
	\sum_{(m_1,m_2) \in \N^{+2}}\mu_{m_1,m_2,0}(f) \left(Z_{m_1}Z_{m_2}-\mathds{1}_{\{m_1=m_2\}}\right) \\
	+ &\sum_{\substack{m_1 \in\N,\\ m_2 \in\N^+}}\mu_{m_1,m_1,m_2}(f) Z_{m_1,m_2} +\sum_{\substack{(m_1,m_2,m_3) \in \N^2\times\N^+: \\ m_1\neq m_2}}\mu_{\bm{m}}(f) Z_{\{m_1,m_2\},m_3}. 
\end{align*}
To prove the convergence of $\sqrt{n}\G_n$, we consider a condition on $\mathcal{F}$ that slightly differs from Assumption \ref{as:vc}-(i).

\begin{hyp}\label{as:vc_degen}
	The class $\mathcal{F}$ admits an envelope $F$ with $PF^2<\infty$ and
	\begin{equation*}
		\int_0^{\infty}\sup_{Q\in\mathcal{Q}}\log N\left(\eta||F||_{Q,2},\mathcal{F},||\cdot||_{Q,2}\right)d\eta<\infty.
	\end{equation*}
\end{hyp}

Assumption~\ref{as:vc_degen} is more stringent than Assumption~\ref{as:vc}-(i). A similar condition was also imposed by \cite{ArconesGine1993} for degenerate U-processes of order 1, see their condition (5.1). 

\begin{thm}\label{thm:degenerate}
		Suppose that $k=2$, Assumptions \ref{as:dgp}-\ref{as:measurability} and \ref{as:vc_degen} hold and $K(f,f)=0$ for all $f\in \mathcal{F}$. Then $\sqrt{n}\mathbb{G}_n$ converges weakly in $\ell^\infty(\mathcal{F})$ to $\G^d$.	
\end{thm}

As with degenerate U-processes \citep[see Section 5 of][]{ArconesGine1993}, the limit process is a Gaussian chaos process. The result is based in particular on a symmetrization lemma and a maximal inequality taylored to these degenerate cases. Specifically, the symmetrized process only includes Rademacher variables at the pair $\{i_1,i_2\}$ level, or products $\eps_{i_1}^{(1)}\eps_{i_1}^{(2)}$ of Rademacher variables. We refer to Lemmas S\ref{lem:sym3} and S\ref{lem:max_ineq_deg_case} in the supplement for more details.

\medskip
Finally, we note that the bootstrap process considered above does not generally converge to $\G^d$.\footnote{The same holds true for  the multiplier bootstrap process considered below.} With i.i.d. data, for instance, one can show that the variance of the bootstrapped mean  converges to $3\V(Y_{i_1,i_2})$. We expect similar phenomena as with U statistics, where the bootstrap is known to fail in degenerate cases \citep{Arcones1992, Arcones1994}. In the close case of separately exchangeable arrays (see Section \ref{sub:unbalanced} below), \cite{menzel2017} shows that a suitable wild bootstrap is consistent for the sample average, whether or not we have degeneracy. Whether such a result generalizes to the empirical process is left for future research.


\subsection{Further results on the bootstrap} 
\label{ssub:further_results_on_the_bootstrap}

Theorem \ref{thm:unifboot} shows convergence of the bootstrap process under conditions on $\mathcal{F}$ that ensure the convergence of the initial process $\G_n$. The following result shows that under moment conditions, convergence of $\G_n$ is actually necessary for the convergence of $\G_n^*$ to a Gaussian process.

\begin{thm}\label{thm:CNboot}
	Suppose that Assumptions \ref{as:dgp}-\ref{as:measurability} hold, $Pf^2<\infty$ for all $f\in\mathcal{F}$ and $\mathcal{F}$ admits an envelope $F$ such that $PF^{1+\delta}<\infty$ for some $\delta>0$. Then, if conditional on $(Y_{\i})_{\i \in \I_k}$ and outer almost surely, the process $\mathbb{G}^{\ast}_n$ converges weakly in $\ell^\infty(\mathcal{F})$ to $\mathbb{G}$, a centered Gaussian process, the process $\mathbb{G}_n$ also converges weakly in $\ell^\infty(\mathcal{F})$ to $\mathbb{G}$.
\end{thm}

Theorem \ref{thm:CNboot} may be seen as a partial extension to jointly exchangeable arrays of Theorem 2.4 in \cite{gine1990}, which, with i.i.d. data, establishes the equivalence between the convergence of the bootstrap process and $PF^2<\infty$ together with convergence of the initial process.

\medskip
With i.i.d. data, several other bootstrap schemes than the multinomial bootstrap are possible: see, e.g., \cite{barbe1995} for an extensive review. The situation is probably no different with jointly exchangeable arrays. To illustrate this, we consider a version of the multiplier bootstrap  adapted to such data \citep[see, e.g.,][for the case of i.i.d. data]{Kosorok2003}. Specifically, let $(\xi_i)_{i=1}^n$ be a sequence of i.i.d. random variables that are centered, have unit variance and are independent from the original data $(Y_{\i})_{\i\in\I_{n,2}}.$ We then consider the following process:
$$\mathbb{G}_n^{m*} : f \mapsto \frac{1}{\sqrt{n}}\sum_{i_1=1}^n\xi_{i_1}\left(\frac{1}{n-1}\sum_{1\leq i_2\neq i_1 \leq n} \left[f(Y_{i_1,i_2})+ f(Y_{i_2,i_1})\right]- 2 \mathbb{P}_nf\right).$$

The next theorem shows the conditional weak convergence of $\mathbb{G}_n^{m*} $ under the same conditions on $\mathcal{F}$ as previously.

\begin{thm}\label{thm:multiplier}
	Suppose that Assumptions \ref{as:dgp}-\ref{as:measurability} and  \ref{as:vc}-(i) hold and $(\xi_i)_{i=1}^n$ is i.i.d. with $\E(\xi_1)=0$, $\V(\xi_1)=1$. Then, conditional on $(Y_{\i})_{\i \in \I_k}$ and outer almost surely, the process $\mathbb{G}^{m\ast}_n$ converges weakly in $\ell^\infty(\mathcal{F})$ to $\mathbb{G}$.
\end{thm}


\subsection{Separately exchangeable arrays}
\label{sub:unbalanced}

Up to now, we have considered cases where the $n$ units that interact stem from the same population. In some cases, however, they do not, because the $k$ populations differ. For instance, we may be interested only in relationships between men and women. In that case, the symmetry condition in Assumption \ref{as:dgp} has to be strengthened:  both the labelling of men and the labelling of women should be irrelevant. This corresponds to so-called separately exchangeable arrays, defined formally in Assumption \ref{as:dgp3} below. Another important motivation for considering separately exchangeable arrays is multiway clustering, namely dependence arising through different dimensions of clustering. For instance, wages of workers may be affected by local shocks or sector-of-activity shocks. In such cases, we observe $Y_{i_1,i_2}$, the wage of a worker in geographical area $i_1$ and sector of activity $i_2$.\footnote{Oftentimes, we actually have several observations per cell, and the number varies from one cell to another. This extension is discussed in Section \ref{sub:heterog} of the supplement.}

\medskip
More generally, we consider in this section random variables $Y_{\i}$ where $\i=(i_1,...,i_k)\in \N^{+k}$, implying that repetitions (e.g. $\i=(1,...,1)$) are allowed. We impose the following condition on these random variables.

\begin{hyp}\label{as:dgp3}
For any $(\pi_1,...,\pi_k)\in \mathfrak{S}(\mathbb{N}^+)^k$, $$(Y_{\i})_{\i \in \N^{+k}}\overset{d}{=}(Y_{\pi_1(i_1),...,\pi_k(i_k)})_{\i \in \N^{+k}}.$$ Moreover, for any $A,B$, disjoint subsets of $\mathbb{N}^+$, $(Y_{\i})_{\i \in A^k}$ is independent of $(Y_{\i})_{\i \in B^k}$.
\end{hyp}

This condition is stronger than Assumption \ref{as:dgp} since it implies in particular equality in distribution for $\pi_1=...=\pi_k$.

\medskip
Let us redefine $\un$ here as $(1,...,1)$ and let $\n=(n_1,...,n_k)$, where $n_j\geq 1$ denotes the number of units observed in population $j$ (or cluster $j$ with multiway clustering). Note that in general,  $n_j \neq n_{j'}$ for $j\neq j'$. The sample at hand is then $(Y_{\i})_{\un\leq \i\leq \n}$, where $\i\geq \i'$ means that $i_j\geq i'_j$ for all $j=1,...,k$. Let $\underline{n}=\min(n_1,...,n_k)$. The empirical measure and empirical process that we consider for separately exchangeable arrays are:
\begin{align*}
\P_{\n}f & = \frac{1}{\prod_{j=1}^k n_j}\sum_{\un \leq \i \leq \n} f(Y_{\i}),\\
\G_{\n}f & = \sqrt{\underline{n}}\left(\P_{\n}f - Pf\right).
\end{align*}

We also consider the ``pigeonhole bootstrap'', suggested by \cite{mccullagh2000resampling} and studied, in the case of the sample mean and for particular models, by \cite{owen2007pigeonhole}. This bootstrap scheme is very close to the one we considered in Section \ref{sec:gen_results} for jointly exchangeable arrays, except that the weights are now independent from one coordinate to another:
\begin{enumerate}
	\item For each $j\in\{1,...,k\}$, $n_j$ elements are sampled with replacement and equal probability in the set $\{1,...,n_j\}$. For each $i_j$ in this set, let $W^j_{i_j}$ denote the number of times $i_j$ is selected this way.
	\item The $k$-tuple $\i=(i_1,...,i_k)$ is then selected $W_{\i}=\prod_{j=1}^k W^j_{i_j}$ times in the bootstrap sample.
\end{enumerate}

The bootstrap process $\G_{\n}^{\ast}$ is thus defined on $\mathcal{F}$ by
$$\G^{\ast}_{\n}f=\sqrt{\underline{n}} \left(\frac{1}{\prod_{j=1}^k n_j}\sum_{\un \leq \i \leq \n}\left(W_{\i} - 1\right)\sum_{\ell=1}^{N_{\i}}f(Y_{\i,\ell})\right).$$

Henceforth, we consider the convergence of $\P_{\n}$, $\G_{\n}$ and $\G^*_{\n}$ as $\underline{n}$ tends to infinity. More precisely, as with multisample U-statistics \citep[see, e.g.][Section 12.2]{vanderVaart2000}, we assume that there is an index $m\in \N^+$, left implicit hereafter, and increasing functions $g_1,...,g_k$ such that for all $j$, $n_j=g_j(m)\rightarrow \infty$ as $m\rightarrow \infty$ (we also assume without loss of generality that for all $m\in\N^+$, $g_j(m+1)>g_j(m)$ for some $j$). The following theorem extends Theorems \ref{thm:unif} and \ref{thm:unifboot} to this set-up.

\begin{thm} \label{thm:unif3}
	Suppose that Assumptions \ref{as:measurability} and \ref{as:dgp3} hold and that for every $j=1,...,k$, there exists $\lambda_j\geq 0$ such that $\underline{n}/n_j\rightarrow \lambda_j\geq 0$. Then:
\begin{enumerate}
	\item If Assumption \ref{as:entropy} holds,
	$\sup_{f\in \mathcal{F}}\left|\P_{\n}f -Pf\right|$ tends to 0 a.s. and in $L^1$.
	\item  If Assumption \ref{as:vc}-(i) holds, the process $\mathbb{G}_n$ converges weakly in $\ell^\infty(\mathcal{F})$ to a centered Gaussian process $\mathbb{G}_\lambda$ on $\mathcal{F}$ as $n$ tends to infinity. Moreover, the covariance kernel $K_\lambda$ of $\mathbb{G}_\lambda$ satisfies:
\begin{equation}
	K_\lambda(f_1,f_2) = \sum_{j=1}^k\lambda_j \Cov\left(f_1(Y_{\un}),f_2(Y_{\deux})\right),
\label{eq:cov_sep}	
\end{equation}
where $\deux$ is the $k$-tuple with 2 in each entry but 1 in entry $j$.
	\item If Assumption \ref{as:vc}-(i) holds, the process $\mathbb{G}^{\ast}_n$ converges weakly to $\mathbb{G}_\lambda$, conditional on $(Y_{\i})_{\i \in \N^{+k}}$ and outer almost surely.
\end{enumerate}
\end{thm}

Theorem \ref{thm:unif3} includes the case where $\lambda_j=0$ for some $j$, corresponding to ``strongly unbalanced'' designs with different rates of convergence to $\infty$ along the different dimensions of the array. In that case, only the dimensions with the slowest rate of convergence contribute to the asymptotic distribution, as can be seen in \eqref{eq:cov_sep}.

\medskip
Because the $(n_j)_{j=1...k}$ are not all equal in general, Theorem \ref{thm:unif3} does not follow directly from Theorem \ref{thm:unif}, even if Assumption \ref{as:dgp3} is stronger than Assumption \ref{as:dgp}. We prove the result by showing a simpler and convenient version of the symmetrization lemma in this setting. We refer to Lemma S\ref{lem:sym2} in the supplement for more details.

\section{Applications to international trade}
\label{sec:application}

Finally, we illustrate the importance of accounting for dependence in real dyadic data, through two applications to international trade data.

\subsection{Evolution of international trade} 
\label{sub:evolution_of_international_trade}

There is a large interest in economics on the evolution of international trade. But before analyzing the causes and consequences of such an evolution, one must check that there is indeed some significant changes. In this first application, we test whether the distribution of exports remains the same between two consecutive years, using Comtrade data on all countries from 2012 to 2018.  We use for that purpose the Kolmogorov-Smirnov (KS) test statistic
$$KS_t=\sup_{u\in\R}\left|\frac{1}{n(n-1)}\sum_{(i_1,i_2)\in\I_{n,2}} \mathds{1}_{\{T_{i_1,i_2,t}\leq u\}} - \mathds{1}_{\{T_{i_1,i_2,t+1}\leq u\}}\right|.$$
where $T_{i_1,i_2,t}$ denotes the trade volume from country $i_1$ to country $i_2$ in year $t$. Let us assume that Assumption \ref{as:dgp} holds, with $Y_{\i}=(T_{\i,t}, T_{\i,t+1})$. Then, under the null hypothesis that the distributions of $T_{\i,t}$ and $T_{\i,t+1}$ are equal, we have, by Theorem \ref{thm:unif}, $\sqrt{n} KS_t	\convD \|\mathbb{G}\|_{\mathcal{F}}$, with $\mathcal{F}=\{f_u(x,y)=\mathds{1}_{\{x\leq u\}}-\mathds{1}_{\{y\leq u\}}\}$. Given the dependence structure both between pairs of countries and across time, the distribution of $\|\mathbb{G}\|_{\mathcal{F}}$ depends on the true data generating process. To estimate it, we rely on the recentered bootstraped test statistic:
$$KS^*_t=\sup_{u\in\R}\left|\frac{1}{n(n-1)}\sum_{(i_1,i_2)\in\I_{n,2}} (W_{\i} -1) \left(\mathds{1}_{\{T_{i_1,i_2,t}\leq u\}}- \mathds{1}_{\{T_{i_1,i_2,t+1}\leq u\}}\right)\right|.$$
We compute the p-value of the test by $\P\left(KS^*_t>KS_t \big| (Y_{\i})_{\i \in\I_{n,k}}\right)$. For the sake of comparison, we also compute  p-values based on alternative forms of dependence that have been considered in applied work on similar data. Specifically, we also assume that the variables $(Y_{\i})_{\i}$ are i.i.d. We then assume pairwise clustering, where $Y_{i_1,i_2}$ and $Y_{i_2,i_1}$ may be dependent, but $Y_{\i}$ and $Y_{\j}$ are independent if $\j$ is not a permutation of $\i$. We also consider one-way clustering according to $i_1$ (and, similarly, according to $i_2$). In this case, $Y_{i_1,i_2}$ and $Y_{i_1,i_3}$ may be dependent, but $Y_{i_1,i_2}$ and $Y_{i'_1,i_3}$ are independent as soon as $i_1\neq i'_1$, whether or not $i_2=i_3$. For each of these cases, we use the bootstrap, but with different bootstrap schemes accounting for these different dependence structures.

\begin{table}[H]
	\begin{center}
		\begin{tabular}{lcccccc}
			Pairs of & KS test & \multicolumn{5}{c}{p-values under different assumptions} \\
			years & statistic & i.i.d. & P.W. cl. & E. cl. & I. cl. & dyadic \\
			\hline
2012-2013 & 0.048 & $<0.001$ & $<0.001$ & $<0.001$ & $<0.001$ & $<0.001$ \\
2013-2014 & 0.018 & $<0.001$ & $<0.001$ & $<0.001$ & 0.026 & 0.038 \\
2014-2015 & 0.022 & $<0.001$ & $<0.001$ & $<0.001$ & 0.005 & 0.007 \\
2015-2016 & 0.002 & 0.44 & 0.391 & 0.377 & 0.951 & 0.998 \\
2016-2017 & 0.012 & $<0.001$ & $<0.001$ & $<0.001$ & 0.215 & 0.254 \\
2017-2018 & 0.045 & $<0.001$ & $<0.001$ & $<0.001$ & $<0.001$ & $<0.001$ \\
			\hline
			\multicolumn{7}{p{300pt}}{{\scriptsize Notes: data from the Comtrade database. ``cl.'', ``E'', ``I'' and ``P.W.''  stand for clustering, exporter, importer and pairwise, respectively. The p-values were obtained with 1,000 bootstrap samples.}}
		\end{tabular}
	\end{center}
	\caption{KS tests of $F_{T_{\i,t}}=F_{T_{\i,t+1}}$ under different dependence assumptions}
	\label{tab:appli1}
\end{table}

\medskip
The results are displayed in Table \ref{tab:appli1}. They suggest significant changes in export volumes in some years but not all. In particular, international trade seems very stable between 2015 and 2017. There is some evidence of changes between 2012 and 2015 but we still do not reject the null hypothesis at the 1\% level for the years 2013-2014. The other columns of the table shows the importance of accounting for dependence along both dimensions. In particular, assuming i.i.d. data or pairwise dependence always leads to a strong rejection of the null, except for 2015-2016.\footnote{\label{foot:iid_test} A concern is that if the data are actually i.i.d. (or, more generally, pairwise dependent), our bootstrap is conservative, which  would explain the discrepancy between the p-values under  pariwise dependence and non-degenerate joint exchangeability.  Using the methodology in Section \ref{sub:test_of_i_i_d_data} of the supplement, we test for pairwise dependence.  For the eight years we consider, the null hypothesis is rejected at all standard levels, with p-values always smaller than $10^{-4}$.} Clustering along exporters also leads to artificially small p-values, in particular for the pairs 2013-2014, 2014-2015 and 2016-2017. In this context, clustering along importers leads to results that are closer to those based on dyadic data.


\subsection{Estimation of a gravity equation} 
\label{sub:estimation_of_a_gravity_equation}

Second, we revisit \cite{silva2006log}, who estimate the so-called gravity equation for international trade. Omitting the year index, this gravity equation states that $T_{i_1,i_2}$ satisfies
\begin{equation}
T_{i_1,i_2} = \exp(\alpha_0) G_{i_1}^{\alpha_1} G_{i_2}^{\alpha_2} D_{i_1,i_2}^{\alpha_3} \exp\left(A_{i_1,i_2}\beta\right) \eta_{i_1,i_2},
\label{eq:gravity}
\end{equation}
where $G_{i}$ denotes country $i$'s GDP, which would correspond to the mass of $i$ in a traditional gravity equation, $D_{i_1,i_2}$ denotes the distance between $i_1$ and $i_2$, $A_{i_1,i_2}$ are additional control variables and $\eta_{i_1,i_2}$ is an unobserved term.

\medskip
To estimate $\theta_0 = (\alpha_0,...,\alpha_3, \beta')'$, 
%
\cite{silva2006log} suggest to use the Poisson pseudo maximum likelihood (PPML for short) estimator $\widehat{\theta}$. The idea, formalized in \cite{gourieroux1984pseudo}, is that with i.i.d data, the PPML estimator is consistent and asymptotically normal for $\theta_0$ even if $T_{\i}$ does not follow a Poisson model, provided that $\E\left[\eta_{\i}|X_{\i}\right]=1$, with $X_{\i}=(1,\ln(G_{i_1}),\ln(G_{i_2}),\ln(D_{\i}),A_{\i})$. This is because the PPML estimator is based on the empirical counterpart of
\begin{equation}\label{eq:moment_PPML}
  \E\left[X'_{\i}\left(T_{\i} - \exp(X_{\i}\theta_0)\right) \right] = 0,
\end{equation}
and this equality holds true if $\E\left[\eta_{\i}|X_{\i}\right]=1$. 

\medskip
Now, assuming as in \cite{silva2006log} that the variables $(Y_{\i})_{\i \in \I_2}$  (with $Y_{\i}=(T_{\i},X_{\i})$) are i.i.d. is restrictive. We suppose instead that Assumption \ref{as:dgp} holds. Then Theorem \ref{thm:Z_est} applies to this setting, implying that $\widehat{\theta}$ is still consistent and asymptotically normal in this case.\footnote{In this case, $\mathcal{H}=\{1,...,\text{dim}(X_{\i})\}$ and $\psi_{\theta,h}(Y_{\i})=X_{h,\i}(T_{\i} - \exp(X_{\i}\theta_0))$. Then the key conditions 2 and 3 in Theorem \ref{thm:Z_est} are satisfied as soon as $\Theta$ is bounded, see e.g. Example 19.7 in \cite{vanderVaart2000}.} Nonetheless, the rates of convergence and asymptotic variance are different in the two cases, resulting in different inference on $\theta_0$.\footnote{The same application has been considered by \cite{graham2018}, who shows,  assuming convergence of a certain sample average, the asymptotic normality of the PPML estimator under the same dependence structure as ours. On the other hand, he neither considers bootstrap-based inference nor proves the consistency of his (asymptotic) variance estimator.}

\medskip
We use the same dataset as \cite{silva2006log}, which covers 136 countries for year 1990, and consider the exact same specification as the one they use in their Table 3. In this specification, the additional control variables $A_{\i}$ include exporter- and importer-level variables, namely their GDP per capita, a dummy variable equal to one if countries are landlocked and a remoteness index, which is the log of GDP-weighted average distance to all other countries. It also includes variables at the pair level, namely dummy variables for contiguity, common language, colonial tie, free-trade agreement and openness. This openness dummy is equal to one if at least one country is part of a preferential trade agreement. We refer to \cite{silva2006log} for additional details.

\medskip
Table \ref{tab:appli} below presents the results. The first column displays the point estimates, which, as expected, are identical to those in \cite{silva2006log}. The other columns display the p-values for the null hypothesis that $\theta_{0j}$, the $j$-th component of $\theta_0$, is equal to 0. We consider the same forms of dependence as with the KS test above. Under joint exchangeability, we compute the p-value $p_j$ for $\theta_{0j}=0$ using $p_j=\mathbb{P}\left(|\widehat{\theta}_j^* - \widehat{\theta}_j|>|\widehat{\theta}_j| \big| (Y_{\i})_{\i \in\I_{n,k}}\right)$. For other forms of dependence, we follow the usual practice of computing the p-values using the asymptotic normality of $\widehat{\theta}_j$ and estimators of the asymptotic variance under these various dependence structures.

\begin{table}[H]
	\begin{center}
		\begin{tabular}{lcccccc}
			& & \multicolumn{5}{c}{p-values under different assumptions} \\
			Variable & Estimator & i.i.d & P.W. cl. & E. cl. & I. cl. & dyadic \\
			\hline
Log(E's GDP)  & 0.732 & $<10^{-3}$ & $<10^{-3}$ & $<10^{-3}$ & $<10^{-3}$ & $<10^{-3}$ \\
Log(I's GDP)  & 0.741 & $<10^{-3}$ & $<10^{-3}$ & $<10^{-3}$ & $<10^{-3}$ & $<10^{-3}$ \\
Log(E's PCGDP)  & 0.157 & 0.003 & $<10^{-3}$ & 0.04 & 0.001 & 0.078 \\
Log(I's PCGDP)  & 0.135 & 0.003 & $<10^{-3}$ & 0.004 & 0.055 & 0.076 \\
Log of distance  & -0.784 & $<10^{-3}$ & $<10^{-3}$ & $<10^{-3}$ & $<10^{-3}$ & $<10^{-3}$ \\
Contiguity   & 0.193 & 0.064 & 0.16 & 0.112 & 0.077 & 0.461 \\
Common-language   & 0.746 & $<10^{-3}$ & $<10^{-3}$ & $<10^{-3}$ & $<10^{-3}$ & 0.056 \\
Colonial-tie   & 0.025 & 0.867 & 0.902 & 0.891 & 0.882 & 0.952 \\
Landlocked E   & -0.863 & $<10^{-3}$ & $<10^{-3}$ & $<10^{-3}$ & $<10^{-3}$ & 0.004 \\
Landlocked I   & -0.696 & $<10^{-3}$ & $<10^{-3}$ & $<10^{-3}$ & $<10^{-3}$ & 0.011 \\
E's remoteness  & 0.66 & $<10^{-3}$ & $<10^{-3}$ & $<10^{-3}$ & $<10^{-3}$ & 0.036 \\
I's remoteness  & 0.562 & $<10^{-3}$ & $<10^{-3}$ & 0.003 & 0.004 & 0.105 \\
P-T agreement   & 0.181 & 0.041 & 0.117 & 0.054 & 0.122 & 0.456 \\
Openness   & -0.107 & 0.416 & 0.522 & 0.498 & 0.453 & 0.771 \\
			\hline
			\multicolumn{7}{p{340pt}}{{\scriptsize Notes: data from \cite{silva2006log}, same specification as in their Table 3. ``cl.'', ``E'', ``I'', ``PCGDP'', ``P-T'', ``P.W.'' stand for clustering, exporter, importer, per capita GPD, preferential-trade and pairwise, respectively. The p-values for the last column were obtained with 1,000 bootstrap samples.}}
		\end{tabular}
	\end{center}
	\caption{Point estimates of $\theta_0$ and p-values of $\theta_{0j}=0$ under different dependence assumptions}
	\label{tab:appli}
\end{table}

\medskip
Using our bootstrap leads to much larger p-values than under the i.i.d. assumption. Only the log of distance and the  log of GDP of the exporter and the importer appear to be significant at the $10^{-3}$ levels, whereas five additional control variables are significant at that level under the i.i.d. assumption. In particular, common language and importer's remoteness are not even significant at the usual 5\% level.\footnote{\label{foot:iid_test2} As in Footnote \ref{foot:iid_test} above, we test for pairwise dependence, to see whether our results could be driven by the fact that our bootstrap is conservative in such cases. We obtain a p-value smaller than $10^{-4}$ and thus reject this hypothesis at all usual levels.} Interestingly, there is also a gap between assuming one-way clustering, either at the exporter or at the importer level, and assuming to have a jointly exchangeable and dissociated array. In the former case, we still have seven variables that are significant at the $10^{-3}$ levels. Confidence intervals, not displayed here, lead to similar conclusions. In particular, compared to the average length of i.i.d.-based 95\% confidence intervals, those based on pairwise clustering are only 8\% wider. Those based on one-way clustering on exporters (resp. importers) are 20\% (resp. 17\%) larger. On the other hand, those based on Assumption \ref{as:dgp} are 136\% wider.


\section{Conclusion} 
\label{sec:conclusion}

While polyadic data are increasingly used in applied work, and empirical researchers routinely account for multiway clustering when computing standard errors, the statistical theory behind these forms of dependence has lagged behind. Following \cite{bickel2009nonparametric} and \cite{menzel2017}, we link these dependence structures to jointly and separately exchangeable arrays. Using representation results for such arrays, we then prove uniform laws of large numbers and central limit theorems. These results imply consistency and asymptotic normality of various nonlinear estimators under such dependence. We also establish the general validity of natural extensions of the standard nonparametric bootstrap to such arrays. Our application shows that using those bootstrap schemes may make a large difference compared to assuming i.i.d. data or clustering along a single dimension, as has often been done.

\medskip
One caveat is that for the bootstrap confidence intervals to be valid, the asymptotic variance of the estimator should be positive. This may not be the case, for instance if the data $(Y_{\i})_{\i\in \I_k}$ are actually i.i.d. Inference based on the wild bootstrap without this positivity condition has been studied for sample averages under multiway clustering by \cite{menzel2017}. How to conduct inference on nonlinear estimators under joint exchangeability or multiway clustering without this positivity condition remains an avenue for future research.


\appendix
\section{Key lemmas} 
\label{sec:proof_lem_sym}

We first state the symmetrisation lemma. Let $(\eps_{A})_{A\subset \mathbb{N}^+}$ denote Rademacher independent variables, independent of $\left(Y_{\i}\right)_{\i \in \I_k}$. Then:

\begin{lem} \label{lem:sym}
	Suppose that Assumptions \ref{as:dgp}-\ref{as:measurability} hold and $P|f|<\infty$ for all $f\in \mathcal{F}$. Then there exist real numbers $C_{1,k},...,C_{k,k}$ depending only on $k$ and $(Y_{\i}^{1})_{\i \in\I_k}$,..., $(Y_{\i}^{k})_{\i \in\I_k}$, jointly exchangeable and dissociated arrays with $Y_{\un}^{j}\overset{d}{=}Y_{\un}$ for all $j\in\{1,...,k\}$, satisfying
	\begin{align*}
	& \E \left[ \Phi \left(\sup_{f \in \mathcal{F}} |\mathbb{P}_{n}f - Pf|\right)\right] \\
	\leq & \frac{1}{k}\sum_{r=1}^k\frac{1}{|\mathcal{E}_r|}\sum_{\e'\in \mathcal{E}_r} \E \left[ \Phi \left(\frac{(n-k)!}{n!}C_{r,k} \sup_{f \in \mathcal{F}}\left|\sum_{\i \in \I_{n,k}}\eps_{\{\i \odot \e' \}^+} f\left(Y_{\i}^r\right)\right|\right)\right],
	\end{align*}
\end{lem}

Though more complicated than its i.i.d. version \citep[see e.g. Lemma 2.3.1 in][] {vanderVaartWellner1996}, it serves the exact same purpose in the proofs of Theorems \ref{thm:unif}-\ref{thm:unifboot}: conditional on the $\left(Y_{\i}^r\right)_{\i \in \I_k}$, the process $f \mapsto \sum_{\i \in \I_{n,k}}\eps_{\{\i \odot \e' \}^+} f\left(Y_{\i}^r\right)$ is sub-Gaussian. In view of the AHK representation, the terms $\eps_{\{\i \odot \e' \}^+}$ could be expected. Given the aforementioned link with U-statistics, Lemma \ref{lem:sym} can also be seen as a generalization of the symmetrization lemma for U-processes for non-degenerate cases, see in particular Theorem 3.5.3 in \cite{delapena1999}.

\medskip
The proof of Lemma~\ref{lem:sym} crucially hinges upon the following decoupling inequality, which may be of independent interest. Hereafter, we let  $\mathcal{A}_r=\{A\subseteq \{1,...n\}: |A|=r\}$.

\begin{lem}\label{lem:coup}
	Let $r\leq k$, $\left(W_{A}\right)_{A\in\mathcal{A}_r}$ be a family of i.i.d. random variables with values in a Polish space $\mathcal{S}$ and $\left(W^{(j)}_{A}\right)_{A\in\mathcal{A}_r}$, $j=1,...,|\mathcal{E}_r|$ be some independent copies of this family. Let $\Phi$ be a non-decreasing convex function from $\mathbb{R}^{+}$ to $\mathbb{R}$ and $\ell$ be a bijection from $\mathcal{E}_r$ to $\{1,...,|\mathcal{E}_r|\}$. Let $\mathcal{H}$  be a pointwise measurable class of functions from $\mathcal{S}^{|\mathcal{E}_r|}\times \I_{n,k}$ to $\mathbb{R}$ such that
$\E\left(\sup_{h\in \mathcal{H}}\left|h\left(\left(W_{\{\i \odot \e\}^+}\right)_{\e \in \mathcal{E}_r},\i\right)\right|\right)<\infty$. Finally, let $L_r= \left(3|\mathcal{E}_r|^{|\mathcal{E}_r|}\right)^{|\mathcal{E}_r|-1}$. Then
	\begin{align*}
& \E\Phi\left(\sup_{h\in \mathcal{H}}\left|\sum_{\i \in \I_{n,k}}h\left(\left(W_{\{\i \odot \e\}^+}\right)_{\e \in \mathcal{E}_r},\i\right)\right|\right) \\
\leq & \E\Phi\left(L_r\sup_{h\in \mathcal{H}}\left|\sum_{\i \in \I_{n,k}}h\left(\left(W^{(\ell(\e))}_{\{\i \odot \e\}^+}\right)_{\e \in \mathcal{E}_r},\i\right)\right|\right).
\end{align*}
\end{lem}

The proof is given in the supplement. This result generalizes the decoupling inequality for $U$-statistics of \cite{de1992decoupling} to our setting. As with $U$-statistics, it is possible to obtain a reverse inequality if $r\in\{1,k-1,k\}$ and $\pi \mapsto h\left(\left(W_{\{\i_\pi \odot \e\}^+}\right)_{\e \in \mathcal{E}_r},\i_\pi\right)$ is constant on $\mathfrak{S}_k$, for all $h\in \mathcal{H}$. With such a reverse inequality, it is possible to replace $Y_{\i}^r$ by $Y_{\i}$ in Lemma \ref{lem:sym}. It is unclear to us, however, whether this reverse inequality still holds if $r\not\in\{1,k-1,k\}$ (implying $k\geq 4$). The key argument for the reverse inequality in \cite{de1992decoupling}  is that by the symmetry condition above, we can replace $h\left(\left(W_{\{\i_\pi \odot \e\}^+}\right)_{\e \in \mathcal{E}_r},\i_\pi\right)$ by an average over $k!$ terms. However,  for the proof to extend to our setting, one would need an average over $|\mathcal{E}_r|!$ terms. This is not possible in general when $|\mathcal{E}_r|>k$, which is the case when $r\not\in\{1,k-1,k\}$.

\medskip
Next, in order to prove the convergence of the empirical process under the bracketing entropy condition (Assumption \ref{as:vc}-(ii)), we establish the following maximal inequality, which is very close to that of \cite{GineNickl2015} for i.i.d. data (see their Lemma 3.5.12). 

\begin{lem}\label{lem:max_ineq}
	Suppose that Assumption \ref{as:dgp} holds. Let $(f_j)_{1\leq j\leq N}$ be  real-valued functions and $\mathcal{F}=\{x\mapsto e f_j(x), e\in\{-1,1\}, j=1,...,N\}$. Then:
$$\E\left[\max_{f \in \mathcal{F}} |\G_n(f)|\right]\leq  2\sqrt{k\log 2N\max_{f\in\mathcal{F}} \V(f(X_1))} + \frac{4k\log 2N\max_{f\in \mathcal{F}} ||f||_{\infty}}{3\sqrt{n}}.$$
\end{lem}

\newpage
\bibliography{biblio}

\newpage


\section*{Supplementary material}

# Supplement to "Empirical Process Results for Exchangeable Arrays"

Laurent Davezies[*]  Xavier D'Haultfœuille [†]  Yannick Guyonvarch [‡]

**Abstract**

This supplement first presents additional extensions to our main results, and in particular the case of multiple observations per $k$-tuple. Second, Monte Carlo simulations illustrate the performance of our bootstrap scheme in finite samples. Third, it displays all the proofs of our results. Section 3 gathers the proofs of the main results in Section 2 of the paper, while Section 4 focuses on the proofs of the extensions. Section 5 collects all the technical lemmas.

## 1 Additional extensions

### 1.1 Multiple observations per $k$-tuple

In some cases, we observe multiple observations for the same $k-$tuple $\boldsymbol{i}$. For instance, in the case of exchanges in a network, we may observe multiple or no such exchanges between $i_1$ and $i_2$. In sport competitions, we may observe $N_{i_1,i_2}$ matches between players $i_1$ and $i_2$, with possibly $N_{i_1,i_2} = 0$. Similarly, in multiway clustering (see Section 3.3 of the main text), we very often have several individuals per "cell" $\boldsymbol{i}$. To deal with this issue, and focusing here for simplicity on the case of jointly exchangeable arrays,[1] we consider that for each $\boldsymbol{i} \in \mathbb{I}_k$, there exists a random variable $N_{\boldsymbol{i}}$ taking values in $\mathbb{N}$ and a sequence $Y_{\boldsymbol{i}} = (Y_{\boldsymbol{i},\ell})_{\ell \geq 1}$, with $Y_{\boldsymbol{i},\ell}$ having support $\mathcal{Y}$, such that we only observe $(N_{\boldsymbol{i}}, (Y_{\boldsymbol{i},\ell})_{1 \leq \ell \leq N_{\boldsymbol{i}}})$. To allow for $N_{\boldsymbol{i}} = 0$, we assume in the following that for any sequence $(a_\ell)_{\ell \geq 1}$, $\sum_{\ell=1}^{0} a_\ell = 0$.

In this set-up, it is often natural to redefine the parameters of interest: if the relevant units of observation are the $N_{\boldsymbol{i}}$ units within each $k-$tuple, then parameters of interest are defined with respect to $\widetilde{P}$ rather than $P$, with

$$\widetilde{P}f = \mathbb{E}\left[\sum_{\ell=1}^{N_{\boldsymbol{1}}} f(Y_{\boldsymbol{1},\ell})\right].$$

In the example of sport matches, this expectation weights equally each match rather than each pair of players and is therefore often more relevant.[2] For instance, the sample average

$$\widehat{\theta} = \frac{\sum_{\boldsymbol{i} \in \mathbb{I}_{n,k}} \sum_{\ell=1}^{N_{\boldsymbol{i}}} Y_{\boldsymbol{i},\ell}}{\sum_{\boldsymbol{i} \in \mathbb{I}_{n,k}} N_{\boldsymbol{i}}}$$

---

[*]CREST-ENSAE, laurent.davezies@ensae.fr

[†]CREST-ENSAE. xavier.dhaultfoeuille@ensae.fr

[‡]CREST-ENSAE. yannick.guyonvarch@ensae.fr

[1]Jointly separable arrays can be treated similarly.

[2]If the unit of interest were the tuple itself, one could consider instead $\widetilde{\widetilde{P}}f = \mathbb{E}\left[\sum_{\ell=1}^{N_{\boldsymbol{1}}} f(Y_{\boldsymbol{1},\ell})/N_{\boldsymbol{1}}\right]$. Our results below on $\widetilde{P}$ directly extend to $\widetilde{\widetilde{P}}$. But note that $\widetilde{\widetilde{P}}$ is not defined in the very common situation where $\mathbb{P}(N_{\boldsymbol{i}} = 0) > 0$.

is an estimator of $\theta_0 = \widetilde{P}(\text{Id})/\widetilde{P}(1)$, where Id denotes the identity function. This parameter also satisfies $\theta_0 = \int y d\widetilde{F}_Y(y)$, with $\widetilde{F}_Y(y) = \widetilde{P}(\mathbb{1}_{\{\cdot \leq y\}})/\widetilde{P}(1)$. Similarly, quantiles would be defined as $\theta_0 = \widetilde{F}_Y^{-1}(\tau)$ for some $\tau \in (0,1)$. More generally, any parameter related to the units within each $k-$tuple is defined with respect to $\widetilde{P}$ rather than $P$.

Accordingly, we study the behavior of $\widetilde{\mathbb{P}}_n$, $\widetilde{\mathbb{G}}_n$ and $\widetilde{\mathbb{G}}_n^*$ defined on $\mathcal{F}$ by:

$$\widetilde{\mathbb{P}}_n f = \frac{(n-k)!}{n!} \sum_{\boldsymbol{i} \in \mathbb{I}_{n,k}} \sum_{\ell=1}^{N_{\boldsymbol{i}}} f(Y_{\boldsymbol{i},\ell}),$$

$$\widetilde{\mathbb{G}}_n f = \sqrt{n}\left(\widetilde{\mathbb{P}}_n(f) - \widetilde{P}f\right),$$

$$\widetilde{\mathbb{G}}_n^* f = \sqrt{n}\frac{(n-k)!}{n!} \sum_{\boldsymbol{i} \in \mathbb{I}_{n,k}} (W_{\boldsymbol{i}} - 1)\sum_{\ell=1}^{N_{\boldsymbol{i}}} f(Y_{\boldsymbol{i},\ell}).$$

The following theorem shows that the previous results extend to this set-up with random $N_{\boldsymbol{i}}$, only up to adaptations of the moment conditions.

**Theorem S1.** *Suppose that Assumption 1 holds with $(N_{\boldsymbol{i}}, Y_{\boldsymbol{i}})$ in place of $Y_{\boldsymbol{i}}$, $\widetilde{P}1 > 0$ and Assumption 2 holds. Then:*

1. *If Assumption 3-(i) holds with $F$ also satisfying $\widetilde{P}F < \infty$, then $\sup_{f \in \mathcal{F}} \left|\widetilde{\mathbb{P}}_n f - \widetilde{P}f\right|$ tends to 0 almost surely and in $L^1$.*

2. *If $\mathbb{E}(N_{\boldsymbol{1}}^2) < \infty$ and Assumption 4-(i) holds with $F$ also satisfying $\mathbb{E}\left(N_{\boldsymbol{1}} \sum_{\ell=1}^{N_{\boldsymbol{1}}} F^2(Y_{\boldsymbol{1},\ell})\right) < \infty$, the process $\widetilde{\mathbb{G}}_n$ converges weakly in $\ell^\infty(\mathcal{F})$ to a centered Gaussian process $\widetilde{\mathbb{G}}$ on $\mathcal{F}$ as $n$ tends to infinity. Moreover, the covariance kernel $\widetilde{K}$ of $\widetilde{\mathbb{G}}$ satisfies:*

$$\widetilde{K}(f_1, f_2) = \frac{1}{(k-1)!^2} \sum_{(\pi,\pi') \in \mathfrak{S}(\{\boldsymbol{1}\}) \times \mathfrak{S}(\{\boldsymbol{1'}\})} \text{Cov}\left(\sum_{\ell=1}^{N_{\pi(\boldsymbol{1})}} f_1(Y_{\pi(\boldsymbol{1}),\ell}), \sum_{\ell=1}^{N_{\pi'(\boldsymbol{1'})}} f_2(Y_{\pi'(\boldsymbol{1'}),\ell})\right).$$

3. *Under the same condition as in 2., the process $\widetilde{\mathbb{G}}_n^*$ converges weakly to $\widetilde{\mathbb{G}}$, conditional on $(Y_{\boldsymbol{i}})_{\boldsymbol{i} \in \mathbb{I}_k}$ and outer almost surely.*

We assume that $(N_{\boldsymbol{i}}, Y_{\boldsymbol{i}})_{\boldsymbol{i} \in \mathbb{I}_k}$, rather than just $(Y_{\boldsymbol{i}})_{\boldsymbol{i} \in \mathbb{I}_k}$, satisfies Assumption 1. Importantly, however, this does not restrict the dependence between $N_{\boldsymbol{i}}$ and $Y_{\boldsymbol{i}}$, or between the $(Y_{\boldsymbol{i},\ell})_\ell$. Hence, conditional on $N_{\boldsymbol{i}}$, the correlation between $Y_{\boldsymbol{i},\ell}$ and $Y_{\boldsymbol{i},\ell'}$ may vary with $N_{\boldsymbol{i}}$, for instance. Note also that even if we focus on $\widetilde{P}$ rather than $P$ here, the conditions on $\mathcal{F}$ remain nearly unchanged, with only modifications of the moment conditions. For uniform LLNs, we simply replace $PF < \infty$ by $\widetilde{P}F < \infty$. For uniform CLTs, instead of replacing $PF^2 < \infty$ by $\widetilde{P}F^2 < \infty$, we require the slightly stronger conditions that $\mathbb{E}(N_{\boldsymbol{1}}^2) < \infty$ and $\mathbb{E}\left(N_{\boldsymbol{1}} \sum_{\ell=1}^{N_{\boldsymbol{1}}} F^2(Y_{\boldsymbol{1},\ell})\right) < \infty$. These conditions are nonetheless equivalent to $\widetilde{P}F^2 < \infty$ when $N_{\boldsymbol{1}}$ is bounded. Note also that with a finite $\mathcal{F}$, our proof would only require $\widetilde{P}F^2 < \infty$.

The proof of Theorem S1 is very similar to those of Theorems 2.1 and 2.2, with one difference. In those theorems, we use the symmetrization lemma to bound the fluctuations of $\mathbb{G}_n$ by a function of

the entropy of the class $\mathcal{F}$. Here, similarly, we bound the fluctuations of $\widetilde{\mathbb{G}}_n$ by a function of the entropy of the class

$$\widetilde{\mathcal{F}} = \left\{ \widetilde{f}(n, y_1, ..., y_n) = \sum_{\ell=1}^n f(y_\ell) : n \in \mathbb{N}, (y_1, ..., y_n) \in \mathcal{Y}^n; f \in \mathcal{F} \right\}.$$

The additional point to prove is that we can control the complexity of $\widetilde{\mathcal{F}}$ under Assumption 4 and the moment conditions above, even if Assumption 4 imposes conditions on $\mathcal{F}$ rather than on $\widetilde{\mathcal{F}}$ directly.

### 1.2 Including "diagonal" elements

In some cases, $Y_i$ may also be defined for $i \in \{1, .., n\}^k \setminus \mathbb{I}_{n,k}$, namely for $k$-tuples featuring identical indices. For instance in international trade, one may also consider sales of national firms in their own countries, corresponding to variables $Y_{i,i}$ for $i \in \mathbb{N}^+$. Let us first consider the case $k = 2$. We still impose Assumption 1, but this time on the array $(Y_{i,j})_{(i,j) \in \mathbb{N}^{+2}}$ instead of $(Y_{i,j})_{(i,j) \in \mathbb{I}_2}$. Then the empirical measures is simply

$$\mathbb{P}_n f = \frac{1}{n^2} \sum_{1 \le i_1, i_2 \le n} f(Y_{i_1, i_2}),$$

whereas the definition of $\mathbb{G}_n$ remains unchanged, with still $Pf = \mathbb{E}[f(Y_{1,2})]$. The bootstrap scheme would remain the same, except that $Y_{i,i}$ would appear with a weight $W_i$, as in the standard bootstrap with i.i.d. data.

To analyse $\mathbb{P}_n$, $\mathbb{G}_n$ and their bootsrap counterpart, we can simply cut them into two parts, the part excluding "diagonal" elements and that including only such elements. Since the $(Y_{i,i})_{i \in \mathbb{N}}$ are i.i.d. under Assumption 1, we can apply uniform LLN and CLT and their bootstrap counterpart to the "diagonal" part. Because we further divide by $1/n$ in $\mathbb{P}_n f$, this part tends to 0 in probability. Hence, Theorems 2.1 and 2.2 directly apply without any modifications to this setting.

With $k > 2$, a similar reasoning holds, except that we have to introduce jointly exchangeable arrays of smaller dimensions $k - j$, with $j \in \{1, ..., k-1\}$ corresponding to the number of repetitions in the $k$-tuples. We can then apply our theorems to those lower dimensional arrays. Because the corresponding averages are multiplied by terms tending to 0, as with $k = 2$, they are asymptotically negligible, and again our theorems apply without any changes.

### 1.3 A statistical test of independence

Let $(Y_{i_1, i_2})_{(i_1, i_2) \in \mathbb{I}_2}$ denote an array of random variables in $\mathbb{R}^d$ satisfying Assumption 1 with $k = 2$. We show in this section how to test that the $(Y_{i_1, i_2})_{(i_1, i_2) \in \mathbb{I}_2, i_1 < i_2}$ are i.i.d. Such a test may be appealing for two reasons. First, standard inference can be conducted on the data under this hypothesis. Second, it may indicate that inference based on our bootstrap is conservative, in case we do not reject the null hypothesis. Note that the restriction $i_1 < i_2$ allows one to have pairwise dependence, i.e. $Y_{i_1, i_2}$ and $Y_{i_2, i_1}$ may still be dependent.

To define the test statistic, let $\overline{Y}$ denote the sample average of the $(Y_{i_1,i_2})_{(i_1,i_2)\in\mathbb{I}_{n,2}}$ and

$$\widehat{\Sigma} = \frac{2}{n(n-1)} \sum_{i_1<i_2} (Y_{i_1,i_2} + Y_{i_2,i_1} - 2\overline{Y})(Y_{i_1,i_2} + Y_{i_2,i_1} - 2\overline{Y})'.$$

Then define $\widehat{Z}_{i_1,i_2} = \widehat{\Sigma}^{-1/2}\left(Y_{i_1,i_2} + Y_{i_2,i_1} - 2\overline{Y}\right)$, with $\widehat{Z}_{i_1,i_2} = (\widehat{Z}^1_{i_1,i_2}, ..., \widehat{Z}^d_{i_1,i_2})'$. The test statistic we consider is

$$T_n = \sum_{\ell=1}^{d} \max\left[0, \left(\frac{6}{n(n-1)(n-2)}\right)^{1/2} \sum_{i_1<i_2<i_3} \widehat{Z}^\ell_{i_1,i_2} \widehat{Z}^\ell_{i_1,i_3}\right]^2. \quad (1.1)$$

Let $F_\ell$ denote the cumulative distribution function of a $\chi^2(\ell)$ distribution, with the convention that $F_0(x) = \mathbb{1}_{x\geq 0}$. Then define $F(x) = \frac{1}{2^d} \sum_{\ell=0}^{d} \binom{d}{\ell} F_\ell(x)$ and $q(\tau) = \inf\{x : F(x) \geq \tau\}$. For testing the null hypothesis at the asymptotic level $1-\alpha$, we consider the test of critical region $W_\alpha = \{T_n > q(1-\alpha)\}$.

The idea behind the test is that if the $(Y_{i_1,i_2})_{(i_1,i_2)\in\mathbb{I}_2,i_1<i_2}$ are i.i.d., $\mathbb{E}\left(Z^\ell_{i_1,i_2} Z^\ell_{i_1,i_3}\right) = 0$ for all $\ell$, with $Z_{i_1,i_2}$ defined as $\widehat{Z}_{i_1,i_2}$ but with $\overline{Y}$ and $\widehat{\Sigma}$ replaced by $\mathbb{E}(Y_{1,2})$ and $\mathbb{V}(Y_{1,2})$, respectively. Thus, intuitively, $T_n$ remains bounded in probability in this case.[3] On the other hand, in non-degenerate cases,

$$\mathbb{E}\left(Z^\ell_{i_1,i_2} Z^\ell_{i_1,i_3}\right) = \mathbb{E}\left[\mathbb{E}(Z^\ell_{i_1,i_2}|U_{i_1})^2\right] > 0, \quad (1.2)$$

where the $(U_{i_1})_{i_1\in\mathbb{N}^+}$ are the random variables appearing in the AHK representation of $(Y_{i_1,i_2})_{(i_1,i_2)\in\mathbb{I}_2}$. As a result, we can expect $T_n$ to tend to infinity. The following proposition formalizes these ideas.[4]

**Proposition S 1.** *Suppose that Assumption 1 holds with $k=2$, $Y_{1,2} = (Y^1_{1,2},...,Y^d_{1,2}) \in \mathbb{R}^d$, $\mathbb{E}(|Y^\ell_{1,2}|^4) < \infty$ for all $\ell \in \{1,...,d\}$ and $\mathbb{V}(Y_{1,2})$ is positive. Then:*

1. *If the $(Y_{i_1,i_2})_{(i_1,i_2)\in\mathbb{I}_2,i_1<i_2}$ are i.i.d., $\lim_{n\to\infty}\Pr(W_\alpha) = \alpha$ for all $\alpha \in (0,1/2)$;*

2. *If (1.2) holds for some $\ell \in \{1,...,d\}$, $\lim_{n\to\infty}\Pr(W_\alpha) = 1$ for all $\alpha \in (0,1)$.*

We apply this test to the trade data we use in Section 4 of the paper. A common concern in our two applications is that if the data are i.i.d., our bootstrap procedure is conservative. This could then explain the discrepancy between the p-values based on our bootstrap and those based on assuming pairwise clustering. In the first application, $Y_{i_1,i_2} = T_{i_1,i_2,t}$ where we recall that $T_{i_1,i_2,t}$ are the exports from country $i_1$ to country $i_2$ during year $t \in \{2012,...,2017\}$. The results are clear-cut: with $T_n > 70$ for all the years we consider, the p-values are always smaller than $10^{-4}$. In the second application, one can show that that the variable that matters for inference on $\theta_0$ is $Y_{i_1,i_2} = X_{i_1,i_2}\left(T_{i_1,i_2} - \exp(X'_{i_1,i_2}\theta_0)\right)$. Estimating $Y_{i_1,i_2}$ by $\widehat{Y}_{i_1,i_2} = X_{i_1,i_2}\left(T_{i_1,i_2} - \exp(X'_{i_1,i_2}\widehat{\theta})\right)$, with $\widehat{\theta}$ the PPML estimator of $\theta_0$,[5] we obtain a large test statistic of $T_n \simeq 276.3$, corresponding once more to a p-value smaller than $10^{-4}$. So in both applications, we reject at all usual levels the null hypothesis above.

---

[3]We may have $\mathbb{E}(Z_{i_1,i_2} Z_{i_1,i_3}) = 0$ in other degenerate cases but we derive below the distribution of $T_n$ only if the $(Y_{i_1,i_2})_{(i_1,i_2)\in\mathbb{I}_2,i_1<i_2}$ are i.i.d.

[4]We impose here the existence of fourth moments. Approximation arguments as those used in the proof of Theorem 3.1 could be used to show the same result under finite second moments only, but at the price of lengthening the proof.

[5]Just as we show in the proof of Proposition S1 that using $\overline{Y}$ and $\widehat{\Sigma}$ instead of $\mathbb{E}(Y_{1,2})$ and $\mathbb{V}(Y_{1,2})$ does not have any effect on the asymptotic distribution of $T_n$, replacing $\theta_0$ by $\widehat{\theta}$ does not affect this asymptotic distribution.

# 2 Monte Carlo simulations

We investigate in this section the finite sample properties of the bootstrap scheme considered above, by studying the coverage probabilities of confidence intervals based on this bootstrap. We consider dyadic data satisfying Assumption 1, with $N_{\boldsymbol{i}} = 1$ for all $\boldsymbol{i} \in \mathbb{I}_2$, and the following dependence structure:

$$Y_{i_1,i_2} = 1 + \mu(\varepsilon_{1i_1} + \varepsilon_{2i_2}) + \sqrt{0.5 - \mu^2}\left(\nu \varepsilon^S_{i_1,i_2} + \sqrt{2-\nu^2}\varepsilon_{i_1,i_2}\right),$$

where the $(\varepsilon_{1i_1}, \varepsilon_{2i_1})_{i_1 \in \mathbb{N}^+}$, $(\varepsilon^S_{i_1,i_2})_{(i_1,i_2) \in \mathbb{I}_2}$ and $(\varepsilon_{i_1,i_2})_{(i_1,i_2) \in \mathbb{I}_2}$ are mutually independent and marginally all standard normal variables. We impose $\text{Corr}(\varepsilon_{1i_1}, \varepsilon_{2i_1}) = 0.8$ and $\varepsilon^S_{i_1,i_2} = \varepsilon^S_{i_2,i_1}$. The parameter $\mu \in [0, 1/\sqrt{2}]$ represents the importance of individual versus pair factors, whereas $\nu \in [0, \sqrt{2}]$ represents the importance of symmetric versus asymmetric shocks. In the baseline scenario, we let $(\mu, \nu) = (\sqrt{0.2}, 1)$. We also consider three other scenarios. In the first, $(\mu, \nu) = (\sqrt{0.2}, 0)$. In the second, $(\mu, \nu) = (0, 1)$, which implies that the limit process is degenerate. In the third, we use exponential(1) distributions instead of standard normal variables, with a Gaussian copula on $(\varepsilon_{1i_1}, \varepsilon_{2i_1})$ and still a correlation parameter of 0.8. Our parameter of interest $\theta_0$ is the median of $Y_{1,2}$, which is thus equal to 1. Hereafter, we study inference on $\theta_0$ based on the empirical median $\widehat{\theta}$, for $n \in \{10, 20, 40, 80\}$.

We inspect the performance of two different confidence intervals. The first is the symmetric interval $[\widehat{\theta} \pm q_{0.95}(|\widehat{\theta}^* - \widehat{\theta}|)]$, where $\widehat{\theta}^*$ denotes the bootstrap counterpart of $\widehat{\theta}$ and $q_\alpha(U)$ denotes the quantile of order $\alpha$ of $U$, conditional on the data $(Y_{\boldsymbol{i}})_{\boldsymbol{i} \in \mathbb{I}_{n,k}}$. The second is the percentile bootstrap interval $[q_{0.025}(\widehat{\theta}^*), q_{0.975}(\widehat{\theta}^*)]$. Given Theorem 2.4, both intervals are asymptotically valid.

Our results are displayed in Table 1. Our two confidence intervals have very good properties in the three non-degenerate scenarios, even for very small sample sizes. They appear to be slightly conservative for small $n$ and with normal variables, but the example of exponential distributions shows that this need not be the rule. In the degenerate scenario where $\mu = 0$, the confidence intervals are very conservative. Finally, the two confidence intervals are very close to each other.

| Scenario | n | Symmetric bootstrap CI | Percentile bootstrap CI |
|---|---|---|---|
| Baseline: | 10 | 0.984 | 0.986 |
| $\mu = \sqrt{0.2}$ | 20 | 0.977 | 0.979 |
| $\nu = 1$ | 40 | 0.969 | 0.971 |
| | 80 | 0.961 | 0.961 |
| Baseline | 10 | 0.98 | 0.983 |
| but $\nu = 0$ | 20 | 0.971 | 0.972 |
| | 40 | 0.965 | 0.968 |
| | 80 | 0.962 | 0.961 |
| Baseline | 10 | 0.996 | 0.997 |
| but $\mu = 0$ | 20 | 0.998 | 0.997 |
| | 40 | 0.999 | 0.998 |
| | 80 | 0.999 | 0.998 |
| Baseline | 10 | 0.939 | 0.935 |
| but exponential | 20 | 0.940 | 0.942 |
| variables | 40 | 0.943 | 0.942 |
| | 80 | 0.945 | 0.944 |

Notes: 5,000 simulations, 200 bootstrap samples for each.

Table 1: Coverage probabilities on the true median (nominal coverage: 95%)

We also consider Kolmogorov-Smirnov (KS) tests for the two-sample problem. Specifically, we are interested in testing whether two variables $Y_1^1$ and $Y_1^2$ have the same distribution, under the assumption that the array of variables $Y_i = (Y_i^1, Y_i^2)$ satisfy Assumption 1. The DGP that we consider for $Y_i^t$ is:

$$Y_{i_1,i_2}^t = \theta(t-1) + \varepsilon_{1i_1}^t + \varepsilon_{2i_2}^t + \varepsilon_{i_1,i_2}^t,$$

where $(\varepsilon_{1i_1}^t)_{i_1}$, $(\varepsilon_{2i_2}^t)_{i_2}$ and $(\varepsilon_{i_1,i_2}^t)_{i_1,i_2}$ are mutually independent. We also suppose that for $j \in \{1,2\}$, $(\varepsilon_{ji_1}^1, \varepsilon_{ji_1}^2) \sim \mathcal{N}(0, \Sigma)$ with $\Sigma_{11} = \Sigma_{12} = 1$, $\Sigma_{12} = 0.5$. Similarly, $(\varepsilon_{i_1,i_2}^1, \varepsilon_{i_1,i_2}^2) \sim \mathcal{N}(0, \Sigma)$. The null hypothesis therefore holds when $\theta = 0$, and not otherwise. To test for the null hypothesis, we rely on the Kolmogorov-Smirnov test. As explained in Section 4.1, the asymptotic distribution of the test statistic depends in a complicated way on the data generating process. Instead of trying to estimate it, we rely on the bootstrap.

Figure 1 shows the power of the bootstrap test as a function of $\theta$ and for different sample sizes, namely $n = 10, 20, 40$ and $80$, as above. The power curves are as expected. In particular, they increase quickly with $n$ when $\theta \neq 0$. The test is slightly conservative for $n \leq 40$ and rejects slightly too much for $n = 80$, but the true level remains close to the nominal for $\theta = 0$.

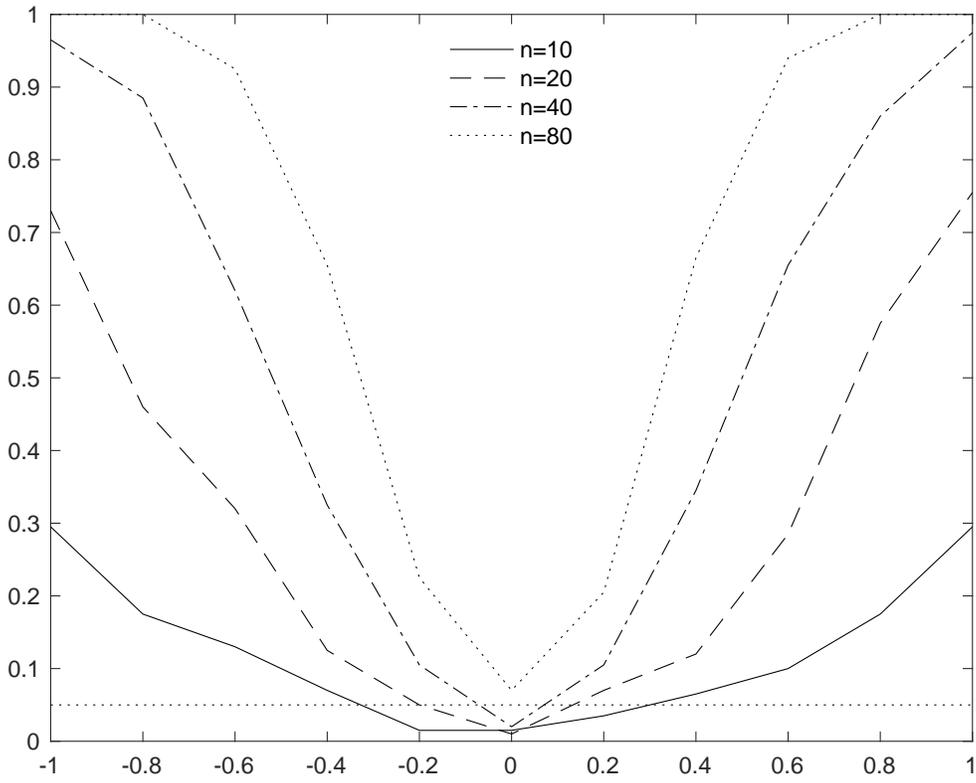

Notes: for each $n$ and $\theta \in \{-1, -0.8, ..., 1\}$, we used 500 samples and 200 bootstrap samples for each.

Figure 1: Power curves of the KS test for different sample sizes

## 3 Proofs of the main results

### 3.1 Notation

To ease the reading, we first summarize the notation we use throughout the proofs. Objects introduced in a single proof are defined therein directly and not reported here. We recall that $k$ denotes the dimension of the array of data. Also, bootstrap counterparts appear with a star.

**Subsets or elements of $\mathbb{N}^k$**

| | |
|---|---|
| $A^+$ | $A \cap (0, \infty)$, for any $A \subset \mathbb{R}$. |
| $\overline{A}$ | $\{\boldsymbol{i} \in A : i_j \neq i_{j'} \text{ if } j \neq j'\}$, for any $A \subset \mathbb{N}^{+k}$. |
| $\overrightarrow{A}$ | $\{\boldsymbol{i} \in \overline{A} : i_j < i_{j'} \text{ if } j < j'\}$ for any $A \subset \mathbb{N}^{+k}$. |
| $|A|$ | the cardinal of $A \subset \mathbb{N}^{+k}$. |
| $\mathfrak{S}(A)$ | The set of permutations on $A$. |
| $\mathfrak{S}_r$ | $\mathfrak{S}(\{1, ..., r\})$ |
| $\mathbb{I}_k$ | $\overline{\mathbb{N}^{+k}}$. |
| $\mathbb{I}_{n,k}$ | $\overline{\{1, ..., n\}^k}$. |
| $\mathcal{E}_r$ | $\{\boldsymbol{e} \in \{0; 1\}^k : \sum_{j=1}^k e_j = r\}$ for $r = 1, ..., k$. |
| $\boldsymbol{i}$ | element of $\mathbb{I}_k$ or $\mathbb{N}^{+k}$, with component $(i_1, ..., i_k)$. |
| $\{\boldsymbol{i}\}$ | the set of distinct elements of $\boldsymbol{i} = (i_1, ..., i_k) \in \mathbb{N}^k$. |
| $\boldsymbol{e}$ | element of $\{0, 1\}^k$. |

7

| | |
|---|---|
| $\boldsymbol{i^e}$ | for $\boldsymbol{i} \in \mathbb{I}_{n,r}$ and $\boldsymbol{e} \in \mathcal{E}_r$, the $k$-dimensional vector with component $i_1$ at the first non-null entry of $\boldsymbol{e}$, $i_2$ at the second non-null entry of $\boldsymbol{e}$ and so on.[6] |
| $\boldsymbol{0}$ | $(0,...,0)$ |
| $\boldsymbol{1}$ | $(1,...,k)$ except in Section 4.4 and Lemmas S2, S6 and S8, where $\boldsymbol{1} = (1,...,1)$. |
| $\boldsymbol{2}_r$ | element of $\mathbb{N}^k$ with 2 at each component but 1 at its $r$th component. |
| $\boldsymbol{i}_\pi$ | $(i_{\pi(1)}, ..., i_{\pi(r)})$, for any $\boldsymbol{i} \in \mathbb{N}^r$ and $\pi \in \mathfrak{S}_r$. |
| $\odot$ | the Hadamard product, i.e. $\boldsymbol{i} \odot \boldsymbol{e} = (i_1 e_1, ..., i_k e_k)$. |

**Sample and random variables**

| | |
|---|---|
| $n$ | Number of units in the population. |
| $\boldsymbol{n}$ | $(n_1, ... n_k)$, with $n_j$ the number of clusters in the $j$-th dimension in Section 4.4. |
| $\Pi_{\boldsymbol{n}}$ | $\prod_{j=1}^k n_j$. |
| $\widetilde{Y}_{\boldsymbol{i}}$ | $(N_{\boldsymbol{i}}, (Y_{\boldsymbol{i},\ell})_{\ell=1...N_{\boldsymbol{i}}})$ (see Section 1 above). |
| $(\varepsilon_A)_{A \in \mathcal{A}}$ | Mutually independent Rademacher random variables (i.e., with values 1 or $-1$ with probability $1/2$), for any set $\mathcal{A}$. |
| $(Y_{\boldsymbol{i}}^r)_{\boldsymbol{i} \in \mathbb{I}_k}$ | jointly exchangeable array defined in Lemma A.1 with marginal distribution $P$. |
| $(\widetilde{Y}_{\boldsymbol{i}}^r)_{\boldsymbol{i} \in \mathbb{I}_k}$ | same as $(Y_{\boldsymbol{i}}^r)_{\boldsymbol{i} \in \mathbb{I}_k}$, but when applying Lemma A.1 to $\widetilde{\mathcal{F}}$ and $(\widetilde{Y}_{\boldsymbol{i}})_{\boldsymbol{i} \in \mathbb{I}_k}$ instead of $\mathcal{F}$ and $(Y_{\boldsymbol{i}})_{\boldsymbol{i} \in \mathbb{I}_k}$. |
| $(\widetilde{Y}_{\boldsymbol{i}}^{r,r'})_{\boldsymbol{i} \in \mathbb{I}_k}$ | same as $(\widetilde{Y}_{\boldsymbol{i}}^r)_{\boldsymbol{i} \in \mathbb{I}_k}$. |
| $\boldsymbol{Y}_{\{\boldsymbol{i}\}}$ | $(Y_{i_1, i_2}, Y_{i_2, i_1})$. |

**Functions and classes of functions**

| | |
|---|---|
| Id | The identity function. |
| $\mathcal{D}$ | $\cup_{n \in \mathbb{N}} (\{n\} \times \mathcal{Y}^n)$. |
| $\mathcal{F}^s$ | $\{g : \exists f \in \mathcal{F} : g(x,y) = [f(x) + f(y)]/2\}$. |
| $\mathcal{F}^2$ | $\{f^2 : f \in \mathcal{F}\}$, for any class of functions $\mathcal{F}$. |
| $\mathcal{F} \times \mathcal{G}$ | $\{(f,g) : f \in \mathcal{F}, g \in \mathcal{G}\}$. |
| $\mathcal{F}_\delta$ | $\left\{h = f_1 - f_2 : (f_1, f_2) \in \mathcal{F} \times \mathcal{F}, \mathbb{E}\left[(f_1(Y_{\boldsymbol{1}}) - f_2(Y_{\boldsymbol{1}}))^2\right] \leq \delta^2\right\}$. |
| $\mathcal{F}_\infty$ | $\{h = f_1 - f_2 : (f_1, f_2) \in \mathcal{F} \times \mathcal{F}\}$. |
| $\widetilde{f}$ | for any function $f$ from $\mathcal{Y}$ to $\mathbb{R}$, the function from $\mathcal{D}$ to $\mathbb{R}$ defined by $\widetilde{f}(n, y_1, ..., y_n) = \sum_{\ell=1}^n f(y_\ell)$. |
| $\widetilde{\mathcal{F}}$ | $\left\{\widetilde{f} : f \in \mathcal{F}\right\}$. $\widetilde{\mathcal{F}}_\delta$ and $\widetilde{\mathcal{F}}_\infty$ are defined similarly. |
| $N(\eta, \mathcal{F}, \|\cdot\|)$ | the minimal number of $\|\cdot\|$-closed balls of radius $\eta$ with centers in $\mathcal{F}$ needed to cover $\mathcal{F}$. If $\|\cdot\|$ is random, $N(\eta, \mathcal{F}, \|\cdot\|)$ denotes the measurable cover of this minimal number. |
| $J_\mathcal{F}(u)$ | $\int_0^u \sup_Q \sqrt{\log N(\eta\|F\|_{Q,2}, \mathcal{F}, \|\cdot\|_{Q,2})} d\eta$, where the supremum is taken over the set of probability measures with finite support. |

---

[6]For instance if $k = 5$, $r = 3$, $\boldsymbol{i} = (6, 9, 2)$ and $\boldsymbol{e} = (0, 1, 1, 0, 1)$, we obtain $\boldsymbol{i^e} = (0, 6, 9, 0, 2)$.

**Probability measures and norms**

Note that we sometimes need to evaluate random variables at some specific value of the probability space. We denote by $\omega$ elements of this probability space $\Omega$.

$$
\begin{array}{ll}
Qf & \int f dQ, \text{ for any probability measure } Q. \\
P & \text{the probability distribution of } Y_i. \\
\mathbb{P}_n, \mathbb{P}'_n & \frac{(n-k)!}{n!} \sum_{\boldsymbol{i}\in\mathbb{I}_{n,k}} \delta_{Y_{\boldsymbol{i}}} \text{ and } \frac{1}{n^k} \sum_{\boldsymbol{i}\in\mathbb{I}_{n,k}} \delta_{Y_{\boldsymbol{i}}}, \text{ respectively.} \\
\mathbb{P}_n^* & \frac{(n-k)!}{n!} \sum_{\boldsymbol{i}\in\mathbb{I}_{n,k}} W_{\boldsymbol{i}} f(Y_{\boldsymbol{i}}), \text{ where } W_{\boldsymbol{i}} \text{ is the bootstrap weight of } \boldsymbol{i}. \\
\mathbb{P}_n^r & \frac{(n-k)!}{n!} \sum_{\boldsymbol{i}\in\mathbb{I}_{n,k}} \delta_{Y_{\boldsymbol{i}}^r} \\
||g||_{\mu,r} & (\int |g|^r d\mu)^{1/r} \text{ for } \mu \text{ a measure and } r \geq 1 \\
||f||_{e,M,1} & \frac{(n-k)!}{n!} \sum_{\boldsymbol{i}\in\overrightarrow{\mathbb{I}_{n,r}}} |\sum_{\pi\in\mathfrak{S}_r} \sum_{\boldsymbol{i}'\in\overline{(\{1,\ldots,n\}\setminus\{\boldsymbol{i}\})^{k-r}}} f\left(Y_{(\boldsymbol{i}_\pi)^e + \boldsymbol{i}'(1-e)}^r\right) \mathbb{1}_{\left\{F\left(Y_{(\boldsymbol{i}_\pi)^e + \boldsymbol{i}'(1-e)}^r\right) \leq M\right\}}|,
\end{array}
$$
with $f\in\mathcal{F}$, $F$ an envelope for $\mathcal{F}$ and $M>0$.

## 3.2 Lemma A.1

We proceed in four steps. First, we obtain an upper bound with a sum of differences that are identically distributed but not independent. Roughly speaking, they are nonetheless "less dependent", as we "decouple" the random variables appearing in the AHK representation (2.1) by introducing independent copies of them (see inequality (3.4) below). In the second step, using a telescopic sum, we further bound our expectation of interest by another one involving sums of differences that are independent, conditional on a suitable $\sigma$-algebra. The third step is the symmetrisation step itself, where Rademacher variables are introduced. The fourth step concludes by combining the previous steps. Note that the key decoupling inequality (3.4) is given separately in Lemma A.2, as it may be of independent interest.

**First step: decoupling**

For any $(r,j)\in\{1,...,k\}\times\mathbb{N}$, let $\left(U_A^{(j)}\right)_{A\subset\mathbb{N}^+:1\leq|A|\leq r}$ and $\left(V_A^{(j)}\right)_{A\subset\mathbb{N}^+:1\leq|A|\leq r}$ denote some independent copies of the $(U_A)_{A\subset\mathbb{N}^+:1\leq|A|\leq r}$. Let $Y_{\boldsymbol{i}}^{(k)} = \tau\left(\left(U_{\{\boldsymbol{i}\odot e\}^+}^{(0)}\right)_{e\in\cup_{j=1}^k \mathcal{E}_j}\right)$ and, for $r<k$,

$$Y_{\boldsymbol{i}}^{(r)} = \tau\left(\left(U_{\{\boldsymbol{i}\odot e\}^+}^{(0)}\right)_{e\in\cup_{j=1}^r \mathcal{E}_j}, \left(V_{\{\boldsymbol{i}\odot e\}^+}^{(0)}\right)_{e\in\cup_{j=r+1}^k \mathcal{E}_j}\right).$$

Because $\mathbb{E}[f(Y_1)] = \mathbb{E}\left[f\left(Y_{\boldsymbol{i}}^{(k)}\right)|Y_{\boldsymbol{i}}^{(0)}\right]$ and $(Y_{\boldsymbol{i}})_{\boldsymbol{i}\in\mathbb{I}_k} \stackrel{d}{=} \left(Y_{\boldsymbol{i}}^{(0)}\right)_{\boldsymbol{i}\in\mathbb{I}_k}$, we obtain, by Jensen's inequality and Lemma S1,

$$\mathbb{E}\left[\Phi\left(\sup_{f\in\mathcal{F}} \left|\frac{(n-k)!}{n!} \sum_{\boldsymbol{i}\in\mathbb{I}_{n,k}} f(Y_{\boldsymbol{i}}) - \mathbb{E}[f(Y_1)]\right|\right)\right]$$
$$\leq \mathbb{E}\left[\Phi\left(\sup_{f\in\mathcal{F}} \left|\frac{(n-k)!}{n!} \sum_{\boldsymbol{i}\in\mathbb{I}_{n,k}} f\left(Y_{\boldsymbol{i}}^{(0)}\right) - f\left(Y_{\boldsymbol{i}}^{(k)}\right)\right|\right)\right]$$
$$\leq \frac{1}{k}\sum_{r=1}^k \mathbb{E}\left[\Phi\left(k\sup_{f\in\mathcal{F}} \left|\frac{(n-k)!}{n!} \sum_{\boldsymbol{i}\in\in\mathbb{I}_{n,k}} f\left(Y_{\boldsymbol{i}}^{(r-1)}\right) - f\left(Y_{\boldsymbol{i}}^{(r)}\right)\right|\right)\right]. \tag{3.1}$$

9

For $\boldsymbol{i} \in \mathbb{R}^k$ and $\pi \in \mathfrak{S}_k$, let $\boldsymbol{i}_\pi = (i_{\pi(1)}, ..., i_{\pi(k)})$. For any $f \in \mathcal{F}$, let also

$$\overline{f}\left(\left(U_{\{i\odot e\}^+}\right)_{e\in\cup_{r=1}^k \mathcal{E}_r}\right) = \frac{1}{k!}\sum_{\pi\in\mathfrak{S}_k} f(Y_{\boldsymbol{i}_\pi}).$$

Note that $\sum_{\boldsymbol{i}\in\mathbb{I}_{n,k}} \overline{f}\left(\left(U_{\{i\odot e\}^+}\right)_{e\in\cup_{r=1}^k \mathcal{E}_r}\right) = \sum_{\boldsymbol{i}\in\mathbb{I}_{n,k}} f(Y_{\boldsymbol{i}})$ and if the components of $\boldsymbol{i}'$ are a permutation of those of $\boldsymbol{i}$ we have

$$\overline{f}\left(\left(U_{\{i\odot e\}^+}\right)_{e\in\cup_{r=1}^k \mathcal{E}_r}\right) = \overline{f}\left(\left(U_{\{i'\odot e\}^+}\right)_{e\in\cup_{r=1}^k \mathcal{E}_r}\right). \tag{3.2}$$

For $r = 1, ..., k$, let $\overline{\mathcal{E}}_r = \cup_{j=r+1}^k \mathcal{E}_j$ and $\underline{\mathcal{E}}_r = \cup_{j=1}^{r-1} \mathcal{E}_r$. Let $\mathcal{U}^r$ be the $\sigma$-algebra generated by the variables $(U^{(0)}_{\{i\odot e\}^+})_{(i,e)\in\mathbb{I}_{n,k}\times\underline{\mathcal{E}}_r}$ and $(V^{(0)}_{\{i\odot e\}^+})_{(i,e)\in\mathbb{I}_{n,k}\times\overline{\mathcal{E}}_r}$. For any $j \in \mathbb{N}$, $\boldsymbol{i} \in \mathbb{I}_{n,k}$ and $e \in \cup_{j'=1}^k \mathcal{E}_{j'}$, let $W^{(j)}_{\{i\odot e\}^+} = \left(U^{(j)}_{\{i\odot e\}^+}, V^{(j)}_{\{i\odot e\}^+}\right)$.

As we will reason conditional on $\mathcal{U}^r$, let us use $\overline{f}_{r,\boldsymbol{i}}(\boldsymbol{w})$ as a shortcut for

$$\overline{f}\left(\left(U^{(0)}_{\{i\odot e\}^+}\right)_{e\in\underline{\mathcal{E}}_r}, \boldsymbol{w}, \left(V^{(0)}_{\{i\odot e\}^+}\right)_{e\in\overline{\mathcal{E}}_r}\right),$$

for any vector $\boldsymbol{w} \in \mathbb{R}^{|\mathcal{E}_r|}$. Let us also define

$$\Delta\overline{f}_{r,\boldsymbol{i}}\left(\left(W^{(0)}_{\{i\odot e\}^+}\right)_{e\in\mathcal{E}_r}\right)$$
$$= k\frac{(n-k)!}{n!}\left[\overline{f}_{r,\boldsymbol{i}}\left(\left(U^{(0)}_{\{i\odot e\}^+}\right)_{e\in\mathcal{E}_r}\right) - \overline{f}_{r,\boldsymbol{i}}\left(\left(V^{(0)}_{\{i\odot e\}^+}\right)_{e\in\mathcal{E}_r}\right)\right].$$

Then, by definition of $Y_{\boldsymbol{i}}^{(r)}$ and $\Delta\overline{f}_{r,\boldsymbol{i}}$,

$$\mathbb{E}\left[\Phi\left(k\sup_{f\in\mathcal{F}}\left|\frac{(n-k)!}{n!}\sum_{\boldsymbol{i}\in\mathbb{I}_{n,k}} f\left(Y_{\boldsymbol{i}}^{(r-1)}\right) - f\left(Y_{\boldsymbol{i}}^{(r)}\right)\right|\right)\bigg|\mathcal{U}^r\right]$$
$$= \mathbb{E}\left[\Phi\left(\sup_{f\in\mathcal{F}}\left|\sum_{\boldsymbol{i}\in\mathbb{I}_{n,k}} \Delta\overline{f}_{r,\boldsymbol{i}}\left(\left(W^{(0)}_{\{i\odot e\}^+}\right)_{e\in\mathcal{E}_r}\right)\right|\right)\bigg|\mathcal{U}^r\right]. \tag{3.3}$$

Remark that the first result in Lemma A.2 applies conditional on $\mathcal{U}^r$. Then, letting $K_{1,r} = \left(3|\mathcal{E}_r||\mathcal{E}_r|\right)^{|\mathcal{E}_r|-1}$ and $\ell$ be an arbitrary bijection from $\mathcal{E}_r$ to $\{1, ..., |\mathcal{E}_r|\}$, we obtain

$$\mathbb{E}\left[\Phi\left(\sup_{f\in\mathcal{F}}\left|\sum_{\boldsymbol{i}\in\mathbb{I}_{n,k}} \Delta\overline{f}_{r,\boldsymbol{i}}\left(\left(W^{(0)}_{\{i\odot e\}^+}\right)_{e\in\mathcal{E}_r}\right)\right|\right)\bigg|\mathcal{U}^r\right]$$
$$\leq \mathbb{E}\left[\Phi\left(K_{1,r}\sup_{f\in\mathcal{F}}\left|\sum_{\boldsymbol{i}\in\mathbb{I}_{n,k}} \Delta\overline{f}_{r,\boldsymbol{i}}\left(\left(W^{(\ell(e))}_{\{i\odot e\}^+}\right)_{e\in\mathcal{E}_r}\right)\right|\right)\bigg|\mathcal{U}^r\right]. \tag{3.4}$$

**Second step: telescoping sum**

Let $\prec$ be a total order on $\mathcal{E}_r$. We note $\boldsymbol{e} \preceq \boldsymbol{e}'$ if $\boldsymbol{e} \prec \boldsymbol{e}'$ or $\boldsymbol{e} = \boldsymbol{e}'$. For every $(\boldsymbol{e}, \boldsymbol{e}') \in \mathcal{E}_r^2$ let

$$\overline{W}^{(\ell,\boldsymbol{e})}_{\{i\odot e'\}^+} = \begin{cases} \left(U^{(\ell(e'))}_{\{i\odot e'\}^+}, U^{(\ell(e'))}_{\{i\odot e'\}^+}\right) & \text{if } \boldsymbol{e}' \prec \boldsymbol{e} \\ \left(V^{(\ell(e'))}_{\{i\odot e'\}^+}, V^{(\ell(e'))}_{\{i\odot e'\}^+}\right) & \text{if } \boldsymbol{e}' \succ \boldsymbol{e} \\ \left(U^{(\ell(e'))}_{\{i\odot e'\}^+}, V^{(\ell(e'))}_{\{i\odot e'\}^+}\right) & \text{if } \boldsymbol{e}' = \boldsymbol{e}. \end{cases}$$

10

Then, for any $e \in \mathcal{E}_r$,

$$\left(\overline{W}^{(\ell,e)}_{\{i\odot e'\}+}\right)_{e'\in\mathcal{E}_r} = \left(U^{(\ell(e'))}_{\{i\odot e'\}+}\mathbb{1}_{\{e'\preceq e\}} + V^{(\ell(e'))}_{\{i\odot e'\}+}\mathbb{1}_{\{e'\succ e\}},\right.$$
$$\left. U^{(\ell(e'))}_{\{i\odot e'\}+}\mathbb{1}_{\{e'\prec e\}} + V^{(\ell(e'))}_{\{i\odot e'\}+}\mathbb{1}_{\{e'\succeq e\}}\right)_{e'\in\mathcal{E}_r}, \quad (3.5)$$

and $\left(\overline{W}^{(\ell,e')}_{\{i\odot e'\}+}\right)_{e'\in\mathcal{E}_r} = \left(W^{(\ell(e'))}_{\{i\odot e'\}+}\right)_{e'\in\mathcal{E}_r}$. $\Delta\overline{f}_{r,i}\left(\left(W^{(\ell(e'))}_{\{i\odot e'\}+}\right)_{e'\in\mathcal{E}_r}\right)$ can be decomposed into the following telescoping sum:

$$\Delta\overline{f}_{r,i}\left(\left(W^{(\ell(e'))}_{\{i\odot e'\}+}\right)_{e'\in\mathcal{E}_r}\right) = \sum_{e\in\mathcal{E}_r} \Delta\overline{f}_{r,i}\left(\left(\overline{W}^{(\ell,e)}_{\{i\odot e'\}+}\right)_{e'\in\mathcal{E}_r}\right).$$

By Lemma S1, we obtain, with $K_{2,r} = |\mathcal{E}_r|K_{1,r}$,

$$\mathbb{E}\left[\Phi\left(K_{1,r}\sup_{f\in\mathcal{F}}\left|\sum_{i\in\mathbb{I}_{n,k}}\Delta\overline{f}_{r,i}\left(\left(W^{\ell(e')}_{\{i\odot e'\}+}\right)_{e'\in\mathcal{E}_r}\right)\right|\right)\Big|\mathcal{U}^r\right]$$
$$\leq \frac{1}{|\mathcal{E}_r|}\sum_{e\in\mathcal{E}_r}\mathbb{E}\left[\Phi\left(K_{2,r}\sup_{f\in\mathcal{F}}\left|\sum_{i\in\mathbb{I}_{n,k}}\Delta\overline{f}_{r,i}\left(\left(\overline{W}^{(\ell,e)}_{\{i\odot e'\}+}\right)_{e'\in\mathcal{E}_r}\right)\right|\right)\Big|\mathcal{U}^r\right]. \quad (3.6)$$

**Third step: symmetrization**

For any $e \in \mathcal{E}_r$, let $\mathcal{U}^r_{\ell,e}$ be the $\sigma$-algebra generated by the same variables as $\mathcal{U}^r$, $(U^{(\ell(e'))}_{\{i\odot e'\}+})_{(i\times e')\in\mathbb{I}_{n,k}\times\underline{\mathcal{E}}_r:e'\prec e}$ and $(V^{(\ell(e'))}_{\{i\odot e'\}+})_{(i,e')\in\mathbb{I}_{n,k}\times\underline{\mathcal{E}}_r:e'\succ e}$. Let $\overrightarrow{\mathbb{I}_{n,k}} = \{(i_1, i_2, ..., i_k) \in \{1, ..., n\}^k : i_1 < i_2 < ... < i_k\} \subset \mathbb{I}_{n,k}$ and $\mathfrak{S}_k$ be the set of permutations of $\{1, ..., k\}$. For any $i = (i_1, ..., i_k) \in \mathbb{N}^k$ and $\pi \in \mathfrak{S}_k$, let $i_\pi$ denote $(i_{\pi(1)}, ..., i_{\pi(k)})$. For any $i \in \mathbb{I}_r$ and $e \in \mathcal{E}_r$, let $i^e$ be the $k$-dimensional vector with component $i_1$ in the first non-null entry of $e$, $i_2$ in the second non-null entry of $e$ and so on. Similarly, for any $i \in \mathbb{I}_{k-r}$ and $e \in \mathcal{E}_r$, let $i^{(1-e)}$ be the $k$-dimensional vector with component $i_1$ at the first null entry of $e$, $i_2$ at the second null entry of $e$ and so on. For instance, if $k = 5$, $r = 3$, $i = (6,9,2)$, $i' = (7,3)$ and $e = (0,1,1,0,1)$, we obtain $i^e = (0,6,9,0,2)$ and $i'^{(1-e)} = (7,0,0,3,0)$.

For every $e \in \mathcal{E}_r$, we have

$$\mathbb{I}_{n,k} = \left\{i^e_\pi + i'^{(1-e)} : i \in \overrightarrow{\mathbb{I}_{n,r}}, \pi \in \mathfrak{S}_r, i' \in \overline{(\{1,...,n\}\setminus\{i\})^{k-r}}\right\}. \quad (3.7)$$

Thus,

$$\sum_{i\in\mathbb{I}_{n,k}}\Delta\overline{f}_{r,i}\left(\left(\overline{W}^{(\ell,e)}_{\{i\odot e'\}+}\right)_{e'\in\mathcal{E}_r}\right)$$
$$= \sum_{i\in\overrightarrow{\mathbb{I}_{n,r}}}\sum_{i'\in\overline{(\{1,...,n\}\setminus\{i\})^{k-r}}}\sum_{\pi\in\mathfrak{S}_r}\Delta\overline{f}_{r,i^e_\pi+i'^{(1-e)}}\left(\left(\overline{W}^{(\ell,e)}_{\{(i^e_\pi+i'^{(1-e)})\odot e'\}+}\right)_{e'\in\mathcal{E}_r}\right).$$

With this new indexation of the sum on $i$ and reasoning conditional on $\mathcal{U}^r_{\ell,e}$, the triple sum above can be rewritten as a sum of $n!/[(n-r)!r!]$ symmetric and independent terms. Hence, it is equal in distribution to

$$\sum_{i\in\overrightarrow{\mathbb{I}_{n,r}}}\varepsilon_{\{i\}}\sum_{i'\in\overline{(\{1,...,n\}\setminus\{i\})^{k-r}}}\sum_{\pi\in\mathfrak{S}_r}\Delta\overline{f}_{r,i^e_\pi+i'^{(1-e)}}\left(\left(\overline{W}^{(\ell,e)}_{\{(i^e_\pi+i'^{(1-e)})\odot e'\}+}\right)_{e'\in\mathcal{E}_r}\right),$$

11

where the $(\varepsilon_A)_{A\subset\{1,\ldots,n\}}$ are i.i.d. Rademacher variables. For every $\boldsymbol{i} \in \overrightarrow{\mathbb{I}_{n,r}}$ and any $\pi \in \mathfrak{S}_r$, we have $\{\boldsymbol{i}\} = \{(\boldsymbol{i}^{\boldsymbol{e}}_\pi + \boldsymbol{i}'^{(\boldsymbol{1}-\boldsymbol{e})}) \odot \boldsymbol{e}\}^+$. Hence, using (3.7) again,

$$\sum_{\boldsymbol{i}\in\overrightarrow{\mathbb{I}_{n,r}}} \varepsilon_{\{\boldsymbol{i}\}} \sum_{\boldsymbol{i}'\in(\{1,\ldots,n\}\setminus\{\boldsymbol{i}\})^{k-r}} \sum_{\pi\in\mathfrak{S}_r} \Delta\overline{f}_{r,\boldsymbol{i}^{\boldsymbol{e}}_\pi+\boldsymbol{i}'^{(\boldsymbol{1}-\boldsymbol{e})}} \left(\left(\overline{W}^{(\ell,\boldsymbol{e})}_{\{(\boldsymbol{i}^{\boldsymbol{e}}_\pi+\boldsymbol{i}'^{(\boldsymbol{1}-\boldsymbol{e})})\odot \boldsymbol{e}'\}^+}\right)_{\boldsymbol{e}'\in\mathcal{E}_r}\right)$$
$$= \sum_{\boldsymbol{i}\in\mathbb{I}_{n,k}} \varepsilon_{\{\boldsymbol{i}\odot\boldsymbol{e}\}^+} \Delta\overline{f}_{r,\boldsymbol{i}} \left(\left(\overline{W}^{(\ell,\boldsymbol{e})}_{\{\boldsymbol{i}\odot\boldsymbol{e}'\}^+}\right)_{\boldsymbol{e}'\in\mathcal{E}_r}\right).$$

Furthermore, for every $\boldsymbol{e} \in \mathcal{E}_r$, by (3.5),

$$\frac{n!}{k(n-k)!}\Delta\overline{f}_{r,\boldsymbol{i}}\left(\left(\overline{W}^{(\ell,\boldsymbol{e})}_{\{\boldsymbol{i}\odot\boldsymbol{e}'\}^+}\right)_{\boldsymbol{e}'\in\mathcal{E}_r}\right)$$
$$=\overline{f}_{r,\boldsymbol{i}}\left(\left(U^{(\ell(\boldsymbol{e}'))}_{\{\boldsymbol{i}\odot\boldsymbol{e}'\}^+}\mathbb{1}_{\{\boldsymbol{e}'\preceq\boldsymbol{e}\}} + V^{(\ell(\boldsymbol{e}'))}_{\{\boldsymbol{i}\odot\boldsymbol{e}'\}^+}\mathbb{1}_{\{\boldsymbol{e}'\succ\boldsymbol{e}\}}\right)_{\boldsymbol{e}\in\mathcal{E}_r}\right)$$
$$-\overline{f}_{r,\boldsymbol{i}}\left(\left(U^{(\ell(\boldsymbol{e}'))}_{\{\boldsymbol{i}\odot\boldsymbol{e}'\}^+}\mathbb{1}_{\{\boldsymbol{e}'\prec\boldsymbol{e}\}} + V^{(\ell(\boldsymbol{e}'))}_{\{\boldsymbol{i}\odot\boldsymbol{e}'\}^+}\mathbb{1}_{\{\boldsymbol{e}'\succeq\boldsymbol{e}\}}\right)_{\boldsymbol{e}\in\mathcal{E}_r}\right).$$

Since for every $(j,j') \in \mathbb{N}^2$, $(U^{(j)}_A)_{A\subseteq\{1,\ldots,n\}}$ and $(V^{(j')}_A)_{A\subseteq\{1,\ldots,n\}}$ are equal in distribution and independent and $(U^{(j)}_A)_{A\subseteq\{1,\ldots,n\}} \perp\!\!\!\perp (U^{(j')}_A)_{A\subseteq\{1,\ldots,n\}}$ whenever $j \neq j'$, we obtain, conditional on $\mathcal{U}^r$,

$$\left(\left(U^{(\ell(\boldsymbol{e}'))}_{\{\boldsymbol{i}\odot\boldsymbol{e}'\}^+}\mathbb{1}_{\{\boldsymbol{e}'\preceq\boldsymbol{e}\}} + V^{(\ell(\boldsymbol{e}'))}_{\{\boldsymbol{i}\odot\boldsymbol{e}'\}^+}\mathbb{1}_{\{\boldsymbol{e}'\succ\boldsymbol{e}\}}\right)_{\boldsymbol{e}'\in\mathcal{E}_r}\right)_{\boldsymbol{i}\in\mathbb{I}_{n,k}}$$
$$\stackrel{d}{=} \left(\left(U^{(\ell(\boldsymbol{e}'))}_{\{\boldsymbol{i}\odot\boldsymbol{e}'\}^+}\mathbb{1}_{\{\boldsymbol{e}'\prec\boldsymbol{e}\}} + V^{(\ell(\boldsymbol{e}'))}_{\{\boldsymbol{i}\odot\boldsymbol{e}'\}^+}\mathbb{1}_{\{\boldsymbol{e}'\succeq\boldsymbol{e}\}}\right)_{\boldsymbol{e}'\in\mathcal{E}_r}\right)_{\boldsymbol{i}\in\mathbb{I}_{n,k}}$$
$$\stackrel{d}{=} \left(\left(U^{(\ell(\boldsymbol{e}'))}_{\{\boldsymbol{i}\odot\boldsymbol{e}'\}^+}\right)_{\boldsymbol{e}'\in\mathcal{E}_r}\right)_{\boldsymbol{i}\in\mathbb{I}_{n,k}}.$$

Then, by independence between $(\varepsilon_A)_{A\subset\mathbb{N}^+:1\leq|A|\leq k}$ and $(U^{(j)}_A, V^{(j)}_A)_{j\in\mathbb{N}, A\subset\mathbb{N}^+:1\leq|A|\leq k}$ and the triangle and Jensen inequalities

$$\frac{1}{|\mathcal{E}_r|}\sum_{\boldsymbol{e}\in\mathcal{E}_r}\mathbb{E}\left[\Phi\left(K_{2,r}\sup_{f\in\mathcal{F}}\left|\sum_{\boldsymbol{i}\in\mathbb{I}_{n,k}}\varepsilon_{\{\boldsymbol{i}\odot\boldsymbol{e}\}^+}\Delta\overline{f}_{r,\boldsymbol{i}}\left(\left(\overline{W}^{(\ell,\boldsymbol{e})}_{\{\boldsymbol{i}\odot\boldsymbol{e}'\}^+}\right)_{\boldsymbol{e}'\in\mathcal{E}_r}\right)\right|\right)\Big|\mathcal{U}^r\right]$$
$$\leq \frac{1}{|\mathcal{E}_r|}\sum_{\boldsymbol{e}\in\mathcal{E}_r}\mathbb{E}\left[\Phi\left(K_{3,r}\sup_{f\in\mathcal{F}}\left|\sum_{\boldsymbol{i}\in\mathbb{I}_{n,k}}\varepsilon_{\{\boldsymbol{i}\odot\boldsymbol{e}\}^+}\overline{f}_{r,\boldsymbol{i}}\left(\left(U^{(\ell(\boldsymbol{e}'))}_{\{\boldsymbol{i}\odot\boldsymbol{e}'\}^+}\right)_{\boldsymbol{e}'\in\mathcal{E}_r}\right)\right|\right)\Big|\mathcal{U}^r\right], \qquad (3.8)$$

where $K_{3,r} = 2k\frac{(n-k)!}{n!}K_{2,r}$.

**Fourth step: conclusion**

Combining Equations (3.1), (3.3), (3.4), (3.6), (3.8) and using the expressions of $K_{1,r}$, $K_{2,r}$ and $K_{3,r}$, we finally obtain

$$\mathbb{E}\left[\Phi\left(\sup_{f\in\mathcal{F}}\left|\frac{(n-k)!}{n!}\sum_{\boldsymbol{i}\in\mathbb{I}_{n,k}}f(Y_{\boldsymbol{i}}) - \mathbb{E}\left[f(Y_{(\boldsymbol{1})})\right]\right|\right)\right]$$
$$\leq \frac{1}{k}\sum_{r=1}^{k}\frac{1}{|\mathcal{E}_r|}\sum_{\boldsymbol{e}\in\mathcal{E}_r}\mathbb{E}\left[\Phi\left(C_{r,k}\frac{(n-k)!}{n!}\sup_{f\in\mathcal{F}}\left|\sum_{\boldsymbol{i}\in\mathbb{I}_{n,k}}\varepsilon_{\{\boldsymbol{i}\odot\boldsymbol{e}'\}^+}f(Y^r_{\boldsymbol{i}})\right|\right)\right],$$

with $C_{r,k} = 2k|\mathcal{E}_r| \left(3|\mathcal{E}_r|^{\mathcal{E}_r}\right)^{|\mathcal{E}_r|-1}$ and

$$Y_{\boldsymbol{i}}^r = \tau\left(\left(U_{\{\boldsymbol{i}\odot e\}+}^{(0)}\right)_{e\in\underline{\mathcal{E}}_r}, \left(U_{\{\boldsymbol{i}\odot e\}+}^{(\ell(e))}\right)_{e\in\mathcal{E}_r}, \left(V_{\{\boldsymbol{i}\odot e\}+}^{(0)}\right)_{e\in\overline{\mathcal{E}}_r}\right).$$

By construction of the $\left(U_A^{(j)}\right)_{A\subset\mathbb{N}^+:1\leq|A|\leq k}$ and $\left(V_A^{(0)}\right)_{A\subset\mathbb{N}^+:r+1\leq|A|\leq k}$, $(Y_{\boldsymbol{i}}^r)_{\boldsymbol{i}\in\mathbb{I}_k}$ is jointly exchangeable and dissociated, with marginal distribution $P$. This concludes the proof.

## 3.3 Lemma A.2

For any $j \in \{1, ..., |\mathcal{E}_r|\}$, let $L_{r,j} = \left(3|\mathcal{E}_r|^{|\mathcal{E}_r|}\right)^{|\mathcal{E}_r|-j}$. We will prove by reverse induction on $j$ that for every function $b$ from $\mathcal{E}_r$ to $\{1, ..., |\mathcal{E}_r|\}$ with $|\mathcal{R}(b)| = j$,

$$\mathbb{E}\Phi\left(\sup_{h\in\mathcal{H}}\left|\sum_{\boldsymbol{i}\in\mathbb{I}_{n,k}} h\left(\left(W_{\{\boldsymbol{i}\odot e\}+}^{(b(e))}\right)_{e\in\mathcal{E}_r}, \boldsymbol{i}\right)\right|\right)$$
$$\leq \mathbb{E}\Phi\left(L_{r,j}\sup_{h\in\mathcal{H}}\left|\sum_{\boldsymbol{i}\in\mathbb{I}_{n,k}} h\left(\left(W_{\{\boldsymbol{i}\odot e\}+}^{(\ell(e))}\right)_{e\in\mathcal{E}_r}, \boldsymbol{i}\right)\right|\right). \quad (3.9)$$

The result follows by considering $j = 1$. (3.9) is in fact an equality when $j = |\mathcal{E}_r|$, so the result holds for the base case. Next, when $b$ is not a bijection, both sides of (3.9) are left unchanged when $b$ is replaced by $\sigma \circ b$ for $\sigma$ a permutation of $\{1, ..., |\mathcal{E}_r|\}$. As a consequence, we can assume without loss of generality that $|b^{-1}(1)| \geq 2$ and $b^{-1}(2) = ... = b^{-1}(|b^{-1}(1)|) = \emptyset$ in the induction step. This induction step is divided into two parts. In the first part, we build an array of random variables $(\widetilde{W}_A^{(e)})_{e\in\mathcal{E}_r, A\in\mathcal{A}_r}$. This array is such that

$$\left(\left(\widetilde{W}_{\{\boldsymbol{i}\odot e\}+}^{(e)}\right)_{e\in b^{-1}(1)}, \left(\widetilde{W}_{\{\boldsymbol{i}\odot e\}+}^{(e)}\right)_{e\notin b^{-1}(1)}\right)_{\boldsymbol{i}\in\mathbb{I}_{n,k}}$$
$$\stackrel{d}{=} \left(\left(W_{\{\boldsymbol{i}\odot e\}+}^{(\ell'(e))}\right)_{e\in b^{-1}(1)}, \left(W_{\{\boldsymbol{i}\odot e\}+}^{(b(e))}\right)_{e\notin b^{-1}(1)}\right)_{\boldsymbol{i}\in\mathbb{I}_{n,k}}, \quad (3.10)$$

with $\ell'$ a bijection from $b^{-1}(1)$ to $\{1, ..., |b^{-1}(1)|\}$. Moreover, it satisfies, for all $\boldsymbol{i} \in \mathbb{I}_{n,k}$,

$$\mathbb{E}\left(h\left(\left(\widetilde{W}_{\{\boldsymbol{i}\odot e\}+}^{(e)}\right)_{e\in\mathcal{E}_r}, \boldsymbol{i}\right)\bigg|\mathcal{W}\right) = \frac{1}{|B(b)|}\sum_{b'\in B(b)} h\left(\left(W_{\{\boldsymbol{i}\odot e\}+}^{(b'(e))}\right)_{e\in\mathcal{E}_r}, \boldsymbol{i}\right), \quad (3.11)$$

where $\mathcal{W}$ denotes the $\sigma$-algebra generated by the $\left(W_A^{(j)}\right)_{A\in\mathcal{A}_r, j=1,...,|\mathcal{E}_r|}$ and

$$B(b) = \{b' : b'(e) = b(e) \text{ if } e \notin b^{-1}(1), b'(e) \in \{1, ..., |b^{-1}(1)|\} \text{ if } e \in b^{-1}(1)\}.$$

In the second part of the induction step, we combine (3.10) and (3.11) with Jensen, convexity and triangle inequalities to get upper bounds on the left-hand side of (3.9).

**First part: construction of the $\widetilde{W}_A^{(e)}$.**

Let $\ell'$ be a bijection from $b^{-1}(1)$ to $\{1, ..., |b^{-1}(1)|\}$ and let $\left(r_A^b\right)_{A\in\mathcal{A}_r}$ be some independent uniform random variables on $\{1, ..., |b^{-1}(1)|\}$. For $(j, l) \in \mathbb{N} \times \mathbb{N}^+$, $\text{rem}(j, l)$ denotes the remainder of the

13

division of $j$ by $l$. For any $(e, A) \in \mathcal{E}_r \times \mathcal{A}_r$, let $\widetilde{W}_A^{(e)} = W_A^{(1+\mathrm{rem}(\ell'(e)+r_A^b,|b^{-1}(1)|))}$ if $e \in b^{-1}(1)$ and $\widetilde{W}_A^{(e)} = W_A^{(b(e))}$ otherwise. Similarly, let $\widehat{W}_A^{(e)} = W_A^{(\ell'(e))}$ if $e \in b^{-1}(1)$ and $\widehat{W}_A^{(e)} = W_A^{(b(e))}$ otherwise. Conditional on $r_A^b$, the function $e \mapsto 1 + \mathrm{rem}\left(\ell'(e) + r_A^b, |b^{-1}(1)|\right)$ is a bijection from $b^{-1}(1)$ to $\{1, ..., |b^{-1}(1)|\}$. It follows that conditional on $r_A^b$, we have

$$\left(\widetilde{W}_A^{(e)}\right)_{e \in \mathcal{E}_r} \stackrel{d}{=} \left(\widehat{W}_A^{(e)}\right)_{e \in \mathcal{E}_r}.$$

Because the right-hand side does not depend on $r_A^b$, the previous equality also holds unconditionally. Independence of the $W_A^{(j)}$s across $A$ ensures

$$\left(\widetilde{W}_A^{(e)}\right)_{e \in \mathcal{E}_r, A \in \mathcal{A}_r} \stackrel{d}{=} \left(\widehat{W}_A^{(e)}\right)_{e \in \mathcal{E}_r, A \in \mathcal{A}_r},$$

or equivalently

$$\left(\widetilde{W}_{\{i \odot e'\}^+}^{(e)}\right)_{e \in \mathcal{E}_r, i \in \mathbb{I}_{n,k}, e' \in \mathcal{E}_r} \stackrel{d}{=} \left(\widehat{W}_{\{i \odot e'\}^+}^{(e)}\right)_{e \in \mathcal{E}_r, i \in \mathbb{I}_{n,k}, e' \in \mathcal{E}_r}.$$

Considering elements such that $e' = e$ in the previous equality yields (3.10).

Next, if $(A_e)_{e \in \mathcal{E}_r}$ is a family of distinct elements of $\mathcal{A}_r$, then uniform distribution and independence of the $r_{A_e}^b$s induces that for every $i \in \mathbb{I}_{n,k}$

$$\mathbb{E}\left(h\left(\left(\widetilde{W}_{A_e}^{(e)}\right)_{e \in \mathcal{E}_r}, i\right) \bigg| \mathcal{W}\right) = \frac{1}{|B(b)|} \sum_{b' \in B(b)} h\left(\left(W_{A_e}^{(b'(e))}\right)_{e \in \mathcal{E}_r}, i\right).$$

For every $i \in \mathbb{I}_{n,k}$, $(\{i \odot e\}^+)_{e \in \mathcal{E}_r}$ is a family of distinct subsets of $\{1, ..., n\}$ of cardinal $r$, so (3.11) follows.

**Second part: upper bound on the LHS of** (3.9)

As $\{2, ..., |b^{-1}(1)|\} \cap \mathcal{R}(b) = \emptyset$, $B(b) \setminus \{b\}$ can be partitioned into two subsets $B_1(b)$ and $B_2(b)$, with

$$B_1(b) = \{b' \in B(b) : |\mathcal{R}(b')| > j = |\mathcal{R}(b)|\},$$
$$B_2(b) = \{b' \in B(b) : b'(e) = m \in \{2, ..., |b^{-1}(1)|\} \,\forall e \in b^{-1}(1)\}.$$

Moreover, $|B_2(b)| = |b^{-1}(1)| - 1$. Let $\mathcal{W}_1$ and $\mathcal{W}_1'$ be the $\sigma$-algebra generated by $\left\{W_A^{(j)}, A \in \mathcal{A}_r, j \in \mathcal{R}(b)\right\}$ and $\left\{W_A^{(j)}, A \in \mathcal{A}_r, j \in \mathcal{R}(b) \setminus \{1\}\right\}$, respectively. The $W_A^{(j)}$s are i.i.d. across $j$. Consequently, for every $b' \in B_2(b)$,

$$\mathbb{E}\left(h\left(\left(W_{\{i \odot e\}^+}^{(b'(e))}\right)_{e \in \mathcal{E}_r}, i\right) \bigg| \mathcal{W}_1\right) = \mathbb{E}\left(h\left(\left(W_{\{i \odot e\}^+}^{(b'(e))}\right)_{e \in \mathcal{E}_r}, i\right) \bigg| \mathcal{W}_1'\right)$$
$$= \mathbb{E}\left(h\left(\left(W_{\{i \odot e\}^+}^{(b(e))}\right)_{e \in \mathcal{E}_r}, i\right) \bigg| \mathcal{W}_1'\right).$$

As a result, using the partition $B(b) = \{b\} \cup B_1(b) \cup B_2(b)$, we obtain

$$h\left(\left(W_{\{i \odot e\}^+}^{(b(e))}\right)_{e \in \mathcal{E}_r}, i\right) = \mathbb{E}\left[\sum_{b' \in B(b)} h\left(\left(W_{\{i \odot e\}^+}^{(b'(e))}\right)_{e \in \mathcal{E}_r}, i\right) \bigg| \mathcal{W}_1\right]$$
$$- \mathbb{E}\left[\sum_{b' \in B_1(b)} h\left(\left(W_{\{i \odot e\}^+}^{(b'(e))}\right)_{e \in \mathcal{E}_r}, i\right) \bigg| \mathcal{W}_1\right]$$
$$- \left(|b^{-1}(1)| - 1\right) \mathbb{E}\left[h\left(\left(W_{\{i \odot e\}^+}^{(b(e))}\right)_{e \in \mathcal{E}_r}, i\right) \bigg| \mathcal{W}_1'\right].$$

Then, by Lemma S1.

$$3\mathbb{E}\Phi\left(\sup_{h\in\mathcal{H}}\left|\sum_{i\in\mathbb{I}_{n,k}}h\left(\left(W_{\{i\odot e\}^+}^{(b(e))}\right)_{e\in\mathcal{E}_r},i\right)\right|\right)$$

$$\leq \mathbb{E}\left[\Phi\left(3\sup_{h\in\mathcal{H}}\left|\sum_{i\in\mathbb{I}_{n,k}}\mathbb{E}\left[\sum_{b'\in B(b)}h\left(\left(W_{\{i\odot e\}^+}^{(b'(e))}\right)_{e\in\mathcal{E}_r},i\right)\bigg|\mathcal{W}_1\right]\right|\right)\right]$$

$$+\mathbb{E}\left[\Phi\left(3\sup_{h\in\mathcal{H}}\left|\sum_{i\in\mathbb{I}_{n,k}}\mathbb{E}\left[\sum_{b'\in B_1(b)}h\left(\left(W_{\{i\odot e\}^+}^{(b'(e))}\right)_{e\in\mathcal{E}_r},i\right)\bigg|\mathcal{W}_1\right]\right|\right)\right]$$

$$+\mathbb{E}\left[\Phi\left(3\left(|b^{-1}(1)|-1\right)\sup_{h\in\mathcal{H}}\left|\sum_{i\in\mathbb{I}_{n,k}}\mathbb{E}\left[h\left(\left(W_{\{i\odot e\}^+}^{(b(e))}\right)_{e\in\mathcal{E}_r},i\right)\bigg|\mathcal{W}_1'\right]\right|\right)\right]. \quad (3.12)$$

Denote by $T_1, T_2$ and $T_3$ the three terms on the RHS and let $\widetilde{b}(e) = \ell'(e)$ if $e \in b^{-1}(1)$ and $\widetilde{b}(e) = b(e)$ otherwise. Then

$$T_1 \leq \mathbb{E}\left[\Phi\left(3\sup_{h\in\mathcal{H}}\left|\sum_{i\in\mathbb{I}_{n,k}}\sum_{b'\in B(b)}h\left(\left(W_{\{i\odot e\}^+}^{(b'(e))}\right)_{e\in\mathcal{E}_r},i\right)\right|\right)\right]$$

$$=\mathbb{E}\left[\Phi\left(3|B(b)|\sup_{h\in\mathcal{H}}\left|\sum_{i\in\mathbb{I}_{n,k}}\mathbb{E}\left(h\left(\left(\widetilde{W}_{\{i\odot e\}^+}^{(e)}\right)_{e\in\mathcal{E}_r},i\right)\bigg|\mathcal{W}\right)\right|\right)\right]$$

$$\leq\mathbb{E}\left[\Phi\left(3|B(b)|\sup_{h\in\mathcal{H}}\left|\sum_{i\in\mathbb{I}_{n,k}}h\left(\left(W_{\{i\odot e\}^+}^{(\widetilde{b}(e))}\right)_{e\in\mathcal{E}_r},i\right)\right|\right)\right]$$

$$\leq\mathbb{E}\left[\Phi\left(3|B(b)|L_{r,j+1}\sup_{h\in\mathcal{H}}\left|\sum_{i\in\mathbb{I}_{n,k}}h\left(\left(W_{\{i\odot e\}^+}^{(\ell(e))}\right)_{e\in\mathcal{E}_r},i\right)\right|\right)\right]. \quad (3.13)$$

The first inequality follows by Jensen's inequality. The first equality is due to (3.11). The second inequality uses Jensen's inequality and (3.10). Finally, (3.13) relies on the induction hypothesis and $|\mathcal{R}(\widetilde{b})| > j$. Similarly,

$$T_2 \leq \frac{1}{|B_1(b)|}\sum_{b'\in B_1(b)}\mathbb{E}\left[\Phi\left(3|B_1(b)|\sup_{h\in\mathcal{H}}\left|\sum_{i\in\mathbb{I}_{n,k}}h\left(\left(W_{\{i\odot e\}^+}^{(b'(e))}\right)_{e\in\mathcal{E}_r},i\right)\right|\right)\right]$$

$$\leq\mathbb{E}\left[\Phi\left(3|B_1(b)|L_{r,j+1}\sup_{h\in\mathcal{H}}\left|\sum_{i\in\mathbb{I}_{n,k}}h\left(\left(W_{\{i\odot e\}^+}^{(\ell(e))}\right)_{e\in\mathcal{E}_r},i\right)\right|\right)\right], \quad (3.14)$$

where the first inequality follows by Jensen's inequality and the second by the induction hypothesis, since $|\mathcal{R}(b')| > j$ for all $b' \in B_1(b)$. Finally, note that for each $i$, all the $\{i \odot e\}^+$s are disjoint so, conditional on $\mathcal{W}_1'$,

$$\left(W_{\{i\odot e\}^+}^{(b(e))}\right)_{e\in\mathcal{E}_r} \stackrel{d}{=} \left(W_{\{i\odot e\}^+}^{(\widetilde{b}(e))}\right)_{e\in\mathcal{E}_r}.$$

As a result,

$$
\begin{aligned}
T_3 =& \mathbb{E}\left[\Phi\left(3\left(|b^{-1}(1)|-1\right)\sup_{h\in\mathcal{H}}\left|\sum_{i\in\mathbb{I}_{n,k}}\mathbb{E}\left[h\left(\left(W^{(\widetilde{b}(e))}_{\{i\odot e\}+}\right)_{e\in\mathcal{E}_r},i\right)\Big|\mathcal{W}'_1\right]\right|\right)\right]\\
\leq& \mathbb{E}\left[\Phi\left(3\left(|b^{-1}(1)|-1\right)L_{r,j+1}\sup_{h\in\mathcal{H}}\left|\sum_{i\in\mathbb{I}_{n,k}}h\left(\left(W^{(\ell(e))}_{\{i\odot e\}+}\right)_{e\in\mathcal{E}_r},i\right)\right|\right)\right],
\end{aligned}
\tag{3.15}
$$

where the inequality follows by Jensen's inequality and the induction hypothesis again. We finally get (3.9) by combining (3.12)-(3.15) with monotonicity of $\Phi$, the expression of $L_{r,j+1}$ and

$$
\max\left(|B(b)|,|B_1(b)|,|b^{-1}(1)-1|\right)\leq|\mathcal{E}_r|^{|\mathcal{E}_r|}.
$$

This concludes the induction step, and thus the proof of the lemma.

## 3.4 Lemma A.3

Let $\lfloor n/k \rfloor$ denote the largest integer smaller or equal to $n/k$. For any array $(A_i)_{i\in\mathbb{I}_{n,k}}$ and any $i \in \mathbb{I}_{n,k}$, we have

$$
A_i = \frac{\sum_{t=1}^{\lfloor n/k \rfloor}\sum_{\pi\in\mathfrak{S}_n} A_{\pi(k(t-1)+1),\ldots,\pi(kt)}\mathbb{1}_{\{\pi(k(t-1)+1)=i_1,\ldots,\pi(kt)=i_k\}}}{\sum_{t=1}^{\lfloor n/k \rfloor}\sum_{\pi\in\mathfrak{S}_n}\mathbb{1}_{\{\pi(k(t-1)+1)=i_1,\ldots,\pi(kt)=i_k\}}}.
$$

The denominator is $\lfloor n/k \rfloor$ times the number of permuations over $\{1,...,n\}$ with $k$ values fixed. It is thus equal to $\lfloor n/k \rfloor(n-k)!$. Hence,

$$
\begin{aligned}
\sum_{i\in\mathbb{I}_{n,k}} A_i &= \frac{1}{\lfloor n/k \rfloor(n-k)!}\sum_{\pi\in\mathfrak{S}_n}\sum_{t=1}^{\lfloor n/k \rfloor} A_{\pi(k(t-1)+1),\ldots,\pi(kt)}\\
&\quad\times\sum_{i\in\mathbb{I}_{n,k}}\mathbb{1}_{\{\pi(k(t-1)+1)=i_1,\ldots,\pi(kt)=i_k\}}\\
&= \frac{1}{\lfloor n/k \rfloor(n-k)!}\sum_{\pi\in\mathfrak{S}_n}\sum_{t=1}^{\lfloor n/k \rfloor} A_{\pi(k(t-1)+1),\ldots,\pi(kt)}.
\end{aligned}
$$

Let $c = \max_{1\leq j\leq N}\|f_j\|_\infty$ and $\sigma^2 = \max_{1\leq j\leq N}\mathbb{V}(f_j(Y_1))$. For any $n \geq k$, let $r_n = (n/\lfloor n/k \rfloor)^{1/2}$. Then:

$$
\frac{\sqrt{n}(n-k)!}{n!}\frac{\sqrt{\lfloor n/k \rfloor}}{2cr_n}\sum_{i\in\mathbb{I}_{n,k}}(g(Y_i)-\mathbb{E}[g(Y_1)]) = \frac{1}{n!}\sum_{\pi\in\mathfrak{S}_n} V_{n,g,\pi},
$$

with

$$
V_{n,g,\pi} = \frac{1}{2c}\sum_{t=1}^{\lfloor n/k \rfloor} g(Y_{\pi(k(t-1)+1),\ldots,\pi(kt)}) - \mathbb{E}[g(Y_1)].
$$

For every $s > 0$, and $g \in \mathcal{F}$ let $\varphi_{n,g}(s) = \mathbb{E}\left[\exp\left(s\frac{1}{n!}\sum_{\pi\in\mathfrak{S}_n} V_{n,g,\pi}\right)\right]$. By convexity of $x\mapsto\exp(sx)$ and joint exchangeability, we get

$$
\varphi_{n,g}(s) \leq \frac{1}{n!}\sum_{\pi\in\mathfrak{S}_n}\mathbb{E}\left[e^{sV_{n,g,\pi}}\right] = \mathbb{E}\left[e^{sV_{n,g,\mathrm{id}}}\right],
$$

16

where id denotes the identity permutation. $V_{n,g,\mathrm{id}}$ is a sum of $\lfloor n/k \rfloor$ i.i.d. random variables. Then, by Theorem 3.1.5 in Giné and Nickl (2015),

$$\mathbb{E}\left[e^{sV_{n,g,\mathrm{id}}}\right] \leq \exp\left(\lfloor n/k\rfloor \frac{\sigma^2}{4c^2}\left(\exp(s)-1-s\right)\right).$$

Hence, Theorem 3.1.10 (b) in Giné and Nickl (2015) entails

$$\mathbb{E}\left[\max_{f\in\mathcal{F}}|\mathbb{G}_n(f)|\right] \leq \frac{2cr_n}{\sqrt{\lfloor n/k\rfloor}}\left(\sqrt{\frac{2\lfloor n/k\rfloor\sigma^2\log 2N}{4c^2}} + \frac{\log 2N}{3}\right)$$

$$\leq 2\sqrt{k\sigma^2\log 2N} + \frac{4ck\log 2N}{3\sqrt{n}},$$

where the last inequality follows by remarking that $r_n \leq \sqrt{2k}$.

## 3.5 Theorem 2.1

### 3.5.1 Uniform law of large numbers

**Convergence in $L^1$ under Assumption 4-(i)** Let $M$ be some arbitrary positive constant. The symmetrization Lemma A.1 applied to the class $\mathcal{G} = \left\{f\mathbb{1}_{\{F\leq M\}}, f\in\mathcal{F}\right\}$ and $\Phi =\mathrm{Id}$ ensures that

$$\mathbb{E}\left[\sup_{\mathcal{F}}|\mathbb{P}_n f - Pf|\right] \leq 2\mathbb{E}\left[F(Y_{\boldsymbol{1}})\mathbb{1}_{\{F(Y_{\boldsymbol{1}})>M\}}\right]$$

$$+ \sum_{r=1}^{k}\sum_{e\in\mathcal{E}_r} K_{r,k}\mathbb{E}\left[\sup_{\mathcal{F}}\left|\frac{(n-k)!}{n!}\sum_{\boldsymbol{i}\in\mathbb{I}_{n,k}}\varepsilon_{\{\boldsymbol{i}\odot e\}^+}f(Y_{\boldsymbol{i}}^r)\mathbb{1}_{\{F(Y_{\boldsymbol{i}}^r)\leq M\}}\right|\right],$$

with $K_{r,k}$ some non negative number depending on $r$ and $k$ only.

For every $(a_{ij})_{i=1\ldots n, j=1\ldots m} \in \mathbb{R}^{nm}$ and independent Rademacher random variables $(\varepsilon_i)_{i=1\ldots n}$, we have (see for instance Lemma 2.3.4 in Giné and Nickl, 2015)

$$\mathbb{E}\left[\max_{j\in\{1,\ldots,m\}}\left|\sum_{i=1}^{n}\varepsilon_i a_{ij}\right|\right] \leq \left[2\log(2m)\max_{j\in\{1,\ldots,m\}}\sum_{i=1}^{n}a_{ij}^2\right]^{1/2}. \quad (3.16)$$

Next, reasoning conditionally on the data, we can consider for every $\eta_1 > 0$ a minimal $\eta_1$-covering of $\mathcal{F}$ for the seminorm $||\cdot||_{e,M,1}$ with closed balls centered in $\mathcal{F}$. This implies

$$\mathbb{E}\left[\sup_{\mathcal{F}}\left|\frac{(n-k)!}{n!}\sum_{\boldsymbol{i}\in\mathbb{I}_{n,k}}\varepsilon_{\{\boldsymbol{i}\odot e\}^+}f(Y_{\boldsymbol{i}}^r)\mathbb{1}_{\{F(Y_{\boldsymbol{i}}^r)\leq M\}}\right|\Big|(Y_{\boldsymbol{i}}^r)_{\boldsymbol{i}\in\mathbb{I}_{n,k}}\right]$$

$$=\mathbb{E}\left[\sup_{\mathcal{F}}\left|\frac{(n-k)!}{n!}\sum_{\overrightarrow{\boldsymbol{i}}\in\mathbb{I}_{n,r}}\varepsilon_{\{\boldsymbol{i}\}}\sum_{\pi\in\mathfrak{S}_r}\sum_{\boldsymbol{i'}\in(\{1,\ldots,n\}\setminus\{\boldsymbol{i}\})^{k-r}}f\left(Y_{(\boldsymbol{i}_\pi)e+\boldsymbol{i'}(\boldsymbol{1}-e)}^r\right)\mathbb{1}_{\left\{F\left(Y_{(\boldsymbol{i}_\pi)e+\boldsymbol{i'}(\boldsymbol{1}-e)}^r\right)\leq M\right\}}\right|\Big|(Y_{\boldsymbol{i}}^r)_{\boldsymbol{i}\in\mathbb{I}_{n,k}}\right]$$

$$\leq M\left(\frac{2\log 2N(\eta_1,\mathcal{F},||\cdot||_{e,M,1})(n-r)!r!}{n!}\right)^{1/2} + \eta_1. \quad (3.17)$$

To obtain the inequality, we apply (3.16) with $m = N(\eta_1,\mathcal{F},||\cdot||_{e,M,1})$ and

$$a_{ij} = \frac{(n-k)!}{n!}\sum_{\pi\in\mathfrak{S}_r}\sum_{\boldsymbol{i'}\in(\{1,\ldots,n\}\setminus\{\boldsymbol{i}\})^{k-r}}f_j\left(Y_{(\boldsymbol{i}_\pi)e+\boldsymbol{i'}(\boldsymbol{1}-e)}^r\right)\mathbb{1}_{\left\{F(Y_{(\boldsymbol{i}_\pi)e+\boldsymbol{i'}(\boldsymbol{1}-e)}^r)\leq M\right\}},$$

where $f_j$ is one of the $N(\eta_1, \mathcal{F}, ||\cdot||_{e,M,1})$ centers of balls needed to cover $\mathcal{F}$. Inequality then (3.17) follows by remarking that

$$\left(\sum_{i=1}^n a_{ij}^2\right)^{1/2} \leq M\binom{n}{r}^{1/2} \frac{(n-k)!}{n!} r! \frac{(n-r)!}{(n-r-(k-r))!} = M\left(\frac{(n-r)!r!}{n!}\right)^{1/2}.$$

Observe that $||g||_{e,M,1} \leq ||g||_{\mathbb{Q}_n^r,1}$. Thus, considering $\eta_1 = \eta||F||_{\mathbb{Q}_n^r,1}$ and using Point 2 of Lemma S12, we have, for every $\eta > 0$,

$$\mathbb{E}\left[\sup_{\mathcal{F}} \left|\frac{(n-k)!}{n!} \sum_{\boldsymbol{i} \in \mathbb{I}_{n,k}} \varepsilon_{\{\boldsymbol{i} \odot \boldsymbol{e}\}^+} f(Y_{\boldsymbol{i}}^r) \mathbb{1}_{\{F(Y_{\boldsymbol{i}}^r) \leq M\}}\right| \Big| (Y_{\boldsymbol{i}}^r)_{\boldsymbol{i} \in \mathbb{I}_{n,k}}\right]$$
$$\leq M\left(\frac{2\log 2 \sup_Q N(\eta||F||_{Q,1}, \mathcal{F}, ||\cdot||_{Q,1})(n-r)!r!}{n!}\right)^{1/2} + \eta||F||_{\mathbb{Q}_n^r,1}.$$

For any $r$ and any $\boldsymbol{i} \in \mathbb{I}_k$, we have $\mathbb{E}(F(Y_{\boldsymbol{i}}^r)) = \mathbb{E}(F(Y_1))$, and next $\mathbb{E}(||F||_{\mathbb{Q}_n^r,1}) = \mathbb{E}(F(Y_1))$. Integration with respect to the distribution of $(Y_{\boldsymbol{i}}^r)_{\boldsymbol{i} \in \mathbb{I}_{n,k}}$ ensures

$$\mathbb{E}\left[\sup_{\mathcal{F}} \left|\frac{(n-k)!}{n!} \sum_{\boldsymbol{i} \in \mathbb{I}_{n,k}} \varepsilon_{\{\boldsymbol{i} \odot \boldsymbol{e}\}^+} f(Y_{\boldsymbol{i}}^r) \mathbb{1}_{\{F(Y_{\boldsymbol{i}}^r) \leq M\}}\right|\right]$$
$$\leq M\left(\frac{2\log 2 \sup_Q N(\eta||F||_{Q,1}, \mathcal{F}, ||\cdot||_{Q,1})(n-r)!r!}{n!}\right)^{1/2} + \eta\mathbb{E}(F(Y_1)).$$

It follows that there exists a constant $K_k'$ such that

$$\mathbb{E}\left[\sup_{\mathcal{F}} |\mathbb{P}_n f - Pf|\right] \leq K_k' \left(\mathbb{E}\left[F(Y_1) \mathbb{1}_{\{F(Y_1) > M\}}\right]\right.$$
$$\left. + M\left(\frac{2\log 2 \sup_Q N(\eta||F||_{Q,1}, \mathcal{F}, ||\cdot||_{Q,1})}{n}\right)^{1/2} + \eta\mathbb{E}(F(Y_1))\right).$$

Picking $M$ and $\eta$ such that $\mathbb{E}\left[F(Y_1) \mathbb{1}_{\{F(Y_1) > M\}}\right] + \eta\mathbb{E}(F(Y_1))$ is small and letting $n$ tend to infinity, we conclude that $\mathbb{E}\left[\sup_{\mathcal{F}} |\mathbb{P}_n f - Pf|\right] = o(1)$.

**Almost-sure convergence under Assumption 4-(i).** Let $\Sigma_n$ the $\sigma$-algebra generated by $\mathcal{H}_n$, the set of functions $g$ from $\mathcal{Y}^{\mathbb{I}_k}$ to $\mathbb{R}$ that are invariant by the action of any permutation $\pi$ on $\mathbb{N}^+$ such that $\pi(j) = j$ for $j \geq n$:

$$g((Y_i)_{i \in \mathbb{I}_k}) = g\left((Y_{(\pi(i_1),...,\pi(i_k))})_{i \in \mathbb{I}_k}\right).$$

Let $h\left((Y_i)_{i \in \mathbb{I}_{n,k}}\right) = \sup_{\mathcal{F}} |\mathbb{P}_n f - Pf|$ and for $l = 1, ..., n+1$, let $\mathbb{P}_{n+1}^{\setminus\{l\}} f = \frac{(n-k)!}{n!} \sum_{\boldsymbol{i} \in \mathbb{I}_{n+1,k}} f(Y_{\boldsymbol{i}}) \mathbb{1}_{\{l \notin \{\boldsymbol{i}\}\}}$. Let $\pi$ denote the transposition on $\mathbb{N}^+$ exchanging $n+1$ and $l$. Exchangeability and the definition of $\mathcal{H}_n$ ensure that

$$\left((Y_i)_{i \in \overline{(\{1,...,n+1\}\setminus\{l\})^k}}, \left(g\left((Y_i)_{i \in \mathbb{I}_k}\right)\right)_{g \in \mathcal{H}_{n+1}}\right) \stackrel{d}{=} \left((Y_{\pi(i)})_{i \in \overline{(\{1,...,n+1\}\setminus\{l\})^k}}, \left(g\left((Y_{\pi(i)})_{i \in \mathbb{I}_k}\right)\right)_{g \in \mathcal{H}_{n+1}}\right)$$
$$\stackrel{\text{a.s.}}{=} \left((Y_i)_{i \in \mathbb{I}_{n,k}}, \left(g\left((Y_i)_{i \in \mathbb{I}_k}\right)\right)_{g \in \mathcal{H}_{n+1}}\right).$$

18

For every $l < n+1$, the above implies that conditional on $\Sigma_{n+1}$, $(Y_i)_{i \in \overline{(\{1,...,n+1\}\setminus\{l\})^k}}$ has the same distribution as $(Y_i)_{i \in \mathbb{I}_{n,k}}$. As a result,

$$\mathbb{E}\left(\sup_{\mathcal{F}}\left|\mathbb{P}_{n+1}^{\setminus\{l\}}f - Pf\right|\Big|\Sigma_{n+1}\right) = \mathbb{E}\left(h((Y_i)_{i \in \overline{(\{1,...,n+1\}\setminus\{l\})^k}})|\Sigma_{n+1}\right)$$
$$= \mathbb{E}\left(h((Y_i)_{i \in \mathbb{I}_{n,k}})|\Sigma_{n+1}\right)$$
$$= \mathbb{E}\left(\sup_{\mathcal{F}}|\mathbb{P}_n f - Pf||\Sigma_{n+1}\right).$$

Because $\sum_{l=1}^{n+1} \mathbb{P}_{n+1}^{\setminus\{l\}} f = \frac{(n-k)!}{n!} \sum_{i \in \mathbb{I}_{n+1,k}} \sum_{l=1}^{n+1} f(Y_i) \mathbf{1}_{\{l \notin \{i\}\}} = \frac{(n+1-k)!}{n!} \sum_{i \in \mathbb{I}_{n+1,k}} f(Y_i)$, we have

$$\frac{1}{n+1}\sum_{l=1}^{n+1} \mathbb{P}_{n+1}^{\setminus\{l\}} f = \mathbb{P}_{n+1}f.$$

The triangle inequality ensures

$$\sup_{\mathcal{F}}|\mathbb{P}_{n+1}f - Pf| \leq \frac{1}{n+1}\sum_{l=1}^{n+1} \sup_{\mathcal{F}}\left|\mathbb{P}_{n+1}^{\setminus\{l\}}f - Pf\right|.$$

Combining the last inequality with $\mathbb{E}\left(\sup_{\mathcal{F}}|\mathbb{P}_{n+1}f - Pf||\Sigma_{n+1}\right) = \sup_{\mathcal{F}}|\mathbb{P}_{n+1}f - Pf|$ yields

$$\sup_{\mathcal{F}}|\mathbb{P}_{n+1}f - Pf| \leq \frac{1}{n+1}\sum_{l=1}^{n+1} \mathbb{E}\left(\sup_{\mathcal{F}}\left|\mathbb{P}_{n+1}^{\setminus\{l\}}f - Pf\right|\Big|\Sigma_{n+1}\right)$$
$$= \mathbb{E}\left(\sup_{\mathcal{F}}|\mathbb{P}_n f - Pf||\Sigma_{n+1}\right).$$

This means that $\sup_{\mathcal{F}}|\mathbb{P}_n f - Pf|$ is a backward submartingale with respect to the decreasing filtration $\Sigma_n$. Hence, by the convergence theorem for backwards submartingale (see, e.g., Theorem 22 of Chapter 24 in Fristedt and Gray, 2013) and its convergence to 0 in $L^1$, $\sup_{\mathcal{F}}|\mathbb{P}_n f - Pf|$ converges almost surely to 0.

**Results under Assumption 4-(ii)** Thanks to the previous almost-sure convergence result applied to $\mathcal{F}$ reduced to a single function, we know that for every $f \in \mathcal{F}$, $\mathbb{P}_n f \xrightarrow{L^1,\text{a.s.}} \mathbb{P}f$. Using this observation, we can replicate the proof of Theorem 2.4.1 in van der Vaart and Wellner (1996).

### 3.5.2 Uniform central limit theorem

We follow the usual strategy here by showing the pointwise convergence, asymptotic equicontinuity and total boundedness of $\mathcal{F}$ (see, e.g., van der Vaart and Wellner, 1996).

**First step: pointwise convergence** Let $(f_1, ..., f_m) \in \mathcal{F} \times ... \times \mathcal{F}$. The Cramer-Wold device ensures the joint asymptotic normality of $(f_1, ..., f_m)$ if the asymptotic normality holds for $f = \sum_{i=1}^m \lambda_i f_i$ for every $(\lambda_1, ..., \lambda_m) \in \mathbb{R}^m$. For $f \in L^2(P)$, $\widehat{\theta} = \frac{(n-k)!}{n!} \sum_{i \in \mathbb{I}_{n,k}} f(Y_i)$ denotes the estimator of $\theta_0 = \mathbb{E}(f(Y_1))$. Theorem A in Silverman (1976) ensures that

$$\sqrt{n}\left(\widehat{\theta} - \theta_0\right) \xrightarrow{d} \mathcal{N}(0, K(f,f)).$$

19

**Second step: asymptotic equicontinuity under Assumption 4-(i)** We have to show that, for every $\epsilon > 0$, $\lim_{\delta \to 0} \limsup_{n \to \infty} \mathbb{P}\left(\sup_{f \in \mathcal{F}_\delta} |\mathbb{G}_n f| > \epsilon\right) = 0$. We show the stronger result that

$$\lim_{\delta \to 0} \limsup_{n \to \infty} \mathbb{E}\left[\sup_{\mathcal{F}_\delta} |\mathbb{G}_n f|\right] = 0.$$

A weighted Rademacher empirical process is sub-Gaussian with respect to the Euclidean norm of the vector of weights. As a result, conditionally on the original data, we can apply Theorem 2.3.6 in Giné and Nickl (2015). This observation implies that for every $r = 1, ..., k$ and $\boldsymbol{e} \in \mathcal{E}_r$,

$$\mathbb{E}\left[\sup_{\mathcal{F}_\delta} \left|\frac{(n-k)!}{n!} \sum_{\boldsymbol{i} \in \mathbb{I}_{n,k}} \varepsilon_{\{\boldsymbol{i} \odot \boldsymbol{e}\}^+} f(Y_{\boldsymbol{i}}^r)\right| \,\Big|\, (Y_{\boldsymbol{i}}^r)_{\boldsymbol{i} \in \mathbb{I}_{n,k}}\right]$$
$$=\mathbb{E}\left[\sup_{\mathcal{F}_\delta} \left|\frac{(n-k)!}{n!} \sum_{\boldsymbol{i} \in \overrightarrow{\mathbb{I}_{n,r}}} \varepsilon_{\{\boldsymbol{i}\}} \sum_{\pi \in \mathfrak{S}_r} \sum_{\boldsymbol{i}' \in \overrightarrow{(\{1,...,n\} \setminus \{\boldsymbol{i}\})^{k-r}}} f\left(Y_{(\boldsymbol{i}_\pi)\boldsymbol{e} + \boldsymbol{i}'(1-\boldsymbol{e})}^r\right)\right| \,\Big|\, (Y_{\boldsymbol{i}}^r)_{\boldsymbol{i} \in \mathbb{I}_{n,k}}\right]$$
$$\leq \frac{4\sqrt{2(n-r)!r!}}{\sqrt{n!}} \int_0^{\sigma_e} \sqrt{\log 2N\left(\varepsilon, \mathcal{F}_\delta, ||\cdot||_{e,2}\right)} d\varepsilon,$$

with

$$||f||_{e,2}^2 = \frac{(n-r)!r!}{n!} \sum_{\boldsymbol{i} \in \overrightarrow{\mathbb{I}_{n,r}}} \left(\frac{(n-k)!}{(n-r)!r!} \sum_{\pi \in \mathfrak{S}_r} \sum_{\boldsymbol{i}' \in \overrightarrow{(\{1,...,n\} \setminus \{\boldsymbol{i}\})^{k-r}}} f\left(Y_{(\boldsymbol{i}_\pi)\boldsymbol{e} + \boldsymbol{i}'(1-\boldsymbol{e})}^r\right)\right)^2$$

and $\sigma_e^2 = \sup_{\mathcal{F}_\delta} ||f||_{e,2}^2$. A convexity argument ensures $||f||_{e,2}^2 \leq ||f||_{\mathbb{P}_n^r,2}^2$. As a result, $N\left(\varepsilon, \mathcal{F}_\delta, ||\cdot||_{e,2}\right) \leq N\left(\varepsilon, \mathcal{F}_\delta, ||\cdot||_{\mathbb{P}_n^r,2}\right)$ and $\sigma_e^2 \leq \sigma_r^2$, with $\sigma_r^2 = \sup_{\mathcal{F}_\delta} ||f||_{\mathbb{P}_n^r,2}^2$. Next for every $r = 1, ..., k$ and $\boldsymbol{e} \in \mathcal{E}_r$:

$$\sqrt{n} \mathbb{E}\left[\sup_{\mathcal{F}_\delta} \left|\frac{(n-k)!}{n!} \sum_{\boldsymbol{i} \in \mathbb{I}_{n,k}} \varepsilon_{\{\boldsymbol{i} \odot \boldsymbol{e}\}^+} f(Y_{\boldsymbol{i}}^r)\right| \,\Big|\, (Y_{\boldsymbol{i}}^r)_{\boldsymbol{i} \in \mathbb{I}_{n,k}}\right]$$
$$\leq 4\sqrt{2k!} \int_0^{\sigma_r} \sqrt{\log 2N\left(\varepsilon, \mathcal{F}_\delta, ||\cdot||_{\mathbb{P}_n^r,2}\right)} d\varepsilon. \tag{3.18}$$

Since $\sqrt{a+b} \leq \sqrt{a} + \sqrt{b}$ and $\mathbb{P}_n^r$ is a (random) probability measure on $\mathcal{Y}$ with finite support for any $r = 1, ..., k$, we obtain

$$\int_0^{\sigma_r} \sqrt{\log 2N\left(\varepsilon, \mathcal{F}_\delta, ||\cdot||_{\mathbb{P}_n^r,2}\right)} d\varepsilon$$
$$\leq \sqrt{\log 2}\sigma_r + ||F||_{\mathbb{P}_n^r,2} \int_0^{\sigma_r/||F||_{\mathbb{P}_n^r,2}} \sup_Q \sqrt{\log N(\eta ||F||_{Q,2}, \mathcal{F}_\delta, ||\cdot||_{Q,2})} d\eta. \tag{3.19}$$

Let $J_{\mathcal{F}_\delta}(u) = \int_0^u \sup_Q \sqrt{\log N(\eta ||F||_{Q,2}, \mathcal{F}_\delta, ||\cdot||_{Q,2})} d\eta$. The functions $x \mapsto \sqrt{x}$ and $(x, y) \mapsto \sqrt{y} J_{\mathcal{F}_\delta}(\sqrt{x}/\sqrt{y})$ are both concave (the latter in view of Point 2 of Lemma S11) and $\mathbb{E}\left(||F||_{\mathbb{P}_n^r,2}^2\right) = \mathbb{E}\left(||F^2||_{\mathbb{P}_n^r,1}\right) = \mathbb{E}(F^2(Y_1))$. Then, by Lemma A.1 applied to the class $\mathcal{F}_\delta$, (3.18)-(3.19) and Jensen's inequality,

$$\mathbb{E}\left[\sup_{\mathcal{F}_\delta} |\mathbb{G}_n f|\right] \lesssim \sum_{r=1}^k \mathbb{E}(\sigma_r^2)^{1/2} + \mathbb{E}(F^2(Y_1))^{1/2} J_{\mathcal{F}_\delta}\left(\frac{\mathbb{E}(\sigma_r^2)^{1/2}}{\mathbb{E}(F^2(Y_1))^{1/2}}\right).$$

Thanks to Points 3 and 4 of Lemmas S12, we further get

$$\mathbb{E}\left[\sup_{\mathcal{F}_\delta}|\mathbb{G}_n f|\right] \lesssim \sum_{r=1}^{k} \mathbb{E}(\sigma_r^2)^{1/2} + \mathbb{E}(F^2(Y_1))^{1/2} J_\mathcal{F}\left(\frac{\mathbb{E}(\sigma_r^2)^{1/2}}{4\mathbb{E}(F^2(Y_1))^{1/2}}\right).$$

As $\lim_{x \to 0} J_\mathcal{F}(x) = 0$, it is sufficient to show that

$$\lim_{\delta \to 0} \limsup_{n \to \infty} \mathbb{E}(\sigma_r^2) = 0, \text{ for every } r = 1, ..., k \tag{3.20}$$

By the triangle inequality and the definition of $\mathcal{F}_\delta$ and $\mathcal{F}_\infty$,

$$\sigma_r^2 = \sup_{\mathcal{F}_\delta}\left|\frac{(n-k)!}{n!}\sum_{\boldsymbol{i}\in\mathbb{I}_{n,k}} f^2(Y_{\boldsymbol{i}}^r)\right| \leq \sup_{\mathcal{F}_\delta}\left|\frac{(n-k)!}{n!}\sum_{\boldsymbol{i}\in\mathbb{I}_{n,k}} f^2(Y_{\boldsymbol{i}}^r) - Pf^2\right| + \delta^2$$

$$\leq \sup_{\mathcal{F}_\infty}\left|\frac{(n-k)!}{n!}\sum_{\boldsymbol{i}\in\mathbb{I}_{n,k}} f^2(Y_{\boldsymbol{i}}^r) - Pf^2\right| + \delta^2.$$

Noting that $4F^2$ is an envelope for $\mathcal{F}_\infty^2$, Point 5 of Lemma S12 yields

$$\sup_Q N\left(\eta\|4F^2\|_{Q,1}, \mathcal{F}_\infty^2, \|\cdot\|_{Q,1}\right) < \infty \text{ for every } \eta > 0.$$

Applying Theorem 2.1.1 to the class $\mathcal{F}_\infty^2$ for the array $(Y_{\boldsymbol{i}}^r)_{\boldsymbol{i}\in\mathbb{I}_k}$, we get

$$\lim_{n\to\infty} \mathbb{E}\left(\sup_{\mathcal{F}_\infty}\left|\frac{(n-k)!}{n!}\sum_{\boldsymbol{i}\in\mathbb{I}_{n,k}} f^2(Y_{\boldsymbol{i}}^r) - Pf^2\right|\right) = 0,$$

and then (3.20) holds.

**Third step: asymptotic equicontinuity under Assumption 4-(ii)** The proof in the i.i.d case is detailed in, e.g., Giné and Nickl (2015), see their Theorem 3.7.38. We simply remark that once the maximal inequality for independent data (cf. Lemma 3.5.12 in Giné and Nickl, 2015) is replaced with Lemma A.3, the proof of Theorem 3.7.38 in Giné and Nickl (2015) remains valid in our setup up to a modification of the constants, which now depend on $k$.

**Fourth step: total boundedness** We start with Assumption 4-(i). Fix $\varepsilon > 0$. The reasoning previously used to control $\sigma_r$ ensures $\lim_{n\to\infty} \mathbb{E}\left(\sup_{\mathcal{F}_\infty}|\mathbb{P}_n f^2 - Pf^2|\right) = 0$. Then we have with probability approaching one and for every $(f_1, f_2) \in \mathcal{F} \times \mathcal{F}$

$$\|f_1 - f_2\|_{P,2}^2 \leq \|f_1 - f_2\|_{\mathbb{P}_n,2}^2 + \varepsilon^2.$$

As a consequence,

$$N(\varepsilon, \mathcal{F}, \|\cdot\|_{P,2}) \leq N\left(\frac{\varepsilon}{\sqrt{2}}, \mathcal{F}, \|\cdot\|_{\mathbb{P}_n,2}\right) + o_p(1)$$

$$\leq \mathbb{1}_{\{\|F\|_{\mathbb{P}_n,2}=0\}} + \sup_Q N\left(\frac{\varepsilon\|F\|_{Q,2}}{\sqrt{2}\|F\|_{\mathbb{P}_n,2}}, \mathcal{F}, \|\cdot\|_{Q,2}\right) \mathbb{1}_{\{\|F\|_{\mathbb{P}_n,2}>0\}} + o_p(1)$$

$$= O_p(1),$$

because $||F||_{\mathbb{P}_n,2}$ converges almost surely to $\mathbb{E}(F^2(Y_1))^{1/2}$. Then $N(\varepsilon, \mathcal{F}, ||\cdot||_{P,2}) < \infty$.

Under Assumption 4-(ii), we use the following relation (see, e.g. van der Vaart and Wellner, 1996, p.84): for every $\varepsilon > 0$, $N(\varepsilon, \mathcal{F}, ||\cdot||_{P,2}) \leq N_{[\,]}(2\varepsilon, \mathcal{F}, ||\cdot||_{P,2})$. Finiteness of $\int_0^\infty \sqrt{\log N_{[\,]}(2\eta, \mathcal{F}, ||\cdot||_{P,2})}d\eta$ is then enough to conclude.

## 3.6 Proposition 2.1

**Sufficient condition**

Let $A_n = n(n-1)\mathbb{E}\left[\sup_{f \in \mathcal{F}} |\mathbb{P}_n f - Pf|\right]$. Remark that $A_n = \mathbb{E}\left[\sup_{\mathcal{F}^s} \left|\sum_{i \in \mathbb{I}_{n,2}} f(\boldsymbol{Y}_{\{i\}}) - \mathbb{E}[f(\boldsymbol{Y}_{\{1\}})]\right|\right]$. Then, by the triangle inequality,

$$A_n \leq \mathbb{E}\left[\sup_{\mathcal{F}^s} \left|\sum_{i \in \mathbb{I}_{n,2}} \mathbb{E}\left[f(\boldsymbol{Y}_{\{i\}})|U_{\{i\}}\right] - \mathbb{E}\left[f(\boldsymbol{Y}_{\{1\}})\right]\right|\right] + \mathbb{E}\left[\sup_{\mathcal{F}^s} \left|\sum_{i \in \mathbb{I}_{n,2}} \left(f(\boldsymbol{Y}_{\{i\}}) - \mathbb{E}\left[f(\boldsymbol{Y}_{\{i\}})|U_{\{i\}}\right]\right)\right|\right]. \tag{3.21}$$

Next, remark that for all $f \in \mathcal{F}^s$,

$$\sum_{i \in \mathbb{I}_{n,2}} \mathbb{E}\left[f(\boldsymbol{Y}_{\{i\}})|U_{\{i\}}\right] - \mathbb{E}\left[f(\boldsymbol{Y}_{\{1\}})\right] = 2 \sum_{i_2 < i_1} \left(\mathbb{E}\left[f(\boldsymbol{Y}_{\{i\}})|U_{\{i\}}\right] - \mathbb{E}\left[f(\boldsymbol{Y}_{\{1\}})\right]\right),$$

which is a sum of $n(n-1)/2$ mutually independent and centered terms. Hence, by the usual symmetrization lemma (see, e.g., Lemma 2.3.1 in van der Vaart and Wellner, 1996),

$$\mathbb{E}\left[\sup_{\mathcal{F}^s} \left|\sum_{i \in \mathbb{I}_{n,2}} \mathbb{E}\left[f(\boldsymbol{Y}_{\{i\}})|U_{\{i\}}\right] - \mathbb{E}\left[f(\boldsymbol{Y}_{\{1\}})\right]\right|\right]$$
$$\leq 4\mathbb{E}\left[\sup_{\mathcal{F}^s} \left|\sum_{i \in \overrightarrow{\mathbb{I}_{n,2}}} \varepsilon_{\{i\}} \mathbb{E}\left[f(\boldsymbol{Y}_{\{i\}})|U_{\{i\}}\right] - \mathbb{E}\left[f(\boldsymbol{Y}_{\{1\}})\right]\right|\right]. \tag{3.22}$$

Let us turn to the second term in the right-hand side of (3.21). By Lemma 2.1, the random map

$$f \mapsto \sum_{i \in \mathbb{I}_{n,2}} \left(f(\boldsymbol{Y}_{\{i\}}) - \mathbb{E}\left[f(\boldsymbol{Y}_{\{i\}})|U_{\{i\}}\right]\right)$$

is, conditional on $(U_{\{i\}})_{i \in \mathbb{I}_{n,2}}$, a nondegenerate and centered (generalized) U-process with symmetric kernels indexed by $\{i\}$. Applying Theorem 3.5.3 in de la Peña and Giné (1999) and their remark 3.5.4 (ii) conditional on $(U_{\{i\}})_{i \in \mathbb{I}_{n,2}}$ and then re-integrating, we obtain

$$\mathbb{E}\left[\sup_{\mathcal{F}^s} \left|\sum_{i \in \mathbb{I}_{n,2}} \left(f(\boldsymbol{Y}_{\{i\}}) - \mathbb{E}\left[f(\boldsymbol{Y}_{\{i\}})|U_{\{i\}}\right]\right)\right|\right]$$
$$\lesssim \mathbb{E}\left[\sup_{\mathcal{F}^s} \left|\sum_{i \in \mathbb{I}_{n,2}} \varepsilon_{i_1} \left(f(\boldsymbol{Y}_{\{i\}}) - \mathbb{E}\left[f(\boldsymbol{Y}_{\{i\}})|U_{\{i\}}\right]\right)\right|\right]$$
$$\lesssim \mathbb{E}\left[\sup_{\mathcal{F}^s} \left|\sum_{i \in \mathbb{I}_{n,2}} \varepsilon_{i_1} f(\boldsymbol{Y}_{\{i\}})\right|\right] + \mathbb{E}\left[\sup_{\mathcal{F}^s} \left|\sum_{i \in \mathbb{I}_{n,2}} \mathbb{E}\left[\varepsilon_{i_1} f(\boldsymbol{Y}_{\{i\}}) | (\varepsilon_{i_1})_{i_1}, (U_{\{i\}})_{i \in \mathbb{I}_{n,2}}\right]\right|\right]$$
$$\lesssim \mathbb{E}\left[\sup_{\mathcal{F}^s} \left|\sum_{i \in \mathbb{I}_{n,2}} \varepsilon_{i_1} f(\boldsymbol{Y}_{\{i\}})\right|\right]. \tag{3.23}$$

(3.21), (3.22) and (3.23) together yield

$$A_n \lesssim \mathbb{E}\left[\sup_{\mathcal{F}^s}\left|\sum_{\overrightarrow{i}\in\mathbb{I}_{n,2}} \varepsilon_{\{i\}}\mathbb{E}\left[f(\boldsymbol{Y}_{\{i\}})|U_{\{i\}}\right]\right|\right] + \mathbb{E}\left[\sup_{\mathcal{F}^s}\left|\sum_{i\in\mathbb{I}_{n,2}} \varepsilon_{i_1} f(\boldsymbol{Y}_{\{i\}})\right|\right]. \tag{3.24}$$

Applying this inequality to $\{f\mathbb{1}_{F\leq M} : f \in \mathcal{F}^s\}$ instead of $\mathcal{F}^s$ and following the same reasoning as in the proof of the uniform LLN in Theorem 2.1 until (3.17), we get for every $\varepsilon > 0$ and $M > 0$

$$\frac{A_n}{n(n-1)} \lesssim \mathbb{E}[F(\boldsymbol{Y_1})\mathbb{1}_{F(\boldsymbol{Y_1})>M}]$$
$$+ M\mathbb{E}\left[\left(\frac{\log N(\varepsilon, \mathcal{F}, ||\cdot||_{M,1,2})}{n^2}\right)^{1/2} + \left(\frac{\log N(\varepsilon, \mathcal{F}, ||\cdot||_{M,1,1})}{n}\right)^{1/2}\right] + \varepsilon, \tag{3.25}$$

with $||f||_{M,1,1} := \frac{1}{n}\sum_{i_1=1}^n \left|\frac{1}{n-1}\sum_{i_2\neq i_1} f(\boldsymbol{Y_i})\mathbb{1}_{\{F(\boldsymbol{Y_i})\leq M\}}\right|$ and
$||f||_{M,1,2} := \frac{1}{n(n-1)}\sum_{\overrightarrow{i}\in\mathbb{I}_{n,2}} \left|\mathbb{E}\left[f(\boldsymbol{Y_i})\mathbb{1}_{\{F(\boldsymbol{Y_i})\leq M\}}|U_{\{i\}}\right]\right|$. Now, mimicking the proof of (5.2.12) on page 230 in de la Peña and Giné (1999), we can see that the condition

$$\max\left(\frac{\log N(\varepsilon, \mathcal{F}, ||\cdot||_{1,2})}{n^2}, \frac{\log N(\varepsilon, \mathcal{F}, ||\cdot||_{1,1})}{n}\right) \xrightarrow{\mathbb{P}} 0 \quad \forall \varepsilon > 0$$

together with $\mathbb{E}[F(\boldsymbol{Y_1})] < \infty$ implies that for every $\varepsilon > 0$ there exists $M_1(\varepsilon) > 0$ such that

$$\limsup_{n\to\infty} \mathbb{E}\left[\left(\frac{\log N(\varepsilon, \mathcal{F}, ||\cdot||_{M_1(\varepsilon),1,2})}{n^2}\right)^{1/2} + \left(\frac{\log N(\varepsilon, \mathcal{F}, ||\cdot||_{M_1(\varepsilon),1,1})}{n}\right)^{1/2}\right] = 0$$

Inequality (3.25) then implies that for every $\varepsilon > 0$, there exists $M_2(\varepsilon) \geq M_1(\varepsilon)$ such that $\limsup_{n\to\infty} A_n/n(n-1) \lesssim \varepsilon$. Hence, $A_n/n(n-1) \to 0$. Finally, to move from convergence in $L^1$ to a.s. convergence, we apply the same backward martingale argument as in Theorem 2.1.

**Necessary condition**

First, $\sup_\mathcal{F} |\mathbb{P}_n f - Pf| \leq \mathbb{P}_n F + PF$. Moreover, $\mathbb{P}_n F$ is uniformly integrable (since $PF < \infty$). Thus, $\sup_\mathcal{F} |\mathbb{P}_n f - Pf|$ is also uniformly integrable, and its almost-sure convergence implies $A_n \xrightarrow{L^1} 0$. We now show the converse of inequalities (3.21) and (3.24). We start by proving that

$$A_n \gtrsim \mathbb{E}\left[\sup_{\mathcal{F}^s}\left|\sum_{\overrightarrow{i}\in\mathbb{I}_{n,2}} \mathbb{E}\left[f(\boldsymbol{Y}_{\{i\}})|U_{\{i\}}\right] - \mathbb{E}\left[f(\boldsymbol{Y}_{\{1\}})\right]\right| + \sup_{\mathcal{F}^s}\left|\sum_{i\in\mathbb{I}_{n,2}}\left(f(\boldsymbol{Y}_{\{i\}}) - \mathbb{E}\left[f(\boldsymbol{Y}_{\{i\}})|U_{\{i\}}\right]\right)\right|\right]. \tag{3.26}$$

To establish (3.26), observe first that almost surely,

$$\sup_{\mathcal{F}^s}\left|\sum_{\overrightarrow{i}\in\mathbb{I}_{n,2}} \mathbb{E}\left[f(\boldsymbol{Y}_{\{i\}})|U_{\{i\}}\right] - \mathbb{E}\left[f(\boldsymbol{Y}_{\{1\}})\right]\right| = \sup_{\mathcal{F}^s}\left|\mathbb{E}\left[\sum_{\overrightarrow{i}\in\mathbb{I}_{n,2}}\left(f(\boldsymbol{Y}_{\{i\}}) - f(\boldsymbol{Y}^*_{\{i\}}))\right)\big|(U_{\{j\}})_{j\in\mathbb{I}_{n,2}}\right]\right|,$$

where $(Y^*_i)_{i\in\mathbb{I}_{n,2}}$ is an independent copy of $(Y_i)_{i\in\mathbb{I}_{n,2}}$. Then, by Jensen's inequality,

$$\mathbb{E}\left[\sup_{\mathcal{F}^s}\left|\sum_{\overrightarrow{i}\in\mathbb{I}_{n,2}} \mathbb{E}\left[f(\boldsymbol{Y}_{\{i\}})|U_{\{i\}}\right] - \mathbb{E}\left[f(\boldsymbol{Y}_{\{1\}})\right]\right|\right] \leq \mathbb{E}\left[\sup_{\mathcal{F}^s}\left|\sum_{\overrightarrow{i}\in\mathbb{I}_{n,2}}\left(f(\boldsymbol{Y}_{\{i\}}) - f(\boldsymbol{Y}^*_{\{i\}}))\right)\right|\right]$$

$$\leq 2\mathbb{E}\left[\sup_{\mathcal{F}^s}\left|\sum_{\overrightarrow{i}\in\mathbb{I}_{n,2}} f(\boldsymbol{Y}_{\{i\}}) - \mathbb{E}\left[f(\boldsymbol{Y}_{\{1\}})\right]\right|\right], \tag{3.27}$$

23

where (3.27) follows by the triangle inequality. We also have, almost surely,

$$\sup_{\mathcal{F}^s} \left| \sum_{i \in \mathbb{I}_{n,2}} \left( f(Y_{\{i\}}) - \mathbb{E}\left[ f(Y_{\{i\}}) | U_{\{i\}} \right] \right) \right| = \sup_{\mathcal{F}^s} \left| \mathbb{E}\left[ \sum_{i \in \mathbb{I}_{n,2}} \left( f(Y_{\{i\}}) - f(Y_{\{i\}}^{**}) \right) | (U_\ell)_{\ell=1}^n, (U_{\{j\}})_{j \in \mathbb{I}_{n,2}} \right] \right|,$$

where $(Y_i^{**})_{i \in \mathbb{I}_{n,2}} = (\tau(V_i, V_j, U_{\{i\}}))_{i \in \mathbb{I}_{n,2}}$, and $(V_\ell)_{\ell=1}^n$ is an independent copy of $(U_\ell)_{\ell=1}^n$. Then, by similar arguments as above,

$$\mathbb{E}\left[ \sup_{\mathcal{F}^s} \left| \sum_{i \in \mathbb{I}_{n,2}} \left( f(Y_{\{i\}}) - \mathbb{E}\left[ f(Y_{\{i\}}) | U_{\{i\}} \right] \right) \right| \right] \leq 2 \mathbb{E}\left[ \sup_{\mathcal{F}^s} \left| \sum_{i \in \mathbb{I}_{n,2}} f(Y_{\{i\}}) - \mathbb{E}\left[ f(Y_{\{1\}}) \right] \right| \right]. \quad (3.28)$$

(3.26) follows by (3.27) and (3.28). Then, by the desymmetrization lemmas for i.i.d. data and U-statistics (see respectively Lemma 2.3.6 in van der Vaart and Wellner, 1996, and Theorem 3.5.3 in de la Peña and Giné, 1999), we obtain

$$\mathbb{E}\left[ \sup_{\mathcal{F}^s} \left| \sum_{i \in \mathbb{I}_{n,2}} f(Y_{\{i\}}) - \mathbb{E}\left[ f(Y_{\{1\}}) \right] \right| \right] \gtrsim \mathbb{E}\left[ \sup_{\mathcal{F}^s} \left| \sum_{i \in \overrightarrow{\mathbb{I}_{n,2}}} \varepsilon_{\{i\}} \left( \mathbb{E}\left[ f(Y_{\{i\}}) | U_{\{i\}} \right] - \mathbb{E}\left[ f(Y_{\{1\}}) \right] \right) \right| \right]$$

$$+ \mathbb{E}\left[ \sup_{\mathcal{F}^s} \left| \sum_{i \in \mathbb{I}_{n,2}} \varepsilon_{i_1} \left( f(Y_{\{i\}}) - \mathbb{E}\left[ f(Y_{\{i\}}) | U_{\{i\}} \right] \right) \right| \right]. \quad (3.29)$$

Up to the conditional expectation, this may be seen as a converse of the symmetrization lemma for $k=2$. Now, let us define

$$B_n = \mathbb{E}\left[ \sup_{\mathcal{F}^s} \left| \frac{1}{n(n-1)} \sum_{i \in \overrightarrow{\mathbb{I}_{n,2}}} \varepsilon_{\{i\}} \mathbb{E}[f(Y_{\{i\}}) | U_{\{i\}}] \right| \Big| (U_{\{i\}})_{i \geq 1} \right],$$

$$C_n = \mathbb{E}\left[ \sup_{\mathcal{F}^s} \left| \frac{1}{n(n-1)} \sum_{i \in \mathbb{I}_{n,2}} \varepsilon_{i_1} f(Y_{\{i\}}) \right| \Big| (U_j)_{j \geq 1} \right],$$

$$\|f\|_{2,1} = \left[ \frac{1}{n} \sum_{i_1=1}^n \left( \frac{1}{n-1} \sum_{i_2 \neq i_1} f(Y_{i_1,i_2}) + f(Y_{i_2,i_1}) \right)^2 \right]^{1/2},$$

$$\|f\|_{2,2} = \left( \frac{1}{n(n-1)} \sum_{1 \leq i_1 < i_2 \leq n} \mathbb{E}\left[ f(Y_{i_1,i_2}) + f(Y_{i_2,i_1}) | U_{\{i_1,i_2\}} \right]^2 \right)^{1/2}.$$

By following the proof of (ii) ⇒ (iii) in Theorem 5.2.2 in de la Peña and Giné (1999), there exists a universal constant $K$ such that for every $\varepsilon > 0$,

$$\varepsilon^2 \log N(\varepsilon, \mathcal{F}, \|\cdot\|_{2,1}) \leq K \frac{n(n-1)}{2} B_n^2 \log\left(2 + B_n^{-1}\right),$$

$$\varepsilon^2 \log N(\varepsilon, \mathcal{F}, \|\cdot\|_{2,2}) \leq K n\, C_n^2 \log\left(2 + C_n^{-1}\right).$$

For every $\varepsilon > 0$, $N(\varepsilon, \mathcal{F}, \|\cdot\|_{1,1}) \leq N(\varepsilon, \mathcal{F}, \|\cdot\|_{2,1})$ and $N(\varepsilon, \mathcal{F}, \|\cdot\|_{1,2}) \leq N(\varepsilon, \mathcal{F}, \|\cdot\|_{2,2})$. Hence,

$$\frac{\log N(\varepsilon, \mathcal{F}, \|\cdot\|_{1,2})}{n(n-1)} \vee \frac{\log N(\varepsilon, \mathcal{F}, \|\cdot\|_{1,1})}{n} \leq \frac{K}{\varepsilon^2} \left( \frac{B_n^2 \log(2 + B_n^{-1})}{2} \vee C_n^2 \log(2 + C_n^{-1}) \right).$$

We now prove that $B_n$ and $C_n$ converge to 0 in $L^1$, and thus in probability. The result will follow by the continuous mapping theorem, since $\lim_{x\to 0} x^2 \log(2+1/x) = 0$. First, by the triangle inequality,

$$\mathbb{E}[B_n] \leq \frac{1}{n(n-1)} \mathbb{E}\left[\sup_{\mathcal{F}^s} \left| \sum_{i \in \overrightarrow{\mathbb{I}_{n,2}}} \varepsilon_{\{i\}} \left(\mathbb{E}\left[f(\boldsymbol{Y}_{\{i\}})|U_{\{i\}}\right] - \mathbb{E}\left[f(\boldsymbol{Y}_{\{1\}})\right]\right) \right| \right]$$

$$+ \frac{1}{n(n-1)} \mathbb{E}\left[\sup_{\mathcal{F}^s} \left| \sum_{i \in \overrightarrow{\mathbb{I}_{n,2}}} \varepsilon_{\{i\}} \mathbb{E}\left[f(\boldsymbol{Y}_{\{1\}})\right] \right| \right].$$

The first term tends to zero by (3.29). By the Cauchy-Schwarz inequality, the second term is bounded by

$$\mathbb{E}\left[F(Y_1)\right] \mathbb{E}\left[\frac{1}{n(n-1)} \left| \sum_{i \in \overrightarrow{\mathbb{I}_{n,2}}} \varepsilon_{\{i\}} \right| \right] \lesssim \frac{1}{\sqrt{n(n-1)}}.$$

Hence, $E(B_n) \to 0$. Turning to $C_n$, we have

$$\mathbb{E}[C_n] \leq \frac{1}{n(n-1)} \left\{ \mathbb{E}\left[\sup_{\mathcal{F}^s} \left| \sum_{i \in \overrightarrow{\mathbb{I}_{n,2}}} \varepsilon_{i_1} \left(f(\boldsymbol{Y}_{\{i\}}) - \mathbb{E}\left[f(\boldsymbol{Y}_{\{i\}})|U_{\{i\}}\right]\right) \right| \right] \right.$$

$$\left. + \mathbb{E}\left[\sup_{\mathcal{F}^s} \left| \sum_{i \in \overrightarrow{\mathbb{I}_{n,2}}} \varepsilon_{i_1} \left(\mathbb{E}\left[f(\boldsymbol{Y}_{\{i\}})|U_{\{i\}}\right] - \mathbb{E}\left[f(\boldsymbol{Y}_{\{1\}})\right]\right) \right| \right] + \mathbb{E}\left[\sup_{\mathcal{F}^s} \left| \sum_{i \in \overrightarrow{\mathbb{I}_{n,2}}} \varepsilon_{i_1} \mathbb{E}\left[f(\boldsymbol{Y}_{\{1\}})\right] \right| \right] \right\}.$$

The first term tends to zero by (3.29). By the same reasoning as above, the third term also tends to 0. Let $D_n$ denote the second term. By the triangle inequality,

$$D_n \leq \frac{1}{n(n-1)} \mathbb{E}\left[\sup_{\mathcal{F}^s} \left| \sum_{i_1=1}^{n} \varepsilon_{i_1} \sum_{i_2 < i_1} \mathbb{E}\left[f(\boldsymbol{Y}_{\{i\}})|U_{\{i\}}\right] - \mathbb{E}\left[f(\boldsymbol{Y}_{\{1\}})\right] \right| \right]$$

$$+ \frac{1}{n(n-1)} \mathbb{E}\left[\sup_{\mathcal{F}^s} \left| \sum_{i_1=1}^{n} \varepsilon_{i_1} \sum_{i_2 > i_1} \mathbb{E}\left[f(\boldsymbol{Y}_{\{i\}})|U_{\{i\}}\right] - \mathbb{E}\left[f(\boldsymbol{Y}_{\{1\}})\right] \right| \right].$$

Let $D_{1n}$ and $D_{2n}$ denote the two expectations on the right-hand side. The two terms are similar, so we only consider $D_{1n}$. The variable in the supremum may be written as $|\sum_{i_1} \varepsilon_{i_1} W_{i_1}|$, where the $(W_{i_1})_{i_1=1...n}$ are mutually independent and centered. Thus, by the desymmetrization lemma for independent variables (see, e.g., Lemma in van der Vaart and Wellner, 1996), Jensen's inequality and $f(\boldsymbol{Y}_{i_1,i_2}) = f(\boldsymbol{Y}_{i_2,i_1})$ for all $f \in \mathcal{F}^s$,

$$D_{1n} \leq \frac{2}{n(n-1)} \mathbb{E}\left[\sup_{\mathcal{F}^s} \left| \sum_{i_1=1}^{n} \sum_{i_2 < i_1} \mathbb{E}\left[f(\boldsymbol{Y}_{\{i\}})|U_{\{i\}}\right] - \mathbb{E}\left[f(\boldsymbol{Y}_{\{1\}})\right] \right| \right]$$

$$\leq \frac{2}{n(n-1)} \mathbb{E}\left[\sup_{\mathcal{F}^s} \left| \sum_{i_1=1}^{n} \sum_{i_2 < i_1} f(\boldsymbol{Y}_{\{i\}}) - \mathbb{E}\left[f(\boldsymbol{Y}_{\{1\}})\right] \right| \right]$$

$$\leq \frac{A_n}{n(n-1)},$$

which tends to 0. The result follows.

## 3.7 Theorem 2.2

We only have to prove the pointwise convergence and the asymptotic equicontinuity, since the total boundedness of $\mathcal{F}$ is proved in Theorem 2.1. In the bootstrap scheme, we sample $n$ units independently in $\{1, ..., n\}$ with replacement and equal probability. Then, for any $i = 1, ..., n$, $i^*$ denotes the $i$-th sampled unit and for $\boldsymbol{i} \in \mathbb{I}_{n,k}$, we let $\boldsymbol{i}^* = (i_1^*, ..., i_k^*)$. Then $\mathbb{P}_n^* f = \frac{(n-k)!}{n!} \sum_{\boldsymbol{i} \in \mathbb{I}_{n,k}} f(Y_{\boldsymbol{i}^*}) \mathbb{1}_{\{\boldsymbol{i}^* \in \mathbb{I}_{n,k}\}}$.

**First step: pointwise convergence**

As in Theorem 2.1, it suffices to prove the result for a single $f$ such that $\mathbb{E}[f(Y_1)^2] < \infty$.

**Substep 1: asymptotic equivalence** Let $\theta = \mathbb{E}(f(Y_1))$, $\theta^* = \frac{(n-k)!}{n!} \sum_{\boldsymbol{i} \in \mathbb{I}_{n,k}} f(Y_{\boldsymbol{i}^*}) \mathbb{1}_{\{\boldsymbol{i}^* \in \mathbb{I}_{n,k}\}}$ and $\widehat{\theta} = \frac{(n-k)!}{n!} \sum_{\boldsymbol{i} \in \mathbb{I}_{n,k}} f(Y_{\boldsymbol{i}})$ its bootstrap counterpart. For $\boldsymbol{i} \in \{1, ..., n\}^k$, let $h(\boldsymbol{i}) = \mathbb{1}_{\{\boldsymbol{i} \in \mathbb{I}_{n,k}\}} \sum_{\pi \in \mathfrak{S}_k} f(Y_{i_\pi})$. We have $\theta^* = \frac{(n-k)!}{n!k!} \sum_{\boldsymbol{i} \in \mathbb{I}_{n,k}} h(\boldsymbol{i}^*)$, $\widehat{\theta} = \frac{(n-k)!}{n!k!} \sum_{\boldsymbol{i} \in \mathbb{I}_{n,k}} h(\boldsymbol{i}) = \frac{(n-k)!}{n!k!} \sum_{\boldsymbol{i} \in \{1,...,n\}^k} h(\boldsymbol{i})$ and $\mathbb{E}\left(\theta^* | (Y_{\boldsymbol{i}})_{\boldsymbol{i} \in \mathbb{I}_k}\right) = \frac{n!}{n^k(n-k)!} \widehat{\theta}$. Let

$$\theta_1^* = \frac{(n-k)!}{n!k!} \sum_{\boldsymbol{j} \in \{1,...,n\}^k} h(j_1^*, j_2, ..., j_k).$$

We have $\mathbb{E}(\theta_1^* | (Y_{\boldsymbol{i}})_{\boldsymbol{i} \in \mathbb{I}_k}) = \widehat{\theta}$. For $(\boldsymbol{i}, \boldsymbol{j}) \in \mathbb{I}_{n,k} \times \{1, ..., n\}^k$, observe that

$$\mathbb{E}\left(h(\boldsymbol{i}^*) h(j_1^*, j_2, ..., j_k) | (Y_{\boldsymbol{i}'})_{\boldsymbol{i}' \in \mathbb{I}_{n,k}}\right)$$
$$= \begin{cases} \frac{n!k!}{n^k(n-k)!} \widehat{\theta} \times \frac{1}{n} \sum_{j=1}^n h(j, j_2, ..., j_k) & \text{if } j_1 \notin \{i_1, ..., i_k\} \\ \frac{1}{n^k} \sum_{i_1=1}^n \left(\sum_{(i_2,...,i_k) \in \{1,...,n\}^{k-1}} h(\boldsymbol{i}) \times h(i_1, j_2, ..., j_k)\right) & \text{otherwise.} \end{cases}$$

Consequently,

$$n\mathbb{E}\left(\theta^* \theta_1^* | (Y_{\boldsymbol{i}'})_{\boldsymbol{i}' \in \mathbb{I}_{n,k}}\right) = n\frac{(n-k)!^2}{n!^2 k!^2}(n-k)\frac{n!}{(n-k)!}\left(\frac{n!k!}{n^k(n-k)!}\widehat{\theta}\frac{1}{n}\sum_{\boldsymbol{j} \in \mathbb{I}_{n,k}} h(\boldsymbol{j})\right)$$
$$+ n\frac{(n-k)!^2}{n!^2 k!^2} k \frac{n!}{(n-k)!}\frac{1}{n^k}\sum_{i_1=1}^n \left(\sum_{(i_2,...,i_k) \in \{1,...,n\}^{k-1}} h(\boldsymbol{i})\right)^2.$$

Hence,

$$n\mathbb{E}\left(\theta^* \theta_1^* | (Y_{\boldsymbol{i}'})_{\boldsymbol{i}' \in \mathbb{I}_{n,k}}\right) = \frac{n-k}{n^k}\frac{n!}{(n-k)!}\widehat{\theta}^2 + \frac{k}{k!^2}\frac{n^k(n-k)!}{n!}\frac{1}{n^{2k-1}}\sum_{i_1=1}^n \left(\sum_{(i_2,...,i_k) \in \{1,...,n\}^{k-1}} h(\boldsymbol{i})\right)^2.$$

Focusing on the last sum, Lemma S7 allows us to conclude that

$$\sum_{i_1=1}^n \left(\sum_{(i_2,...,i_k) \in \{1,...,n\}^{k-1}} h(\boldsymbol{i})\right)^2$$
$$= \sum_{\boldsymbol{j} \in \{1,...,n\}^{2k-1}} h(j_1, ..., j_k) h(j_1, j_{k+1}, ..., j_{2k-1})$$
$$= \sum_{c=0}^{k-1} \binom{k-1}{c}^2 \left(n^{2k-1-c} \mathbb{E}\left[h(1, ..., k) h(1, ..., 1+c, k+1, ..., 2k-c-1)\right] + o_{\text{a.s.}}(n^{2k-1-c})\right).$$

26

As $\frac{n^k(n-k)!}{n!}$ converges to 1, the quantity $\frac{n^k(n-k)!}{n!} \frac{1}{n^{2k-1}} \sum_{i_1=1}^{n} \left( \sum_{(i_2,...,i_k) \in \{1,...,n\}^{k-1}} h(\boldsymbol{i}) \right)^2$ converges almost surely to $\mathbb{E}\left(h(\boldsymbol{1})h(\boldsymbol{1}')\right)$.

Combining the exchangeability assumption, symmetry of $h$ and a combinatorial argument (see the proof of Theorem 12.3 in van der Vaart, 2000), we obtain

$$n\mathbb{E}\left(\theta^{*2} | (Y_{\boldsymbol{i}'})_{\boldsymbol{i}' \in \mathbb{I}_{n,k}}\right)$$
$$= n \frac{(n-k)!^2}{n!^2 k!^2} \sum_{\boldsymbol{i} \in \mathbb{I}_{n,k}} \sum_{\boldsymbol{j} \in \mathbb{I}_{n,k}} \mathbb{E}\left(h(\boldsymbol{i}^*)h(\boldsymbol{j}^*) | (Y_{\boldsymbol{i}'})_{\boldsymbol{i}' \in \mathbb{I}_{n,k}}\right)$$
$$= n \frac{(n-k)!}{n! k!^2} \mathbb{E}\left(h(1^*, ..., k^*)^2 | (Y_{\boldsymbol{i}'})_{\boldsymbol{i}' \in \mathbb{I}_{n,k}}\right)$$
$$= n \frac{(n-k)!^2}{n!^2} \sum_{l=0}^{k} \binom{n}{k}\binom{k}{l}\binom{n-k}{k-l} \mathbb{E}\left[h(1^*, ..., k^*)h(1^*, ..., l^*, (k+1)^*, ..., (2k-l)^*) | (Y_{\boldsymbol{i}'})_{\boldsymbol{i}' \in \mathbb{I}_{n,k}}\right].$$

When $l = 0$

$$n \frac{(n-k)!^2}{n!^2} \binom{n}{k}\binom{k}{l}\binom{n-k}{k-l} \mathbb{E}\left[h(1^*, ..., k^*)h(1^*, ..., l^*, (k+1)^*, ..., (2k-l)^*) | (Y_{\boldsymbol{i}'})_{\boldsymbol{i}' \in \mathbb{I}_{n,k}}\right]$$
$$= \frac{n(n-k)!^2 n!}{n!^2 k!^2 (n-2k)!} \left(\frac{1}{n^k} \sum_{\boldsymbol{i} \in \{1,...,n\}^k} h(\boldsymbol{i})\right)^2 = n \frac{n!}{n^{2k}(n-2k)!} \widehat{\theta}^2.$$

For every $l = 1, ..., k$,

$$\mathbb{E}\left[h(1^*, ..., (k-1)^*, k^*)h(1^*, ..., l^*, (k+1)^*, ..., (2k-l)^*) | (Y_{\boldsymbol{i}'})_{\boldsymbol{i}' \in \mathbb{I}_{n,k}}\right]$$
$$= \frac{1}{n^l} \sum_{\boldsymbol{i} \in \{1,...,n\}^l} \left(\frac{1}{n^{k-l}} \sum_{\boldsymbol{j} \in \{1,...,n\}^{k-l}} h(i_1, ..., i_l, j_1, ..., j_{k-l})\right)^2$$
$$= \frac{1}{n^{2k-l}} \sum_{\boldsymbol{j} \in \{1,...,n\}^{2k-l}} h(j_1, ..., j_k) h(j_1, ..., j_l, j_{k+1}, ..., j_{2k-l})$$
$$= \frac{1}{n^{2k-l}} \sum_{c=0}^{k-l} \binom{k-l}{c}^2 \left(n^{2k-l-c} \mathbb{E}\left[h(1, ..., k)h(1, ..., l+c, k+1, ..., 2k-c-l)\right] + o_{\text{a.s.}}(n^{2k-l-c})\right)$$
$$= \mathbb{E}\left[h(1, ..., k)h(1, ..., l, k+1, ..., 2k-l)\right] + o_{\text{a.s.}}(1),$$

using Lemma S7 once more. As $n \frac{(n-k)!^2}{n!^2} \binom{n}{k}\binom{k}{l}\binom{n-k}{k-l} = O(n^{1-k+k-l}) = o(1)$ for every $l \geq 2$ and $n \frac{(n-k)!^2}{n!^2} \binom{n}{k}\binom{k}{1}\binom{n-k}{k-1} = \frac{k^2}{k!^2} + o(1)$, we get

$$n\mathbb{E}\left(\theta^{*2} | (Y_{\boldsymbol{i}'})_{\boldsymbol{i}' \in \overline{\{1,...,n\}}}^{k}\right) = n \frac{n!}{n^{2k}(n-2k)!} \widehat{\theta}^2 + \frac{k^2}{k!^2} \mathbb{E}\left[h(1, ..., k)h(1, ..., l, k+1, ..., 2k-1)\right] + o_{\text{a.s.}}(1).$$

27

We also have
$$n\mathbb{E}\left(\theta_1^{*2}|(Y_{\boldsymbol{i'}})_{\boldsymbol{i'}\in\mathbb{I}_{n,k}}\right) = n\frac{(n-k)!^2}{n!^2 k!^2}\sum_{\boldsymbol{i}\in\{1,...,n\}^k}\sum_{\boldsymbol{j}\in\{1,...,n\}^k}\mathbb{E}\left(h(i_1^*,i_2,...,i_k)h(j_1^*,j_2,...,j_k)|(Y_{\boldsymbol{i'}})_{\boldsymbol{i'}\in\mathbb{I}_{n,k}}\right)$$
$$= n\frac{(n-k)!^2}{n!^2 k!^2}n\frac{1}{n}\sum_{i_1=1}^{n}\left(\sum_{(i_2,...,i_k)\in\{1,...,n\}^{k-1}}h(\boldsymbol{i})\right)^2$$
$$+ n\frac{(n-k)!^2}{n!^2 k!^2}\frac{n(n-1)}{n^2}\left(\sum_{\boldsymbol{i}\in\{1,...,n\}^k}h(\boldsymbol{i})\right)^2$$
$$= n\frac{(n-k)!^2}{n!^2 k!^2}\sum_{i_1=1}^{n}\left(\sum_{(i_2,...,i_k)\in\{1,...,n\}^{k-1}}h(\boldsymbol{i})\right)^2 + (n-1)\widehat{\theta}^2.$$

It follows that
$$\mathbb{E}\left(n\left((\theta^* - \widehat{\theta}) - k(\theta_1^* - \widehat{\theta})\right)^2 |(Y_{\boldsymbol{i'}})_{\boldsymbol{i'}\in\mathbb{I}_{n,k}}\right)$$
$$= n\mathbb{E}\left(\theta^{*2}|(Y_{\boldsymbol{i'}})_{\boldsymbol{i'}\in\mathbb{I}_{n,k}}\right) + nk^2\mathbb{E}\left(\theta_1^{*2}|(Y_{\boldsymbol{i'}})_{\boldsymbol{i'}\in\mathbb{I}_{n,k}}\right) + n(k-1)^2\widehat{\theta}^2$$
$$- 2kn\mathbb{E}\left(\theta^*\theta_1^*|(Y_{\boldsymbol{i'}})_{\boldsymbol{i'}\in\mathbb{I}_{n,k}}\right) + 2n(k-1)\frac{n!}{n^k(n-k)!}\widehat{\theta}^2 - 2n(k-1)k\widehat{\theta}^2$$
$$= n\widehat{\theta}^2\left(\frac{n!}{n^{2k}(n-2k)!} + k^2\frac{(n-1)}{n} + (k-1)^2 + \left(2(k-1) - 2k\frac{n-k}{n}\right)\frac{n!}{n^k(n-k)!} - 2(k-1)k\right)$$
$$+ \left(\frac{k^2}{k!^2} + \frac{k^2}{k!^2} - 2\frac{k^2}{k!^2}\right)\mathbb{E}(h(1,...,k)h(\boldsymbol{1'})) + R,$$

with $R \xrightarrow{a.s.} 0$ and $\widehat{\theta}^2 \xrightarrow{a.s.} \theta_0^2$. Moreover $\frac{n!}{n^{2k}(n-2k)!} = 1 - \frac{1}{n}(k(2k-1)) + O(n^{-2})$, $\frac{n!}{n^k(n-k)!} = 1 - \frac{1}{n}\left(\frac{k(k-1)}{2}\right) + O(n^{-2})$. Next
$$\left(\frac{n!}{n^{2k}(n-2k)!} + k^2\frac{(n-1)}{n} + (k-1)^2 + \left(2(k-1) - 2k\frac{n-k}{n}\right)\frac{n!}{n^k(n-k)!} - 2(k-1)k\right)$$
$$= \left(\frac{n!}{n^{2k}(n-2k)!} + k^2\frac{(n-1)}{n} + (k-1)^2 + 2\left(\frac{k^2}{n} - 1\right)\frac{n!}{n^k(n-k)!} - 2(k-1)k\right)$$
$$= 1 + k^2 + (k-1)^2 - 2 - 2k^2 + 2k + \frac{1}{n}\left(k - 2k^2 - k^2 + 2k^2 + k(k-1)\right) + O(n^{-2})$$
$$= O(n^{-2}).$$

We have thus proved that
$$\mathbb{E}\left(n\left((\theta^* - \widehat{\theta}) - k(\theta_1^* - \widehat{\theta})\right)^2 |(Y_{\boldsymbol{i'}})_{\boldsymbol{i'}\in\mathbb{I}_{n,k}}\right) \xrightarrow{a.s.} 0.$$

Characterization of the convergence in distribution for the bootstrap using the bounded-Lipschitz metric ensures that it is sufficient to prove the asymptotic normality of $\sqrt{n}k(\theta_1^* - \widehat{\theta})$. Indeed if $L$ is a random variable whose distribution is the limit distribution of $\sqrt{n}k(\theta^* - \widehat{\theta})$ we have:
$$\sup_{h\in BL_1(\mathbb{R})}\left|\mathbb{E}\left(h(\sqrt{n}(\theta^* - \widehat{\theta}))|(Y_{\boldsymbol{i'}})_{\boldsymbol{i'}\in\mathbb{I}_{n,k}}\right) - \mathbb{E}(h(L))\right|$$
$$\leq \sup_{h\in BL_1(\mathbb{R})}\left|\mathbb{E}\left(h(\sqrt{n}k(\theta_1^* - \widehat{\theta}))|(Y_{\boldsymbol{i'}})_{\boldsymbol{i'}\in\mathbb{I}_{n,k}}\right) - \mathbb{E}(h(L))\right|$$
$$+ \mathbb{E}\left(\left|\sqrt{n}((\theta^* - \widehat{\theta}) - k(\theta_1^* - \widehat{\theta}))\right||(Y_{\boldsymbol{i'}})_{\boldsymbol{i'}\in\mathbb{I}_{n,k}}\right).$$

Then, by Levy criterion, pointwise convergence follows if we prove that for every $t \in \mathbb{R}$

$$\left|\mathbb{E}\left(\exp\left(\mathfrak{i}t\sqrt{n}k(\theta_1^* - \widehat{\theta})\right)|(Y_{i'})_{i' \in \mathbb{I}_{n,k}}\right) - \mathbb{E}\left(\exp(\mathfrak{i}tL)\right)\right| = o_{\text{a.s.}}(1), \tag{3.30}$$

where $\mathfrak{i}^2 = -1$. The next two substeps are devoted to proving this result.

**Substep 2:** $\lim_n \mathbb{E}\left(\left|\mathbb{E}\left(\exp\left(\mathfrak{i}t\sqrt{n}k!(\theta_1^* - \widehat{\theta})\right)|(Y_{i'})_{i' \in \mathbb{I}_{n,k}}\right) - e^{-t^2 \mathbb{V}\mathbb{E}(h(\mathbf{1})|U_{\{1\}})/2}\right|\right) = 0.$ Let us define

$$a_{n,i} = \frac{(n-k)!}{(n-1)!} \sum_{(i_2, ..., i_k) \in \mathbb{I}_{n,k-1}} h(i, i_2, ..., i_k),$$

$$a_{n,i}^* = \frac{(n-k)!}{(n-1)!} \sum_{(i_2, ..., i_k) \in \mathbb{I}_{n,k-1}} h(i^*, i_2, ..., i_k).$$

Given the sampling procedure in the bootstrap we have $\mathbb{E}(g(a_{n,i}^*)|(Y_{i'})_{i' \in \mathbb{I}_{n,k}}) = \frac{1}{n}\sum_{i=1}^n g(a_{n,i})$. Furthermore, $(a_{n,i}^*)_{i=1,...,n}$ forms an i.i.d. sequence conditional on $(Y_{i'})_{i' \in \mathbb{I}_{n,k}}$. Let $Z_n = \sqrt{n}k!(\theta_1^* - \widehat{\theta})$. Remark that $Z_n = \sum_{i=1}^n z_{n,i}/\sqrt{n}$, with $z_{n,i} = a_{n,i}^* - \frac{1}{n}\sum_{i'=1}^n a_{n,i'}$. We have for every $\epsilon > 0$ and $t \in \mathbb{R}$,

$$\left|\exp\left(\frac{\mathfrak{i}tz_{n,i}}{\sqrt{n}}\right) - \left(1 + \frac{\mathfrak{i}tz_{n,i}}{\sqrt{n}} - \frac{t^2 z_{n,i}^2}{2n}\right)\right| \leq \min\left(\frac{|t^3 z_{n,i}|^3}{\sqrt{n}^3}, \frac{t^2 z_{n,i}^2}{n}\right)$$

$$\leq \frac{|t^3 z_{n,i}^3|}{\sqrt{n}^3} \mathbb{1}_{\{|z_{n,i}| < \epsilon\sqrt{n}\}} + \frac{t^2 z_{n,i}^2}{n} \mathbb{1}_{\{|z_{n,i}| > \epsilon\sqrt{n}\}}$$

$$\leq \left(\epsilon|t|^3 + t^2 \mathbb{1}_{\{|z_{n,i}| > \epsilon\sqrt{n}\}}\right) \frac{z_{n,i}^2}{n}.$$

Let $V_n = \mathbb{E}(z_{n,i}^2|(Y_{i'})_{i' \in \mathbb{I}_{n,k}}) = \frac{1}{n}\sum_{i=1}^n a_{n,i}^2 - \left(\frac{1}{n}\sum_{i=1}^n a_{n,i}\right)^2$ and $V = \mathbb{V}\mathbb{E}(h(\mathbf{1})|U_{\{1\}})$. Lemma S7 and the fact that $(h(\mathbf{i}))_{\mathbf{i} \in \mathbb{I}_k}$ is $k$ jointly exchangeable and dissociated allow us to claim that

$$V_n \xrightarrow{L^1, \text{a.s.}} \mathbb{E}[h(1, ..., k)h(\mathbf{1}')] - \mathbb{E}[h(1, ..., k)]^2$$

$$= \mathbb{E}\left[\mathbb{E}\left[h(1, ..., k)h(\mathbf{1}')|U_{\{1\}}\right]\right] - \mathbb{E}\left[\mathbb{E}\left[h(1, ..., k)|U_{\{1\}}\right]\right]^2$$

$$= \mathbb{E}\left[\mathbb{E}\left[h(1, ..., k)|U_{\{1\}}\right]^2\right] - \mathbb{E}\left[\mathbb{E}\left[h(1, ..., k)|U_{\{1\}}\right]\right]^2 = V,$$

where the last equality can be recovered thanks to Assumption 1 and the almost-sure representation of $(h(\mathbf{i}))_{\mathbf{i} \in \mathbb{I}_k}$. As $\mathbb{E}(z_{n,i}|(Y_{i'})_{i' \in \mathbb{I}_{n,k}}) = 0$, we deduce from the triangle inequality that

$$\left|\mathbb{E}\left(\exp\left(\frac{\mathfrak{i}tz_{n,i}}{\sqrt{n}}\right)\Big|(Y_{i'})_{i' \in \mathbb{I}_{n,k}}\right) - \left(1 - \frac{t^2 V_n}{2n}\right)\right| \leq \epsilon|t|^3 \frac{V_n}{n} + \frac{t^2}{n}\mathbb{E}\left(z_{n,i}^2 \mathbb{1}_{\{|z_{n,i}| > \epsilon\sqrt{n}\}}|(Y_{i'})_{i' \in \mathbb{I}_{n,k}}\right),$$

and then

$$\left|\mathbb{E}\left(\exp\left(\frac{\mathfrak{i}tz_{n,i}}{\sqrt{n}}\right)\Big|(Y_{i'})_{i' \in \mathbb{I}_{n,k}}\right) - \left(1 - \frac{t^2 V}{2n}\right)\right|$$

$$\leq \epsilon|t|^3 \frac{V_n}{n} + \frac{t^2}{n}\mathbb{E}\left(z_{n,i}^2 \mathbb{1}_{\{|z_{n,i}| > \epsilon\sqrt{n}\}}|(Y_{i'})_{i' \in \mathbb{I}_{n,k}}\right) + \frac{t^2}{2n}|V_n - V|.$$

Because $|\prod_{i=1}^n a_i - \prod_{i=1}^n b_i| \leq \sum_{i=1}^n |a_i - b_i|$ if $\max_{i=1,\ldots,n} \max(|a_i|, |b_i|) \leq 1$ and since the $(z_{n,i})_{i=1\ldots n}$ are i.i.d. conditional on the data, we obtain

$$\left| \mathbb{E}\left(\exp(itZ_n)|(Y_{\boldsymbol{i'}})_{\boldsymbol{i'} \in \mathbb{I}_{n,k}}\right) - \exp\left(-\frac{t^2 V}{2}\right) \right|$$
$$\leq \epsilon |t|^3 V_n + t^2 \mathbb{E}\left(z_{n,1}^2 \mathbb{1}_{\{|z_{n,1}| > \epsilon \sqrt{n}\}} | (Y_{\boldsymbol{i'}})_{\boldsymbol{i'} \in \mathbb{I}_{n,k}}\right) + t^2 |V_n - V| + \left|\exp\left(-\frac{t^2 V}{2}\right) - \left(1 - \frac{t^2 V}{2n}\right)^n\right|.$$

A convexity argument and the Cauchy-Schwarz inequality ensure $z_{n,1}^2 \leq 2a_{n,1}^{*2} + 2\left(\frac{1}{n}\sum_{i'=1}^n a_{n,i'}\right)^2 \leq 2a_{n,1}^{*2} + 2\frac{(n-k)!}{n!}\sum_{\boldsymbol{i} \in \mathbb{I}_{n,k}} h(\boldsymbol{i})^2$. This implies

$$\mathbb{E}\left(z_{n,1}^2 \mathbb{1}_{\{|z_{n,1}| > \epsilon \sqrt{n}\}}\right)$$
$$\leq 2\mathbb{E}\left[\mathbb{E}\left(a_{n,1}^{*2} \mathbb{1}_{\{a_{n,1}^{*2} > \epsilon^2 n/4\}} | (Y_{\boldsymbol{i'}})_{\boldsymbol{i'} \in \mathbb{I}_{n,k}}\right)\right] + 2\mathbb{E}\left[\mathbb{E}\left(a_{n,1}^{*2} | (Y_{\boldsymbol{i'}})_{\boldsymbol{i'} \in \mathbb{I}_{n,k}}\right) \mathbb{1}_{\left\{\frac{(n-k)!}{n!}\sum_{\boldsymbol{i} \in \mathbb{I}_{n,k}} h(\boldsymbol{i})^2 > \epsilon^2 n/4\right\}}\right]$$
$$+ 2\mathbb{E}\left[\frac{(n-k)!}{n!} \sum_{\boldsymbol{i} \in \mathbb{I}_{n,k}} h(\boldsymbol{i})^2 \mathbb{E}\left(\mathbb{1}_{\{a_{n,1}^{*2} > \epsilon^2 n/4\}} | (Y_{\boldsymbol{i'}})_{\boldsymbol{i'} \in \mathbb{I}_{n,k}}\right)\right]$$
$$+ 2\mathbb{E}\left[\frac{(n-k)!}{n!} \sum_{\boldsymbol{i} \in \mathbb{I}_{n,k}} h(\boldsymbol{i})^2 \mathbb{1}_{\left\{\frac{(n-k)!}{n!}\sum_{\boldsymbol{i} \in \mathbb{I}_{n,k}} h(\boldsymbol{i})^2 > \epsilon^2 n/4\right\}}\right]$$
$$\leq 2\mathbb{E}\left[a_{n,1}^2 \mathbb{1}_{\{a_{n,1}^2 > \epsilon^2 n/4\}}\right] + 2\mathbb{E}\left[h^2(\boldsymbol{1}) \mathbb{1}_{\left\{\frac{(n-k)!}{n!}\sum_{\boldsymbol{i} \in \mathbb{I}_{n,k}} h^2(\boldsymbol{i}) > \epsilon^2 n/4\right\}}\right]$$
$$+ 2\left(1 - \frac{k}{n}\right)\mathbb{E}\left[h^2(2,3,\ldots,k+1) \mathbb{1}_{\{a_{n,1}^2 > \epsilon^2 n/4\}}\right] + 2\frac{k}{n}\mathbb{E}\left[h^2(\boldsymbol{1}) \mathbb{1}_{\{a_{n,1}^2 > \epsilon^2 n/4\}}\right]$$
$$+ 2\mathbb{E}\left[h^2(\boldsymbol{1}) \mathbb{1}_{\left\{\frac{(n-k)!}{n!}\sum_{\boldsymbol{i} \in \mathbb{I}_{n,k}} h^2(\boldsymbol{i}) > \epsilon^2 n/4\right\}}\right]. \tag{3.31}$$

Conditional on $U_{\{1\}}$, $(h(1, i_2, \ldots, i_k))_{(i_2, \ldots, i_k) \in \overline{(\mathbb{N}^+ \setminus \{1\})^{k-1}}}$ is a jointly exchangeable and dissociated array of dimension $k - 1$. Hence $a_{n,1} \xrightarrow{\text{a.s.}} \mathbb{E}\left(h(\boldsymbol{1})|U_{\{1\}}\right)$. Furthermore, $\frac{(n-k)!}{n!}\sum_{\boldsymbol{i} \in \mathbb{I}_{n,k}} h^2(\boldsymbol{i}) \xrightarrow{\text{a.s.}} \mathbb{E}(h^2(1, \ldots, k))$. As a result, all the indicator functions on the right-hand side of the last inequality in (3.31) tend to 0 almost surely. The dominated convergence theorem also ensures that $\mathbb{E}\left(z_{n,1}^2 \mathbb{1}_{\{|z_{n,1}| > \epsilon \sqrt{n}\}}\right) \to 0$ for every $\epsilon > 0$. Further, $\mathbb{E}(|V_n - V|) \to 0$ and $\left|\exp\left(-\frac{t^2 V}{2}\right) - \left(1 - \frac{t^2 V}{2n}\right)^n\right|$ converges almost surely to 0 and is bounded. As a consequence,

$$\limsup_n \mathbb{E}\left(\left|\mathbb{E}(\exp(itZ_n)|(Y_{\boldsymbol{i'}})_{\boldsymbol{i'} \in \mathbb{I}_{n,k}}) - e^{-t^2 V/2}\right|\right) \leq \epsilon |t|^3.$$

Since $\epsilon$ could be chosen arbitrarily small, we finally get

$$\lim_n \mathbb{E}\left(\left|\mathbb{E}(\exp(itZ_n)|(Y_{\boldsymbol{i'}})_{\boldsymbol{i'} \in \mathbb{I}_{n,k}}) - e^{-t^2 V/2}\right|\right) = 0.$$

**Substep 3: conclusion on the almost-sure weak convergence of the bootstrap mean**
We finally prove the almost-sure convergence of $\mathbb{E}(\exp(itZ_n)|(Y_{\boldsymbol{i'}})_{\boldsymbol{i'} \in \mathbb{I}_{n,k}})$, not only its convergence in $L^1$ as above. Recall that $V = \mathbb{V}\mathbb{E}(h(\boldsymbol{1})|U_{\{1\}})$ with $U$ stemming from the AHK representation of

$h(i)$. We have

$$\mathbb{E}(Z_n^2|(Y_i)_{i\in\mathbb{I}_k}) = \frac{1}{n}\sum_{i=1}^n \mathbb{V}\left(a_{n,i}^*|(Y_i)_{i\in\mathbb{I}_k}\right) = V_n \xrightarrow{a.s.} V.$$

Given $(Y_i)_{i\in\mathbb{I}_k}$, $Z_n$ is bounded in probability: for every $\varepsilon \in (0,1)$, considering

$$\eta((Y_i)_{i\in\mathbb{I}_k}) = \frac{\sup_n \mathbb{E}\left(Z_n^2|(Y_i)_{i\in\mathbb{I}_k}\right)}{\varepsilon},$$

we have $\mathbb{P}\left(Z_n^2 \geq \eta((Y_i)_{i\in\mathbb{I}_k})|(Y_i)_{i\in\mathbb{I}_k}\right) \leq \varepsilon$ by Markov's inequality. Given $(Y_i)_{i\in\mathbb{I}_k}$, every subsequence $Z_{\sigma(n)}$ admits a further subsequence $Z_{\sigma'\circ\sigma(n)}$ that converges in distribution to $L_{\sigma'\circ\sigma}$, by Prohorov's Theorem. By Levy's criterion for weak convergence, this means that there is a set $\Omega'$ of probability one, independent of $\sigma'$ and $\sigma$, such that for every $\omega \in \Omega'$, $\mathbb{E}\left(e^{itZ_{\sigma'\circ\sigma(n)}}|(Y_i)_{i\in\mathbb{I}_k} = (Y_i(\omega))_{i\in\mathbb{I}_k}\right)$ converges to $\mathbb{E}(e^{itL_{\sigma'\circ\sigma}}|(Y_i)_{i\in\mathbb{I}_k} = (Y_i(\omega))_{i\in\mathbb{I}_k})$ for every $t \in \mathbb{R}$. Note that $L_{\sigma'\circ\sigma}$ could depend on $(Y_i)_{i\in\mathbb{I}_k}$. We can now write

$$\mathbb{E}\left[\left|\mathbb{E}[e^{itL_{\sigma'\circ\sigma}}|(Y_i)_{i\in\mathbb{I}_k}] - \exp(-t^2V/2)\right|\right]$$
$$\leq \mathbb{E}\left[\left|\mathbb{E}[e^{itL_{\sigma'\circ\sigma}}|(Y_i)_{i\in\mathbb{I}_k}] - \mathbb{E}[e^{itZ_{\sigma'\circ\sigma(n)}}|(Y_i)_{i\in\mathbb{I}_k}]\right|\right] + \mathbb{E}\left[\left|\mathbb{E}[e^{itZ_{\sigma'\circ\sigma(n)}}|(Y_i)_{i\in\mathbb{I}_k}] - \exp(-t^2V/2)\right|\right].$$

The first term on the right-hand side converges to 0 by dominated convergence. The second term converges to 0 by the result proved in the second substep. We finally have that almost surely, $\mathbb{E}\left[e^{itL_{\sigma'\circ\sigma}}|(Y_i)_{i\in\mathbb{I}_k}\right] = \exp(-t^2V/2)$ for every $t \in \mathbb{R}$, every subsequence $\sigma$ and some subsequence $\sigma'$. From Urysohn's subsequence principle (see Tao, 2011, Section 2.1.17, Pages 185-186), this means that almost surely, $Z_n$ converges in distribution conditionally on $(Y_i)_{i\in\mathbb{I}_k}$ to $\mathcal{N}(0, V)$. We conclude that (3.30) holds with $L \sim \mathcal{N}\left(0, \frac{k^2}{k!^2}V\right)$.

**Second step: Asymptotic equicontinuity**

Let $\mathcal{F}_\delta = \{f = f_1 - f_2 : (f_1, f_2) \in \mathcal{F} \times \mathcal{F}, \mathbb{E}\left(f^2(Y_1)\right) \leq \delta^2\}$. We have to show

$$\lim_{\delta \to 0}\limsup_{n\to\infty} \mathbb{E}\left(\sup_{\mathcal{F}_\delta} |\mathbb{G}_n^*(f)||(Y_i)_{i\in\mathbb{I}_k}\right) \stackrel{a.s.}{=} 0.$$

Let $N^* = \frac{(n-k)!}{n!}\sum_{i\in\mathbb{I}_{n,k}} \mathbb{1}_{\{i^*\in\mathbb{I}_{n,k}\}}$. Note that $\mathbb{E}[\mathbb{P}_n^*f|(Y_i)_{i\in\mathbb{I}_k}] = \mathbb{P}_n'f = \frac{1}{n^k}\sum_{i\in\mathbb{I}_{n,k}} f(Y_i)$. By independence of the $i^*$ with $(Y_i)_{i\in\mathbb{I}_k}$, we have:

$$\mathbb{E}\left[\sup_{f\in\mathcal{F}_\delta} |\mathbb{G}_n^*f|\bigg|(Y_i)_{i\in\mathbb{I}_k}\right]$$
$$\leq \mathbb{E}\left[\sup_{f\in\mathcal{F}_\delta} \sqrt{n}\left|\mathbb{P}_n^*f - \mathbb{P}_n'f\right|\bigg|(Y_i)_{i\in\mathbb{I}_k}\right] + \sqrt{n}\left(1 - \frac{n!}{n^k(n-k)!}\right)\frac{(n-k)!}{n!}\sum_{i\in\mathbb{I}_{n,k}} F(Y_i)$$
$$\leq \mathbb{E}\left[\sup_{f\in\mathcal{F}_\delta} \sqrt{n}\left|\mathbb{P}_n^*f - \mathbb{P}_n'f\right|\bigg|(Y_i)_{i\in\mathbb{I}_k}\right] + \frac{(n-k)!}{n!}\sum_{i\in\mathbb{I}_{n,k}} F(Y_i) \times o(1)$$

Because $\frac{(n-k)!}{n!}\sum_{i\in\mathbb{I}_{n,k}} F(Y_i) \xrightarrow{a.s.} \mathbb{E}(F(Y_1))$, we only have to show that

$$\limsup_{n\to\infty} \mathbb{E}\left[\sup_{f\in\mathcal{F}_\delta} \sqrt{n}\left|\mathbb{P}_n^*f - \mathbb{P}_n'f\right|\bigg|(Y_i)_{i\in\mathbb{I}_k}\right] \xrightarrow{a.s.} 0 \text{ as } \delta \to 0.$$

Using the symmetrization step of Lemma S5, we can write that for some constant $C_k$ that depends on $k$ only

$$\mathbb{E}\left[\sup_{f \in \mathcal{F}_\delta} \sqrt{n} \left|\mathbb{P}_n^* f - \mathbb{P}_n' f\right| \Big| (Y_i)_{i \in \mathbb{I}_k}\right]$$

$$\leq k C_k \sqrt{n} \mathbb{E}\left[\sup_{f \in \mathcal{F}_\delta} \left|\frac{1}{n} \sum_{i_1=1}^n \varepsilon_{\{i_1\}} \frac{(n-k)!}{(n-1)!} \sum_{(i_2,\ldots,i_k):\boldsymbol{i} \in \mathbb{I}_{n,k}} f(Y_{\boldsymbol{i}^*}) \mathbb{1}_{\{\boldsymbol{i}^* \in \mathbb{I}_{n,k}\}}\right| \Big| (Y_i)_{i \in \mathbb{I}_k}, N^* > 0\right] \mathbb{P}(N^* > 0).$$

We have

$$\mathbb{E}\left[\sup_{f \in \mathcal{F}_\delta} \left|\frac{1}{n} \sum_{i_1=1}^n \varepsilon_{\{i_1\}} \frac{(n-k)!}{(n-1)!} \sum_{(i_2,\ldots,i_k):\boldsymbol{i} \in \mathbb{I}_{n,k}} f(Y_{\boldsymbol{i}^*}) \mathbb{1}_{\{\boldsymbol{i}^* \in \mathbb{I}_{n,k}\}}\right| \Big| (Y_i)_{i \in \mathbb{I}_k}, (\boldsymbol{i}^*)_{i \in \mathbb{I}_{n,k}}, N^* > 0\right]$$

$$\leq \frac{4\sqrt{2}}{\sqrt{n}} \int_0^{\sigma_{1,2}} \sqrt{\log 2 N(\varepsilon, \mathcal{F}_\delta, ||\cdot||_{1,2}^*)} d\varepsilon,$$

for $||f||_{1,2}^{2*} = \frac{1}{n} \sum_{i_1=1}^n \left(\frac{(n-k)!}{(n-1)!} \sum_{(i_2,\ldots,i_k):\boldsymbol{i} \in \mathbb{I}_{n,k}} f(Y_{\boldsymbol{i}^*}) \mathbb{1}_{\{\boldsymbol{i}^* \in \mathbb{I}_{n,k}\}}\right)^2$ and $\sigma_{1,2}^2 = \sup_{\mathcal{F}_\delta} ||f||_{1,2}^{*2}$. We now reason conditional on $N^* > 0$. The Cauchy-Schwarz inequality ensures $||f||_{1,2}^{*2} \leq N^* ||f||_{\mathbb{P}_n^*,2}^{*2}$ for $||f||_{\mathbb{P}_n^*,2}^{*2} = N^{*-1} \frac{(n-k)!}{n!} \sum_{i \in \mathbb{I}_{n,k}} f^2(Y_{\boldsymbol{i}^*}) \mathbb{1}_{\{\boldsymbol{i}^* \in \mathbb{I}_{n,k}\}}$. It follows that (see Point 1 of Lemma S12) $\sigma_{1,2}^2 \leq \sigma_n^{*2} = \sup_{\mathcal{F}_\delta} N^* ||f||_{\mathbb{P}_n^*,2}^{*2}$ and

$$N(\varepsilon, \mathcal{F}_\delta, ||\cdot||_{1,2}^*) \leq N(\varepsilon, \mathcal{F}_\delta, N^{*1/2} ||\cdot||_{\mathbb{P}_n^*,2}^*) \leq N(\varepsilon N^{*-1/2}, \mathcal{F}_\delta, ||\cdot||_{\mathbb{P}_n^*,2}^*).$$

Monotonicity of the integral, Points 3-4 of Lemma S12 and $\sqrt{a+b} \leq \sqrt{a} + \sqrt{b}$ entail

$$\mathbb{E}\left[\sup_{f \in \mathcal{F}_\delta} |\mathbb{G}_n^* f| \Big| (Y_i)_{i \in \mathbb{I}_k}\right]$$

$$\leq K_k' \mathbb{E}\left[\sigma_n^* + \int_0^{\sigma_n^*} \sqrt{\log N\left(4\varepsilon N^{*-1/2}, \mathcal{F}, ||\cdot||_{\mathbb{P}_n^*,2}^*\right)} d\varepsilon \Big| (Y_i)_{i \in \mathbb{I}_k}, N^* > 0\right] \mathbb{P}(N^* > 0),$$

for some constant $K_k'$ depending only on $k$. Furthermore, when $N^* > 0$ the following holds:

$$\int_0^{\sigma_n^*} \sqrt{\log N\left(4\varepsilon N^{*-1/2}, \mathcal{F}, ||\cdot||_{\mathbb{P}_n^*,2}^*\right)} d\varepsilon$$

$$= \int_0^{\sigma_n^*} \sqrt{\log N\left(\varepsilon ||F||_{\mathbb{P}_n^*,2}^* / (4 N^{*1/2} ||F||_{\mathbb{P}_n^*,2}^*), \mathcal{F}, ||\cdot||_{\mathbb{P}_n^*,2}^*\right)} d\varepsilon$$

$$= 4 N^{*1/2} ||F||_{\mathbb{P}_n^*,2}^* \int_0^{\sigma_n^*/(4 N^{*1/2} ||F||_{\mathbb{P}_n^*,2}^*)} \sqrt{\log N\left(\varepsilon ||F||_{\mathbb{P}_n^*,2}^*, \mathcal{F}, ||\cdot||_{\mathbb{P}_n^*,2}^*\right)} d\varepsilon$$

$$\leq 4 \sqrt{N^* ||F||_{\mathbb{P}_n^*,2}^{*2}} J_{\mathcal{F}}\left(\frac{\sqrt{\sigma_n^{*2}}}{4\sqrt{N^* ||F||_{\mathbb{P}_n^*,2}^{*2}}}\right).$$

This, Lemma S11, the facts that $\mathbb{E}(\sigma_n^{*2} | (Y_i)_{i \in \mathbb{I}_k}, N^* > 0) = \mathbb{E}(\sigma_n^{*2} | (Y_i)_{i \in \mathbb{I}_k}) / \mathbb{P}(N^* > 0)$, $\mathbb{E}(N^* ||F||_{\mathbb{P}_n^*,2}^{*2} | (Y_i)_{i \in \mathbb{I}_k}, N^* > 0) = \frac{1}{n^k} \sum_{i \in \mathbb{I}_{n,k}} F^2(Y_i) / \mathbb{P}(N^* > 0)$ and Jensen's inequality thus ensure

$$\mathbb{E}\left[\sup_{f \in \mathcal{F}_\delta} |\mathbb{G}_n^* f| \Big| (Y_i)_{i \in \mathbb{I}_k}\right]$$

$$\leq K_k' \left(\mathbb{E}\left(\sigma_n^{*2} | (Y_i)_{i \in \mathbb{I}_k}\right)^{1/2} + \left(\frac{1}{n^k} \sum_{i \in \mathbb{I}_{n,k}} F^2(Y_i)\right)^{1/2} J_{\mathcal{F}}\left(\frac{\mathbb{E}\left(\sigma_n^{*2} | (Y_i)_{i \in \mathbb{I}_k}\right)^{1/2}}{4 \left(\frac{1}{n^k} \sum_{i \in \mathbb{I}_{n,k}} F^2(Y_i)\right)^{1/2}}\right) \sqrt{\mathbb{P}(N^* > 0)}\right).$$

Since $\frac{1}{n^k} \sum_{i \in \mathbb{I}_{n,k}} F^2(Y_i) \xrightarrow{a.s.} \mathbb{E}\left(F^2(Y_1)\right)$, we only have to show that

$$\limsup_{n \to \infty} \mathbb{E}\left(\sigma_n^{*2} | (Y_i)_{i \in \mathbb{I}_k}\right) \xrightarrow{a.s.} 0. \tag{3.32}$$

We have

$$\sigma_n^{*2} = \sup_{\mathcal{F}_\delta} |\mathbb{P}_n^* f^2| \leq \sup_{\mathcal{F}_\delta} \left|\mathbb{P}_n^* f^2 - \mathbb{P}_n f^2\right| + \sup_{\mathcal{F}_\delta} \left|\mathbb{P}_n f^2 - P f^2\right| + \delta^2$$

$$\leq \sup_{\mathcal{F}_\infty} \left|\mathbb{P}_n^* f^2 - \mathbb{P}_n f^2\right| + \sup_{\mathcal{F}_\infty} \left|\mathbb{P}_n f^2 - P f^2\right| + \delta^2.$$

Point 5 of Lemma S12 implies that $\sup_Q N\left(\eta \|4F^2\|_{Q,1}, \mathcal{F}_\infty^2, \|\cdot\|_{Q,1}\right) < \infty$ for every $\eta > 0$. Theorem 2.1 and Lemma S5 imply

$$\mathbb{E}\left(\sup_{\mathcal{F}_\infty} \left|\mathbb{P}_n^* f^2 - \mathbb{P}_n f^2\right| \Big| (Y_i)_{i \in \mathbb{I}_k}\right) \xrightarrow{a.s.} 0 \quad \text{and} \quad \sup_{\mathcal{F}_\infty} \left|\mathbb{P}_n f^2 - P f^2\right| \xrightarrow{a.s.} 0,$$

which finally yields (3.32). The result follows.

## 3.8 Theorem 2.3

The proof is the same as that of Theorem 13.4 in Kosorok (2006), with one change only: we have to check that $\mathbb{G}$, the limit of $\theta \mapsto \sqrt{n}(\Psi_n(\theta) - \Psi(\theta))$, is continuous. Given the kernel of $\mathbb{G}$, it suffices to check that for all $(\pi, \pi') \in \mathfrak{S}(\{\mathbf{1}\}) \times \mathfrak{S}(\{\mathbf{1}'\})$,

$$\sup_{h \in \mathcal{H}} \left|\text{Cov}\left([\psi_{\theta,h} - \psi_{\theta_0,h}](Y_{\pi(\mathbf{1})}), [\psi_{\theta,h} - \psi_{\theta_0,h}](Y_{\pi'(\mathbf{1}')})\right)\right| \to 0. \tag{3.33}$$

By Cauchy-Schwarz's inequality and joint exchangeability, this covariance is smaller than

$$\mathbb{E}\left\{[\psi_{\theta,h} - \psi_{\theta_0,h}]^2 (Y_{\pi(\mathbf{1})})\right\} = P(\psi_{\theta,h} - \psi_{\theta_0,h})^2.$$

Therefore, Condition 4 ensures that (3.33) holds. The result follows.

## 3.9 Theorem 2.4

The first result follows by Theorem 2.1.2 because the class $\{u \mapsto 1\{u \leq y\} : y \in \mathbb{R}^p\}$ is pointwise measurable and satisfies Assumption 4. The second point follows directly from Point 1 and the functional delta method, see e.g. Theorem 20.8 in van der Vaart (2000). Finally, Point 3 follows from Theorem 2.2 and the functional delta method for the bootstrap, see e.g. Theorem 23.9 in van der Vaart (2000).

# 4 Proofs of the extensions

## 4.1 Theorem 3.1

For random variables (respectively vectors, matrices) $X$ and $L$ indexed by a set $\mathcal{I}$, $(X_i)_{i \in \mathcal{I}} \overset{d}{\leadsto} (L_i)_{i \in \mathcal{I}}$ means that $(X_i)_{i \in \mathcal{I}'}$ converges weakly to $(L_i)_{i \in \mathcal{I}'}$ for any finite subset $\mathcal{I}'$ of $\mathcal{I}$. This differs from $(X_i)_{i \in \mathcal{I}} \xrightarrow{d} (L_i)_{i \in \mathcal{I}}$, which means that $\mathcal{I}$, $(X_i)_{i \in \mathcal{I}} \overset{d}{\leadsto} (L_i)_{i \in \mathcal{I}}$ but also that $(X_i)_{i \in \mathcal{I}}$ is asymptotically tight. Asymptotic tightness is implied by asymptotic equicontinuity and total boundedness of the process $(X_i)_{i \in \mathcal{I}}$.

### 4.1.1 Pointwise convergence

Because $\mathbb{G}f = 0$, we have $\mathbb{E}(f(Y_{1,2})|U_1) = \mathbb{E}(f(Y_{1,2}|U_2) = \mathbb{E}(f(Y_{1,2}))$, and next $\mu_{\boldsymbol{m}} = 0$ for any $\boldsymbol{m} \in \mathcal{M}_0 = \{\boldsymbol{m} \in \mathbb{N}^3 : \min(m_1, m_2) = 0, \max(m_1, m_2) > 0, m_3 = 0\}$. So $\mu_{\boldsymbol{m}} \neq 0$ only if $\boldsymbol{m} \in \mathcal{M}_1 \cup \mathcal{M}_2 \cup \mathcal{M}_3 = \{\boldsymbol{m} \in \mathbb{N}^3 : \min(m_1, m_2) > 0; m_3 = 0\} \cup \{\boldsymbol{m} \in \mathbb{N}^3 : m_1 = m_2, m_3 > 0\} \cup \{\boldsymbol{m} \in \mathbb{N}^3 : m_1 \neq m_2, m_3 > 0\}$. Let $\mathcal{M} = \mathcal{M}_0 \cup \mathcal{M}_1 \cup \mathcal{M}_2 \cup \mathcal{M}_3$, because $\mathbb{V}(f(Y_{1,2})) < \infty$, we also have $\sum_{\boldsymbol{m} \in \mathcal{M}} \mu_{\boldsymbol{m}}^2(f) < \infty$.

We have:

$$\frac{n}{n(n-1)} \sum_{\boldsymbol{i} \in \mathbb{I}_{n,2}} \psi_{m_1}(U_{i_1})\psi_{m_2}(U_{i_2})$$
$$= \frac{n}{n-1}\left(\left(\frac{1}{\sqrt{n}}\sum_{i=1}^n \psi_{m_1}(U_i)\right)\left(\frac{1}{\sqrt{n}}\sum_{i=1}^n \psi_{m_2}(U_i)\right) - \frac{1}{n}\sum_{i=1}^n \psi_{m_1}(U_i)\psi_{m_2}(U_i)\right).$$

It follows from the law of large numbers, the multivariate central limit theorem and the continuous mapping theorem that

$$\left(\frac{n}{n(n-1)} \sum_{\boldsymbol{i} \in \mathbb{I}_{n,2}} \psi_{m_1}(U_{i_1})\psi_{m_2}(U_{i_2})\right)_{\boldsymbol{m} \in \mathcal{M}_1} \stackrel{d}{\rightsquigarrow} \left(Z_{m_1}Z_{m_2} - \mathbb{1}_{\{m_1 = m_2\}}\right)_{\boldsymbol{m} \in \mathcal{M}_1} \tag{4.1}$$

By almost-sure convergence of the sample mean of jointly exchangeable arrays (Eagleson and Weber, 1978), ergodicity of dissociated arrays (Kallenberg, 2005), the independence of the $U_i$ and the orthogonality of the $\psi_m$, we obtain

$$\frac{1}{n(n-1)} \sum_{\boldsymbol{i} \in \mathbb{I}_{n,2}} \psi_{m_1}(U_{i_1})\psi_{m_2}(U_{i_2})\psi_{m'_1}(U_{i_1})\psi_{m'_2}(U_{i_2}) \xrightarrow{\text{a.s.}} \mathbb{1}_{\{m_1 = m'_1, m_2 = m'_2\}}. \tag{4.2}$$

We have:

$$\left(\frac{n}{n(n-1)} \sum_{\boldsymbol{i} \in \mathbb{I}_{n,2}} \psi_{m_1}(U_{i_1})\psi_{m_2}(U_{i_2})\psi_{m_3}(U_{\{i_1,i_2\}})\right)_{\boldsymbol{m} \in \mathcal{M}_2 \cup \mathcal{M}_3}$$
$$= \left(\frac{n}{\sqrt{n(n-1)}}\frac{1}{\sqrt{n(n-1)}} \sum_{1 \leq i_1 < i_2 \leq n} [\psi_{m_1}(U_{i_1})\psi_{m_2}(U_{i_2}) + \psi_{m_1}(U_{i_2})\psi_{m_2}(U_{i_1})]\psi_{m_3}(U_{\{i_1,i_2\}})\right)_{\boldsymbol{m} \in \mathcal{M}_2 \cup \mathcal{M}_3}$$

Note that the $\psi_m$ are uniformly bounded, so for any $\eta > 0$,

$$\mathbb{1}_{\{|(\psi_{m_1}(U_{i_1})\psi_{m_2}(U_{i_2}) + \psi_{m_1}(U_{i_2})\psi_{m_2}(U_{i_1}))\psi_{m_3}(U_{\{i_1,i_2\}})| > \sqrt{n(n-1)}\eta\}} = 0$$

for $\sqrt{n(n-1)}\eta > 2^{3/2}$. Then conditionally on $(U_i)_{i \geq 1}$, the Lindeberg-Feller theorem and the almost-sure convergence (4.2) imply:

$$\left(\begin{array}{c}\left(\frac{1}{\sqrt{n(n-1)}}\sum_{\boldsymbol{i} \in \mathbb{I}_{n,2}} \psi_{m_1}(U_{i_1})\psi_{m_2}(U_{i_2})\psi_{m_3}(U_{\{i_1,i_2\}})\right)_{\boldsymbol{m} \in \mathcal{M}_2} \\ \left(\frac{1}{\sqrt{n(n-1)}}\sum_{\boldsymbol{i} \in \mathbb{I}_{n,2}} \psi_{m_1}(U_{i_1})\psi_{m_2}(U_{i_2})\psi_{m_3}(U_{\{i_1,i_2\}})\right)_{\boldsymbol{m} \in \mathcal{M}_3}\end{array}\right) \stackrel{d}{\rightsquigarrow} \left(\begin{array}{c}(Z_{m_1,m_3})_{\boldsymbol{m} \in \mathcal{M}_2} \\ (Z_{\{m_1;m_2\},m_3})_{\boldsymbol{m} \in \mathcal{M}_3}\end{array}\right)$$

Considering $V_{\boldsymbol{m}} = \frac{n}{n(n-1)} \sum_{\boldsymbol{i} \in \mathbb{I}_{n,2}} \psi_{m_1}(U_{i_1}) \psi_{m_2}(U_{i_2}) \psi_{m_3}(U_{\{i_1,i_2\}})$ for $\boldsymbol{m} \in \mathcal{M}_1 \cup \mathcal{M}_2 \cup \mathcal{M}_3$, the previous convergence means that for any $(t_{\boldsymbol{m}}) \in \mathbb{R}^{\mathcal{M}}$ and any finite subset $\mathcal{M}'_2 \subset \mathcal{M}_2$ and $\mathcal{M}'_3 \subset \mathcal{M}_3$:

$$\mathbb{E}\left[\exp\left(i \sum_{\boldsymbol{m} \in \mathcal{M}'_2 \cup \mathcal{M}'_3} t_{\boldsymbol{m}} V_{\boldsymbol{m}}\right) \bigg| (U_i)_{i \geq 1}\right] \xrightarrow{a.s.} \exp\left(-\frac{1}{2} \sum_{\boldsymbol{m} \in \mathcal{M}'_2 \cup \mathcal{M}'_3} t_{\boldsymbol{m}}^2\right)$$

Let $\mathcal{M}' = \mathcal{M}'_1 \cup \mathcal{M}'_2 \cup \mathcal{M}'_3$ a finite subset of $\mathcal{M}$ with $\mathcal{M}'_1 \subset \mathcal{M}_1$, $\mathcal{M}'_2 \subset \mathcal{M}_2$ and $\mathcal{M}'_3 \subset \mathcal{M}_3$, triangle inequality ensures that

$$\left|\mathbb{E}\left[\exp\left(i \sum_{\boldsymbol{m} \in \mathcal{M}'} t_{\boldsymbol{m}} V_{\boldsymbol{m}}\right)\right] - \exp\left(-\frac{1}{2} \sum_{\boldsymbol{m} \in \mathcal{M}'_2 \cup \mathcal{M}'_3} t_{\boldsymbol{m}}^2\right) \mathbb{E}\left[\exp\left(i \sum_{\boldsymbol{m} \in \mathcal{M}'_1} t_{\boldsymbol{m}} V_{\boldsymbol{m}}\right)\right]\right|$$

$$\leq \mathbb{E}\left[\left|\exp\left(i \sum_{\boldsymbol{m} \in \mathcal{M}'_1} t_{\boldsymbol{m}} V_{\boldsymbol{m}}\right)\right| \times \left|\mathbb{E}\left[\exp\left(i \sum_{\boldsymbol{m} \in \mathcal{M}'_2 \cup \mathcal{M}'_3} t_{\boldsymbol{m}} V_{\boldsymbol{m}}\right) \bigg| (U_i)_{i \geq 1}\right] - \exp\left(-\frac{1}{2} \sum_{\boldsymbol{m} \in \mathcal{M}'_2 \cup \mathcal{M}'_3} t_{\boldsymbol{m}}^2\right)\right|\right]$$

$$\leq \mathbb{E}\left[\left|\mathbb{E}\left[\exp\left(i \sum_{\boldsymbol{m} \in \mathcal{M}'_2 \cup \mathcal{M}'_3} t_{\boldsymbol{m}} V_{\boldsymbol{m}}\right) \bigg| (U_i)_{i \geq 1}\right] - \exp\left(-\frac{1}{2} \sum_{\boldsymbol{m} \in \mathcal{M}'_2 \cup \mathcal{M}'_3} t_{\boldsymbol{m}}^2\right)\right|\right],$$

which tends to 0 by the dominated convergence theorem. Convergence (4.1) implies

$$\lim_n \mathbb{E}\left[\exp\left(i \sum_{\boldsymbol{m} \in \mathcal{M}'_1} t_{\boldsymbol{m}} V_{\boldsymbol{m}}\right)\right] = \mathbb{E}\left[\exp\left(i \left(\sum_{\boldsymbol{m} \in \mathcal{M}'_1} t_{\boldsymbol{m}}(Z_{m_1} Z_{m_2} - \mathbb{1}_{\{m_1 = m_2\}})\right)\right)\right],$$

and again by dominated convergence and the Levy theorem:

$$\begin{pmatrix} \left(\frac{n}{n(n-1)} \sum_{\boldsymbol{i} \in \mathbb{I}_{n,2}} \psi_{m_1}(U_{i_1}) \psi_{m_2}(U_{i_2})\right)_{\boldsymbol{m} \in \mathcal{M}_1} \\ \left(\frac{n}{n(n-1)} \sum_{\boldsymbol{i} \in \mathbb{I}_{n,2}} \psi_{m_1}(U_{i_1}) \psi_{m_2}(U_{i_2}) \psi_{m_3}(U_{\{i_1,i_2\}})\right)_{\boldsymbol{m} \in \mathcal{M}_2} \\ \left(\frac{n}{n(n-1)} \sum_{\boldsymbol{i} \in \mathbb{I}_{n,2}} \psi_{m_1}(U_{i_1}) \psi_{m_2}(U_{i_2}) \psi_{m_3}(U_{\{i_1,i_2\}})\right)_{\boldsymbol{m} \in \mathcal{M}_3} \end{pmatrix} \xrightarrow{d} \begin{pmatrix} \left(Z_{m_1} Z_{m_2} - \mathbb{1}_{\{m_1 = m_2\}}\right)_{\boldsymbol{m} \in \mathcal{M}_1} \\ (Z_{m_1, m_3})_{\boldsymbol{m} \in \mathcal{M}_2} \\ \left(Z_{\{m_1, m_2\}, m_3}\right)_{\boldsymbol{m} \in \mathcal{M}_3} \end{pmatrix}$$
(4.3)

To finish the proof of pointwise convergence, we use an approximation argument. For $\overline{m} \in \mathbb{N}^+ \cup \infty$, let $\mathcal{M}(\overline{m}) = \{\boldsymbol{m} \in \{0, ..., \overline{m}\}^3 : \boldsymbol{m} \neq (0, 0, 0)\}$ and $\mathcal{M}_j(\overline{m}) = \mathcal{M}_j \cap \mathcal{M}(\overline{m})$, for $j = 1, 2, 3$. Let

$$R_{\overline{m}} = \frac{n}{n(n-1)} \sum_{\boldsymbol{i} \in \mathbb{I}_{n,2}} \sum_{\boldsymbol{m} \in \bigcup_{j=1,2,3} \mathcal{M}_j(\overline{m})} \mu_{\boldsymbol{m}}(f) \psi_{m_1}(U_{i_1}) \psi_{m_2}(U_{i_2}) \psi_{m_3}(U_{\{i_1,i_2\}}),$$

so that $R_\infty = n^{1/2} \mathbb{G}_n f$. Similarly, let

$$\mathbb{G}^d_{\overline{m}}(f) = \sum_{\boldsymbol{m} \in \mathcal{M}_1(\overline{m})} \mu_{\boldsymbol{m}}(Z_{m_1} Z_{m_2} - \mathbb{1}_{\{m_1 = m_2\}}) + \sum_{\boldsymbol{m} \in \mathcal{M}_2(\overline{m})} \mu_{\boldsymbol{m}} Z_{m_1, m_3} + \sum_{\boldsymbol{m} \in \mathcal{M}_3(\overline{m})} \mu_{\boldsymbol{m}} Z_{\{m_1, m_2\}, m_3},$$

so that $\mathbb{G}^d_\infty(f) = \mathbb{G}^d(f)$. For any Lipschitz function $\varphi$ from $(\mathbb{R}, ||\cdot||_1)$ to $(\mathbb{R}, ||\cdot||_1)$ with Lipschitz coefficient $C_\varphi$, we have:

$$\left|\mathbb{E}\left[\varphi(R_\infty)\right] - \mathbb{E}\left[\varphi(\mathbb{G}^d_\infty(f))\right]\right|$$
$$\leq C_\varphi \mathbb{E}\left[|R_\infty - R_{\overline{m}}|\right] + \left|\mathbb{E}\left[\varphi(R_{\overline{m}})\right] - \mathbb{E}\left[\varphi(\mathbb{G}^d_{\overline{m}}(f))\right]\right| + C_\varphi \mathbb{E}\left[\left|\mathbb{G}^d_\infty(f) - \mathbb{G}^d_{\overline{m}}(f)\right|\right]$$
$$\leq C_\varphi \sqrt{\mathbb{V}(R_\infty - R_{\overline{m}})} + \left|\mathbb{E}\left[\varphi(R_{\overline{m}})\right] - \mathbb{E}\left[\varphi(\mathbb{G}^d_{\overline{m}}(f))\right]\right| + C_\varphi \sqrt{\mathbb{V}(\mathbb{G}^d_\infty(f) - \mathbb{G}^d_{\overline{m}}(f))} \quad (4.4)$$

35

Now, remark that

$$\mathbb{V}(\mathbb{G}_\infty^d(f) - \mathbb{G}_{\overline{m}}^d(f)) = \sum_{m_1<m_2}(\mu_{m_1,m_2,0}(f)+\mu_{m_2,m_1,0}(f))^2 \mathbb{1}_{\{\overline{m}<m_2\}} + 2\sum_{m\geq\overline{m}}\mu_{m,m,0}^2(f)$$
$$+ \sum_{\boldsymbol{m}\in\bigcup_{j=2,3}\mathcal{M}_j\setminus\mathcal{M}_j(\overline{m})} \mu_{\boldsymbol{m}}^2$$
$$\leq 2 \sum_{\boldsymbol{m}\in\mathcal{M}\setminus\mathcal{M}(\overline{m})} \mu_{\boldsymbol{m}}^2(f). \tag{4.5}$$

Moreover, because,

$$\mathbb{E}(\psi_{m_1}(U_{i_1})\psi_{m_2}(U_{i_2})\psi_{m_3}(U_{\{i_1,i_2\}})\psi_{m_1'}(U_{i_1'})\psi_{m_2'}(U_{i_2'})\psi_{m_3'}(U_{\{i_1',i_2'\}}))$$
$$=\mathbb{1}_{\{m_3=m_3',\}}\left(\mathbb{1}_{\{i_1=i_1',i_2=i_2'\}}\mathbb{1}_{\{m_1=m_1',m_2=m_2'\}} + \mathbb{1}_{\{i_1=i_2',i_2=i_1'\}}\mathbb{1}_{\{m_1=m_2',m_2=m_1'\}}\right)$$

we have

$$\mathbb{V}(R_\infty - R_{\overline{m}}) = \frac{n}{n-1}\mathbb{V}(\mathbb{G}_\infty^d(f) - \mathbb{G}_{\overline{m}}^d(f)) \leq 4\sum_{\boldsymbol{m}\in\mathcal{M}\setminus\mathcal{M}(\overline{m})}\mu_{\boldsymbol{m}}^2(f). \tag{4.6}$$

By (4.3), $\lim_{n\to\infty}\left|\mathbb{E}\left[\varphi(R_{\overline{m}})\right] - \mathbb{E}\left[\varphi(\mathbb{G}_{\overline{m}}^d(f))\right]\right| = 0$ for any fixed $\overline{m}$. Next, by choosing $\overline{m}$ so that $\sum_{\boldsymbol{m}\in\mathcal{M}\setminus\mathcal{M}(\overline{m})}\mu_{\boldsymbol{m}}^2$ is arbitrarily small, (4.4)-(4.6) ensure that $\mathbb{E}\left[\varphi(R_\infty)\right] \to \mathbb{E}\left[\varphi(\mathbb{G}_\infty^d(f))\right]$. Therefore, $n^{1/2}\mathbb{G}_n(f) \xrightarrow{d} \mathbb{G}_\infty^d(f)$.

### 4.1.2 Asymptotic equicontinuity

We want to prove $\lim_{\delta\to 0}\limsup_{n\to\infty}\mathbb{E}\left[\sup_{f\in\mathcal{F}_\delta}\left|\frac{n}{n(n-1)}\sum_{\boldsymbol{i}\in\mathbb{I}_{n,2}}f(Y_{\boldsymbol{i}}) - \mathbb{E}[f(Y_{\boldsymbol{1}})]\right|\right] = 0$. We first use Lemma S3:

$$\mathbb{E}\left[\sup_{f\in\mathcal{F}_\delta}\left|\frac{n}{n(n-1)}\sum_{\boldsymbol{i}\in\mathbb{I}_{n,2}}f(Y_{\boldsymbol{i}}) - \mathbb{E}[f(Y_{\boldsymbol{1}})]\right|\right]$$
$$\lesssim \underbrace{\mathbb{E}\left[\sup_{f\in\mathcal{F}_\delta}\left|\frac{n}{n(n-1)}\sum_{\boldsymbol{i}\in\mathbb{I}_{n,2}}\varepsilon_{\{\boldsymbol{i}\}}f(Y_{\boldsymbol{i}})\right|\right]}_{=:A_1(n,\delta)} + \underbrace{\mathbb{E}\left[\sup_{f\in\mathcal{F}_\delta^s}\left|\frac{n}{n(n-1)}\sum_{\boldsymbol{i}\in\mathbb{I}_{n,2}}\varepsilon_{i_1}^1\varepsilon_{i_2}^2 f(Y_{\{\boldsymbol{i}\}}^1)\right|\right]}_{=:A_2(n,\delta)}. \tag{4.7}$$

To control $A_1(n,\delta)$, we remark that in the proof of Theorem 2.1.2 (second step in Section 3.5.2), the following result is given (setting $k$ and $r$ equal to 2)

$$\mathbb{E}\left[\sup_{f\in\mathcal{F}_\delta}\left|\frac{1}{n(n-1)}\sum_{\boldsymbol{i}\in\mathbb{I}_{n,2}}\varepsilon_{\{\boldsymbol{i}\}}f(Y_{\boldsymbol{i}})\right|\right] \leq \frac{8}{\sqrt{n(n-1)}}\mathbb{E}\left[\int_0^{\sigma_2}\sqrt{\log 2N\left(\eta,\mathcal{F}_\delta,||\cdot||_{\mathbb{P}_n^2,2}\right)}d\eta\right],$$

where $\sigma_2^2 := \sup_{f\in\mathcal{F}_\delta}||f||_{\mathbb{P}_n^2,2}$ and $||\cdot||_{\mathbb{P}_n^2,2}$ is defined in Section 3.1. As a result,

$$\mathbb{E}\left[\sup_{f\in\mathcal{F}_\delta}\left|\frac{n}{n(n-1)}\sum_{\boldsymbol{i}\in\mathbb{I}_{n,2}}\varepsilon_{\{\boldsymbol{i}\}}f(Y_{\boldsymbol{i}})\right|\right] \leq 8\sqrt{2}\mathbb{E}\left[\int_0^{\sigma_2}\sqrt{\log 2N\left(\eta,\mathcal{F}_\delta,||\cdot||_{\mathbb{P}_n^2,2}\right)}d\eta\right].$$

It is then shown in the second step in Section 3.5.2 that

$$\lim_{\delta\to 0}\limsup_{n\to\infty}\mathbb{E}\left[\int_0^{\sigma_2}\sqrt{\log 2N\left(\eta,\mathcal{F}_\delta,||\cdot||_{\mathbb{P}_n^2,2}\right)}d\eta\right] = 0$$

36

under Assumption 4-(i), and thus also if $\int_0^\infty \sup_{Q \in \mathcal{Q}} \log N(\eta\|F\|_{Q,2}, \mathcal{F}, \|\cdot\|_{Q,2})\, d\eta < \infty$. Hence,

$$\lim_{\delta \to 0} \limsup_{n \to \infty} A_1(n, \delta) = 0. \tag{4.8}$$

To control $A_2(n, \delta)$, remark first that by Assumption 2,

$$\mathbb{E}\left[\sup_{f \in \mathcal{F}_\delta^s} \left|\sum_{i \in \overrightarrow{\mathbb{I}_{n,2}}} \varepsilon_{i_1}^1 \varepsilon_{i_2}^2 f(Y_{\{i\}}^1)\right| \,\Big|\, (Y_i^1)_{i \in \mathbb{I}_2}\right] = \mathbb{E}\left[\sup_{f \in \mathcal{G}_\delta^s} \left|\sum_{i \in \overrightarrow{\mathbb{I}_{n,2}}} \varepsilon_{i_1}^1 \varepsilon_{i_2}^2 f(Y_{\{i\}}^1)\right| \,\Big|\, (Y_i^1)_{i \in \mathbb{I}_2}\right],$$

where $\mathcal{G}_\delta^s$ is a countable subclass of $\mathcal{F}_\delta^s$. Let $\|f\| := \sqrt{\frac{2}{n(n-1)} \sum_{i \in \overrightarrow{\mathbb{I}_{n,2}}} f(Y_{\{i\}}^1)^2}$ and remark that $T := \left\{\left(\sqrt{\frac{2}{n(n-1)}} f(Y_{\{i\}}^1)\right)_{i \in \overrightarrow{\mathbb{I}_{n,2}}} : f \in \mathcal{G}_\delta^s\right\}$ is a countable subset of $\mathbb{R}^{n(n-1)/2}$ with diameter

$$D := \sup_{(f_1, f_2) \in \mathcal{G}_\delta^s \times \mathcal{G}_\delta^s} \left\|\sqrt{\frac{2}{n(n-1)}}(f_1 - f_2)\right\|_T = \sup_{(f_1, f_2) \in \mathcal{F}_\delta^s \times \mathcal{F}_\delta^s} \|f_1 - f_2\|,$$

where $\|\cdot\|_T$ is defined in Lemma S10. Then, by this lemma,

$$\mathbb{E}\left[\sup_{f \in \mathcal{F}_\delta^s} \left|\sqrt{\frac{2}{n(n-1)}} \sum_{i \in \overrightarrow{\mathbb{I}_{n,2}}} \varepsilon_{i_1}^1 \varepsilon_{i_2}^2 f(Y_{\{i\}}^1)\right| \,\Big|\, (Y_i^1)_{i \in \mathbb{I}_2}\right] \lesssim \int_0^{D/2} \log N(\eta, \mathcal{G}_\delta^s, \|\cdot\|)\, d\eta.$$

Let $\mathbb{Q}_n := \frac{1}{n(n-1)} \sum_{i \in \mathbb{I}_{n,2}} \delta_{\{Y_i^1\}}$. By the triangle inequality used twice, we have

$$\frac{D}{2} \leq \sup_{g \in \mathcal{G}_\delta^s} \|g\|$$

$$\leq \sqrt{\frac{1}{2n(n-1)}} \sup_{f \in \mathcal{F}_\delta} \left\{\sqrt{\sum_{i \in \overrightarrow{\mathbb{I}_{n,2}}} f(Y_{i_1, i_2}^1)^2} + \sqrt{\sum_{i \in \overrightarrow{\mathbb{I}_{n,2}}} f(Y_{i_2, i_1}^1)^2}\right\}$$

$$\leq \sqrt{\frac{1}{2n(n-1)}} \sup_{f \in \mathcal{F}_\delta} \left\{\sqrt{\sum_{i \in \mathbb{I}_{n,2}} f(Y_{i_1, i_2}^1)^2} + \sqrt{\sum_{i \in \mathbb{I}_{n,2}} f(Y_{i_2, i_1}^1)^2}\right\}$$

$$\leq \sqrt{2} \sup_{f \in \mathcal{F}_\delta} \|f\|_{\mathbb{Q}_n, 2}.$$

In the same fashion, we can prove that for every $(g_1, g_2) \in \mathcal{G}_\delta^s \times \mathcal{G}_\delta^s$, there exist $(f_1, f_2) \in \mathcal{F}_\delta \times \mathcal{F}_\delta$ such that $\|g_1 - g_2\| \leq \sqrt{2}\|f_1 - f_2\|_{\mathbb{Q}_n, 2}$, $g_1(x, y) = (f_1(x) + f_1(y))/2$ and $g_2(x, y) = (f_2(x) + f_2(y))/2$. As a result, for every $\varepsilon > 0$, $N(\varepsilon, \mathcal{G}_\delta^s, \|\cdot\|) \leq N(\varepsilon/\sqrt{2}, \mathcal{F}_\delta, \|\cdot\|_{\mathbb{Q}_n, 2})$. We get

$$\int_0^{D/2} \log N(\eta, \mathcal{G}_\delta^s, \|\cdot\|)\, d\eta \leq \int_0^{\sqrt{2}\sigma} \log N(\eta, \mathcal{F}_\delta, \|\cdot\|_{\mathbb{Q}_n, 2})\, d\eta$$

$$\leq \int_0^{\sqrt{2}\sigma/\|2F\|_{\mathbb{Q}_n, 2}} \sup_Q \log N(\eta\|2F\|_{Q, 2}, \mathcal{F}_\delta, \|\cdot\|_{Q, 2})\, d\eta,$$

where $\sigma := \sup_{f \in \mathcal{F}_\delta} \|f\|_{\mathbb{Q}_n, 2}$. Integrating over $(Y_i^1)_{i \in \mathbb{I}_2}$ and reasoning as in the end of the second step in Section 3.5.2, we obtain

$$\lim_{\delta \to 0} \limsup_{n \to \infty} \mathbb{E}\left[\sup_{f \in \mathcal{F}_\delta^s} \left|\sqrt{\frac{2}{n(n-1)}} \sum_{i \in \overrightarrow{\mathbb{I}_{n,2}}} \varepsilon_{i_1}^1 \varepsilon_{i_2}^2 f(Y_{\{i\}}^1)\right|\right] = 0.$$

This in turn implies that $\lim_{\delta \to 0} \limsup_{n \to \infty} A_2(n, \delta) = 0$. This last result, combined with (4.7) and (4.8) is enough to conclude.

## 4.2 Theorem 3.2

Hereafter, we let $K$ denote the covariance kernel of $\mathbb{G}$, $\rho(f_1, f_2) = \sqrt{K(f_1 - f_2, f_1 - f_2)}$ and $\nu(f_1, f_2) = \sqrt{K(f_1 - f_2, f_1 - f_2) + [\mathbb{E}(f_1 - f_2)]^2}$. We also introduce $\mathcal{F}_\delta = \{g = f_1 - f_2 : (f_1, f_2) \in \mathcal{F} \times \mathcal{F}, \nu(f_1, f_2) \leq \delta\}$ and note that for all $(f_1, f_2) \in \mathcal{F} \times \mathcal{F}$,

$$\nu(f_1, f_2) = \sqrt{\mathbb{E}\left[\mathbb{E}((f_1 - f_2)(Y_{1,2}) + (f_1 - f_2)(Y_{2,1})|U_1)^2\right]}.$$

Finally, we define the auxiliary bootstrap process

$$\mathbb{G}_n^{*\prime}(f) = \sqrt{n}\left(\frac{1}{n(n-1)} \sum_{1 \leq i,j \leq n} \frac{f(Y_{i,j^*}) + f(Y_{j^*,i})}{2} \mathbb{1}_{\{i \neq j^*\}} - \frac{1}{n(n-1)} \sum_{i,j \in \mathbb{I}_{n,2}} f(Y_{i,j})\right).$$

We prove the result in three steps. We first prove pointwise convergence of $\mathbb{G}_n f$ and total boundedness with respect to $\nu$. Next, we prove the convergence of the process $\mathbb{G}_n^{*\prime}$. Finally, we show the asymptotic equicontinuity of $\mathbb{G}_n$ with respect to $\nu$.

**1. Pointwise convergence and total boundedness with respect to $\nu$.**

First, by the pointwise convergence established in the proof of Theorem 2.1 and since $\mathbb{E}(f^2(Y_{1,2})) < \infty$, $\mathbb{G}_n f$ converges to $\mathbb{G} f$ for all $f \in \mathcal{F}$. Second, because the process $\mathbb{G}_n^*$ converges weakly in $\ell^\infty(\mathcal{F})$ to $\mathbb{G}$, $(\mathcal{F}, \rho)$ is totally bounded. The set $\{f : |\mathbb{E}(f)| \leq 2\mathbb{E}(F)\}$ is also totally bounded for the semimetric $\rho'(f_1, f_2) = |\mathbb{E}(f_1 - f_2)|$. Now, let $(B_i^1)_{i=1,\ldots,n_1}$ (resp. $(B_{i'}^2)_{i'=1,\ldots,n_2}$) denotes a $\delta/2\sqrt{2}$-covering of $(\mathcal{F}, \rho)$ (resp. of $(\mathcal{F}, \rho')$). Each non-empty set $B_i^1 \cap B_{i'}^2$ is included in a $\nu$-ball of radius $\delta$ and $\mathcal{F} \subset \bigcup_{i,i'}(B_i^1 \cap B_{i'}^2)$. This ensures that $(\mathcal{F}, \nu)$ is totally bounded.

**2. Asymptotic equicontinuity of the process $\mathbb{G}_n^{*\prime}$.**

Let $(1^{**}, \ldots, n^{**})$ be an independent copy of $(1^*, \ldots, n^*)$ and let

$$\mathbb{G}_n^{*\prime\prime}(f) = \sqrt{n}\left(\frac{1}{n(n-1)} \sum_{1 \leq i,j \leq n} \frac{f(Y_{i^{**},j^*}) + f(Y_{j^*,i^{**}})}{2} \mathbb{1}_{\{i^{**} \neq j^*\}} - \frac{1}{n(n-1)} \sum_{1 \leq i,j \leq n} f(Y_{i,j})\mathbb{1}_{\{i \neq j\}}\right).$$

Conditional on $(Y_{i,j})_{i,j \in \mathbb{I}_2}$, $\mathbb{G}_n^*$ is a $U$-process on the class $\{g : \{1, \ldots, n\} \to \mathbb{R} : g(i,j) = (f(Y_{i,j}) + f(Y_{j,i}))/2, f \in \mathcal{F}\}$, while $\mathbb{G}_n^{*\prime\prime}$ corresponds to its decoupled version. Then by Theorem 3.1.1 of de la Peña and Giné (1999),

$$\mathbb{E}\left[\sup_\mathcal{F} |\mathbb{G}_n^{*\prime\prime}(f)|\bigg|(Y_{i,j})_{i,j \in \mathbb{I}_2}\right] \lesssim \mathbb{E}\left[\sup_\mathcal{F} |\mathbb{G}_n^*(f)|\bigg|(Y_{i,j})_{i,j \in \mathbb{I}_2}\right].$$

Moreover, we have :

$$\mathbb{E}\left[\sup_\mathcal{F} |\mathbb{G}_n^{*\prime\prime}(f)|\bigg|(Y_{i,j})_{i,j \in \mathbb{I}_2}, (1^*, \ldots, n^*)\right]$$

$$= \frac{1}{n^n} \sum_{1 \leq i_1, \ldots, i_n \leq n} \sup_\mathcal{F} \left|\sqrt{n}\left(\frac{1}{n(n-1)} \sum_{1 \leq j',j \leq n} \frac{f(Y_{i_{j'},j^*}) + f(Y_{j^*,i_{j'}})}{2} \mathbb{1}_{\{i_{j'} \neq j^*\}} - \frac{1}{n(n-1)} \sum_{1 \leq i,j \leq n} f(Y_{i,j})\mathbb{1}_{\{i \neq j\}}\right)\right|$$

$$\geq \sup_\mathcal{F} \left|\sqrt{n}\left(\frac{1}{n(n-1)} \frac{1}{n^n} \sum_{1 \leq i_1, \ldots, i_n \leq n} \sum_{1 \leq j',j \leq n} \frac{f(Y_{i_{j'},j^*}) + f(Y_{j^*,i_{j'}})}{2} \mathbb{1}_{\{i_{j'} \neq j^*\}} - \frac{1}{n(n-1)} \sum_{1 \leq i,j \leq n} f(Y_{i,j})\mathbb{1}_{\{i \neq j\}}\right)\right|$$

$$= \sup_\mathcal{F} |\mathbb{G}_n^{*\prime}(f)|.$$

38

Integrating over $(1^*, ..., n^*)$, we then obtain $\mathbb{E}\left[\sup_{\mathcal{F}} |\mathbb{G}_n^{*'}(f)| \Big| (Y_{i,j})_{i,j \in \mathbb{I}_2}\right] \lesssim \mathbb{E}\left[\sup_{\mathcal{F}} |\mathbb{G}_n^*(f)| \Big| (Y_{i,j})_{i,j \in \mathbb{I}_2}\right]$. The asymptotic equicontinuity of $\mathbb{G}_n^{*'}$ follows from that of $\mathbb{G}_n^*$.

### 3. Asymptotic equicontinuity of $\mathbb{G}_n$ with respect to $\nu$.

The idea is to prove asymptotic equicontinuity for symmetrized processes and then exploit the symmetrization lemma A.1. First, by Lemma S4, we have

$$\frac{(1-e^{-1})(1-e^{-1/2})}{\sqrt{2}} \mathbb{E}\left(\sup_{\mathcal{F}_\delta} \left| \frac{\sqrt{n}}{n(n-1)} \sum_{i,j \in \mathbb{I}_{n,2}} \varepsilon_{\{i,j\}} \frac{f(Y_{i,j}) + f(Y_{j,i})}{2} \right| \Big| (Y_{i,j})_{i,j \in \mathbb{I}_2}\right)$$

$$\leq \mathbb{E}\left(\sup_{\mathcal{F}_\delta} \left| \frac{\sqrt{n}}{n(n-1)} \sum_{i,j=1}^n \varepsilon_{\{i,j\}} \frac{f(Y_{i^*,j^*}) + f(Y_{j^*,i^*})}{2} \mathbb{1}_{\{i^* \neq j^*\}} \right| \Big| (Y_{i,j})_{i,j \in \mathbb{I}_2}\right)$$

$$= \frac{1}{n(n-1)} \sum_{i',j'=1}^n \mathbb{E}\left(\sup_{\mathcal{F}_\delta} \left| \frac{\sqrt{n}}{n^2} \sum_{i,j=1}^n \varepsilon_{\{i,j\}} \frac{f(Y_{i',j'}) + f(Y_{j',i'})}{2} \mathbb{1}_{\{i' \neq j'\}} \right| \Big| (Y_{i,j})_{i,j \in \mathbb{I}_2}\right)$$

$$\leq \frac{1}{n(n-1)} \sum_{i',j'=1}^n \sup_{\mathcal{F}_\delta} \left| \frac{f(Y_{i',j'}) + f(Y_{j',i'})}{2} \mathbb{1}_{\{i \neq j'\}} \right| n^{-3/2} \mathbb{E}\left| \sum_{1 \leq i,j \leq n} \varepsilon_{\{i,j\}} \right|$$

$$\leq \frac{2}{\sqrt{n}} \frac{1}{n(n-1)} \sum_{i,j \in \mathbb{I}_{n,2}} F(Y_{i,j}). \tag{4.9}$$

Besides, using $\mathbb{1}_{\{a>b\}} \leq (a/b)^\delta$ and convexity of $u \mapsto u^{1+\delta}$,

$$\mathbb{E}\left(\frac{1}{n(n-1)} \sum_{i,j \in \mathbb{I}_{n,2}} F(Y_{i,j}) \mathbb{1}_{\left\{\frac{1}{n(n-1)} \sum_{i,j \in \mathbb{I}_{n,2}} F(Y_{i,j}) > M\right\}}\right) \leq M^{-\delta} \mathbb{E}\left[\left(\frac{1}{n(n-1)} \sum_{i,j \in \mathbb{I}_{n,2}} F(Y_{i,j})\right)^{1+\delta}\right]$$

$$\leq M^{-\delta} \mathbb{E}(F^{1+\delta}(Y_{1,2})),$$

with $\mathbb{E}(F^{1+\delta}(Y_{1,2})) < \infty$. It follows that $\sum_{i,j \in \mathbb{I}_{n,2}} F(Y_{i,j})/n(n-1)$ is uniform integrable, namely

$$\lim_{M \to \infty} \sup_n \mathbb{E}\left(\frac{1}{n(n-1)} \sum_{i,j \in \mathbb{I}_{n,2}} F(Y_{i,j}) \mathbb{1}_{\left\{\frac{1}{n(n-1)} \sum_{i,j \in \mathbb{I}_{n,2}} F(Y_{i,j}) > M\right\}}\right) = 0.$$

By (4.9) and monotonicity of $y \mapsto y \mathbb{1}_{\{y > M\}}$,

$$\mathbb{E}\left(\sup_{\mathcal{F}_\delta} \left| \frac{\sqrt{n}}{n(n-1)} \sum_{i,j \in \mathbb{I}_{n,2}} \varepsilon_{\{i,j\}} \frac{f(Y_{i,j}) + f(Y_{j,i})}{2} \right| \Big| (Y_{i,j})_{i,j \in \mathbb{I}_2}\right)$$

is an uniformly integrable sequence as well. Then, Fatou's inequality for uniformly integrable random variables (see (see Shiryaev, 2007, Remark 2 in Chapter 1), together with (4.9), imply

$$\limsup_n \mathbb{E}\left(\sup_{\mathcal{F}_\delta} \left| \frac{\sqrt{n}}{n(n-1)} \sum_{i,j \in \mathbb{I}_{n,2}} \varepsilon_{\{i,j\}} \frac{f(Y_{i,j}) + f(Y_{j,i})}{2} \right|\right)$$

$$\leq \mathbb{E}\left[\limsup_n \mathbb{E}\left(\sup_{\mathcal{F}_\delta} \left| \frac{\sqrt{n}}{n(n-1)} \sum_{i,j \in \mathbb{I}_{n,2}} \varepsilon_{\{i,j\}} \frac{f(Y_{i,j}) + f(Y_{j,i})}{2} \right| \Big| (Y_{i,j})_{i,j \in \mathbb{I}_2}\right)\right]$$

$$= 0. \tag{4.10}$$

Next, let $(\varepsilon_j)_{j \in \mathbb{N}^+}$ denote independent Rademacher variables that are independent of the data. By the first inequality in Proposition 2.2 of Giné and Zinn (1990),

$$\frac{1-e^{-1}}{\sqrt{2}} \mathbb{E}\left(\sup_{\mathcal{F}_\delta} \left| \frac{\sqrt{n}}{n(n-1)} \sum_{i,j \in \mathbb{I}_{n,2}} \varepsilon_j \frac{f(Y_{i,j}) + f(Y_{j,i})}{2} \right| \middle| (Y_{i,j})_{i,j \in \mathbb{I}_2} \right)$$

$$= \frac{1-e^{-1}}{\sqrt{2}} \mathbb{E}\left(\sup_{\mathcal{F}_\delta} \left| \frac{\sqrt{n}}{n} \sum_{j=1}^n \frac{\varepsilon_j}{n-1} \sum_{i=1}^n \frac{f(Y_{i,j}) + f(Y_{j,i})}{2} \mathbb{1}_{\{i \neq j\}} \right| \middle| (Y_{i,j})_{i,j \in \mathbb{I}_2} \right)$$

$$\leq \mathbb{E}\left(\sup_{\mathcal{F}_\delta} \left| \frac{\sqrt{n}}{n} \sum_{j=1}^n \frac{\varepsilon_j}{n-1} \sum_{i=1}^n \frac{f(Y_{i,j^*}) + f(Y_{j^*,i})}{2} \mathbb{1}_{\{i \neq j^*\}} \right| \middle| (Y_{i,j})_{i,j \in \mathbb{I}_2} \right)$$

$$= \frac{1}{n} \sum_{j'=1}^n \mathbb{E}\left(\sup_{\mathcal{F}_\delta} \left| \frac{\sqrt{n}}{n} \sum_{j=1}^n \frac{\varepsilon_j}{n-1} \sum_{i=1}^n \frac{f(Y_{i,j'}) + f(Y_{j',i})}{2} \mathbb{1}_{\{i \neq j'\}} \right| \middle| (Y_{i,j})_{i,j \in \mathbb{I}_2} \right)$$

$$\leq \frac{1}{n} \sum_{j'=1}^n \sup_{\mathcal{F}_\delta} \left| \frac{1}{n-1} \sum_{i=1}^n \frac{f(Y_{i,j'}) + f(Y_{j',i})}{2} \mathbb{1}_{\{i \neq j'\}} \right|$$

$$\leq \frac{2}{n(n-1)} \sum_{i,j \in \mathbb{I}_{n,2}} F(Y_{i,j}). \tag{4.11}$$

As above, this inequality and monotonicity of $y \mapsto y \mathbb{1}_{\{y > M\}}$ imply that

$$\mathbb{E}\left(\sup_{\mathcal{F}_\delta} \left| \frac{\sqrt{n}}{n(n-1)} \sum_{i,j \in \mathbb{I}_{n,2}} \varepsilon_j \frac{f(Y_{i,j}) + f(Y_{j,i})}{2} \right| \middle| (Y_{i,j})_{i,j \in \mathbb{I}_2} \right)$$

is an uniformly integrable sequence. Then, by Fatou's inequality again,

$$\limsup_n \mathbb{E}\left(\sup_{\mathcal{F}_\delta} \left| \frac{\sqrt{n}}{n(n-1)} \sum_{i,j \in \mathbb{I}_{n,2}} \varepsilon_j \frac{f(Y_{i,j}) + f(Y_{j,i})}{2} \right| \right)$$

$$\leq \mathbb{E}\left( \limsup_n \mathbb{E}\left(\sup_{\mathcal{F}_\delta} \left| \frac{\sqrt{n}}{n(n-1)} \sum_{i,j \in \mathbb{I}_{n,2}} \varepsilon_j \frac{f(Y_{i,j}) + f(Y_{j,i})}{2} \right| \middle| (Y_{i,j})_{i,j \in \mathbb{I}_2} \right) \right). \tag{4.12}$$

Now, we refine (4.11). Using again the first inequality in Proposition 2.2 of Giné and Zinn (1990), the triangle inequality and the symmetrization lemma for independent random variable, we have:

$$\frac{1-e^{-1}}{\sqrt{2}} \mathbb{E}\left(\sup_{\mathcal{F}_\delta} \left| \frac{\sqrt{n}}{n(n-1)} \sum_{i,j \in \mathbb{I}_{n,2}} \varepsilon_j \frac{f(Y_{i,j}) + f(Y_{j,i})}{2} \right| \middle| (Y_{i,j})_{i,j \in \mathbb{I}_2} \right)$$

$$\leq \mathbb{E}\left(\sup_{\mathcal{F}_\delta} \left| \sqrt{n} \left( \sum_{j=1}^n \varepsilon_j \left[ \frac{1}{n(n-1)} \sum_{1 \leq i \leq n} \frac{f(Y_{i,j^*}) + f(Y_{j^*,i})}{2} \mathbb{1}_{\{i \neq j^*\}} - \frac{1}{n(n-1)} \sum_{i'',j'' \in \mathbb{I}_{n,2}} f(Y_{i'',j''}) \right] \right) \right|$$

$$+ \sup_{\mathcal{F}_\delta} \left| \frac{1}{\sqrt{n}} \sum_{j=1}^n \varepsilon_j \frac{1}{n(n-1)} \sum_{i'',j'' \in \mathbb{I}_{n,2}} f(Y_{i'',j''}) \right| \middle| (Y_{i,j})_{i,j \in \mathbb{I}_2} \right)$$

$$\leq 2\mathbb{E}\left(\sup_{\mathcal{F}_\delta} \left| \sqrt{n} \left( \frac{1}{n(n-1)} \sum_{1 \leq i,j \leq n} \frac{f(Y_{i,j^*}) + f(Y_{j^*,i})}{2} \mathbb{1}_{\{i \neq j^*\}} - \frac{1}{n(n-1)} \sum_{i'',j'' \in \mathbb{I}_{n,2}} f(Y_{i'',j''}) \right) \right| \middle| (Y_{i,j})_{i,j \in \mathbb{I}_2} \right)$$

$$+ \sup_{\mathcal{F}_\delta} \left| \frac{1}{n(n-1)} \sum_{i,j \in \mathbb{I}_{n,2}} f(Y_{i,j}) \right|.$$

40

Combined with (4.12) and Step 2, this entails

$$\limsup_n \mathbb{E}\left(\sup_{\mathcal{F}_\delta}\left|\frac{\sqrt{n}}{n(n-1)}\sum_{i,j\in\mathbb{I}_{n,2}}\varepsilon_j\frac{f(Y_{i,j})+f(Y_{j,i})}{2}\right|\right)$$

$$\leq \frac{\sqrt{2}}{1-e^{-1}}\left\{2\mathbb{E}\left[\limsup_n \mathbb{E}\left(\sup_{\mathcal{F}_\delta}|\mathbb{G}_n^{*\prime}(f)|\,\Big|(Y_{i,j})_{i,j\in\mathbb{I}_2}\right)\right]+\mathbb{E}\left[\limsup_n \sup_{\mathcal{F}_\delta}\left|\frac{1}{n(n-1)}\sum_{i,j\in\mathbb{I}_{n,2}}f(Y_{i,j})\right|\right]\right\}$$

$$\leq \frac{\sqrt{2}}{1-e^{-1}}\left\{2\mathbb{E}\left[\sup_{\mathcal{F}_\delta}|\mathbb{G}(f)|\right]+\mathbb{E}\left[\limsup_n \sup_{\mathcal{F}_\delta}\left|\frac{1}{n(n-1)}\sum_{i,j\in\mathbb{I}_{n,2}}f(Y_{i,j})\right|\right]\right\}. \qquad (4.13)$$

Now, (4.11) implies that

$$\limsup_n \mathbb{E}\left(\sup_{\mathcal{F}_\delta}\left|\frac{1}{n(n-1)}\sum_{i,j\in\mathbb{I}_{n,2}}\varepsilon_j\frac{f(Y_{i,j})+f(Y_{j,i})}{2}\right|\right)=0.$$

Together with (4.10), this implies, by the symmetrization lemma A.1 applied to the array $Z_{i,j}=(Y_{i,j},Y_{j,i})$ and the class $\{g(z_{i,j})=(f(y_{i,j})+f(y_{j,i}))/2:f\in\mathcal{F}_\delta\}$), that

$$\sup_{\mathcal{F}_\delta}\left|\frac{1}{n(n-1)}\sum_{i,j\in\mathbb{I}_{n,2}}\frac{f(Y_{i,j})+f(Y_{j,i})}{2}-\mathbb{E}\left[\frac{f(Y_{1,2})+f(Y_{2,1})}{2}\right]\right|\xrightarrow{L^1}0.$$

By a backward submartingale argument, convergence also holds almost surely. Hence, by the triangle and Cauchy-Schwarz inequalities, we have, almost surely,

$$\limsup_n \sup_{\mathcal{F}_\delta}\left|\frac{1}{n(n-1)}\sum_{i,j\in\mathbb{I}_{n,2}}f(Y_{i,j})\right|\leq \sup_{\mathcal{F}_\delta}\left|\mathbb{E}\left[\frac{f(Y_{1,2})+f(Y_{2,1})}{2}\right]\right|$$

$$\leq \sup_{\mathcal{F}_\delta}\sqrt{\mathbb{E}\left(\mathbb{E}(f(Y_{1,2})+f(Y_{2,1})|U_1)^2\right)}/2$$

$$\leq \delta/2.$$

Plugging in this inequality in (4.13), we obtain

$$\lim_{\delta\to 0}\limsup_n \mathbb{E}\left(\sup_{\mathcal{F}_\delta}\left|\frac{\sqrt{n}}{n(n-1)}\sum_{i,j\in\mathbb{I}_{n,2}}\varepsilon_j\frac{f(Y_{i,j})+f(Y_{j,i})}{2}\right|\right)=0.$$

Combined with (4.10), this implies, by the symmetrization lemma A.1,

$$\lim_{\delta\to 0}\limsup_n \mathbb{E}\left(\sup_{\mathcal{F}_\delta}\left|\frac{\sqrt{n}}{n(n-1)}\sum_{i,j\in\mathbb{I}_{n,2}}f(Y_{i,j})\right|\right)=0.$$

### 4.3 Theorem 3.3

**First step: pointwise convergence**

With a slight abuse of notation, we assimilate $f\in\mathcal{F}^s$ with $f\in\mathcal{F}$. Let $V=\mathbb{V}(\mathbb{E}(2f(\mathbf{Y}_{\{1\}})|U_1))$. We first establish that for all $(f,t)\in\mathcal{F}\times\mathbb{R}$,

$$\mathbb{E}\left[e^{it\mathbb{G}_n^{m*}f}|(Y_{\boldsymbol{i}})_{\boldsymbol{i}\in\mathbb{I}_2}\right]\xrightarrow{L^1}e^{-t^2V/2}. \qquad (4.14)$$

41

Let $z_{ni_1} := 2\xi_{i_1}\left(\frac{1}{n-1}\sum_{i_2 \neq i_1} f(\mathbf{Y}_{\{i\}}) - \mathbb{P}_n f\right)$, $V_{ni_1} := \mathbb{V}(z_{ni_1}|(Y_i)_{i\in\mathbb{I}_2})$. We have

$$\left|\mathbb{E}\left[e^{it\mathbb{G}_n^{m*}f}|(Y_i)_{i\in\mathbb{I}_2}\right] - e^{-t^2V/2}\right|$$

$$= \left|\prod_{i_1=1}^n \mathbb{E}\left[e^{itz_{ni_1}/\sqrt{n}}|(Y_i)_{i\in\mathbb{I}_2}\right] - e^{-t^2V/2}\right|$$

$$\leq \left|e^{-t^2\sum_{i_1=1}^n V_{ni_1}/2n} - e^{-t^2V/2}\right| + \left|\prod_{i_1=1}^n \mathbb{E}\left[e^{itz_{ni_1}/\sqrt{n}}|(Y_i)_{i\in\mathbb{I}_2}\right] - \prod_{i_1=1}^n e^{-(t^2V_{ni_1})/2n}\right|$$

$$\leq \frac{t^2}{2}\underbrace{\left|\frac{1}{n}\sum_{i_1=1}^n V_{ni_1} - V\right|}_{=:A_{n1}} + \sum_{i_1=1}^n \left|\mathbb{E}\left[e^{itz_{ni_1}/\sqrt{n}}|(Y_i)_{i\in\mathbb{I}_2}\right] - e^{-(t^2V_{ni_1})/2n}\right|$$

$$\leq \frac{t^2 A_{n1}}{2} + \underbrace{\sum_{i_1=1}^n \left|\exp\left(-t^2\frac{V_{ni_1}}{2n}\right) - \left(1 - \frac{t^2 V_{ni_1}}{2n}\right)\right|}_{=:A_{n2}} + \underbrace{\sum_{i_1=1}^n \left|\mathbb{E}\left[e^{itz_{ni_1}/\sqrt{n}}|(Y_i)_{i\in\mathbb{I}_2}\right] - \left(1 - \frac{t^2 V_{ni_1}}{2n}\right)\right|}_{=:A_{n3}}.$$

To obtain the second inequality, we use $|\exp(a) - \exp(b)| \leq |a-b|$ for all $a,b \leq 0$ and

$$\left|\prod_{i=1}^n a_i - \prod_{i=1}^n b_i\right| \leq \sum_i |a_i - b_i|$$

which holds for all positive $(a_i, b_i)_{i=1...n} \in [0,1]^{2n}$. Now, by convergence of sample means for jointly exchangeable arrays of dimensions 2 and 3,

$$\frac{1}{n}\sum_{i_1=1}^n V_{ni_1} = \frac{1}{n}\sum_{i_1=1}^n \left(\frac{2}{n-1}\sum_{i_2 \neq i_1} f(\mathbf{Y}_{\{i_1,i_2\}})\right)^2 - \left(\frac{2}{n(n-1)}\sum_{\mathbf{i}\in\mathbb{I}_{n,2}} f(\mathbf{Y}_{\mathbf{i}})\right)^2 \quad (4.15)$$

$$= \frac{1}{(n-1)}\mathbb{P}_n f^2 + \frac{1}{n(n-1)^2}\sum_{\mathbf{i}\in\mathbb{I}_{n,3}} f(\mathbf{Y}_{\{i_1,i_2\}})f(\mathbf{Y}_{\{i_1,i_3\}}) - \frac{2}{n(n-1)}\sum_{\mathbf{i}\in\mathbb{I}_{n,2}} f(\mathbf{Y}_{\mathbf{i}})$$

$$\xrightarrow{L^1} \mathbb{E}\left[f(\mathbf{Y}_{\{\mathbf{1}\}})f(\mathbf{Y}_{\{1,3\}})\right] - \mathbb{E}(2f(\mathbf{Y}_{\{\mathbf{1}\}}))^2 = V.$$

Hence, $A_{n1} \xrightarrow{L^1} 0$.

Next, let us consider $A_{n2}$. For all $x > \delta$, we have $1 - e^{-x} = cx$ for some $c \in [0,1]$, so that $|e^{-x} - (1-x)| \leq x$. For all $0 < x \leq \delta$, $e^{-x} = 1 - x + (x^2/2)e^{-x'}$ for some $x' \in [0,1]$. Thus, $|e^{-x} - (1-x)| \leq \delta x/2$. Hence, $|e^{-x} - (1-x)| \lesssim [\delta x + \mathbb{1}\{x > \delta\}x]$. Therefore,

$$\left|\exp\left(-t^2\frac{V_{ni_1}}{2n}\right) - \left(1 - \frac{t^2 V_{ni_1}}{2n}\right)\right| \lesssim \delta \frac{t^2 V_{ni_1}}{2n} + \mathbb{1}\left\{\frac{t^2 V_{ni_1}}{2n} > \delta\right\}\frac{t^2 V_{ni_1}}{2n}.$$

Hence, taking $\delta = \varepsilon/\sqrt{n}$, for any $\varepsilon > 0$, we get

$$\mathbb{E}[A_{n2}] \lesssim t^2\left\{\frac{\varepsilon}{\sqrt{n}}\mathbb{E}[V_{n1}] + \mathbb{E}\left[V_{n1}\mathbb{1}\{V_{n1} > \varepsilon\sqrt{n}\}\right]\right\}. \quad (4.16)$$

By definition of $V_{n1}$ and joint exchangeability,

$$\mathbb{E}\left[V_{n1}\mathbb{1}\left\{V_{n1} > \varepsilon\sqrt{n}\right\}\right]$$

$$\leq \mathbb{E}\left[\left(\frac{4}{(n-1)^2}\sum_{i_2=2}^{n}f^2(\boldsymbol{Y}_{\{1,i_2\}}) + \frac{4}{(n-1)^2}\sum_{2\leq i_2\neq i_3\leq n}f(\boldsymbol{Y}_{\{1,i_2\}})f(\boldsymbol{Y}_{\{1,i_3\}})\right)\mathbb{1}\left\{V_{n1}>\varepsilon\sqrt{n}\right\}\right]$$

$$\leq \frac{1}{(n-1)}\mathbb{E}\left[4f^2(\boldsymbol{Y}_{\{1\}})\mathbb{1}\left\{V_{n1}>\varepsilon\sqrt{n}\right\}\right] + 4\mathbb{E}\left[\left|f(\boldsymbol{Y}_{\{1\}})f(\boldsymbol{Y}_{\{1,3\}})\right|\mathbb{1}\left\{V_{n1}>\varepsilon\sqrt{n}\right\}\right]. \quad (4.17)$$

By (4.15), $V_{n1}$ involves two sample means. Both converge almost surely, the first by the strong law of large numbers for exchangeable sequences (see e.g. Kingman, 1978). Thus, $V_{n1}/\sqrt{n} \xrightarrow{\text{a.s.}} 0$ and by dominated convergence applied to the right-hand side of (4.17),

$$\mathbb{E}\left[V_{n1}\mathbb{1}\left\{V_{n1}>\varepsilon\sqrt{n}\right\}\right] \to 0. \quad (4.18)$$

Using this, $\mathbb{E}[V_{n1}] \to V$ and (4.16), we get $A_{n2} \xrightarrow{L^1} 0$.

Finally let us turn to $A_{n3}$. Taylor expansions ensure that $\exp(\mathrm{i}tx) = 1 + \mathrm{i}tx - t^2x^2\exp(\mathrm{i}tx^+)$ for some $x^+$ in $[0,x]$ or $[x,0]$, and $\exp(\mathrm{i}tx) = 1 + \mathrm{i}tx - t^2x^2\exp(\mathrm{i}tx) - (\mathrm{i}t^3/6)\exp(\mathrm{i}tx^*)x^3$ for some $x^*$ in $[0,x]$ or $[x,0]$. Using the first for $|x|>\delta$ and the second for $|x|\leq\delta$ yields

$$\exp(\mathrm{i}tx) = 1 + \mathrm{i}tx - \frac{1}{2}t^2x^2 + \frac{\mathbb{1}_{\{|x|>\delta\}}}{2}t^2x^2(1-\exp(\mathrm{i}tx^+)) - \frac{\mathrm{i}\mathbb{1}_{\{|x|\leq\delta\}}}{6}t^3x^3\exp(\mathrm{i}tx^*)$$

for some $x^+, x^*$ in $[0,x]$ or $[x,0]$. Thus, using $|1-\exp(\mathrm{i}tx)|\leq 2$ and $\mathbb{1}_{\{|z_{ni_1}|\leq\delta\}}z_{ni_1}^3 \leq \delta z_{ni_1}^2$, we obtain

$$\left|\mathbb{E}\left[e^{\mathrm{i}tz_{ni_1}/\sqrt{n}}|(Y_i)_{i\in\mathbb{I}_2}\right] - \left(1 - \frac{t^2 V_{ni_1}}{2n}\right)\right| = \left|\frac{t^2}{n}\mathbb{E}\left[\mathbb{1}_{\{|z_{ni_1}|>\delta\}}z_{ni_1}^2(1-\exp(\mathrm{i}tz_{ni_1}^+))|(Y_i)_{i\in\mathbb{I}_2}\right]\right.$$

$$\left. - \frac{\mathrm{i}t^3}{6n^{3/2}}E\left[\mathbb{1}_{\{|z_{ni_1}|\leq\delta\}}z_{ni_1}^3\exp(\mathrm{i}tz_{ni_1}^*))|(Y_i)_{i\in\mathbb{I}_2}\right]\right|$$

$$\lesssim \frac{t^2}{n}\mathbb{E}\left[\mathbb{1}_{\{|z_{ni_1}|>\delta\}}z_{ni_1}^2|(Y_i)_{i\in\mathbb{I}_2}\right] + \frac{|t|^3\delta}{n^{3/2}}V_{ni_1}.$$

Hence, taking $\delta = \varepsilon\sqrt{n}$, we get

$$A_{n3} \lesssim \frac{t^2}{n}\sum_{i_1=1}^{n}\mathbb{E}\left[\mathbb{1}_{\{|z_{ni_1}|>\varepsilon\sqrt{n}\}}z_{ni_1}^2|(Y_i)_{i\in\mathbb{I}_2}\right] + \frac{|t|^3\varepsilon}{n}\sum_{i_1=1}^{n}V_{ni_1}.$$

Hence,

$$\mathbb{E}[A_{n3}] \lesssim t^2\mathbb{E}\left[\mathbb{1}_{\{|z_{n1}|>\varepsilon\sqrt{n}\}}z_{n1}^2\right] + |t|^3\varepsilon\mathbb{E}[V_{n1}].$$

Because $\mathbb{E}[V_{n1}]\to V$, the second term can be made arbitrarily small by choosing $\varepsilon$ appropriately. Further, by definition of $z_{n1}$ and $V_{n1}$,

$$\mathbb{1}_{\{|z_{n1}|>\varepsilon\sqrt{n}\}} \leq \mathbb{1}_{\{V_{n1}>\varepsilon^2\sqrt{n}\}} + \mathbb{1}_{\{|\xi_1|>n^{1/4}\}}.$$

Thus, by independence between $V_{n1}$ and $\xi_1$,

$$\mathbb{E}\left[\mathbb{1}_{\{|z_{n1}|>\varepsilon\sqrt{n}\}}z_{n1}^2\right] \leq E[\xi_1^2]\mathbb{E}\left[\mathbb{1}_{\{|V_{n1}|>\varepsilon^2\sqrt{n}\}}V_{n1}\right] + E[\mathbb{1}_{\{|\xi_1|>n^{1/4}\}}\xi_1^2]\mathbb{E}[V_{n1}].$$

Thus,
$$\mathbb{E}\left[\mathbb{1}_{\{|z_{n1}|>\varepsilon\sqrt{n}\}}z_{n1}^2\right] \leq E[\xi_1^2]\mathbb{E}\left[\mathbb{1}_{\{|V_{n1}|>\varepsilon^2\sqrt{n}\}}V_{n1}\right] + E[\mathbb{1}_{\{|\xi_1|>n^{1/4}\}}\xi_1^2]\mathbb{E}\left[V_{n1}\right].$$

The first term tends to zero by (4.18). The second term tends to 0 by the dominated convergence theorem. Hence, the left-hand side converges to 0, showing that $A_{n3} \xrightarrow{L^1} 0$. Thus, (4.14) holds.

Now, the same reasoning on $A_{n1}$ as above but replacing convergences in $L^1$ by almost-sure convergences show that
$$\mathbb{E}\left[(\mathbb{G}_n^{m*}f)^2 \,|\, (Y_{\boldsymbol{i}})_{\boldsymbol{i}\in\mathbb{I}_2}\right] = \frac{1}{n}\sum_{i_1}^n V_{ni_1} \xrightarrow{\text{a.s.}} V.$$

Then, by the same argument as in the Substep 3 of the proof of pointwise convergence in Theorem 2.2, convergence in (4.14) also holds almost surely:
$$\mathbb{E}\left[e^{it\mathbb{G}_n^{m*}f}|(Y_{\boldsymbol{i}})_{\boldsymbol{i}\in\mathbb{I}_2}\right] \xrightarrow{\text{a.s.}} e^{-t^2 V/2}.$$

Hence, conditional on $(Y_{\boldsymbol{i}})_{\boldsymbol{i}\in\mathbb{I}_2}$ and almost surely, $\mathbb{G}_n^{m*}f$ converges in distribution to $Z \sim \mathcal{N}(0,V)$.

**Second step: Asymptotic equicontinuity**

We want to prove that almost surely, $\lim_{\delta\to 0}\limsup_{n\to\infty}\mathbb{E}[\sup_{f\in\mathcal{F}_\delta}|\mathbb{G}_n^*f|\,|(Y_{\boldsymbol{i}})_{\boldsymbol{i}\in\mathbb{I}_2}] = 0$. By the triangle inequality, it suffices to control separately
$$\mathbb{E}\left[\sup_{f\in\mathcal{F}_\delta}\left|\frac{1}{\sqrt{n}}\sum_{i_1=1}^n \xi_{i_1}\left(\frac{1}{n-1}\sum_{1\leq i_2\neq i_1\leq n} f(Y_{i_1,i_2}) - \mathbb{P}_n f\right)\right|\bigg|(Y_{\boldsymbol{i}})_{\boldsymbol{i}\in\mathbb{I}_2}\right]$$

and $\mathbb{E}\left[\sup_{f\in\mathcal{F}_\delta}\left|\frac{1}{\sqrt{n}}\sum_{i_1=1}^n \xi_{i_1}\left(\frac{1}{n-1}\sum_{1\leq i_2\neq i_1\leq n} f(Y_{i_2,i_1}) - \mathbb{P}_n f\right)\right|\big|(Y_{\boldsymbol{i}})_{\boldsymbol{i}\in\mathbb{I}_2}\right]$. The two terms can be controled in a similar fashion, so we only prove the result for the first term. Now,

$$\mathbb{E}\left[\sup_{f\in\mathcal{F}_\delta}\left|\frac{1}{\sqrt{n}}\sum_{i_1=1}^n \xi_{i_1}\left(\frac{1}{n-1}\sum_{1\leq i_2\neq i_1\leq n} f(Y_{i_1,i_2}) - \mathbb{P}_n f\right)\right|\bigg|(Y_{\boldsymbol{i}})_{\boldsymbol{i}\in\mathbb{I}_2}\right]$$
$$\leq \sup_{f\in\mathcal{F}_\delta}|\mathbb{P}_n f| \times \mathbb{E}\left[\left|\frac{1}{\sqrt{n}}\sum_{i_1=1}^n \xi_{i_1}\right|\right] + \mathbb{E}\left[\sup_{f\in\mathcal{F}_\delta}\left|\frac{1}{\sqrt{n}}\sum_{i_1=1}^n \xi_{i_1}\frac{1}{n-1}\sum_{1\leq i_2\neq i_1\leq n} f(Y_{i_1,i_2})\right|\bigg|(Y_{\boldsymbol{i}})_{\boldsymbol{i}\in\mathbb{I}_2}\right]$$
$$=:A_1 + A_2.$$

Because the $(\xi_i)_{i=1...n}$ are i.i.d. and standardized, $\mathbb{E}\left[\left|\frac{1}{\sqrt{n}}\sum_{i_1=1}^n \xi_{i_1}\right|\right] \leq \sqrt{\mathbb{E}[\xi_1^2]} = 1$. Then, using the triangle and Jensen inequalities,
$$A_1 \leq \sup_{f\in\mathcal{F}_\delta}|\mathbb{P}_n f| \leq \sup_{f\in\mathcal{F}_\delta}|\mathbb{P}_n f - \mathbb{P}f| + \sup_{f\in\mathcal{F}_\delta}|\mathbb{P}f| \leq 2\sup_{f\in\mathcal{F}}|\mathbb{P}_n f - \mathbb{P}f| + \sup_{f\in\mathcal{F}_\delta}\sqrt{\mathbb{P}f^2}.$$

We know that $\sup_{f\in\mathcal{F}}|\mathbb{P}_n f - \mathbb{P}f| = o_{\text{a.s.}}(1)$ and $\sup_{f\in\mathcal{F}_\delta}\sqrt{\mathbb{P}f^2} \leq \delta$ by construction. As a result, for every $\omega$ in a set of probability one and every $\delta$, there exists $n_{\omega,\delta}$ such that for every $n \geq n_{\omega,\delta}$
$$A_1(\omega) := \sup_{f\in\mathcal{F}_\delta}\left|\frac{1}{n(n-1)}\sum_{\boldsymbol{i}\in\mathbb{I}_{n,2}} f(Y_{\boldsymbol{i}}(\omega))\right| \times \mathbb{E}\left[\left|\frac{1}{\sqrt{n}}\sum_{i_1=1}^n \xi_{i_1}\right|\right] \leq 2\delta. \tag{4.19}$$

44

We now control $A_2$. Conditional on $(Y_i)_{i \in \mathbb{I}_2}$, we are dealing with the supremum of a centered empirical process over n i.n.i.d terms. We obtain

$$A_2 \leq 4\sqrt{2} \mathbb{E}\left[\int_0^{\sigma_{1,2}} \sqrt{\log 2N\left(\varepsilon, \mathcal{F}_\delta, ||\cdot||_{1,2}^m\right)} d\varepsilon \Big| (Y_i)_{i \in \mathbb{I}_2}\right],$$

where $||f||_{1,2}^{m2} := \frac{1}{n}\sum_{i_1=1}^n \xi_{i_1}^2 \left(\frac{1}{n-1}\sum_{1 \leq i_2 \neq i_1 \leq n} f(Y_{i_1,i_2})\right)^2$ and $\sigma_{1,2}^2 := \sup_{f \in \mathcal{F}_\delta} ||f||_{1,2}^{m2}$. Let $N^m := \frac{1}{n}\sum_{i_1=1}^n \xi_{i_1}^2$. By convexity, $||f||_{1,2}^{m2} \leq N^m ||f||_{2,2}^{m2}$, where $||f||_{2,2}^{m2} := \frac{1}{n(n-1)N^m}\sum_{i \in \mathbb{I}_{n,2}} \xi_{i_1}^2 f(Y_i)^2$. We deduce that $\sigma_{1,2}^2 \leq \sigma_{2,2}^2 := N^m \sup_{f \in \mathcal{F}_\delta} ||f||_{2,2}^{m2}$.

Following the proof of Theorem 2.2 and acknowledging that conditional on $(\xi_{i_1})_{i_1=1}^n$ and $(Y_i)_{i \in \mathbb{I}_2}$, $\sum_{i \in \mathbb{I}_{n,2}} \xi_{i_1}^2 \delta_{\{Y_i\}}/[n(n-1)N^m]$ is a probability measure, we obtain

$$\mathbb{E}\left[\sup_{f \in \mathcal{F}_\delta} \left|\frac{1}{\sqrt{n}}\sum_{i_1=1}^n \varepsilon_{i_1} \xi_{i_1} \frac{1}{n-1}\sum_{1 \leq i_2 \neq i_1 \leq n} f(Y_{i_1,i_2})\right| \Big| (Y_i)_{i \in \mathbb{I}_2}\right]$$
$$\lesssim \left\{\mathbb{E}[\sigma_{2,2}^2|(Y_i)_{i \in \mathbb{I}_2}]^{1/2} + 4\left(\mathbb{P}_n F^2\right)^{1/2} J_{\mathcal{F}}\left(\frac{\mathbb{E}[\sigma_{2,2}^2|(Y_i)_{i \in \mathbb{I}_2}]^{1/2}}{4(\mathbb{P}_n F^2)^{1/2}}\right)\right\}. \quad (4.20)$$

Since $\mathbb{P}_n F^2 \xrightarrow{a.s.} PF^2 > 0$, we only have to control $\mathbb{E}[\sigma_{2,2}^2|(Y_i)_{i \in \mathbb{I}_2}]$. By the triangle inequality and definition of $\mathcal{F}_\delta$,

$$\mathbb{E}[\sigma_{2,2}^2|(Y_i)_{i \in \mathbb{I}_2}] \leq \mathbb{E}\left[\sup_{f \in \mathcal{F}_\infty}\left|\frac{1}{n}\sum_{i_1=1}^n \xi_{i_1}^2 \frac{1}{n-1}\sum_{1 \leq i_2 \neq i_1 \leq n} f(Y_{i_1,i_2})^2 - \mathbb{P}_n f^2\right| \Big| (Y_i)_{i \in \mathbb{I}_2}\right]$$
$$+ \sup_{f \in \mathcal{F}_\infty} |\mathbb{P}_n f^2 - \mathbb{P} f^2| + \delta^2. \quad (4.21)$$

It is shown in the proof of Theorem 2.2 that $\sup_{f \in \mathcal{F}_\infty} |\mathbb{P}_n f^2 - \mathbb{P} f^2| = o_{a.s.}(1)$. We turn to the first term in (4.21). Conditional on $(Y_i)_{i \in \mathbb{I}_2}$, $\left((\xi_{i_1}^2 - 1)\frac{1}{n-1}\sum_{1 \leq i_2 \neq i_1 \leq n} f(Y_{i_1,i_2})^2\right)_{i_1=1}^n$ is a centered i.n.i.d. sequence. Then, by standard truncation, symmetrization arguments and Lemma 2.3.4 in Giné and Nickl (2015), we obtain, for every possibly random $\eta_1 > 0$,

$$\mathbb{E}\left[\sup_{f \in \mathcal{F}_\infty}\left|\frac{1}{n}\sum_{i_1=1}^n \xi_{i_1}^2 \frac{1}{n-1}\sum_{1 \leq i_2 \neq i_1 \leq n} f(Y_{i_1,i_2})^2 - \mathbb{P}_n f^2\right| \Big| (Y_i)_{i \in \mathbb{I}_2}\right]$$
$$\leq 2\mathbb{E}\left[M\left(\frac{2\log 2N(\eta_1, \mathcal{F}_\infty^2, ||\cdot||_{M,1}^m)}{n}\right)^{1/2} + \eta_1 \Big| (Y_i)_{i \in \mathbb{I}_2}, N^m > 0\right] \mathbb{P}(N^m > 0)$$
$$+ \frac{2}{n(n-1)}\sum_{i \in \mathbb{I}_{n,2}} F(Y_i)^2 \mathbb{E}\left[\mathbb{1}_{\{\xi_1^2 F(Y_i)^2 > M\}} \Big| Y_i\right]$$
$$=: 2(A_3 + A_4),$$

where $||f||_{M,1}^m := \frac{1}{n}\sum_{i_1=1}^n \xi_{i_1}^2 \left|\frac{1}{n-1}\sum_{1 \leq i_2 \neq i_1 \leq n} f(Y_{i_1,i_2}) \mathbb{1}_{\{F(Y_{i_1,i_2})^2 \leq M\}}\right|$. Moreover, $||f||_{M,1}^m \leq N^m ||f||_1^m$ where $||f||_1^m := \frac{1}{nN^m}\sum_{i_1=1}^n \xi_{i_1}^2 \frac{1}{n-1}\sum_{1 \leq i_2 \neq i_1 \leq n} |f(Y_{i_1,i_2})|$. Picking $\eta_1 = \eta N^m ||F^2||_1^m$ for some positive constant $\eta$, we arrive at

$$A_3 \leq M\left(\frac{2\log 2\sup_Q N(\eta ||F^2||_{Q,1}, \mathcal{F}_\infty^2, ||\cdot||_{Q,1})}{n}\right)^{1/2} + \eta \mathbb{P}_n F^2.$$

45

Lemma S12.5 enables us to write

$$\mathbb{E}\left[\sup_{f\in\mathcal{F}_\infty}\left|\frac{1}{n}\sum_{i_1=1}^{n}\xi_{i_1}^2\frac{1}{n-1}\sum_{1\leq i_2\neq i_1\leq n}f(Y_{i_1,i_2})^2-\mathbb{P}_n f^2\right|\Big|(Y_i)_{i\in\mathbb{I}_2}\right]$$
$$\leq 2\left\{M\left(\frac{2\log 2\sup_Q N(\eta||F||_{Q,2}/8,\mathcal{F},||\cdot||_{Q,2})^2}{n}\right)^{1/2}+\eta\mathbb{P}_n F^2+A_4\right\}.$$

For every $M>0$, $A_4$ converges a.s. to $\mathbb{E}\left[\xi_1^2 F(Y_1)^2\mathbb{1}_{\{\xi_1^2 F(Y_1)^2>M\}}\right]$. By combining this with (4.21) and the end of the proof of Lemma S5, we conclude that $\limsup_{n\to\infty}\mathbb{E}[\sigma_{2,2}^2|(Y_i)_{i\in\mathbb{I}_2}]\leq\delta^2$ a.s. Asymptotic equicontinuity then follows from (4.19) and (4.20).

### 4.4 Theorem 3.4

Recall that $n_1,...,n_k$ are all indexed by an index $m$, though we most often leave this dependence implicit hereafter. They also satisfy, as $m\to\infty$, $\underline{n}=\min(n_1,...,n_k)\to\infty$ and $\underline{n}/n_k\to\lambda_j$.

#### 4.4.1 Uniform law of large numbers

The triangle inequality and the symmetrization Lemma S2 for the class $\mathcal{G}=\left\{f\mathbb{1}_{\{F\leq M\}}:f\in\mathcal{F}\right\}$ and $\Phi=\mathrm{Id}$ ensure that for every $M>0$

$$\mathbb{E}\left[\sup_{f\in\mathcal{F}}|\mathbb{P}_n f-Pf|\right]\leq 2\mathbb{E}\left[F(Y_1)\mathbb{1}_{\{F(Y_1)>M\}}\right]+2\sum_{e\in\cup_{r=1}^k\mathcal{E}_r}\mathbb{E}\left[\sup_{f\in\mathcal{F}}\left|\frac{1}{\Pi_n}\sum_{1\leq i\leq n}\varepsilon_{i\odot e}f(Y_i)\mathbb{1}_{\{F(Y_i)\leq M\}}\right|\right].$$

For every $e\in\cup_{j=1}^k\mathcal{E}_j$, let

$$||f||_{e,M,1}=\frac{1}{\Pi_n}\sum_{e\leq c\leq n\odot e}\left|\sum_{1-e\leq c'\leq n\odot(1-e)}f(Y_i)\mathbb{1}_{\{F(Y_{c+c'})\leq M\}}\right|.$$

Using the same steps as in Part 1 of the proof of Theorem 2.1, we get for every $e\in\cup_{j=1}^k\mathcal{E}_j$, every $M>0$ and every possibly random $\eta_1\geq 0$,

$$\mathbb{E}\left[\sup_{f\in\mathcal{F}}\left|\frac{1}{\Pi_n}\sum_{1\leq i\leq n}\varepsilon_{i\odot e}f(Y_i)\mathbb{1}_{\{F(Y_i)\leq M\}}\right|\right]\leq\mathbb{E}\left[\sqrt{2\log 2N(\eta_1,\mathcal{F},||\cdot||_{e,M,1})}M\frac{1}{\sqrt{\prod_{j:e_j=1}n_j}}+\eta_1\right].$$

Observe that $||f||_{e,M,1}\leq ||f||_{\mathbb{Q}_n,1}:=\frac{1}{\Pi_n}\sum_{1\leq i\leq n}|f(Y_i)|$. Letting $\eta_1=\eta||F||_{\mathbb{Q}_n,1}$, we can follow the proof of Point 1 in Theorem S1 to conclude that $\mathbb{E}\left[\sup_\mathcal{F}|\mathbb{P}_n f-Pf|\right]$ tends to 0 as $m\to\infty$.

We now turn to proving almost-sure convergence. Let $\Sigma_n$ be the $\sigma$-algebra generated by $\mathcal{H}_n$ the set of functions $g$ from $\mathcal{D}^{\mathbb{N}^{+k}}$ to $\mathbb{R}$ that are invariant by the action of any $(\pi_1,...,\pi_k)$, with $\pi_r$ any permutation on $\mathbb{N}^+$ such that $\pi_r(j)=j$ if $j\geq n_r$ for $r=1,...,k$:

$$g\left((Y_i)_{i\in\mathbb{N}^{+k}}\right)=g\left(\left(Y_{\pi_1(i_1),...,\pi_k(i_k)}\right)_{i\in\mathbb{N}^{+k}}\right).$$

For every $n'\geq n$, $n'\neq n$, let $\mathbb{J}_{n,n'}=\mathbb{I}_{n'_1,n'_1-n_1}\times...\times\mathbb{I}_{n'_k,n'_k-n_k}$. Then, for every $q=(q_1,...,q_k)\in\mathbb{J}_{n,n'}$, let

$$\mathbb{P}_{n,n'}^q f=\frac{1}{\Pi_n}\sum_{1\leq i\leq n'}f(Y_i)\mathbb{1}_{\{i_1\notin\{q_1\},...,i_k\notin\{q_k\}\}}.$$

46

We observe that for every $\boldsymbol{n}, \boldsymbol{n'}, q$,

$$\mathbb{E}\left(\sup_{f\in\mathcal{F}}\left|\mathbb{P}^q_{\boldsymbol{n},\boldsymbol{n'}}f - Pf\right|\big|\Sigma_{\boldsymbol{n'}}\right) = \mathbb{E}\left(\sup_{f\in\mathcal{F}}|\mathbb{P}_{\boldsymbol{n}}f - Pf|\big|\Sigma_{\boldsymbol{n'}}\right).$$

Moreover,

$$\sum_{q\in\mathbb{J}_{\boldsymbol{n},\boldsymbol{n'}}}\mathbb{P}^q_{\boldsymbol{n},\boldsymbol{n'}}f = \frac{1}{\Pi_{\boldsymbol{n}}}\sum_{1\leq i\leq \boldsymbol{n'}}f(Y_i)\sum_{q\in\mathbb{J}_{\boldsymbol{n},\boldsymbol{n'}}}\mathbb{1}_{\{i_1\notin\{q_1\},\ldots,i_k\notin\{q_k\}\}}$$

$$=\prod_{j=1}^{k}\frac{(n'_j-1)!}{n_j!}\sum_{1\leq i\leq \boldsymbol{n'}}f(Y_i).$$

and next, $\mathbb{P}_{\boldsymbol{n'}}f = \left(\prod_{j=1}^{k}\frac{n_j!}{n'_j!}\right)\sum_{q\in\mathbb{J}_{\boldsymbol{n},\boldsymbol{n'}}}\mathbb{P}^q_{\boldsymbol{n},\boldsymbol{n'}}f = \frac{1}{|\mathbb{J}_{\boldsymbol{n},\boldsymbol{n'}}|}\sum_{q\in\mathbb{J}_{\boldsymbol{n},\boldsymbol{n'}}}\mathbb{P}^q_{\boldsymbol{n},\boldsymbol{n'}}f$. Furthermore,

$$\sup_{f\in\mathcal{F}}|\mathbb{P}_{\boldsymbol{n'}}f - Pf| = \mathbb{E}\left(\sup_{f\in\mathcal{F}}|\mathbb{P}_{\boldsymbol{n'}}f - Pf|\big|\Sigma_{\boldsymbol{n'}}\right).$$

This last equality, combined with those just above and the triangle inequality give

$$\sup_{f\in\mathcal{F}}|\mathbb{P}_{\boldsymbol{n'}}f - Pf| \leq \frac{1}{|\mathbb{J}_{\boldsymbol{n},\boldsymbol{n'}}|}\sum_{q\in\mathbb{J}_{\boldsymbol{n},\boldsymbol{n'}}}\mathbb{E}\left(\sup_{f\in\mathcal{F}}\left|\mathbb{P}^q_{\boldsymbol{n},\boldsymbol{n'}}f - Pf\right|\big|\Sigma_{\boldsymbol{n'}}\right)$$

$$=\mathbb{E}\left(\sup_{f\in\mathcal{F}}|\mathbb{P}_{\boldsymbol{n}}f - Pf|\big|\Sigma_{\boldsymbol{n'}}\right).$$

Then considering $\boldsymbol{n} = (n_1(m), \ldots, n_k(m))$ and $\boldsymbol{n'} = (n_1(m+1), \ldots, n_k(m+1))$, we deduce from the almost-sure convergence of backwards submartingales that $\sup_{f\in\mathcal{F}}|\mathbb{P}_{\boldsymbol{n'}}f - Pf|$ converges almost surely to 0 when $m$ tends to infinity.

### 4.4.2 Uniform central limit theorem

**First step: pointwise weak convergence** To prove the pointwise weak convergence, the line of reasoning is the same as what we resorted to in the first step of the proof of Theorem 2.1.2: for every $f \in \mathcal{F}$, we need to find a suitable $L_2$-approximation of $\mathbb{G}_{\boldsymbol{n}}f$, denoted $H_1 f$, i.e. as $m \to \infty$ $H_1 f$ must satisfy $\mathbb{E}\left[|\mathbb{G}_{\boldsymbol{n}}f - H_1 f|^2\right] = o(1)$ and $H_1 f \xrightarrow{d} \mathcal{N}(0, K(f,f))$. We pick $H_1 f = \sum_{e\in\mathcal{E}_1}\sum_{1\leq i\leq \boldsymbol{n}}\mathbb{E}\left[\mathbb{G}_{\boldsymbol{n}}f|U_{i\odot e}\right]$, where $(U_{i\odot e})_{1\leq i\leq \boldsymbol{n}, e\in\mathcal{E}_1}$ are i.i.d terms that appear in the AHK representation of $(Y_i)_{1\leq i\leq \boldsymbol{n}}$. Let $\boldsymbol{i}^r$ be a vector with all its entries equal to one except the $r$-th one, which is equal to $i_r$. The AHK representation ensures

$$H_1(f) = \sum_{e\in\mathcal{E}_1}\sum_{1\leq i\leq \boldsymbol{n}}\mathbb{E}\left[\mathbb{G}_{\boldsymbol{n}}f|U_{i\odot e}\right]$$

$$= \sum_{r=1}^{k}\frac{\sqrt{\boldsymbol{n}}}{n_r}\sum_{i_r=1}^{n_r}\left(\mathbb{E}\left[f(Y_{\boldsymbol{i}^r})|U_{i_r}\right] - \mathbb{E}[f(Y_{\boldsymbol{1}})]\right)$$

$$\xrightarrow{d}\mathcal{N}(0, K(f,f)).$$

The convergence in distribution comes from the standard central limit theorem applied for each $e \in \mathcal{E}_1$ separately, the mutual independence of terms across $e \in \mathcal{E}_1$ in the previous expression and the fact that $\sqrt{\boldsymbol{n}/n_r} \to \sqrt{\lambda_r}$.

To conclude that $\mathbb{G}_{\boldsymbol{n}} f \xrightarrow{d} \mathcal{N}(0, K(f,f))$ as $m \to \infty$, we rely on the weak convergence of $H_1 f$ and Section C.2.1 in Davezies et al. (2018). The main step there amounts to showing that $\lim_{m\to\infty} \mathbb{V}(H_1 f)/\mathbb{V}(\mathbb{G}_{\boldsymbol{n}} f) = 1$.

**Second step: asymptotic equicontinuity** Following the same reasoning as in the proof of Part 2 of Theorem 2.1, with the symmetrization lemma S2 instead of Lemma A.1, we have

$$\mathbb{E}\left[\sup_{f \in \mathcal{F}_\delta} |\mathbb{G}_{\boldsymbol{n}} f|\right] = \mathbb{E}\left[\sup_{f \in \mathcal{F}_\delta} |\mathbb{G}_{\boldsymbol{n}} f|\right] \lesssim \mathbb{E}\left(\int_0^{\sigma_n} \sqrt{\log 2N\left(\varepsilon, \mathcal{F}_\delta, \|\cdot\|_{\mu_{\boldsymbol{n}},2}\right)} d\varepsilon\right),$$

where $\mu_{\boldsymbol{n}} = \frac{1}{\Pi_{\boldsymbol{n}}} \sum_{\mathbf{1} \le \boldsymbol{i} \le \boldsymbol{n}} \delta_{Y_{\boldsymbol{i}}}$. $\|f\|_{\mu_{\boldsymbol{n}},2}^2$ and $\sigma_n^2$ are defined in the same way as in the proof of Part 2 of Theorem 2.1 (with $\mu_{\boldsymbol{n}}$ instead of $\mu_n$). Still following this proof, we obtain

$$\mathbb{E}\left[\sup_{f \in \mathcal{F}_\delta} |\mathbb{G}_{\boldsymbol{n}} f|\right] \lesssim \mathbb{E}\left(\sigma_n^2\right)^{1/2} + \mathbb{E}\left(F^2(Y_{\mathbf{1}})\right)^{1/2} J_{\mathcal{F}}\left(\frac{\mathbb{E}\left(\sigma_n^2\right)^{1/2}}{4\mathbb{E}\left(F^2(Y_{\mathbf{1}})\right)^{1/2}}\right).$$

Recalling that $\mathbb{E}\left(\sigma_n^2\right) \le \mathbb{E}\left[\sup_{f \in \mathcal{F}_\infty} |\mathbb{P}_{\boldsymbol{n}} f^2 - \mathbb{P} f^2|\right] + \delta^2$, we can follow the end of the asymptotic equicontinuity proof of Part 2 of Theorem 2.1 with obvious minor changes to conclude.

**Third step: total boundedness** We refer to the proof of Theorem 2.1.

### 4.4.3 Convergence of the bootstrap process

As previously, we only have to prove the pointwise convergence and the asymptotic equicontinuity.

**First step: pointwise convergence** Let $\boldsymbol{i}^* = (i_1^*, ..., i_k^*)$ denote the cell obtained by sampling $i_j^*$ with replacement in $1, ..., n_j$ for every $j = 1, ..., k$.
We have the almost-sure representation

$$\boldsymbol{i}^* = (F_{n_1}^{-1}[U^*_{(i_1,0,...,0)}], ..., F_{n_k}^{-1}[U^*_{(0,...,0,i_k)}]),$$

with $(U^*_A)_{A \in \mathbb{N}^k}$ a family of i.i.d. uniform random variables and $F_{n_j}^{-1}$ the quantile function of the discrete uniform distribution on $\{1, ..., n_j\}$. Conditional on the data $(Y_{\boldsymbol{i}})_{\boldsymbol{i} \in \mathbb{N}^{*k}}$, we can thus follow an approach similar to the one we used in the jointly exchangeable case. Let $H_1^* f = \sum_{e \in \mathcal{E}_1} \sum_{\mathbf{1} \le \boldsymbol{i} \le \boldsymbol{n}} \mathbb{E}\left[\mathbb{G}_{\boldsymbol{n}}^* f | (Y_{\boldsymbol{i}})_{\boldsymbol{i} \in \mathbb{N}^{*k}}, U^*_{\boldsymbol{i} \odot e}\right]$ and $h(\boldsymbol{i}) = f(Y_{\boldsymbol{i}})$. $H_1^* f$ can also be written

$$\sqrt{\underline{n}} \sum_{r=1}^k \left(\frac{1}{\Pi_{\boldsymbol{n}}} \sum_{\mathbf{1} \le \boldsymbol{i} \le \boldsymbol{n}} h(i_1, ..., i_{r-1}, i_r^*, i_{r+1}, ..., i_k) - \mathbb{P}_{\boldsymbol{n}} f\right).$$

We first show that $\mathbb{E}\left[(\mathbb{G}_{\boldsymbol{n}}^* f - H_1^* f)^2 | (Y_{\boldsymbol{i}})_{\boldsymbol{i} \in \mathbb{N}^{*k}}\right] = o_{\text{a.s.}}(1)$. Expanding the square in the previous

48

formula gives

$$
\mathbb{E}\left[(\mathbb{G}_{\boldsymbol{n}}^*f - H_1^*f)^2 \,|\, (Y_{\boldsymbol{i}})_{\boldsymbol{i}\in\mathbb{N}^{+k}}\right]
$$

$$
=\underline{n}\bigg\{\mathbb{E}\bigg[\bigg(\sum_{r=1}^k \frac{1}{\Pi_{\boldsymbol{n}}} \sum_{\mathbf{1}\leq \boldsymbol{i}\leq \boldsymbol{n}} h(i_1,...,i_{r-1},i_r^*,i_{r+1},...,i_k)\bigg)^2 \bigg| (Y_{\boldsymbol{i}})_{\boldsymbol{i}\in\mathbb{N}^{+k}}\bigg]
$$

$$
- 2\mathbb{E}\bigg[\bigg(\sum_{r=1}^k \frac{1}{\Pi_{\boldsymbol{n}}} \sum_{\mathbf{1}\leq \boldsymbol{i}\leq \boldsymbol{n}} h(i_1,...,i_{r-1},i_r^*,i_{r+1},...,i_k)\bigg) \mathbb{P}_{\boldsymbol{n}}^*f \,\bigg|\, (Y_{\boldsymbol{i}})_{\boldsymbol{i}\in\mathbb{N}^{+k}}\bigg]
$$

$$
+ \mathbb{E}\left[(\mathbb{P}_{\boldsymbol{n}}^*f)^2 \,|\, (Y_{\boldsymbol{i}})_{\boldsymbol{i}\in\mathbb{N}^{+k}}\right] - (k-1)^2 \,(\mathbb{P}_{\boldsymbol{n}}f)^2\bigg\}.
$$

Let $A_{\boldsymbol{n}} = \sum_{r=1}^k \frac{1}{\Pi_{\boldsymbol{n}}^2} \sum_{\substack{\mathbf{1}\leq \boldsymbol{i},\boldsymbol{i}'\leq \boldsymbol{n} \\ i_r = i'_r}} h(\boldsymbol{i})h(\boldsymbol{i}')$. We can show

$$
\mathbb{E}\left[\bigg(\sum_{r=1}^k \frac{1}{\Pi_{\boldsymbol{n}}} \sum_{\mathbf{1}\leq \boldsymbol{i}\leq \boldsymbol{n}} h(i_1,...,i_{r-1},i_r^*,i_{r+1},...,i_k)\bigg)^2 \,\bigg|\, (Y_{\boldsymbol{i}})_{\boldsymbol{i}\in\mathbb{N}^{+k}}\right]
$$

$$
= (\mathbb{P}_{\boldsymbol{n}}f)^2 \left(\sum_{r=1}^k \frac{(n_r-1)}{n_r} + k(k-1)\right) + A_{\boldsymbol{n}},
$$

$$
\mathbb{E}\left[\bigg(\sum_{r=1}^k \frac{1}{\Pi_{\boldsymbol{n}}} \sum_{\mathbf{1}\leq \boldsymbol{i}\leq \boldsymbol{n}} h(i_1,...,i_{r-1},i_r^*,i_{r+1},...,i_k)\bigg) \mathbb{P}_{\boldsymbol{n}}^*f \,\bigg|\, (Y_{\boldsymbol{i}})_{\boldsymbol{i}\in\mathbb{N}^{+k}}\right]
$$

$$
= (\mathbb{P}_{\boldsymbol{n}}f)^2 \sum_{r=1}^k \frac{(n_r-1)}{n_r} + A_{\boldsymbol{n}},
$$

and $\mathbb{E}\left[(\mathbb{P}_{\boldsymbol{n}}^*f)^2 \,|\, (Y_{\boldsymbol{i}})_{\boldsymbol{i}\in\mathbb{N}^{+k}}\right] = \frac{\prod_{j=1}^k (n_j-1)}{\Pi_{\boldsymbol{n}}} (\mathbb{P}_{\boldsymbol{n}}f)^2 + B_{\boldsymbol{n}}$, where

$$
B_{\boldsymbol{n}} = \frac{1}{\Pi_{\boldsymbol{n}}} \sum_{r=1}^k \sum_{\boldsymbol{e}\in\mathcal{E}_r} \frac{\prod_{1\leq j\leq k: e_j=0}(n_j-1)}{\prod_{1\leq j\leq k: e_j=1} n_j \left(\prod_{1\leq j\leq k: e_j=0} n_j\right)^2} \sum_{\substack{\mathbf{1}\leq \boldsymbol{i},\boldsymbol{i}'\leq \boldsymbol{n} \\ i_j = i'_j \,\forall j: e_j=1}} h(\boldsymbol{i})h(\boldsymbol{i}').
$$

For every $\boldsymbol{e} \in \cup_{r=2}^k \mathcal{E}_r$, we can write the following decomposition

$$
\sum_{\substack{\mathbf{1}\leq \boldsymbol{i},\boldsymbol{i}'\leq \boldsymbol{n} \\ i_j = i'_j \,\forall j: e_j=1}} h(\boldsymbol{i})h(\boldsymbol{i}') = \sum_{\substack{\boldsymbol{e}'\in\cup_{r=1}^k \mathcal{E}_r \\ e'_j=1 \text{ if } e_j=1}} \sum_{(\boldsymbol{i},\boldsymbol{i}')\in\mathcal{I}_{\boldsymbol{n},\boldsymbol{e}'}} h(\boldsymbol{i})h(\boldsymbol{i}'),
$$

with $\mathcal{I}_{\boldsymbol{n},\boldsymbol{e}'} = \{(\boldsymbol{i},\boldsymbol{i}') : \mathbf{1}\leq \boldsymbol{i},\boldsymbol{i}' \leq \boldsymbol{n},\, i_r = i'_r \text{ if } e'_r = 1 \text{ and } i_r \neq i'_r \text{ otherwise}\}$. Applying Lemma S8, we conclude that for every $\boldsymbol{e} \in \cup_{r=2}^k \mathcal{E}_r$,

$$
\sum_{\substack{\mathbf{1}\leq \boldsymbol{i},\boldsymbol{i}'\leq \boldsymbol{n} \\ i_j = i'_j \,\forall j: e_j=1}} h(\boldsymbol{i})h(\boldsymbol{i}') = O_{\text{a.s.}}\left(\Pi_{\boldsymbol{n}} \prod_{1\leq j\leq k: e_j=0}(n_j-1)\right),
$$

$$
B_{\boldsymbol{n}} = \sum_{r=1}^k \frac{n_r \prod_{1\leq j\leq k: j\neq r}(n_j-1)}{\Pi_{\boldsymbol{n}}} \frac{1}{\Pi_{\boldsymbol{n}}^2} \sum_{\substack{\mathbf{1}\leq \boldsymbol{i},\boldsymbol{i}'\leq \boldsymbol{n} \\ i_r = i'_r}} h(\boldsymbol{i})h(\boldsymbol{i}') + O_{\text{a.s.}}(\underline{n}^{-2}).
$$

49

By combining all those elements, we obtain

$$\mathbb{E}\left[(\mathbb{G}_n^* f - H_1^* f)^2 \,|\, (Y_{\boldsymbol{i}})_{\boldsymbol{i}\in\mathbb{N}^{+k}}\right] = \underline{n}\left\{\frac{1}{\Pi_{\boldsymbol{n}}^2}\sum_{r=1}^k \left(\frac{n_r \prod_{1\leq j\leq k: j\neq r}(n_j-1)}{\Pi_{\boldsymbol{n}}} - 1\right)\sum_{\substack{1\leq \boldsymbol{i}\leq \boldsymbol{n}\\1\leq \boldsymbol{i}'\leq \boldsymbol{n}\\i_r=i'_r}} h(\boldsymbol{i})h(\boldsymbol{i}')\right.$$

$$\left. + (\mathbb{P}_{\boldsymbol{n}}f)^2 \left(\frac{\prod_{1\leq j\leq k: j\neq r}(n_j-1)}{\Pi_{\boldsymbol{n}}} - 1 + \sum_{r=1}^k \frac{1}{n_r}\right) + O_{\text{a.s.}}\left(\underline{n}^{-2}\right)\right\}.$$

Noting that $\frac{n_r \prod_{1\leq j\leq k: j\neq r}(n_j-1)}{\Pi_{\boldsymbol{n}}} - 1 = o(1)$, $\frac{\prod_{1\leq j\leq k: j\neq r}(n_j-1)}{\Pi_{\boldsymbol{n}}} - 1 + \sum_{r=1}^k \frac{1}{n_r} = O\left(\underline{n}^{-2}\right)$ and $\frac{1}{\Pi_{\boldsymbol{n}}^2}\sum_{\substack{1\leq \boldsymbol{i}\leq \boldsymbol{n}\\1\leq \boldsymbol{i}'\leq \boldsymbol{n}\\i_r=i'_r}} h(\boldsymbol{i})h(\boldsymbol{i}') = O_{\text{a.s.}}(\underline{n}^{-1})$, again by Lemma S8, we conclude that

$$\mathbb{E}\left[(\mathbb{G}_n^* f - H_1^* f)^2 \,|\, (Y_{\boldsymbol{i}})_{\boldsymbol{i}\in\mathbb{N}^{+k}}\right] = o_{\text{a.s.}}(1).$$

To prove the asymptotic normality of $H_1 f$ conditional on $(Y_{\boldsymbol{i}})_{\boldsymbol{i}\in\mathbb{N}^{+k}}$, we remark that

$$H_1 f = \sum_{r=1}^k \sqrt{\frac{n}{n_r}} \sum_{i_r=1}^{n_r} \frac{z^*_{m,r,i_r}}{\sqrt{n_r}},$$

where $z^*_{m,r,i_r} = \frac{1}{\prod_{1\leq j\leq k: j\neq r} n_j}\sum_{i_j=1,\ldots,n_j, \forall j\neq r}(h(i_1,\ldots,i_{r-1},i_r^*,i_{r+1},\ldots,i_k) - \mathbb{P}_{\boldsymbol{n}}f).$ , For every $r = 1,\ldots,k$, $\left(z^*_{m,r,i_r}\right)_{i_r=1\ldots n_r}$ is an i.i.d. sequence of centered random variables conditional on $(Y_{\boldsymbol{i}})_{\boldsymbol{i}\in\mathbb{N}^{+k}}$ with a distribution that depends on $m$. Since

$$\mathbb{V}\left(z^*_{m,r,1} \,|\, (Y_{\boldsymbol{i}})_{\boldsymbol{i}\in\mathbb{N}^{+k}}\right) = \frac{1}{n_r \prod_{1\leq j\leq k: j\neq r} n_j^2}\sum_{\substack{1\leq \boldsymbol{i},\boldsymbol{i}'\leq \boldsymbol{n}\\i_r=i'_r}} h(\boldsymbol{i})h(\boldsymbol{i}') - (\mathbb{P}_{\boldsymbol{n}}f)^2,$$

we can conclude thanks to Point 1 of this theorem and Lemma S8 that $\mathbb{V}\left(z^*_{m,r,1} \,|\, (Y_{\boldsymbol{i}})_{\boldsymbol{i}\in\mathbb{N}^{+k}}\right) \xrightarrow{\text{a.s.}} \mathbb{E}[h(\boldsymbol{1})h(\boldsymbol{2}_r)] - \mathbb{E}[h(\boldsymbol{1})]^2 = \text{Cov}(h(\boldsymbol{1}), h(\boldsymbol{2}_r)) = V_r$. It is not difficult to see that arguments similar to those of substeps 2 and 3 of Section 3.7 apply. Then, for every $r = 1,\ldots,k$ and every $t\in\mathbb{R}$,

$$\mathbb{E}\left[\exp\left(\mathrm{i}t\sum_{i_r=1}^{n_r}\frac{z^*_{m,r,i_r}}{\sqrt{n_r}}\right)\,|\,(Y_{\boldsymbol{i}})_{\boldsymbol{i}\in\mathbb{N}^{+k}}\right] \xrightarrow{\text{a.s.}} \exp\left(-\frac{t^2 V_r}{2}\right).$$

The continuous mapping theorem, the fact that $\frac{n}{n_r}\to\lambda_r$ and the mutual independence between the $k$ sequences $\left(z^*_{m,r,i_r}\right)_{i_r=1\ldots n_r}$ ($r=1,\ldots,k$) conditional on the data imply that

$$\mathbb{E}\left[\exp\left(\mathrm{i}t H_1 f\right)\,|\,(Y_{\boldsymbol{i}})_{\boldsymbol{i}\in\mathbb{N}^{+k}}\right] = \prod_{r=1}^k \mathbb{E}\left[\exp\left(\mathrm{i}\sqrt{\frac{n}{n_r}}t\sum_{i_r=1}^{n_r}\frac{z^*_{m,r,i_r}}{\sqrt{n_r}}\right)\,|\,(Y_{\boldsymbol{i}})_{\boldsymbol{i}\in\mathbb{N}^{+k}}\right]$$

$$\xrightarrow{\text{a.s.}} \exp\left(-\frac{t^2 \sum_{r=1}^k \lambda_r V_r}{2}\right).$$

The result follows.

**Second step: asymptotic equicontinuity**   First, we have
$$(i^*)_{1\leq i\leq n} = \left(F_{n_1}^{-1}[U^*_{(i_1,0,...,0)}], ..., F_{n_k}^{-1}[U^*_{(0,...,0,i_k)}]\right)_{1\leq i\leq n}.$$
This representation ensures that the symmetrization Lemma S2 for the class $\mathcal{F}_\delta$ and $\Phi=\mathrm{Id}$ is valid. We notice that the representation is simplified as only terms associated with $e \in \mathcal{E}_1$ appear. This implies that the telescoping argument in the proof of Lemma S2 only has to be undertaken over $\mathcal{E}_1$. The following symmetrization inequality thus holds:
$$\mathbb{E}\left[\sup_{f\in\mathcal{F}_\delta}|\mathbb{G}_n^* f|\,\middle|\,(Y_{i'})_{i'\geq 1}\right] \leq 2\sum_{e\in\mathcal{E}_1}\mathbb{E}\left[\sup_{f\in\mathcal{F}_\delta}\left|\frac{1}{\sqrt{\Pi_n}}\sum_{1\leq i\leq n}\varepsilon_{i\odot e}f(Y_{i^*})\right|\,\middle|\,(Y_{i'})_{i'\geq 1}\right].$$
For every $e\in\mathcal{E}_1$, let $r_e$ be the position of the unique non-null element of $e$. This allows us to define
$$\|f\|_{e,2}^* = \frac{1}{n_{r_e}}\sum_{i_{r_e}=1}^{n_{r_e}}\left[\frac{1}{\Pi_{j\neq r_e}n_j}\sum_{(i_1,...,i_{r_e-1},i_{r_e+1},...,i_k:1\leq i\leq n)}f(Y_{i^*})\right]^2,$$
and $\sigma_{n,e}^* = \sup_{f\in\mathcal{F}_\delta}\|f\|_{e,2}^*$. Then, by Theorem 2.3.6 in Giné and Nickl (2015), we obtain
$$\mathbb{E}\left[\sup_{f\in\mathcal{F}_\delta}|\mathbb{G}_n^* f|\,\middle|\,(Y_{i'})_{i'\geq 1}\right]$$
$$\leq 8\sqrt{2}\sum_{e\in\mathcal{E}_1}\frac{1}{\sqrt{n_{r_e}}}\mathbb{E}\left[\sqrt{\log 2}\,\sigma_{n,e}^* + \int_0^{\sigma_{n,e}^*}\sqrt{\log N(\varepsilon,\mathcal{F}_\delta,\|\cdot\|_{e,2}^*)}d\varepsilon\,\middle|\,(Y_{i'})_{i'\geq 1}\right].$$
By a convexity argument, we have, for every $e\in\mathcal{E}_1$ $\|f\|_{e,2}^* \leq \|f\|_{\mathbb{Q}_n^*,2}^*$, with $\|f\|_{\mathbb{Q}_n^*,2}^{*2} = \frac{1}{\Pi_n}\sum_{1\leq i\leq n}f(Y_{i^*})^2$. We also have $\sigma_{n,e}^{*2} \leq \sigma_n^{*2}$, with $\sigma_n^{*2} = \sup_{f\in\mathcal{F}_\delta}\frac{1}{\Pi_n}\sum_{1\leq i\leq n}(f(Y_{i^*}))^2$. Then, using Points 1-4 of Lemma S12 and reasoning as in Theorem 2.2, we get
$$\mathbb{E}\left[\sup_{f\in\mathcal{F}_\delta}|\mathbb{G}_n^* f|\,\middle|\,(Y_{i'})_{i'\geq 1}\right] \leq 8\sqrt{2}k\left\{\sqrt{\log 2}\sqrt{\mathbb{E}\left[\sigma_n^{*2}|(Y_{i'})_{i'\geq 1}\right]}\right.$$
$$\left.+4\sqrt{\frac{1}{\Pi_n}\sum_{1\leq i\leq n}F^2(Y_i)}J_\mathcal{F}\left(\frac{\sqrt{\mathbb{E}\left[\sigma_n^{*2}|(Y_{i'})_{i'\geq 1}\right]}}{4\sqrt{\frac{1}{\Pi_n}\sum_{1\leq i\leq n}F^2(Y_i)}}\right)\right\}.$$
We have $\frac{1}{\Pi_n}\sum_{1\leq i\leq n}F^2(Y_i) \xrightarrow{a.s.} \mathbb{E}(F^2(Y_1)) > 0$ and
$$\sigma_n^{*2} = \sup_{f\in\mathcal{F}_\delta}|\mathbb{P}_n^* f^2| \leq \sup_{f\in\mathcal{F}_\infty}\left|\mathbb{P}_n^* f^2 - \mathbb{P}_n f^2\right| + \sup_{f\in\mathcal{F}_\infty}\left|\mathbb{P}_n f^2 - Pf^2\right| + \delta^2.$$
Moreover, we have shown in the proof of Point 2 that $\sup_{f\in\mathcal{F}_\infty}|\mathbb{P}_n f^2 - Pf^2| \xrightarrow{a.s.} 0$. Thus, it suffices to show
$$\mathbb{E}\left(\sup_{f\in\mathcal{F}_\infty}\left|\mathbb{P}_n^* f^2 - \mathbb{P}_n f^2\right|\,\middle|\,(Y_i)_{i\geq 1}\right) \xrightarrow{a.s.} 0.$$
The symmetrization argument we used to control $\mathbb{E}\left[\sup_{f\in\mathcal{F}_\delta}|\mathbb{G}_n^* f|\,\middle|\,(Y_{i'})_{i'\geq 1}\right]$ still applies and gives
$$\mathbb{E}\left[\sup_{f\in\mathcal{F}_\infty}\left|\mathbb{P}_n^* f^2 - \mathbb{P}_n f^2\right|\,\middle|\,(Y_i)_{i\geq 1}\right] \leq 4\frac{1}{\Pi_n}\sum_{1\leq i\in n}(F(Y_i))^2 \mathbb{1}_{\{(F(Y_i))^2 > M\}}$$
$$+2\sum_{e\in\mathcal{E}_1}\mathbb{E}\left[\sup_{f\in\mathcal{F}_\infty}\left|\frac{1}{\Pi_n}\sum_{1\leq i\in n}\varepsilon_{i\odot e}(f(Y_{i^*}))^2\mathbb{1}_{\{(F(Y_{i^*}))^2\leq M\}}\right|\,\middle|\,(Y_i)_{i\geq 1}\right].$$

Using the seminorm

$$||g||^*_{e,M,1} = \frac{1}{\Pi_n} \sum_{e \leq c \leq n \odot e} \left| \sum_{1-e \leq c' \leq n \odot (1-e)} g(Y_{(c+c')^*}) \right|$$

and reasoning as in Theorem 2.2, we obtain

$$\sum_{e \in \mathcal{E}_1} \mathbb{E}\left[\sup_{f \in \mathcal{F}_\infty} \left|\frac{1}{\Pi_n} \sum_{1 \leq i \leq n} \varepsilon_{i \odot e}(f(Y_{i^*}))^2 \mathbb{1}_{\{(F(Y_{i^*}))^2 \leq M\}}\right| \Big| (Y_i)_{i \geq 1}\right]$$
$$\leq 4k\sqrt{2\log 2 \sup_Q N^2(\eta ||F||_{Q,2}, \mathcal{F}, ||\cdot||_{Q,2})} M \frac{1}{\sqrt{\underline{n}}} + 8k\eta \frac{1}{\Pi_n} \sum_{1 \leq i \leq n} F^2(Y_i).$$

This is enough to conclude that $\mathbb{E}\left(\sup_{f \in \mathcal{F}_\infty} |\mathbb{P}_n^* f^2 - \mathbb{P}_n f^2| \,|\, (Y_i)_{i \geq 1}\right) \xrightarrow{\text{a.s.}} 0$. The result follows.

## 4.5 Theorem S1

### 4.5.1 Uniform law of large numbers

We remark that $\sup_{f \in \mathcal{F}} |\widetilde{\mathbb{P}}_n f - \widetilde{P}f| = \sup_{\widetilde{f} \in \widetilde{\mathcal{F}}} |\mathbb{P}_n \widetilde{f} - P\widetilde{f}|$. Following the same reasoning as in the proof of Theorem 2.1, for every positive $M$ and $\eta_1$ (with $\eta_1$ possibly random) and some constants $K_{r,k}$, there exists a jointly exchangeable and dissociated array $(\widetilde{Y}_i^r)_{i \in \mathbb{I}_k} = (N_i^r, (Y_{i,\ell}^r)_{\ell \geq 1})_{i \in \mathbb{I}_k}$ such that $\widetilde{Y}_i^r \stackrel{d}{=} (N_i, (Y_{i,\ell})_{\ell \geq 1})$ for all $i \in \mathbb{I}_{n,k}$ and

$$\mathbb{E}\left[\sup_{\widetilde{f} \in \widetilde{\mathcal{F}}} |\mathbb{P}_n \widetilde{f} - P\widetilde{f}|\right]$$
$$\leq \mathbb{E}\left[\widetilde{F}(\widetilde{Y}_1) \mathbb{1}_{\{\widetilde{F}(\widetilde{Y}_1) > M\}}\right]$$
$$+ \sum_{r=1}^{k} \sum_{e \in \mathcal{E}_r} K_{r,k} \mathbb{E}\left[\sup_{\mathcal{F}} \left|\frac{(n-k)!}{n!} \sum_{i \in \mathbb{I}_{n,k}} \varepsilon_{\{i \odot e\}^+} \widetilde{f}(\widetilde{Y}_i^r) \mathbb{1}_{\{\widetilde{F}(\widetilde{Y}_i^r) \leq M\}}\right| \Big| \overline{N_1^r} > 0\right] \mathbb{P}\left(\overline{N_1^r} > 0\right)$$
$$\leq \mathbb{E}\left[\widetilde{F}(\widetilde{Y}_1) \mathbb{1}_{\{\widetilde{F}(\widetilde{Y}_1) > M\}}\right]$$
$$+ \sum_{r=1}^{k} \sum_{e \in \mathcal{E}_r} K_{r,k} \mathbb{E}\left[\sqrt{2\log 2N(\eta_1, \widetilde{\mathcal{F}}, ||\cdot||_{e,M,1})} M \frac{\sqrt{(n-r)!r!}}{\sqrt{n!}} + \eta_1 \Big| \overline{N_1^r} > 0\right] \mathbb{P}\left(\overline{N_1^r} > 0\right),$$

where $\overline{N_p^r} = \frac{(n-k)!}{n!} \sum_{i \in \mathbb{I}_{n,k}} (N_i^r)^p$. Moreover, $||\widetilde{f}||_{e,M,1} \leq \overline{N_1^r} ||f||_{\mathbb{Q}_n^r,1}$ with

$$\mathbb{Q}_n^r = \frac{1}{\sum_{i \in \mathbb{I}_{n,k}} N_i^r} \sum_{i \in \mathbb{I}_{n,k}} \sum_{\ell=1}^{N_i^r} \delta_{Y_{i,\ell}^r}.$$

Letting $\eta_1 = \eta \overline{N_1^r} ||F||_{\mathbb{Q}_n^r,1}$ for an arbitrary $\eta > 0$, we have $N(\eta_1, \widetilde{\mathcal{F}}, ||\cdot||_{e,M,1}) \leq N(\eta_1, \mathcal{F}, \overline{N_1^r} \cdot ||\cdot||_{\mathbb{Q}_n^r,1}) = N(\overline{N_1^r}^{-1} \eta_1, \mathcal{F}, ||\cdot||_{\mathbb{Q}_n^r,1})$ whenever $\overline{N_1^r} > 0$. Combining this insight with the fact that

$\mathbb{E}\left[||\widetilde{F}||_{\mathbb{Q}_n^r,1}|\overline{N_1^r} > 0\right] = \mathbb{E}\left[\widetilde{F}(\widetilde{Y}_{\mathbf{1}})\right]/\mathbb{P}\left(\overline{N_1^r} > 0\right)$, we get

$$\mathbb{E}\left[\sup_{\widetilde{\mathcal{F}}}\left|\mathbb{P}_n\widetilde{f} - P\widetilde{f}\right|\right]$$
$$\leq \mathbb{E}\left[\widetilde{F}(\widetilde{Y}_{\mathbf{1}})\mathbb{1}_{\{\widetilde{F}(\widetilde{Y}_{\mathbf{1}})>M\}}\right] + \sum_{r=1}^{k}\sum_{e\in\mathcal{E}_r} K_{r,k}\sqrt{2\log 2\sup_Q N(\eta||F||_{Q,1}, \mathcal{F}, ||\cdot||_{Q,1})} M \frac{\sqrt{(n-r)!r!}}{\sqrt{n!}}$$
$$+ \eta \sum_{r=1}^{k}\sum_{e\in\mathcal{E}_r} K_{r,k}\mathbb{E}\left[\widetilde{F}(\widetilde{Y}_{\mathbf{1}})\right].$$

Considering $M$ sufficiently large and $\eta$ sufficiently small and next $n$ tending to $\infty$ we deduce that $\mathbb{E}\left[\sup_{\widetilde{\mathcal{F}}}|\mathbb{P}_n\widetilde{f} - P\widetilde{f}|\right]$ tends to 0 as $n \to \infty$.

Let $\Sigma_n$ be the $\sigma$-algebra generated by $\mathcal{H}_n$ the set of functions $g$ from $\mathcal{D}^{\mathbb{I}_k}$ to $\mathbb{R}$ that are invariant by the action of any permutation $\pi$ on $\mathbb{N}^+$ such that $\pi(j) = j$ for $j \geq n$:

$$g\left((\widetilde{Y}_i)_{i\in\mathbb{I}_k}\right) = g\left((\widetilde{Y}_{\pi(i)})_{i\in\mathbb{I}_k}\right).$$

Following the same reasoning as in the proof of Theorem 2.1, we conclude that $\left(\sup_{\widetilde{\mathcal{F}}}|\mathbb{P}_n\widetilde{f} - P\widetilde{f}|, \Sigma_n\right)_{n\geq 1}$ is a backwards submartingale ensuring the almost-sure convergence of $\sup_{\widetilde{\mathcal{F}}}|\mathbb{P}_n\widetilde{f} - P\widetilde{f}|$.

### 4.5.2 Uniform central limit theorem

The pointwise weak convergence is ensured by the first step of the proof of Theorem 2.1.2 applied to the class $\widetilde{\mathcal{F}}$ because for every $f \in \mathcal{F}$ we have $\mathbb{E}\left[\left(\sum_{\ell=1}^{N_1} f(Y_{\mathbf{1},\ell})\right)^2\right] < \infty$. We just have to show the asymptotic equicontinuity and total boundedness of $\widetilde{\mathcal{F}}$.

Reasoning as in the proof of Theorem 2.1, we get

$$\mathbb{E}\left[\sup_{f\in\mathcal{F}_\delta}\left|\widetilde{\mathbb{G}}_n f\right|\right] = \mathbb{E}\left[\sup_{\widetilde{f}\in\widetilde{\mathcal{F}}_\delta}\left|\mathbb{G}_n\widetilde{f}\right|\right] \lesssim \sum_{r=1}^{k}\mathbb{E}\left(\int_0^{\widetilde{\sigma}_n^r}\sqrt{\log 2N\left(\varepsilon, \widetilde{\mathcal{F}}_\delta, ||\cdot||_{\mu_n^r,2}\right)}d\varepsilon\right),$$

with $\mu_n^r = \frac{(n-k)!}{n!}\sum_{i\in\mathbb{I}_{n,k}}\delta_{\{(N_i^r, (Y_{i,\ell}^r)_{N_i^r\geq\ell\geq 1})\}}$ and $(\widetilde{\sigma}_n^r)^2 = \sup_{\widetilde{\mathcal{F}}_\delta}||\widetilde{f}||_{\mu_n^r,2}^2$. If $\overline{N_2^r} = 0$, we remark that $\int_0^{\widetilde{\sigma}_n^r}\sqrt{\log 2N\left(\varepsilon, \widetilde{\mathcal{F}}_\delta, ||\cdot||_{\mu_n^r,2}\right)}d\varepsilon = 0$. As a result, we can write

$$\mathbb{E}\left[\sup_{f\in\mathcal{F}_\delta}\left|\widetilde{\mathbb{G}}_n f\right|\right] \lesssim \sum_{r=1}^{k}\mathbb{E}\left(\int_0^{\widetilde{\sigma}_n^r}\sqrt{\log 2N\left(\varepsilon, \widetilde{\mathcal{F}}_\delta, ||\cdot||_{\mu_n^r,2}\right)}d\varepsilon\bigg|\overline{N_2^r} > 0\right)\mathbb{P}\left(\overline{N_2^r} > 0\right).$$

Reasoning conditional on $\overline{N_2^r} > 0$, we let $\mathbb{Q}_n^r = \frac{1}{\sum_{i\in\mathbb{I}_{n,k}}(N_i^r)^2}\sum_{i\in\mathbb{I}_{n,k}} N_i^r \sum_{\ell=1}^{N_i^r}\delta_{\{Y_{i,\ell}^r\}}$. For every $f \in \mathcal{F}_\delta$ and $\widetilde{f}$ the corresponding element in $\widetilde{\mathcal{F}}_\delta$, we have by the Cauchy-Schwarz inequality

$$||\widetilde{f}||_{\mu_n^r,2}^2 \leq \overline{N_2^r}||f||_{\mathbb{Q}_n^r,2}^2, \tag{4.22}$$

and next $N(\varepsilon, \widetilde{\mathcal{F}}_\delta, ||\cdot||_{\mu_n^r,2}) \leq N(\varepsilon, \mathcal{F}_\delta, \overline{N_2^r}^{1/2}||\cdot||_{\mathbb{Q}_n^r,2})$. Moreover, Points 1, 3 and 4 of Lemma S12 ensure that $N(\varepsilon, \widetilde{\mathcal{F}}_\delta, ||\cdot||_{\mu_n^r,2}) \leq N^2(\varepsilon/4\overline{N_2^r}^{1/2}, \mathcal{F}, ||\cdot||_{\mathbb{Q}_n^r,2})$. The inequality $\sqrt{a+b} \leq \sqrt{a} + \sqrt{b}$,

Lemma S11, the fact that $\mathbb{E}\left[\tilde{\sigma}_n^r | \overline{N_2^r} > 0\right] = \mathbb{E}\left[\tilde{\sigma}_n^r\right]/\mathbb{P}\left(\overline{N_2^r} > 0\right)$, $\mathbb{E}\left[N_1^r \sum_{\ell=1}^{N_1^r} F^2(Y_{1,\ell}^r) | \overline{N_2^r} > 0\right] = \mathbb{E}\left[N_1 \sum_{\ell=1}^{N_1} F^2(Y_{1,\ell})\right]/\mathbb{P}\left(\overline{N_2^r} > 0\right)$ and Jensen's inequality imply

$$\mathbb{E}\left[\sup_{f \in \mathcal{F}_\delta} |\widetilde{\mathbb{G}}_n f|\right] \lesssim \sum_{r=1}^k \mathbb{E}\left[(\tilde{\sigma}_n^r)^2\right]^{1/2} + \mathbb{E}\left[N_1 \sum_{\ell=1}^{N_1} F^2(Y_{1,\ell})\right]^{1/2} J_{\mathcal{F}}\left(\frac{\mathbb{E}\left[(\tilde{\sigma}_n^r)^2\right]^{1/2}}{4\mathbb{E}\left(N_1 \sum_{\ell=1}^{N_1} F^2(Y_{1,\ell})\right)^{1/2}}\right).$$

To prove asymptotic equicontinuity, we now follow the end of the second step of the proof of Theorem 2.1.2, starting at (3.20). We can thus claim that it is sufficient to show for every $r = 1, ..., k$

$$\lim_{n \to \infty} \mathbb{E}\left[\sup_{f \in \widetilde{\mathcal{F}}_\infty} |\mu_n^r f^2 - P f^2|\right] = 0.$$

To prove this, we replace Theorem 2.1.1 with Theorem S1.1 and adapt the "change of measure" step in the proof of the latter in the spirit of (4.22). For every positive $M$ and $\eta$, we arrive at

$$\mathbb{E}\left[\sup_{\widetilde{\mathcal{F}}_\infty} |\mu_n^r f^2 - P f^2|\right] \lesssim \mathbb{E}\left[\left(\widetilde{F}(\widetilde{Y}_1)\right)^2 \mathbb{1}_{\left\{\left(\widetilde{F}(\widetilde{Y}_1)\right)^2 > M\right\}}\right]$$
$$+ \sqrt{\frac{\log 2 \sup_Q N^2(\eta \|F\|_{Q,2}, \mathcal{F}, \|\cdot\|_{Q,2})}{n}} M + \eta \mathbb{E}\left[N_1 \sum_{\ell=1}^{N_1} F^2(Y_{1,\ell})\right].$$

Then, by choosing $M$ large enough, $\eta$ small enough and letting $n$ tend to infinity, we deduce that $\mathbb{E}\left[\sup_{\widetilde{\mathcal{F}}_\infty} |\mu_n^r f^2 - P f^2|\right] \to 0$ for every $r = 1, ..., k$.

To conclude the proof of weak convergence, we have to verify total boundedness. By the Markov inequality, we have just shown $\sup_{\widetilde{\mathcal{F}}_\infty} |\mu_n^r f^2 - P f^2| = o_p(1)$ for $r = 1, ..., k$. Fixing $r$, this entails that for every $\varepsilon > 0$ there exists $R_\varepsilon = o_p(1)$ such that for every pair $(f_1, f_2) \in \mathcal{F} \times \mathcal{F}$

$$\mathbb{E}\left[\left(\widetilde{f_1}(\widetilde{Y}_1) - \widetilde{f_2}(\widetilde{Y}_1)\right)^2\right] \leq \|\widetilde{f_1} - \widetilde{f_2}\|_{\mu_n^r, 2}^2 + R_\varepsilon.$$

For every $c > 1$, by definition of covering numbers

$$N(c\varepsilon, \widetilde{\mathcal{F}}, \|\cdot\|_{P,2}) \leq N\left(\varepsilon, \widetilde{\mathcal{F}}, \|\cdot\|_{\mu_n^r, 2}\right) + o_p(1).$$

If $\overline{N_2^r} \|F\|_{\mathbb{Q}_n^r, 2}^2 > 0$, let $U = \varepsilon/(2\overline{N_2^r}^{1/2} \|F\|_{\mathbb{Q}_n^r, 2})$. We have $\overline{N_2^r} \|F\|_{\mathbb{Q}_n^r, 2}^2 \xrightarrow{a.s.} \mathbb{E}\left(N_1 \sum_{\ell=1}^{N_1} F^2(Y_{1,\ell})\right) > 0$. Starting from the last inequality, we obtain, for every $\varepsilon > 0$,

$$N(\varepsilon, \widetilde{\mathcal{F}}, \|\cdot\|_{P,2}) \leq N\left(\frac{\varepsilon}{2}, \widetilde{\mathcal{F}}, \|\cdot\|_{\mu_n^r, 2}\right) + o_p(1)$$
$$\leq N\left(\frac{\varepsilon}{2}, \mathcal{F}, \overline{N_2^r}^{1/2} \|\cdot\|_{\mathbb{Q}_n^r, 2}\right) + o_p(1)$$
$$= N\left(U \|F\|_{\mathbb{Q}_n^r, 2}, \mathcal{F}, \|\cdot\|_{\mathbb{Q}_n^r, 2}\right) \mathbb{1}_{\{\overline{N_2^r} \|F\|_{\mathbb{Q}_n^r, 2}^2 > 0\}} + \mathbb{1}_{\{\overline{N_2^r} \|F\|_{\mathbb{Q}_n^r, 2}^2 = 0\}} + o_p(1)$$
$$\leq \sup_Q N\left(U \|F\|_{Q,2}, \mathcal{F}, \|\cdot\|_{Q,2}\right) \mathbb{1}_{\{\overline{N_2^r} \|F\|_{\mathbb{Q}_n^r, 2}^2 > 0\}} + o_p(1)$$
$$< \infty,$$

where the second inequality is a consequence of the Cauchy-Schwarz inequality and the equality on the third line is a consequence of Point 1 of Lemma S12. Hence, total boundedness holds.

54

### 4.5.3 Convergence of the bootstrap process

The triangle inequality ensures that for every $f \in \mathcal{F}$, we have $\mathbb{E}\left[\widetilde{f}(\widetilde{Y}_1)^2\right] \leq \mathbb{E}\left[\widetilde{F}(\widetilde{Y}_1)^2\right] < \infty$. The pointwise weak convergence thus follows from Theorem 2.2 applied to a finite class. The total boundedness of $(\widetilde{\mathcal{F}}, ||\cdot||_{P,2})$ has already been proved (see the proof of Theorem S1.2). As a result, there only remains to show asymptotic equicontinuity.

The proof follows closely that of Theorem 2.2. Under the moment condition on the envelope given in the statement of Theorem S1.2, it is sufficient to prove

$$\lim_{\delta \to 0} \limsup_{n \to \infty} \sqrt{n} \mathbb{E}\left[\sup_{f \in \mathcal{F}_\delta} \left|\widetilde{\mathbb{P}}_n^* f - \widetilde{\mathbb{P}}_n' f\right| \Big| (\widetilde{Y}_i)_{i \in \mathbb{I}_k}\right] \stackrel{\text{a.s.}}{=} 0,$$

where $\widetilde{\mathbb{P}}_n' f = \frac{1}{n^k} \sum_{i \in \mathbb{I}_{n,k}} \sum_{\ell=1}^{N_i} f(Y_{i,\ell})$.

Let $\overline{N_2^*} = \frac{(n-k)!}{n!} \sum_{i \in \mathbb{I}_{n,k}} N_{i^*}^2 \mathbb{1}_{\{i^* \in \mathbb{I}_{n,k}\}}$. Following the start of the asmptotic equicontinuity proof of Theorem 2.2, we have

$$\sqrt{n}\mathbb{E}\left[\sup_{f \in \mathcal{F}_\delta} \left|\widetilde{\mathbb{P}}_n^* f - \widetilde{\mathbb{P}}_n' f\right| \Big| (\widetilde{Y}_i)_{i \in \mathbb{I}_k}\right] \lesssim \mathbb{E}\left[\int_0^{\sigma_{1,2}^*} \sqrt{\log 2N(\varepsilon, \widetilde{\mathcal{F}}_\delta, ||\cdot||_{1,2}^*)} d\varepsilon \mathbb{1}_{\{\overline{N_2^*}>0\}} \Big| (\widetilde{Y}_i)_{i \in \mathbb{I}_k}\right],$$

for $||\widetilde{f}||_{1,2}^{*2} = \frac{1}{n} \sum_{i_1=1}^n \left(\frac{(n-k)!}{(n-1)!} \sum_{(i_2,\ldots,i_k): i \in \mathbb{I}_{n,k}} \widetilde{f}(\widetilde{Y}_{i^*}) \mathbb{1}_{\{i^* \in \mathbb{I}_{n,k}\}}\right)^2$ and $\sigma_{1,2}^{*2} = \sup_{\widetilde{\mathcal{F}}_\delta} ||\widetilde{f}||_{1,2}^{*2}$. The Cauchy-Schwarz inequality ensures that for every $f \in \mathcal{F}_\delta$, $||\widetilde{f}||_{1,2}^{*2} \leq \overline{N_2^*} ||f||_{\mathbb{Q}_n^*,2}^{*2}$, with

$$||f||_{\mathbb{Q}_n^*,2}^{*2} = \frac{1}{\sum_{i \in \mathbb{I}_{n,k}} N_{i^*}^2} \sum_{i \in \mathbb{I}_{n,k}} N_{i^*} \sum_{\ell=1}^{N_{i^*}} f^2(Y_{i^*,\ell}) \mathbb{1}_{\{i^* \in \mathbb{I}_{n,k}\}}.$$

It follows from Point 1 of Lemma S12 that

$$N(\varepsilon, \widetilde{\mathcal{F}}_\delta, ||\cdot||_{1,2}) \leq N(\varepsilon, \mathcal{F}_\delta, \overline{N_2^*}^{1/2} ||\cdot||_{\mathbb{Q}_n^*,2}^*) \leq N(\varepsilon \overline{N_2^*}^{-1/2}, \mathcal{F}_\delta, ||\cdot||_{\mathbb{Q}_n^*,2}^*).$$

The Cauchy-Schwarz inequality also implies

$$\sigma_{1,2}^{*2} \leq \widetilde{\sigma}_n^{*2} = \sup_{\widetilde{f} \in \widetilde{\mathcal{F}}_\delta} \frac{(n-k)!}{n!} \sum_{i \in \mathbb{I}_{n,k}} \left(\widetilde{f}(\widetilde{Y}_{i^*})\right)^2 \mathbb{1}_{\{i^* \in \mathbb{I}_{n,k}\}} = \sup_{\widetilde{f} \in \widetilde{\mathcal{F}}_\delta} |\mathbb{P}_n^* \widetilde{f}^2|.$$

Following again the proof of Theorem 2.2, we can write

$$\mathbb{E}\left[\sup_{f \in \mathcal{F}_\delta} \left|\widetilde{\mathbb{G}}_n^* f\right| \Big| (\widetilde{Y}_i)_{i \in \mathbb{I}_k}\right] \lesssim \mathbb{E}\left(\widetilde{\sigma}_n^{*2} | (\widetilde{Y}_i)_{i \in \mathbb{I}_k}\right)^{1/2}$$
$$+ \left(\frac{1}{n^k} \sum_{i \in \mathbb{I}_{n,k}} N_i \sum_{\ell=1}^{N_i} F^2(Y_{i,\ell})\right)^{1/2} J_{\mathcal{F}}\left(\frac{\mathbb{E}\left(\widetilde{\sigma}_n^{*2} | (\widetilde{Y}_i)_{i \in \mathbb{I}_k}\right)^{1/2}}{4\left(\frac{1}{n^k} \sum_{i \in \mathbb{I}_{n,k}} N_i \sum_{\ell=1}^{N_i} F^2(Y_{i,\ell})\right)^{1/2}}\right) \sqrt{A_n},$$

where $A_n = \mathbb{P}\left(\overline{N_2^*} > 0 | (\widetilde{Y}_i)_{i \in \mathbb{I}_k}\right)$.

Since $\frac{1}{n^k} \sum_{i \in \mathbb{I}_{n,k}} N_i \sum_{\ell=1}^{N_i} F^2(Y_{i,\ell}) \stackrel{\text{a.s.}}{\longrightarrow} \mathbb{E}\left(N_1 \sum_{\ell=1}^{N_1} F^2(Y_{1,\ell})\right) > 0$ and $A_n \leq 1$, we only have to show that

$$\lim_{\delta \to 0} \limsup_{n \to \infty} \mathbb{E}\left(\widetilde{\sigma}_n^{*2} | (\widetilde{Y}_i)_{i \in \mathbb{I}_k}\right) \stackrel{\text{a.s.}}{\longrightarrow} 0|.$$

We have:

$$\widetilde{\sigma}_n^{*2} = \sup_{\widetilde{f}\in\widetilde{\mathcal{F}}_\delta} |\mathbb{P}_n^* \widetilde{f}^2|$$

$$\leq \sup_{\widetilde{f}\in\widetilde{\mathcal{F}}_\delta} \left|\mathbb{P}_n^* \widetilde{f}^2 - \frac{n!}{n^k(n-k)!}\mathbb{P}_n \widetilde{f}^2\right| + \frac{n!}{n^k(n-k)!}\left(\sup_{\widetilde{f}\in\widetilde{\mathcal{F}}_\delta}\left|\mathbb{P}_n \widetilde{f}^2 - P\widetilde{f}^2\right| + \delta^2\right)$$

$$\leq \sup_{\widetilde{f}\in\widetilde{\mathcal{F}}_\infty} \left|\mathbb{P}_n^* \widetilde{f}^2 - \frac{n!}{n^k(n-k)!}\mathbb{P}_n \widetilde{f}^2\right| + \sup_{\widetilde{f}\in\widetilde{\mathcal{F}}_\infty}\left|\mathbb{P}_n \widetilde{f}^2 - P\widetilde{f}^2\right| + \delta^2$$

In the proof of Theorem S1.2, we have shown that $\sup_{\widetilde{f}\in\widetilde{\mathcal{F}}_\infty} |\mu_n^r \widetilde{f}^2 - P\widetilde{f}^2|$ converges in $L^1$ to 0 for every $r = 1, ..., k$. A similar proof can be used to claim that $\sup_{\widetilde{f}\in\widetilde{\mathcal{F}}_\infty} |\mathbb{P}_n \widetilde{f}^2 - P\widetilde{f}^2|$ converges in $L^1$ to 0. A backward submartingale argument used in the proof of Theorem 2.1.1 ensures that this convergence is almost sure. Because $\frac{n!}{n^k(n-k)!}$ tends to 1, it is sufficient to show that

$$\mathbb{E}\left(\sup_{\widetilde{f}\in\widetilde{\mathcal{F}}_\infty}\left|\mathbb{P}_n^* \widetilde{f}^2 - \frac{n!}{n^k(n-k)!}\mathbb{P}_n \widetilde{f}^2\right| \Big| (\widetilde{Y}_i)_{i\in\mathbb{I}_k}\right) \xrightarrow{a.s.} 0.$$

To do so, one simply has to mimic the proof of Lemma S5 with just one change: we need to upper bound covering numbers over the class $\widetilde{\mathcal{F}}_\infty^2$ for some random $L_1$ pseudometric using Assumption 4-(i) (which is an assumption on $\mathcal{F}$). This can be achieved thanks to Points 1, 2 and 4 of Lemma S12. The rest of the proof of Lemma S5 is left unchanged, up to notational change.

### 4.6 Proof of Proposition S1

**First part.**

Let $\widehat{S}_{ni_1}^\ell = \sum_{i_2=i_1+1}^{n-1} \sum_{i_3=i_2+1}^{n} \widehat{Z}_{i_1,i_2}^\ell \widehat{Z}_{i_1,i_3}^\ell$, $S_{ni_1}^\ell = \sum_{i_2=i_1+1}^{n-1} \sum_{i_3=i_2+1}^{n} Z_{i_1,i_2}^\ell Z_{i_1,i_3}^\ell$,

$$\widehat{N}_n^\ell = \left(\frac{6}{n(n-1)(n-2)}\right)^{1/2} \left(\sum_{i_1=1}^{n-2} \widehat{S}_{ni_1}^\ell\right), \quad N_n^\ell = \left(\frac{6}{n(n-1)(n-2)}\right)^{1/2} \left(\sum_{i_1=1}^{n-2} S_{ni_1}^\ell\right),$$

$\widehat{N}_n = (\widehat{N}_n^1, ..., \widehat{N}_n^d)$ and $N_n = (N_n^1, ..., N_n^d)$. Then $T_n = \sum_{\ell=1}^d \max(0, \widehat{N}_n^\ell)^2$.

First, we show that
$$N_n \xrightarrow{d} \mathcal{N}(0, \mathrm{Id}_d). \tag{4.23}$$

By the Cramer-Wold device, it suffices to show that for any $t \in \mathbb{R}^d$, $t'N_n \xrightarrow{d} \mathcal{N}(0, t't)$. We have

$$t'N_n = \left(\frac{6}{n(n-1)(n-2)}\right)^{1/2} \left(\sum_{i_1=1}^{n-2} t'S_{ni_1}\right),$$

with $S_{ni_1} = (S_{ni_1}^1, ..., S_{ni_1}^d)$. Moreover, using $\mathbb{V}(Z_{i_1,i_2}) = \mathrm{Id}_d$ and independence of the $(Z_{i_1,i_2})_{(i_1,i_2)\in I_2, i_1<i_2}$

under the the null hypothesis, we get

$$\mathbb{V}(t'S_{ni_1}) = \mathbb{V}\left(\sum_{\ell=1}^{d} t_\ell \sum_{i_2=i_1+1}^{n-1} \sum_{i_3=i_2+1}^{n} Z_{i_1,i_2}^{\ell} Z_{i_1,i_3}^{\ell}\right)$$

$$= \sum_{i_2=i_1+1}^{n-1} \sum_{i_3=i_2+1}^{n} \sum_{\ell=1}^{d} t_\ell \sum_{\ell'=1}^{d} t_{\ell'} \text{Cov}\left(Z_{i_1,i_2}^{\ell} Z_{i_1,i_3}^{\ell}, Z_{i_1,i_2}^{\ell'} Z_{i_1,i_3}^{\ell'}\right)$$

$$= \frac{(n-i_1)(n-i_1-1)}{2} t't.$$

Hence,

$$\frac{6}{n(n-1)(n-2)} \sum_{i_1=1}^{n-2} \mathbb{V}(t'S_{ni_1}) \to t't. \tag{4.24}$$

Moreover, under the null hypothesis, the $(S_{ni_1})_{i_1=1\ldots n}$ are independent. Then, by Lyapunov's CLT, $t'N_n \xrightarrow{d} \mathcal{N}(0, t't)$ provided that

$$\frac{\sum_{i_1=1}^{n-2} \mathbb{E}\left[(t'S_{ni_1})^4\right]}{\left(\sum_{i_1=1}^{n-2} \mathbb{V}(t'S_{ni_1})\right)^2} \to 0. \tag{4.25}$$

We have

$$\sum_{i_1=1}^{n-2} \mathbb{E}\left[(t'S_{ni_1})^4\right]$$

$$= \sum_{i_1=1}^{n-2} \mathbb{E}\left[\left(\sum_{\ell=1}^{d} t_\ell \sum_{i_2=i_1+1}^{n-1} \sum_{i_3=i_2+1}^{n} Z_{i_1,i_2}^{\ell} Z_{i_1,i_3}^{\ell}\right)^4\right]$$

$$\leq d^3 \sum_{\ell=1}^{d} |t_\ell|^4 \sum_{i_1=1}^{n-2} \mathbb{E}\left[\left(\sum_{i_2=i_1+1}^{n-1} \sum_{i_3=i_2+1}^{n} Z_{i_1,i_2}^{\ell} Z_{i_1,i_3}^{\ell}\right)^4\right]$$

$$\leq d^4 \sup_{\ell=1,\ldots d} |t_\ell|^4 \sum_{i_1=1}^{n-2} \sum_{i_1 < i_2^{(1)} < i_3^{(1)} \leq n} \sum_{i_1 < i_2^{(2)} < i_3^{(2)} \leq n} \sum_{i_1 < i_2^{(3)} < i_3^{(3)} \leq n} \sum_{i_1 < i_2^{(4)} < i_3^{(4)} \leq n} \mathbb{E}\left[\prod_{j=1}^{4} Z_{i_1,i_2^{(j)}}^{\ell} Z_{i_1,i_3^{(j)}}^{\ell}\right]$$

$$= O(n^5),$$

where the last equality holds because for any $\ell = 1, \ldots, d$, $\mathbb{E}\left[\prod_{j=1}^{4} Z_{i_1,i_2^{(j)}}^{\ell} Z_{i_1,i_3^{(j)}}^{\ell}\right] = 0$ if $\{i_2^{(j)}, i_3^{(j)}\} \not\subset \{i_2^{(j')}, i_3^{(j')} : j' \neq j\}$ for at least one $j \in \{1, \ldots, 4\}$. Combined with (4.24), this implies (4.25). Hence, (4.23) holds.

Next, we prove that

$$\widehat{N}_n^{\ell} - N_n^{\ell} = o_p(1) \quad \forall \ell \in \{1, \ldots, d\}. \tag{4.26}$$

Let $\delta_{i_1,i_2}^{\ell} = \widehat{Z}_{i_1,i_2}^{\ell} - Z_{i_1,i_2}^{\ell}$. We have

$$\widehat{N}_n^{\ell} - N_n^{\ell} = \left(\frac{6}{n(n-1)(n-2)}\right)^{1/2} \left[\sum_{i_1=1}^{n-2} \sum_{i_2=i_1+1}^{n-1} \sum_{i_3=i_2+1}^{n} \delta_{i_1,i_2}^{\ell} \left(\delta_{i_1,i_3}^{\ell} + Z_{i_1,i_3}^{\ell}\right) + Z_{i_1,i_2}^{\ell} \delta_{i_1,i_3}^{\ell}\right]. \tag{4.27}$$

Moreover, with $\widehat{s}_{ij}$ (resp. $s_{ij}$) denoting the $(i,j)$-th term of $\widehat{\Sigma}^{-1/2}$ (resp. $\mathbb{V}(Y_{1,2})^{-1/2}$),

$$\delta^\ell_{i_1,i_2} = \sum_{j=1}^d \widehat{s}_{\ell j} 2\left(\mathbb{E}(Y^j_{1,2}) - \overline{Y^j}\right) + (\widehat{s}_{\ell j} - s_{\ell j})\left(Y^j_{i_1,i_2} + Y^j_{i_2,i_1} - 2\mathbb{E}(Y^j_{1,2})\right)$$

$$:= \sum_{j=1}^d \delta^\ell_{j,i_1,i_2}.$$

In the first term of (4.27), we then have a double sum over $(j,j') \in \{1,...,d\}^2$. We prove that each of these terms is an $o_p(1)$. We have

$$\sum_{i_1<i_2<i_3} \delta^\ell_{j,i_1,i_2} \delta^\ell_{j',i_1,i_3}$$

$$= \frac{n(n-1)(n-2)}{6} R^\ell_{1nj} R^\ell_{1nj'} + R^\ell_{1nj'} R^\ell_{2nj} \sum_{i_1<i_2<i_3}\left(Y^j_{i_1,i_2} + Y^j_{i_2,i_1} - 2\mathbb{E}(Y^j_{1,2})\right)$$

$$+ R^\ell_{1nj} R^\ell_{2nj'} \sum_{i_1<i_2<i_3}\left(Y^{j'}_{i_1,i_3} + Y^{j'}_{i_3,i_1} - 2\mathbb{E}(Y^{j'}_{1,2})\right)$$

$$+ R^\ell_{2nj} R^\ell_{2nj'} \sum_{i_1<i_2<i_3}\left(Y^j_{i_1,i_2} + Y^j_{i_2,i_1} - 2\mathbb{E}(Y^j_{1,2})\right)\left(Y^{j'}_{i_1,i_3} + Y^{j'}_{i_3,i_1} - 2\mathbb{E}(Y^{j'}_{1,2})\right)$$

$$= \frac{n(n-1)(n-2)}{6} R^\ell_{1nj} R^\ell_{1nj'} + R^\ell_{1nj'} R^\ell_{2nj} \sum_{i_1<i_2}(n-i_2)\left(Y^j_{i_1,i_2} + Y^j_{i_2,i_1} - 2\mathbb{E}(Y^j_{1,2})\right)$$

$$+ R^\ell_{1nj} R^\ell_{2nj'} \sum_{\substack{1\leq i_1\leq n-2 \\ i_1+1\leq i_3\leq n}}(i_3-i_1-1)\left(Y^{j'}_{i_1,i_3} + Y^{j'}_{i_3,i_1} - 2\mathbb{E}(Y^{j'}_{1,2})\right)$$

$$+ R^\ell_{2nj} R^\ell_{2nj'} \sum_{i_1<i_2<i_3}\left(Y^j_{i_1,i_2} + Y^j_{i_2,i_1} - 2\mathbb{E}(Y^j_{1,2})\right)\left(Y^{j'}_{i_1,i_3} + Y^{j'}_{i_3,i_1} - 2\mathbb{E}(Y^{j'}_{1,2})\right), \quad (4.28)$$

with $R^\ell_{1nj} = \widehat{s}_{\ell j} 2\left(\mathbb{E}(Y^j_{1,2}) - \overline{Y^j}\right)$ and $R^\ell_{2nj} = \widehat{s}_{\ell j} - s_{\ell j}$. Under the null hypothesis, the $(Y^j_{i_1,i_2} + Y^j_{i_2,i_1})_{(i_1,i_2)\in I_2, i_1<i_2}$ are i.i.d. Then, by the usual LLN and CLT, $R^\ell_{1nj} = O_p(n^{-1})$ and the first term on the right-hand side is an $O_p(n)$. Also, $R^\ell_{2nj} = o_p(1)$. The second moments of the first and second sum in (4.28) are $O(n^4)$. Then the second and third terms of (4.28) are $o_p(n)$. Moreover, the fourth term $T_{4n}$ of (4.28) satisfies

$$T_{4n} = R^\ell_{2nj} R^\ell_{2nj'} \sum_{i_1=1}^{n-2}\left[\sum_{i_2=i_1+1}^{n-1}\sum_{i_3=i_2+1}^n \left(Y^j_{i_1,i_2} + Y^j_{i_2,i_1} - 2\mathbb{E}(Y^j_{1,2})\right)\left(Y^{j'}_{i_1,i_3} + Y^{j'}_{i_3,i_1} - 2\mathbb{E}(Y^{j'}_{1,2})\right)\right],$$

where the $n-2$ terms indexed by $i_1$ are independent with variance of order $n^2$ uniformly in $i_1 \in \{1,...,n-2\}$. Thus, $T_{4n} = R^\ell_{2nj} R^\ell_{2nj'} O_p(n^{3/2}) = o_p(n^{3/2})$ and finally,

$$\left(\frac{6}{n(n-1)(n-2)}\right)^{1/2} \sum_{i_1<i_2<i_3} \delta^\ell_{j,i_1,i_2} \delta^\ell_{j',i_1,i_3} = o_p(1).$$

Using a similar reasoning, we get

$$\left(\frac{6}{n(n-1)(n-2)}\right)^{1/2} \sum_{i_1<i_2<i_3} \delta^\ell_{j,i_1,i_2} Z^\ell_{i_1,i_3} = o_p(1), \quad \left(\frac{6}{n(n-1)(n-2)}\right)^{1/2} \sum_{i_1<i_2<i_3} Z^\ell_{i_1,i_2} \delta^\ell_{j,i_1,i_3} = o_p(1).$$

Then, in view of (4.27), (4.26) follows.

Now, by (4.23), (4.26) and the continuous mapping theorem,

$$T_n \xrightarrow{d} \sum_{\ell=1}^{d} \max(0, V_\ell)^2, \qquad (4.29)$$

where $(V_1, ..., V_d) \sim \mathcal{N}(0, \text{Id}_d)$. Let $D_\ell = \mathbb{1}_{V_\ell > 0}$ and $E_\ell = V_\ell^2$. Then $\max(0, V_\ell)^2 = D_\ell E_\ell$ and $(D_1 E_1, ..., D_d E_d)$ are mutually independent with $D_\ell \sim$ Bernoulli$(1/2)$ and $E_\ell \sim \chi^2(1)$. Moreover,

$$\sum_{\ell=1}^{d} \max(0, V_\ell)^2 | D_1, ..., D_n \sim \chi^2 \left( \sum_{\ell=1}^{d} D_\ell \right).$$

Because $\sum_{\ell=1}^{d} D_\ell \sim$ Binomial$(n, 1/2)$, the cumulative distribution function of $\sum_{\ell=1}^{d} \max(0, V_\ell)^2$ is equal to $F$, as defined in the proposition. Finally, observe that $\alpha < 1/2$, $F(0) \leq 1/2$ and $F$ is continuous and strictly increasing on $(0, \infty)$. Thus, $F$ is continuous at $q(1-\alpha)$ and $F(q(1-\alpha)) = 1-\alpha$. The result follows by definition of $W_\alpha$ and (4.29).

**Second part.**

let $\ell$ be such that $\mathbb{E}(Z^\ell_{i_1, i_2} Z^\ell_{i_1, i_3}) > 0$. We have $T_n \geq \max(0, \widehat{N}^\ell_n)^2$ so it suffices to prove that $\widehat{N}^\ell_n \xrightarrow{\mathbb{P}} \infty$. First, we show that for some $C > 0$,

$$\mathbb{P}\left( n^{-3/2} N^\ell_n > C \right) \to 1. \qquad (4.30)$$

For that purpose, we introduce the 3-dimensional array $(A_{\boldsymbol{i}})_{\boldsymbol{i} \in \mathbb{I}_3}$ by $A_{i_1, i_2, i_3} = Z^\ell_{i_1, i_2} Z^\ell_{i_1, i_3}$ if $i_1 < i_2 < i_3$, and $A_{\pi(\boldsymbol{i})} = A_{\boldsymbol{i}}$ for all $\pi \in \mathfrak{S}(\{\boldsymbol{i}\})$. Then $(A_{\boldsymbol{i}})_{\boldsymbol{i} \in \mathbb{I}_3}$ satisfies Assumption 1. Therefore, by the LLN for jointly exchangeable and dissociated arrays,

$$\frac{6}{n(n-1)(n-2)} \sum_{i_1 < i_2 < i_3} Z^\ell_{i_1, i_2} Z^\ell_{i_1, i_3} = \frac{1}{n(n-1)(n-2)} \sum_{\boldsymbol{i} \in \mathbb{I}_{n,3}} A_{\boldsymbol{i}} \xrightarrow{\mathbb{P}} \mathbb{E}(Z^\ell_{i_1, i_2} Z^\ell_{i_1, i_3}) > 0.$$

Hence, (4.30) holds. Now, we consider the remainder term $\widehat{N}^\ell_n - N^\ell_n$. We follow the same strategy as above. We first use (4.27) and then (4.28). By the LLN and CLT for jointly exchangeable, dissociated arrays, $R^\ell_{1nj} = O_p(n^{-1/2})$ and $R^\ell_{2nj} = o_p(1)$. By the same CLT, the first two sums are $O_p(n^{5/2})$, whereas by the LLN again, the third sum is an $O_p(n^3)$. So at the end,

$$\sum_{i_1 < i_2 < i_3} \delta^\ell_{j, i_1, i_2} \delta^\ell_{j', i_1, i_3} = O_p(n^2) + 2 O_p(n^{-1/2}) o_p(1) O_p(n^{5/2}) + o_p(1) O_p(n^3) = o_p(n^3).$$

We obtain a similar result for $\sum_{i_1 < i_2 < i_3} \delta^\ell_{j, i_1, i_2} Z^\ell_{i_1, i_3}$ and $\sum_{i_1 < i_2 < i_3} Z^\ell_{i_1, i_2} \delta^\ell_{j, i_1, i_3}$. So at the end, we get

$$\widehat{N}^\ell_n - N^\ell_n = o_p(n^{3/2}).$$

Combined with (4.30), this proves that $\widehat{N}^\ell_n \xrightarrow{\mathbb{P}} \infty$. The result follows.

## 5 Technical lemmas

### 5.1 Results related to the symmetrization lemma

Below, $\Phi$ denotes a non-decreasing convex function from $\mathbb{R}^+$ to $\mathbb{R}$.

**Lemma S1** (A useful inequality). *Let $m \in \mathbb{N}^+$ and $(X_1, ..., X_m)$ be any random variables with values in $\mathcal{X}$ and $\mathcal{H}$ be a pointwise measurable class of functions from $\mathcal{X}$ to $\mathbb{R}$. Then*

$$\mathbb{E}\Phi\left[\sup_{h\in\mathcal{H}}\left|\sum_{j=1}^{m}h(X_j)\right|\right] \leq \frac{1}{m}\sum_{j=1}^{m}\mathbb{E}\Phi\left[m\sup_{h\in\mathcal{H}}|h(X_j)|\right].$$

**Lemma S2** (Symmetrization, separately exchangeable, unbalanced and dissociated arrays).
*Let $k \in \mathbb{N}^+$, $\boldsymbol{n} = (n_1, ..., n_k) \in \mathbb{N}^{+k}$ and $(Y_{\boldsymbol{i}})_{\boldsymbol{1}\leq\boldsymbol{i}\leq\boldsymbol{n}}$ a family of random variables with values in a Polish space, such that*

$$(Y_{\boldsymbol{i}})_{\boldsymbol{1}\leq\boldsymbol{i}\leq\boldsymbol{n}} \stackrel{a.s.}{=} \left(\tau\left((U_{\boldsymbol{i}\odot e})_{e\in\cup_{r=1}^{k}\mathcal{E}_r}\right)\right)_{\boldsymbol{1}\leq\boldsymbol{i}\leq\boldsymbol{n}}$$

*for $(U_A)_{A\in\mathbb{N}^k}$ a family of i.i.d. real random variables and some measurable function $\tau$. Let $\mathcal{G}$ a pointwise measurable class of integrable functions of $Y_{\boldsymbol{1}}$. We have*

$$\mathbb{E}\left[\Phi\left(\sup_{g\in\mathcal{G}}\left|\frac{1}{\Pi_{\boldsymbol{n}}}\sum_{\boldsymbol{1}\leq\boldsymbol{i}\leq\boldsymbol{n}}g(Y_{\boldsymbol{i}}) - \mathbb{E}[g(Y_{\boldsymbol{1}})]\right|\right)\right]$$
$$\leq \frac{1}{2^k-1}\sum_{e\in\cup_{r=1}^{k}\mathcal{E}_r}\mathbb{E}\left[\Phi\left(2(2^k-1)\sup_{g\in\mathcal{G}}\left|\frac{1}{\Pi_{\boldsymbol{n}}}\sum_{\boldsymbol{1}\leq\boldsymbol{i}\leq\boldsymbol{n}}\varepsilon_{\boldsymbol{i}\odot e}g(Y_{\boldsymbol{i}})\right|\right)\right],$$

*where $(\varepsilon_A)_{A\in\mathbb{N}^k}$ are i.i.d. Rademacher variables, independent of $(Y_{\boldsymbol{i}})_{\boldsymbol{i}\in\mathbb{N}^{+k}}$.*

**Lemma S3** (Symmetrization in degenerate cases, jointly exchangeable, balanced and dissociated arrays).
*Suppose that $k=2$, Assumptions 1-2 and 4-(i) hold and $\mathbb{G}f = 0$ for all $f \in \mathcal{F}$. Then there exists $(Y_{\boldsymbol{i}}^1)_{\boldsymbol{i}\in\mathbb{I}_2}$, a jointly exchangeable and dissociated array with $Y_{\boldsymbol{1}}^1 \stackrel{d}{=} Y_{\boldsymbol{1}}$, satisfying*

$$\mathbb{E}\left[\sup_{f\in\mathcal{F}}\left|\frac{1}{n(n-1)}\sum_{\boldsymbol{i}\in\mathbb{I}_{n,2}}f(Y_{\boldsymbol{i}}) - \mathbb{E}[f(Y_{\boldsymbol{1}})]\right|\right]$$
$$\lesssim \mathbb{E}\left[\sup_{f\in\mathcal{F}}\left|\frac{1}{n(n-1)}\sum_{\boldsymbol{i}\in\mathbb{I}_{n,2}}\varepsilon_{\{\boldsymbol{i}\}}f(Y_{\boldsymbol{i}})\right|\right] + \mathbb{E}\left[\sup_{f\in\mathcal{F}}\left|\frac{1}{n(n-1)}\sum_{\boldsymbol{i}\in\mathbb{I}_{n,2}}\varepsilon_{i_1}^1\varepsilon_{i_2}^2 f(Y_{\boldsymbol{i}}^1)\right|\right],$$

*with $(\varepsilon_A)_{A\in\mathbb{N}^2}$, $(\varepsilon_i^j)_{i\geq 1, j\in\{1,2\}}$ mutually independent arrays of i.i.d. Rademacher variables, also independent of $(Y_{\boldsymbol{i}})_{\boldsymbol{i}\in\mathbb{I}_2}$ and $(Y_{\boldsymbol{i}}^1)_{\boldsymbol{i}\in\mathbb{I}_2}$.*

**Lemma S4** (Partial extension of Prop. 2.2 in Giné and Zinn, 1990).
*Let $n \geq 2$, $(\varepsilon_{\{i,j\}})_{1\leq i<j\leq n}$ be Rademacher independent variables and $1^*, ..., n^*$ be i.i.d. variables, uniformly distributed on $\{1, ..., n\}$ and independent of $(\varepsilon_{\{i,j\}})_{1\leq i<j\leq n}$. Let $(x_{\{i,j\}})_{1\leq i<j\leq n}$ be a non-random array of size $n \times n$ with components in a Banach space of norm $||\cdot||_B$ and such that $x_{\{i,i\}} = x_{\{i\}} = 0$. Then*

$$\mathbb{E}\left(\left\|\sum_{1\leq i,j\leq n}\varepsilon_{\{i,j\}}x_{\{i^*,j^*\}}\right\|_B\right) \geq \frac{(1-e^{-1})(1-e^{-1/2})}{\sqrt{2}}\mathbb{E}\left(\left\|\sum_{1\leq i,j\leq n}\varepsilon_{\{i,j\}}x_{\{i,j\}}\right\|_B\right).$$

### 5.1.1 Proof of Lemma S1

By the triangle inequality and properties of the supremum,

$$\sup_{h \in \mathcal{H}} \left| \sum_{j=1}^{m} h(X_j) \right| \leq \frac{1}{m} \sum_{j=1}^{m} m \sup_{h \in \mathcal{H}} |h(X_j)|.$$

The result follows by monotonicity and convexity of $\Phi$.

### 5.1.2 Proof of Lemma S2

The proof is much simpler than that of Lemma A.1 because there is much more invariance in separately exchangeable arrays than in jointly exchangeable ones. Consequently the decoupling and recoupling steps used in the proof of Lemma A.1 are not necessary.

To get the result, we introduce $\left(U_A^{(1)}\right)_{A \in \mathbb{N}^k}$ which is an independent copy of $(U_A)_{A \in \mathbb{N}^k}$. We assume without loss of generality that the last argument of $\tau$ is $U_{i \odot \mathbf{1}} = U_i$. On the set $\cup_{l=1}^{k} \mathcal{E}_l$, $\prec$ is the strict total order used (implicitly) to enumerate the arguments of $\tau$ in the statement of the Lemma. We extend this order to $\cup_{l=0}^{k} \mathcal{E}_l$ considering that $\mathbf{0} \prec \boldsymbol{e} \preceq \mathbf{1}$ for every $\boldsymbol{e} \in \cup_{l=1}^{k} \mathcal{E}_l$. For every $(\boldsymbol{e}, \boldsymbol{e}') \in \left(\cup_{l=0}^{k} \mathcal{E}_l\right)^2$, we write $\boldsymbol{e} \preceq \boldsymbol{e}'$ if $\boldsymbol{e} \prec \boldsymbol{e}'$ or $\boldsymbol{e} = \boldsymbol{e}'$. We also let $Y_i^{(\boldsymbol{e})} = \tau\left(\left(U_{i \odot \boldsymbol{e}}^{(1)}\right)_{\mathbf{0} \prec \boldsymbol{e}' \preceq \boldsymbol{e}}, (U_{i \odot \boldsymbol{e}'})_{\boldsymbol{e} \prec \boldsymbol{e}' \preceq \mathbf{1}}\right)$ for every $\boldsymbol{e} \in \cup_{l=1}^{k} \mathcal{E}_l$ (hence $Y_i = Y_i^{(\mathbf{0})}$). Convexity of $\Phi$ then implies

$$\mathbb{E}\left[\Phi\left(\sup_{g \in \mathcal{G}} \left|\frac{1}{\Pi_n} \sum_{1 \leq i \leq n} g(Y_i) - \mathbb{E}[g(Y_\mathbf{1})]\right|\right)\right]$$

$$\leq \mathbb{E}\left[\Phi\left(\sup_{g \in \mathcal{G}} \left|\frac{1}{\Pi_n} \sum_{1 \leq i \leq n} g\left(Y_i^{(\mathbf{0})}\right) - g\left(Y_i^{(\mathbf{1})}\right)\right|\right)\right]$$

$$= \mathbb{E}\left[\Phi\left(\sup_{g \in \mathcal{G}} \left|\frac{1}{\Pi_n} \sum_{1 \leq i \leq n} \sum_{\mathbf{0} \prec \boldsymbol{e} \preceq \mathbf{1}} g\left(Y_i^{(\boldsymbol{e}_{prec})}\right) - g\left(Y_i^{(\boldsymbol{e})}\right)\right|\right)\right]$$

$$\leq \frac{1}{2^k - 1} \sum_{\mathbf{0} \prec \boldsymbol{e} \preceq \mathbf{1}} \mathbb{E}\left[\Phi\left((2^k - 1) \sup_{g \in \mathcal{G}} \left|\frac{1}{\Pi_n} \sum_{1 \leq i \leq n} g\left(Y_i^{(\boldsymbol{e}_{prec})}\right) - g\left(Y_i^{(\boldsymbol{e})}\right)\right|\right)\right]$$

$$= \frac{1}{2^k - 1} \sum_{\boldsymbol{e} \in \cup_{l=1}^{k} \mathcal{E}_l} \mathbb{E}\left[\Phi\left((2^k - 1) \sup_{g \in \mathcal{G}} \left|\frac{1}{\Pi_n} \sum_{1 \leq i \leq n} g\left(Y_i^{(\boldsymbol{e}_{prec})}\right) - g\left(Y_i^{(\boldsymbol{e})}\right)\right|\right)\right],$$

with $\boldsymbol{e}_{prec}$ the element that precedes $\boldsymbol{e}$ for the strict total order $\prec$. For every $\boldsymbol{e} \in \cup_{l=1}^{k} \mathcal{E}_l$, note that

$$\sum_{1 \leq i \leq n} g\left(Y_i^{(\boldsymbol{e}_{prec})}\right) - g\left(Y_i^{(\boldsymbol{e})}\right)$$

$$= \sum_{\boldsymbol{e} \leq \boldsymbol{c} \leq n \odot \boldsymbol{e}} \sum_{\mathbf{1} - \boldsymbol{e} \leq \boldsymbol{c}' \leq n \odot (\mathbf{1} - \boldsymbol{e})} g\left(Y_{\boldsymbol{c}+\boldsymbol{c}'}^{(\boldsymbol{e}_{prec})}\right) - g\left(Y_{\boldsymbol{c}+\boldsymbol{c}'}^{(\boldsymbol{e})}\right).$$

Furthermore,

$$\left(\sum_{\mathbf{1} - \boldsymbol{e} \leq \boldsymbol{c}' \leq n \odot (\mathbf{1} - \boldsymbol{e})} g\left(Y_{\boldsymbol{c}+\boldsymbol{c}'}^{(\boldsymbol{e}_{prec})}\right) - g\left(Y_{\boldsymbol{c}+\boldsymbol{c}'}^{(\boldsymbol{e})}\right)\right)_{\boldsymbol{e} \leq \boldsymbol{c} \leq n \odot \boldsymbol{e}}$$

61

is an array of independent and symmetric random variables conditional on $\left(\left(U^{(1)}_{i \odot e}\right)_{\mathbf{0} \prec e' \prec e}, (U_{i \odot e})_{e \prec e' \preceq \mathbf{1}}\right)$. Standard symmetrization arguments (see for instance van der Vaart and Wellner, 1996, Lemma 2.3.1 in the i.i.d. case) entail

$$\mathbb{E}\left[\Phi\left((2^k-1)\sup_{g\in\mathcal{G}}\left|\frac{1}{\Pi_n}\sum_{e\leq c\leq n\odot e}\sum_{\mathbf{1}-e\leq c'\leq n\odot(\mathbf{1}-e)}g\left(Y^{(e_{prec})}_{c+c'}\right)-g\left(Y^{(e)}_{c+c'}\right)\right|\right)\right]$$

$$\leq \mathbb{E}\left[\Phi\left(2(2^k-1)\sup_{g\in\mathcal{G}}\left|\frac{1}{\Pi_n}\sum_{e\leq c\leq n\odot e}\varepsilon_c \sum_{\mathbf{1}-e\leq c'\leq n\odot(\mathbf{1}-e)}g\left(Y_{c+c'}\right)\right|\right)\right]$$

$$= \mathbb{E}\left[\Phi\left(2(2^k-1)\sup_{g\in\mathcal{G}}\left|\frac{1}{\Pi_n}\sum_{\mathbf{1}\leq i\leq n}\varepsilon_{i\odot e}g\left(Y_i\right)\right|\right)\right].$$

### 5.1.3 Proof of Lemma S3

Recall that under Assumption 1, $(Y_i)_{i\in\mathbb{I}_2} = (\tau(U_{i_1},U_{i_2},U_{\{i\}}))_{i\in\mathbb{I}_2}$ for some $\tau$ and i.i.d random variables $(U_A)_{A\subset\mathbb{N}^+,1\leq|A|\leq 2}$. Let $(V_i)_{i\geq 1}$, $(U^j_i)_{i\geq 1,j\in\{1,2\}}$ and $(V^j_i)_{i\geq 1,j\in\{1,2\}}$ be independent copies of $(U_i)_{i\geq 1}$, also independent from $(U_{\{i\}})_{i\in\mathbb{I}_2}$. Let also $Y^1_i = \tau(V_{i_1},V_{i_2},U_{\{i\}})$. By the triangle inequality used twice, Jensen's inequality and a standard symmetrization argument for sums of independent variables,

$$\mathbb{E}\left[\sup_{f\in\mathcal{F}}\left|\frac{1}{n(n-1)}\sum_{i\in\mathbb{I}_{n,2}}f(Y_i)-\mathbb{E}[f(Y_\mathbf{1})]\right|\right]$$

$$\leq \mathbb{E}\left[\sup_{f\in\mathcal{F}}\left|\frac{1}{n(n-1)}\sum_{i\in\mathbb{I}_{n,2}}\left(f(\tau(V_{i_1},V_{i_2},V_{\{i\}}))-f(\tau(V_{i_1},V_{i_2},U_{\{i\}}))\right)\right|\right]$$

$$+ \mathbb{E}\left[\sup_{f\in\mathcal{F}}\left|\frac{1}{n(n-1)}\sum_{i\in\mathbb{I}_{n,2}}\left(f(\tau(U_{i_1},U_{i_2},U_{\{i\}}))-f(\tau(V_{i_1},V_{i_2},U_{\{i\}}))\right)\right|\right]$$

$$\leq 2\mathbb{E}\left[\sup_{f\in\mathcal{F}}\left|\frac{1}{n(n-1)}\sum_{i\in\mathbb{I}_{n,2}}\epsilon_{\{i\}}f(Y_i)\right|\right]+\mathbb{E}\left[\sup_{f\in\mathcal{F}^s}\left|\frac{1}{n(n-1)}\sum_{i\in\mathbb{I}_{n,2}}\left(f(\boldsymbol{Y}_{\{i\}})-f(\boldsymbol{Y}^1_{\{i\}})\right)\right|\right]. \quad (5.1)$$

We now bound the second term on the right-hand side. We apply Theorem 3.5.3 and Remark 3.5.4 in de la Peña and Giné (1999), with $r=1$ and conditionally on $(U_{\{i\}})_{i\in\mathbb{I}_2}$. After re-integrating, this yields

$$\mathbb{E}\left[\sup_{f\in\mathcal{F}^s}\left|\frac{1}{n(n-1)}\sum_{i\in\mathbb{I}_{n,2}}\left(f(\boldsymbol{Y}_{\{i\}})-f(\boldsymbol{Y}^1_{\{i\}})\right)\right|\right]$$

$$\lesssim \mathbb{E}\left[\sup_{f\in\mathcal{F}^s}\left|\frac{1}{n(n-1)}\sum_{i\in\mathbb{I}_{n,2}}\varepsilon^1_{i_1}\left(f(\boldsymbol{Y}^2_{\{i\}})-f(\boldsymbol{Y}^3_{\{i\}})\right)\right|\right]$$

$$\lesssim \mathbb{E}\left[\sup_{f\in\mathcal{F}^s}\left|\frac{1}{n(n-1)}\sum_{i\in\mathbb{I}_{n,2}}\varepsilon^1_{i_1}\left(f(\boldsymbol{Y}^2_{\{i\}})-\mathbb{E}[f(\boldsymbol{Y}_\mathbf{1})]\right)\right|\right], \quad (5.2)$$

where $(Y_{\boldsymbol{i}}^2, Y_{\boldsymbol{i}}^3)_{\boldsymbol{i}\in\mathbb{I}_2} = \left(\tau(U_{i_1}^1, U_{i_2}^2, U_{\{\boldsymbol{i}\}}), \tau(V_{i_1}^1, V_{i_2}^2, U_{\{\boldsymbol{i}\}})\right)_{\boldsymbol{i}\in\mathbb{I}_2}$. The triangle inequality was used on the last line. Next, let $Y_{\boldsymbol{i}}^4 = \tau(U_{i_1}^1, V_{i_2}^2, V_{\{\boldsymbol{i}\}})$. Since $\mathbb{G}f = 0$ for all $f \in \mathcal{F}$, we have, for all $f \in \mathcal{F}^s$,

$$\mathbb{E}\left[f\left(\boldsymbol{Y}_{\{\boldsymbol{i}\}}^4\right)\Big|U_{i_1}^1\right] = \mathbb{E}\left[f\left(\tau(U_{i_1}^1, U_{i_2}^1, U_{\{\boldsymbol{i}\}})\right)\Big|U_{i_1}^1\right] = \mathbb{E}\left[f\left(\boldsymbol{Y_1}\right)\right].$$

As a result, by Jensen's inequality

$$\mathbb{E}\left[\sup_{f\in\mathcal{F}^s}\left|\frac{1}{n(n-1)}\sum_{\boldsymbol{i}\in\mathbb{I}_{n,2}}\varepsilon_{i_1}^1\left(f(\boldsymbol{Y}_{\{\boldsymbol{i}\}}^2) - \mathbb{E}[f(\boldsymbol{Y_1})]\right)\right|\right]$$

$$=\mathbb{E}\left[\sup_{f\in\mathcal{F}^s}\left|\frac{1}{n(n-1)}\sum_{\boldsymbol{i}\in\mathbb{I}_{n,2}}\mathbb{E}\left[\varepsilon_{i_1}^1\left(f(\boldsymbol{Y}_{\{\boldsymbol{i}\}}^2) - f(\boldsymbol{Y}_{\boldsymbol{i}}^4)\right)\Big|(\varepsilon_{i_1})_{i_1}, (U_{i_1}^1)_{i_1}, (U_{i_2}^2)_{i_2}, (U_{\{\boldsymbol{i}\}})_{\boldsymbol{i}\in\mathbb{I}_{n_2}}\right]\right|\right]$$

$$\leq\mathbb{E}\left[\sup_{f\in\mathcal{F}^s}\left|\frac{1}{n(n-1)}\sum_{\boldsymbol{i}\in\mathbb{I}_{n,2}}\varepsilon_{i_1}^1\left(f(\boldsymbol{Y}_{\{\boldsymbol{i}\}}^2) - f(\boldsymbol{Y}_{\{\boldsymbol{i}\}}^4)\right)\right|\right]$$

$$=2\mathbb{E}\left[\sup_{f\in\mathcal{F}^s}\left|\frac{1}{n(n-1)}\sum_{\boldsymbol{i}\in\overrightarrow{\mathbb{I}_{n,2}}}\varepsilon_{i_1}^1\left(f(\boldsymbol{Y}_{\{\boldsymbol{i}\}}^2) - f(\boldsymbol{Y}_{\{\boldsymbol{i}\}}^4)\right)\right|\right]. \tag{5.3}$$

Conditional on $\left(U_{i_1}^1, \varepsilon_{i_1}^1\right)_{i_1=1}^n$, the variables $\sum_{i_2+1\leq i_1\leq n}\varepsilon_{i_1}^1\left(f(\boldsymbol{Y}_{\{\boldsymbol{i}\}}^2) - f(\boldsymbol{Y}_{\{\boldsymbol{i}\}}^4)\right)/(n(n-1))$ indexed by $i_2$ are mutually independent and centered, for every $f \in \mathcal{F}^s$. Then, by the symmetrization lemma for independent variables and the triangle inequality, we have

$$\mathbb{E}\left[\sup_{f\in\mathcal{F}^s}\left|\frac{1}{n(n-1)}\sum_{\boldsymbol{i}\in\overrightarrow{\mathbb{I}_{n,2}}}\varepsilon_{i_1}^1\left(f(\boldsymbol{Y}_{\{\boldsymbol{i}\}}^2) - f(\boldsymbol{Y}_{\{\boldsymbol{i}\}}^4)\right)\right|\right] \leq \mathbb{E}\left[\sup_{f\in\mathcal{F}^s}\left|\frac{1}{n(n-1)}\sum_{\boldsymbol{i}\in\overrightarrow{\mathbb{I}_{n,2}}}\varepsilon_{i_1}^1\varepsilon_{i_2}^2\left(f(\boldsymbol{Y}_{\{\boldsymbol{i}\}}^2) - f(\boldsymbol{Y}_{\{\boldsymbol{i}\}}^4)\right)\right|\right]$$

$$\leq 2\mathbb{E}\left[\sup_{f\in\mathcal{F}^s}\left|\frac{1}{n(n-1)}\sum_{\boldsymbol{i}\in\overrightarrow{\mathbb{I}_{n,2}}}\varepsilon_{i_1}^1\varepsilon_{i_2}^2 f(\boldsymbol{Y}_{\{\boldsymbol{i}\}}^2)\right|\right]$$

$$=\mathbb{E}\left[\sup_{f\in\mathcal{F}}\left|\frac{1}{n(n-1)}\sum_{\boldsymbol{i}\in\mathbb{I}_{n,2}}\varepsilon_{i_1}^1\varepsilon_{i_2}^2 f(Y_{\boldsymbol{i}}^2)\right|\right]. \tag{5.4}$$

The result follows by combining (5.1)-(5.4).

### 5.1.4 Proof of Lemma S4

Let $(\varepsilon'_{\{i,j\},\{i',j'\}})_{1\leq i,i',j,j'\leq n}$ be independent Rademacher variables, independent of the variables defined in the lemma. Conditionally on $(1^*, ..., n^*)$, and next unconditionally, we have:

$$\left(\varepsilon_{\{i,j\}}x_{\{i^*,j^*\}}\right)_{i,j=1,...,n} \stackrel{d}{=} \left(\varepsilon'_{\{i,j\},\{i^*,j^*\}}x_{\{i^*,j^*\}}\right)_{i,j=1,...,n}.$$

Let $(e_{\{i,j\}})_{1\leq i,j\leq n} \in \{-1,1\}^{n(n+1)/2}$. Conditionally on $1^*, ..., n^*$ and next unconditionally,

$$\left(e_{\{i,j\}}\sum_{1\leq i',j'\leq n}\varepsilon'_{\{i,j\},\{i',j'\}}\mathbb{1}_{\{\{i^*,j^*\}=\{i',j'\}\}}\right)_{1\leq i,j\leq n} \stackrel{d}{=} \left(\sum_{1\leq i',j'\leq n}\varepsilon'_{\{i,j\},\{i',j'\}}\mathbb{1}_{\{\{i^*,j^*\}=\{i',j'\}\}}\right)_{1\leq i,j\leq n}.$$

Let $(\varepsilon''_{\{i,j\}})_{1\leq i,j\leq n}$ denote independent Rademacher variables, independent of all other variables. Using the previous equalities in distribution and Jensen's inequality, we obtain:

$$\mathbb{E}\left(\left\|\sum_{1\leq i,j\leq n}\varepsilon_{\{i,j\}}x_{\{i^*,j^*\}}\right\|_B\right)$$

$$=\mathbb{E}\left(\left\|\sum_{1\leq i,j\leq n}\varepsilon'_{\{i,j\},\{i^*,j^*\}}x_{\{i^*,j^*\}}\right\|_B\right)$$

$$=\mathbb{E}\left(\left\|\sum_{1\leq i,j\leq n}\sum_{1\leq i',j'\leq n}\varepsilon'_{\{i,j\},\{i',j'\}}x_{\{i',j'\}}\mathbb{1}_{\{\{i^*,j^*\}=\{i',j'\}\}}\right\|_B\right)$$

$$=\mathbb{E}\left(\left\|\sum_{1\leq i',j'\leq n}x_{\{i',j'\}}\sum_{1\leq i,j\leq n}\varepsilon'_{\{i,j\},\{i',j'\}}\mathbb{1}_{\{\{i^*,j^*\}=\{i',j'\}\}}\mathbb{1}_{\{i'\neq j'\}}\right\|_B\right)$$

$$=\mathbb{E}\left(\left\|\sum_{1\leq i',j'\leq n}\varepsilon''_{\{i',j'\}}x_{\{i',j'\}}\sum_{1\leq i,j\leq n}\varepsilon'_{\{i,j\},\{i',j'\}}\mathbb{1}_{\{\{i^*,j^*\}=\{i',j'\}\}}\mathbb{1}_{\{i'\neq j'\}}\right\|_B\right)$$

$$=\mathbb{E}\left(\left\|\sum_{1\leq i',j'\leq n}\varepsilon''_{\{i',j'\}}x_{\{i',j'\}}\left|\sum_{1\leq i,j\leq n}\varepsilon'_{\{i,j\},\{i',j'\}}\mathbb{1}_{\{\{i^*,j^*\}=\{i',j'\}\}}\mathbb{1}_{\{i'\neq j'\}}\right|\right\|_B\right)$$

$$\mathbb{E}\left(\left\|\sum_{1\leq i,j\leq n}\varepsilon_{\{i,j\}}x_{\{i^*,j^*\}}\right\|_B\right)$$

$$\geq\mathbb{E}\left(\left\|\sum_{1\leq i',j'\leq n}\varepsilon''_{\{i',j'\}}x_{\{i',j'\}}\mathbb{E}\left(\left|\sum_{1\leq i,j\leq n}\varepsilon'_{\{i,j\},\{i',j'\}}\mathbb{1}_{\{\{i^*,j^*\}=\{i',j'\}\}}\right|\,\Big|\,\varepsilon''\right)\right\|_B\right)$$

$$=\mathbb{E}\left(\left\|\sum_{1\leq i',j'\leq n}\varepsilon''_{\{i',j'\}}x_{\{i',j'\}}\mathbb{E}\left(\left|\sum_{1\leq i,j\leq n}\varepsilon'_{\{i,j\},\{i',j'\}}\mathbb{1}_{\{\{i^*,j^*\}=\{i',j'\}\}}\mathbb{1}_{\{i'\neq j'\}}\right|\right)\right\|_B\right)$$

$$=\mathbb{E}\left(\left\|\sum_{1\leq i',j'\leq n}\varepsilon''_{\{i',j'\}}x_{\{i',j'\}}\mathbb{E}\left(\left|\sum_{1\leq i,j\leq n}\varepsilon'_{\{i,j\},\{1,2\}}\mathbb{1}_{\{\{i^*,j^*\}=\{1,2\}\}}\right|\right)\right\|_B\right)$$

$$=\mathbb{E}\left(\left\|\sum_{1\leq i',j'\leq n}\varepsilon''_{\{i',j'\}}x_{\{i',j'\}}\right\|_B\right)\mathbb{E}\left(\left|\sum_{1\leq i,j\leq n}\varepsilon'_{\{i,j\},\{1,2\}}\mathbb{1}_{\{\{i^*,j^*\}=\{1,2\}\}}\right|\right). \quad (5.5)$$

Next, the Khintchine and Markov inequalities yield:

$$\mathbb{E}\left(\left|\sum_{1\leq i,j\leq n}\varepsilon'_{\{i,j\},\{1,2\}}\mathbb{1}_{\{\{i^*,j^*\}=\{1,2\}\}}\right|\right)\geq\frac{1}{\sqrt{2}}\mathbb{E}\left[\left(\sum_{1\leq i,j\leq n}\mathbb{1}_{\{\{i^*,j^*\}=\{1,2\}\}}\right)^{1/2}\right]$$

$$\geq\frac{1}{\sqrt{2}}\mathbb{P}\left(\sum_{1\leq i<j\leq n}\mathbb{1}_{\{\{i^*,j^*\}=\{1,2\}\}}\geq 1\right)$$

$$=\frac{1}{\sqrt{2}}\mathbb{P}\left(\bigcup_{1\leq i<j\leq n}\{\{i^*,j^*\}=\{1,2\}\}\right)$$

$$=\frac{1}{\sqrt{2}}\mathbb{P}\left(\exists i\leq n:i^*=1\right)\mathbb{P}\left(\exists j\leq n:j^*=2|\exists i\leq n:i^*=1\right)$$

Hence,

$$\mathbb{E}\left(\left|\sum_{1\leq i,j\leq n} \varepsilon'_{\{i,j\},\{1,2\}} \mathbb{1}_{\{\{i^*,j^*\}=\{1,2\}\}}\right|\right) \geq \frac{1}{\sqrt{2}}\left(1-\left(1-\frac{1}{n}\right)^n\right)\left(1-\left(1-\frac{1}{n}\right)^{n-1}\right).$$

Then, because $1-x \leq e^{-x}$ and $n \geq 2$,

$$\mathbb{E}\left(\left|\sum_{1\leq i,j\leq n} \varepsilon'_{\{i,j\},\{1,2\}} \mathbb{1}_{\{\{i^*,j^*\}=\{1,2\}\}}\right|\right) \geq \frac{1}{\sqrt{2}}\left(1-e^{-1}\right)\left(1-e^{-1/2}\right).$$

The result follows from this and (5.5).

## 5.2 Results related to laws of large numbers

**Lemma S5.** *Under Assumptions 1-3, $\mathbb{E}\left(\sup_{\mathcal{F}} |\mathbb{P}_n^* f - \mathbb{P}_n f| \, | (Y_i)_{i\in \mathbb{I}_k}\right) \xrightarrow{a.s.} 0$.*

**Lemma S6.** *Suppose that Assumptions 2, 3 and 6 hold and $\boldsymbol{n} = (n_1(m),...,n_k(m))$ satisfies $n_j(m) \to \infty$ as $m \to \infty$ for $j = 1,...,k$. Then $\mathbb{E}\left(\sup_{\mathcal{F}} |\mathbb{P}_{\boldsymbol{n}}^* f - \mathbb{P}_{\boldsymbol{n}} f| \, | (Y_{\boldsymbol{i}})_{1\leq \boldsymbol{i} \leq \boldsymbol{n}}\right) \xrightarrow{a.s.} 0$ as $m \to \infty$.*

**Lemma S7** (Control of sums of quadratic terms).
*Let $h(\boldsymbol{i}) = \mathbb{1}_{\{\boldsymbol{i}\in \mathbb{I}_{n,k}\}} \sum_{\pi\in \mathfrak{S}_k} Y_{\boldsymbol{i}_\pi}$. If Assumption 1 holds and $\mathbb{E}[Y_{\boldsymbol{1}}^2] < \infty$, then, for every $j = 0,...,k$,*

$$\sum_{\boldsymbol{i}\in\{1,...,n\}^{2k-j}} h(i_1,...,i_k)h(i_1,...,i_l,i_{k+1},...,i_{2k-j})$$
$$= \sum_{c=0}^{k-j} \binom{k-j}{c}^2 \left(n^{2k-j-c}\mathbb{E}[h(1,...,k)h(1,...,j+c,k+1,...,2k-c-j)] + o_{a.s.}(n^{2k-j-c})\right).$$

**Lemma S8** (Control of sums of quadratic terms under separate exchangeability).
*Suppose Assumption 6 holds, $\mathbb{E}[Y_{\boldsymbol{1}}^2] < \infty$ and $\boldsymbol{n} = (n_1(m),...,n_k(m)) \in \mathbb{N}^{+k}$ satisfies $n_j(m) \to \infty$ when $m \to \infty$ for every $j = 1,...,k$. Then for every $\boldsymbol{e} \in \cup_{r=1}^k \mathcal{E}_r$*

$$\frac{1}{\prod_{r=1}^k n_r(n_r-1)^{1-e_r}} \sum_{(\boldsymbol{i},\boldsymbol{i}')\in \mathcal{I}_{\boldsymbol{n},\boldsymbol{e}}} Y_{\boldsymbol{i}} Y_{\boldsymbol{i}'} \xrightarrow{a.s.} \mathbb{E}[Y_{\boldsymbol{1}} Y_{\boldsymbol{b_e}}],$$

*where $\boldsymbol{b_e}$ is a k-dimensional vector such that its j-th entry is equal to 1 if $e_j = 1$ and 2 otherwise and $\mathcal{I}_{\boldsymbol{n},\boldsymbol{e}} = \left\{(\boldsymbol{i},\boldsymbol{i}') : \boldsymbol{1} \leq \boldsymbol{i},\boldsymbol{i}' \leq \boldsymbol{n}, \, i_j = i'_j \text{ if } e_j = 1, \, i_j \neq i'_j \text{ otherwise}\right\}$.*

### 5.2.1 Proof of Lemma S5

Let $i^*$ the $i$th index sampled with replacement in $\{1,...,n\}$. The $i^*$s are distributed as $i^* \overset{i.i.d}{\sim} \mathcal{U}_{\{1,...,n\}}$. For every $\boldsymbol{i} = (i_1,...,i_k) \in \mathbb{I}_{n,k}$, $\boldsymbol{i}^*$ stands for $(i_1^*,...,i_k^*)$. Conditional on the data and for every $f \in \mathcal{F}$, $\mathbb{P}_n^* f = \frac{(n-k)!}{n!} \sum_{\boldsymbol{i}\in \mathbb{I}_{n,k}} f(Y_{\boldsymbol{i}^*}) \mathbb{1}_{\{\boldsymbol{i}^* \in \mathbb{I}_{n,k}\}}$. We remark $\mathbb{E}\left(f(Y_{\boldsymbol{i}^*}) \mathbb{1}_{\{\boldsymbol{i}^* \in \mathbb{I}_{n,k}\}} | (Y_{\boldsymbol{i}})_{\boldsymbol{i}\in \mathbb{I}_k}\right) = \mathbb{P}'_n f = \mathbb{E}[\mathbb{P}_n^* f | (Y_{\boldsymbol{i}})_{\boldsymbol{i}\in \mathbb{I}_k}]$.

Conditional on $(Y_{\boldsymbol{i}})_{\boldsymbol{i}\in \mathbb{I}_k}$, $\frac{(n-k)!}{n!} \sum_{\boldsymbol{i}\in \mathbb{I}_{n,k}} f(Y_{\boldsymbol{i}^*}) \mathbb{1}_{\{\boldsymbol{i}^* \in \mathbb{I}_{n,k}\}}$ is a U-statistics since $f(Y_{\boldsymbol{i}^*}) \mathbb{1}_{\{\boldsymbol{i}^* \in \mathbb{I}_{n,k}\}}$ admits a representation $f(\tau(U_{i_1},...,U_{i_k})) \mathbb{1}_{\{(U_{i_1},...,U_{i_k})\in \mathbb{I}_{n,k}\}}$ for i.i.d. $U_i = i^*$. We also have that $\frac{(n-k)!}{n!} \sum_{\boldsymbol{i}\in \mathbb{I}_{n,k}} f(Y_{\boldsymbol{i}^*}) \mathbb{1}_{\{\boldsymbol{i}^* \in \mathbb{I}_{n,k}\}} = \frac{(n-k)!}{n!} \sum_{\boldsymbol{i}\in \mathbb{I}_{n,k}} h(\boldsymbol{i}^*)$ with $h : \boldsymbol{i} \mapsto \frac{1}{k!} \sum_{\pi\in \mathfrak{S}_k} f(Y_{\boldsymbol{i}_\pi}) \mathbb{1}_{\{\boldsymbol{i}_\pi \in \mathbb{I}_{n,k}\}}$. As a

65

result, the inequality proved on page 1508 in Arcones and Giné (1993) is valid with their $f$ replaced with $h$ (in particular, the sixth inequality on the latter page is true as $h$ is symmetric in its arguments and $h(\cdot)$ does not depend on $\boldsymbol{i}$). Then, for some $C_k > 0$ depending on $k$ only,

$$\mathbb{E}\left[\sup_{f\in\mathcal{F}_\delta}|\mathbb{P}_n^*f - \mathbb{P}'_nf|\,\Big|(Y_{\boldsymbol{i}})_{\boldsymbol{i}\in\mathbb{I}_k}\right] \leq kC_k\mathbb{E}\left[\sup_{f\in\mathcal{F}_\delta}\left|\frac{(n-k)!}{n!}\sum_{\boldsymbol{i}\in\mathbb{I}_{n,k}}\varepsilon_{\{i_1\}}f(Y_{\boldsymbol{i}^*})\mathbb{1}_{\{\boldsymbol{i}^*\in\mathbb{I}_{n,k}\}}\right|\,\Big|(Y_{\boldsymbol{i}})_{\boldsymbol{i}\in\mathbb{I}_k}\right].$$

Let $N^* = \frac{(n-k)!}{n!}\sum_{\boldsymbol{i}\in\mathbb{I}_{n,k}}\mathbb{1}_{\{\boldsymbol{i}^*\in\mathbb{I}_{n,k}\}}$. If $N^* = 0$, we sample fewer than $k$ different units in the bootstrap. In that case, the supremum of the Rademacher process is always equal to 0. As a result,

$$\mathbb{E}\left[\sup_{f\in\mathcal{F}}\left|\frac{(n-k)!}{n!}\sum_{\boldsymbol{i}\in\mathbb{I}_{n,k}}\varepsilon_{\{i_1\}}f(Y_{\boldsymbol{i}^*})\mathbb{1}_{\{F(Y_{\boldsymbol{i}^*})\leq M\}}\mathbb{1}_{\{\boldsymbol{i}^*\in\mathbb{I}_{n,k}\}}\right|\,\Big|(Y_{\boldsymbol{i}})_{\boldsymbol{i}\in\mathbb{I}_k}\right]$$
$$=\mathbb{E}\left[\sup_{f\in\mathcal{F}}\left|\frac{(n-k)!}{n!}\sum_{\boldsymbol{i}\in\mathbb{I}_{n,k}}\varepsilon_{\{i_1\}}f(Y_{\boldsymbol{i}^*})\mathbb{1}_{\{F(Y_{\boldsymbol{i}^*})\leq M\}}\mathbb{1}_{\{\boldsymbol{i}^*\in\mathbb{I}_{n,k}\}}\right|\,\Big|(Y_{\boldsymbol{i}})_{\boldsymbol{i}\in\mathbb{I}_k}, N^* > 0\right]\mathbb{P}(N^* > 0).$$

We now adapt the steps of the proof of Theorem 2.1. Conditional on $\left((Y_{\boldsymbol{i}})_{\boldsymbol{i}\in\mathbb{I}_k},(\boldsymbol{i}^*)_{\boldsymbol{i}\in\mathbb{I}_{n,k}}\right)$ and $N^* > 0$, we can consider for every $\eta_1 > 0$ and every $\boldsymbol{e}\in\mathcal{E}_1$ a minimal $\eta_1$-covering of $\mathcal{F}$ for the seminorm

$$||g||_{M,1}^* = \frac{(n-k)!}{n!}\sum_{i_1=1}^n\left|\sum_{(i_2,\ldots,i_k):\boldsymbol{i}\in\mathbb{I}_{n,k}}g(Y_{\boldsymbol{i}^*})\mathbb{1}_{\{F(Y_{\boldsymbol{i}^*})\leq M\}}\mathbb{1}_{\{\boldsymbol{i}^*\in\mathbb{I}_{n,k}\}}\right|$$

with balls centered in $\mathcal{F}$. This implies

$$\mathbb{E}\left[\sup_{\mathcal{F}}\left|\frac{(n-k)!}{n!}\sum_{\boldsymbol{i}\in\mathbb{I}_{n,k}}\varepsilon_{\{\boldsymbol{i}\odot\boldsymbol{e}\}+}f(Y_{\boldsymbol{i}^*})\mathbb{1}_{\{F(Y_{\boldsymbol{i}^*})\leq M\}}\mathbb{1}_{\{\boldsymbol{i}^*\in\mathbb{I}_{n,k}\}}\right|\,\Big|(Y_{\boldsymbol{i}})_{\boldsymbol{i}\in\mathbb{I}_k},(\boldsymbol{i}^*)_{\boldsymbol{i}\in\mathbb{I}_{n,k}},N^* > 0\right]$$
$$\leq\sqrt{2\log 2N\left(\eta_1,\mathcal{F},||\cdot||_{M,1}^*\right)}M\frac{1}{\sqrt{n}}+\eta_1.$$

Remark that $||g||_{M,1}^* \leq N^*||g||_{\mathbb{Q}_n,1}^*$ where $||g||_{\mathbb{Q}_n,1}^* = N^{*-1}\frac{(n-k)!}{n!}\sum_{\boldsymbol{i}\in\mathbb{I}_{n,k}}|g(Y_{\boldsymbol{i}^*})|\mathbb{1}_{\{\boldsymbol{i}^*\in\mathbb{I}_{n,k}\}}$, for $\mathbb{Q}_n = N^{*-1}\frac{(n-k)!}{n!}\sum_{\boldsymbol{i}\in\mathbb{I}_{n,k}}\delta_{\{Y_{\boldsymbol{i}^*}\}}\mathbb{1}_{\{\boldsymbol{i}^*\in\mathbb{I}_{n,k}\}}$ a (random) probability measure with finite support on $\mathcal{Y}$ that is well-defined when $N^* > 0$. Then, for every $\eta > 0$, letting $\eta_1 = \eta N^*||F||_{\mathbb{Q}_n,1}^*$ and using Point 2 of Lemma S12 and Point 1 of Lemma S12,

$$\mathbb{E}\left[\sup_{\mathcal{F}}\left|\frac{(n-k)!}{n!}\sum_{\boldsymbol{i}\in\mathbb{I}_{n,k}}\varepsilon_{\{i_1\}}f(Y_{\boldsymbol{i}^*})\mathbb{1}_{\{F(Y_{\boldsymbol{i}^*})\leq M\}}\mathbb{1}_{\{\boldsymbol{i}^*\in\mathbb{I}_{n,k}\}}\right|\,\Big|(Y_{\boldsymbol{i}})_{\boldsymbol{i}\in\mathbb{I}_k},(\boldsymbol{i}^*)_{\boldsymbol{i}\in\mathbb{I}_{n,k}},N^* > 0\right]$$
$$\leq\sqrt{2\log 2\sup_Q N\left(\eta||F||_{Q,1},\mathcal{F},||\cdot||_{Q,1}\right)}M\frac{1}{\sqrt{n}}+\eta N^*||F||_{\mathbb{Q}_n,1}^*.$$

Integration with respect to $(\boldsymbol{i}^*)_{\boldsymbol{i}\in\mathbb{I}_{n,k}}|(Y_{\boldsymbol{i}})_{\boldsymbol{i}\in\mathbb{I}_k}, N^* > 0$ combined with the fact that

$\mathbb{E}\left[N^*||F||^*_{\mathbb{Q}_{n,1}}|(Y_i)_{i\in\mathbb{I}_k}, N^*>0\right] = \mathbb{E}\left[N^*||F||^*_{\mathbb{Q}_{n,1}}|(Y_i)_{i\in\mathbb{I}_k}\right]/\mathbb{P}(N^*>0)$ leads to

$$\mathbb{E}\left[\sup_{\mathcal{F}}\left|\frac{(n-k)!}{n!}\sum_{i\in\mathbb{I}_{n,k}}\varepsilon_{\{i_1\}}f(Y_{i^*})\mathbb{1}_{\{F(Y_{i^*})\leq M\}}\mathbb{1}_{\{i^*\in\mathbb{I}_{n,k}\}}\right|\Big|(Y_i)_{i\in\mathbb{I}_k}\right]$$

$$\leq \sqrt{2\log 2\sup_Q N(\eta||F||_{Q,1},\mathcal{F},||\cdot||_{Q,1})}M\frac{1}{\sqrt{n}} + \eta\frac{(n-k)!}{n!}\sum_{i\in\mathbb{I}_{n,k}}\mathbb{E}\left[F(Y_{i^*})\mathbb{1}_{\{i^*\in\mathbb{I}_{n,k}\}}|(Y_i)_{i\in\mathbb{I}_k}\right]$$

$$=\sqrt{2\log 2\sup_Q N(\eta||F||_{Q,1},\mathcal{F},||\cdot||_{Q,1})}M\frac{1}{\sqrt{n}} + \eta\frac{1}{n^k}\sum_{i\in\mathbb{I}_{n,k}} F(Y_i).$$

By almost-sure convergence of the sample mean of jointly exchangeable, dissociated arrays, we can choose $\eta$ such that for $n$ large enough, the right-hand side is arbitrary small. Hence,

$$\mathbb{E}\left[\sup_{f\in\mathcal{F}}|\mathbb{P}_n^*f - \mathbb{P}_n'f|\,\Big|(Y_i)_{i\in\mathbb{I}_k}\right] \xrightarrow{\text{a.s.}} 0. \tag{5.6}$$

Finally, the triangle inequality enables us to write

$$\mathbb{E}\left[\sup_{f\in\mathcal{F}}|\mathbb{P}_n^*f - \mathbb{P}_n f|\,\Big|(Y_i)_{i\in\mathbb{I}_k}\right]$$

$$\leq \mathbb{E}\left[\sup_{f\in\mathcal{F}}\left|\left(\frac{(n-k)!}{n!} - \frac{1}{n^k}\right)\sum_{i\in\mathbb{I}_{n,k}} f(Y_i)\right|\,\Big|(Y_i)_{i\in\mathbb{I}_k}\right] + \mathbb{E}\left[\sup_{f\in\mathcal{F}}|\mathbb{P}_n^*f - \mathbb{P}_n'f|\,\Big|(Y_i)_{i\in\mathbb{I}_k}\right]$$

$$\leq \left(1 - \frac{n!}{n^k(n-k)!}\right)\frac{(n-k)!}{n!}\sum_{i\in\mathbb{I}_{n,k}} F(Y_i) + \mathbb{E}\left[\sup_{f\in\mathcal{F}}|\mathbb{P}_n^*f - \mathbb{P}_n'f|\,\Big|(Y_i)_{i\in\mathbb{I}_k}\right].$$

Using (5.6) and $\frac{n!}{n^k(n-k)!} \to 1$, we conclude that $\mathbb{E}\left[\sup_{f\in\mathcal{F}}|\mathbb{P}_n^*f - \mathbb{P}_n f|\,\big|(Y_i)_{i\in\mathbb{I}_k}\right] \xrightarrow{\text{a.s.}} 0$.

### 5.2.2 Proof of Lemma S6

For every $j = 1,...,k$, let $i_j^*$ the $i_j$-th index sampled with replacement in $[1;n_j]$. The $i_j^*$s are distributed as $i_j^* \overset{i.i.d}{\sim} \mathcal{U}_{[1;n_j]}$ and the $k$ sequences $(i_1^*)_{i_1=1}^{n_1},...,(i_k^*)_{i_k=1}^{n_k}$ are also mutually independent. For every $\mathbf{1} \leq \mathbf{i} \leq \mathbf{n}$, $\mathbf{i}^*$ denotes $(i_1^*,...,i_k^*)$. Conditional on the data and for every $f \in \mathcal{F}$, $\mathbb{P}_{\mathbf{n}}^*f = \frac{1}{\Pi_{\mathbf{n}}}\sum_{\mathbf{1}\leq\mathbf{i}\leq\mathbf{n}} h(\mathbf{i}^*)$ with $h(\mathbf{i}^*) = f(Y_{\mathbf{i}^*})$. We have: $\mathbb{E}\left[\mathbb{P}_{\mathbf{n}}^*f|(Y_{\mathbf{i}})_{\mathbf{i}\in\mathbb{N}^{+k}}\right] = \mathbb{P}_{\mathbf{n}}f$. Note that conditional on $(Y_{\mathbf{i}})_{\mathbf{i}\in\mathbb{N}^{+k}}$, $(\mathbf{i}^*)_{\mathbf{i}\in\mathbb{I}_{n,k}}$ is a family of random vectors that admit a representation $\mathbf{i}^* = \tau((U_{\mathbf{i}\odot e})_{e\in\mathcal{E}_1})$ with $(U_{\mathbf{i}})_{\mathbf{0}\leq\mathbf{i}\leq\mathbf{n}}$ i.i.d. random variables (consider $\tau: (u_1,...,u_k) \in [0,1]^k \mapsto (\lceil n_1\times u_1\rceil,...,\lceil n_k\times u_k\rceil)$ where $\lceil\cdot\rceil$ denotes the ceiling function and $U_{\mathbf{i}} \sim \mathcal{U}_{[0,1]}$). As a result, conditionally on the data, Lemma S2 applies to $\widetilde{Y}_{\mathbf{i}} = \mathbf{i}^*$, $\mathcal{G} = \{h: h(\mathbf{i}^*) = f(Y_{\mathbf{i}^*}), f\in\mathcal{F}\}$ and $\Phi = \text{Id}$. Moreover, because only terms involving $e \in \mathcal{E}_1$ appear in the representation of $\mathbf{i}^*$, a simplification of the proof of Lemma S2 leads to the following inequality:

$$\mathbb{E}\left[\sup_{f\in\mathcal{F}}|\mathbb{P}_{\mathbf{n}}^*f - \mathbb{P}_{\mathbf{n}}f|\,\Big|(Y_{\mathbf{i}})_{\mathbf{i}\in\mathbb{N}^{+k}}\right] \leq \frac{2}{\Pi_{\mathbf{n}}}\sum_{\mathbf{i}=1}^{\mathbf{n}} F(Y_{\mathbf{i}})\mathbb{1}_{\{F(Y_{\mathbf{i}})>M\}}$$

$$+ 2\sum_{e\in\mathcal{E}_1}\mathbb{E}\left[\sup_{f\in\mathcal{F}}\left|\frac{1}{\Pi_{\mathbf{n}}}\sum_{\mathbf{1}\leq\mathbf{i}\leq\mathbf{n}}\varepsilon_{\mathbf{i}\odot e}f(Y_{\mathbf{i}^*})\mathbb{1}_{\{F(Y_{\mathbf{i}^*})\leq M\}}\right|\,\Big|(Y_{\mathbf{i}})_{\mathbf{i}\in\mathbb{N}^{+k}}\right].$$

The rest of the proof is similar to that of $\sup_{\mathcal{F}} |\mathbb{P}_n f - Pf| \xrightarrow{L^1} 0$: in fact, with $||h||_{e,M,1}$ redefined as $||h||_{e,M,1} = \frac{1}{\Pi_n} \sum_{e \leq c \leq n \odot e} \left| \sum_{1-e \leq c' \leq n \odot (1-e)} f(Y_{i^*}) \mathbb{1}_{\{F(Y_{i^*}) \leq M\}} \right|$, we have for every $e \in \mathcal{E}_1$, $M > 0$ and $\eta_1 \geq 0$, possibly random,

$$\mathbb{E}\left[\sup_{f \in \mathcal{F}} \left| \frac{1}{\Pi_n} \sum_{1 \leq i \leq n} \varepsilon_{i \odot e} f(Y_{i^*}) \mathbb{1}_{\{F(Y_{i^*}) \leq M\}} \right| \bigg| (Y_i)_{i \in \mathbb{N}^{+k}} \right]$$

$$\leq \mathbb{E}\left[\sup_{f \in \mathcal{F}} \left| \frac{1}{\Pi_n} \sum_{1 \leq i \leq n} \varepsilon_{i \odot e} f(Y_{i^*}) \mathbb{1}_{\{F(Y_{i^*}) \leq M\}} \right| \bigg| (Y_i)_{i \in \mathbb{N}^{+k}} \right]$$

$$\leq \mathbb{E}\left[\sqrt{2 \log 2 N(\eta_1, \mathcal{F}, || \cdot ||_{e,M,1})} M \frac{1}{\sqrt{n}} + \eta_1 \bigg| (Y_i)_{i \in \mathbb{N}^{+k}} \right].$$

### 5.2.3 Proof of Lemma S7

By definition of $h(\cdot)$, we have

$$\sum_{i \in \{1,...,n\}^{2k-j}} h(i_1, ..., i_k) h(i_1, ..., i_j, i_{k+1}, ..., i_{2k-j})$$

$$= \sum_{i \in \{1,...,n\}^j} \sum_{i' \in (\{1,...,n\} \setminus \{i\})^{k-j}} \sum_{i'' \in (\{1,...,n\} \setminus \{i\})^{k-j}} h(i, i') h(i, i'')$$

$$= \sum_{c=0}^{k-j} \binom{k-j}{c}^2 \sum_{i \in \{1,...,n\}^{j+c}} \sum_{i' \in (\{1,...,n\} \setminus \{i\})^{k-j-c}} \sum_{i'' \in (\{1,...,n\} \setminus (\{i\} \cup \{i'\}))^{k-j-c}} h(i, i') h(i, i'').$$

Since $h$ is invariant by permutation of its entries, the last equality holds by distinguishing between cases depending on the number of common values in the vectors $(i_{j+1}, ..., i_k)$ and $(i_{k+1}, ..., i_{2k-j})$. As $(h(i))_{i \in \mathbb{I}_{n,k}}$ is a $k$-dimensional jointly exchangeable array,

$$\left( h(i, i') h(i, i'') \right)_{i \in \overline{\{1,...,n\}^{j+c}}, i' \in \overline{(\{1,...,n\} \setminus \{i\})^{k-j-c}}, i'' \in \overline{(\{1,...,n\} \setminus (\{i\} \cup \{i'\}))^{k-j-c}}}$$

is a $(2k - j - c)$-dimensional jointly exchangeable array. Moreover $\mathbb{E}(Y_1^2) < \infty$ ensures that $\mathbb{E}(|h(1, ..., k) h(1, ..., j+c, k+1, ..., 2k-j-c)|) < \infty$ so that Theorem 2.1 can be applied to a class $\mathcal{F}$ reduced to the identity function. The equivalence $\frac{n!}{(n-(2k-j-c))!} \sim n^{2k-j-c}$ concludes the proof.

### 5.2.4 Proof of Lemma S8

Let $\Sigma_{m,e}$ the $\sigma$-algebra generated by the set of functions $g$ from $\mathcal{D}^{\mathbb{N}^{+k}} \times \mathcal{D}^{\mathbb{N}^{+k}}$ to $\mathbb{R}$ such that:

$$g((Y_i, Y_{i'})_{(i,i') \in \mathcal{I}_{n,e}}) = g((Y_{\pi_1(i_1),...,\pi_k(i_k)}, Y_{\pi_1(i'_1),...,\pi_k(i'_k)})_{(i,i') \in \mathcal{I}_{n,e}}),$$

for every set of permutations $\pi_1, ..., \pi_k$ such that for every $r = 1, ..., k$, $\pi_r(i) = i$ if $i \geq n_r$. Let

$$W_m = \frac{1}{\prod_{r=1}^{k} n_r (n_r - 1)^{(1-e_r)}} \sum_{(i,i') \in \mathcal{I}_{n,e}} Y_i Y_{i'}.$$

By construction, we have for every $\underline{n} \in \mathbb{N}^+$, $W_m = \mathbb{E}[W_m | \Sigma_{m,e}] = \mathbb{E}[Y_1 Y_{b_e} | \Sigma_{m,e}]$. Furthermore, $\Sigma_{m,e} \supseteq \Sigma_{m+1,e}$ so that

$$\mathbb{E}[W_m | \Sigma_{m+1,e}] = \mathbb{E}[\mathbb{E}[Y_1 Y_{b_e} | \Sigma_{m,e}] | \Sigma_{m+1,e}] = \mathbb{E}[Y_1 Y_{b_e} | \Sigma_{m+1,e}] = W_{m+1}.$$

As a result, we can conclude that $(W_m, \Sigma_{m,e})_{m \geq 1}$ is a backward martingale. From this follows that $W_m \xrightarrow{\text{a.s.}} \mathbb{E}[Y_1 Y_{b_e} | \Sigma_{\infty,e}]$ where $\Sigma_{\infty,e} = \cap_{m \geq 1} \Sigma_{m,e}$. By the dissociation assumption, this sigma-algebra is trivial (see Lemma 7.35 in Kallenberg, 2005), hence $W_m \xrightarrow{\text{a.s.}} \mathbb{E}[Y_1 Y_{b_e}]$.

## 5.3 Contraction and maximal inequalities in degenerate cases

**Lemma S9** (A Kinchine-Kahane inequality).
*For every $1 < p < q < \infty$ and every bounded subset $T$ of $\mathbb{R}^{n(n-1)/2}$*

$$\mathbb{E}\left[\sup_{t \in T} \left| \frac{1}{n(n-1)} \sum_{\boldsymbol{i} \in \overrightarrow{\mathbb{I}_{n,2}}} \varepsilon_{i_1}^1 \varepsilon_{i_2}^2 t_{\boldsymbol{i}} \right|^q \right]^{1/q} \leq \frac{q-1}{p-1} \mathbb{E}\left[\sup_{t \in T} \left| \frac{1}{n(n-1)} \sum_{\boldsymbol{i} \in \overrightarrow{\mathbb{I}_{n,2}}} \varepsilon_{i_1}^1 \varepsilon_{i_2}^2 t_{\boldsymbol{i}} \right|^p \right]^{1/p}.$$

**Lemma S10** (A maximal inequality).
*Let $T$ be a countable subset of $\mathbb{R}^{n(n-1)/2}$ that contains the null vector and for every $t \in T$, $\|t\|_T^2 := \mathbb{E}\left[\left| \frac{1}{n-1} \sum_{\boldsymbol{i} \in \overrightarrow{\mathbb{I}_{n,2}}} \varepsilon_{i_1}^1 \varepsilon_{i_2}^2 t_{\boldsymbol{i}} \right|^2\right] = \frac{1}{(n-1)^2} \sum_{\boldsymbol{i} \in \overrightarrow{\mathbb{I}_{n,2}}} t_{\boldsymbol{i}}^2$. Let also $D := \sup_{(t_1,t_2) \in T \times T} \|t_1 - t_2\|_T$ be the diameter of $T$ for $\|\cdot\|_T$. There exists a constant $K$ such that*

$$\mathbb{E}\left[\sup_{t \in T} \left| \frac{1}{n-1} \sum_{\boldsymbol{i} \in \overrightarrow{\mathbb{I}_{n,2}}} \varepsilon_{i_1}^1 \varepsilon_{i_2}^2 t_{\boldsymbol{i}} \right| \right] \leq K \int_0^{D/2} \log N(\eta, T, \|\cdot\|_T) d\eta.$$

### 5.3.1 Proof of Lemma S9

The reasoning is the same as that at the end of the proof of Theorem 3.2.1 in de la Peña and Giné (1999) and is therefore omitted.

### 5.3.2 Proof of Lemma S10

Let $A_n(t) := \frac{1}{n-1} \sum_{\boldsymbol{i} \in \overrightarrow{\mathbb{I}_{n,2}}} \varepsilon_{i_1}^1 \varepsilon_{i_2}^2 t_{\boldsymbol{i}}$. We first prove that for every $m \geq 1$ and every set of elements of $T$ of cardinality $m$

$$\mathbb{E}\left[\max_{1 \leq j \leq m} A_n(t_j)\right] \leq e \max_{1 \leq j \leq m} \sqrt{\mathbb{E}[A_n(t_j)^2]} \log m + \frac{\max_{1 \leq j \leq m} \sqrt{\mathbb{E}[A_n(t_j)^2]}}{\sqrt{\pi}}$$
$$+ \frac{2 \max_{1 \leq j \leq m} \sqrt{\mathbb{E}[A_n(t_j)^2]}}{\sqrt{\pi}} \sqrt{e \log m + \frac{\max_{1 \leq j \leq m} \sqrt{\mathbb{E}[A_n(t_j)^2]}}{\sqrt{\pi}}}. \quad (5.7)$$

Using Lemma S9 with $T = \{t_j\}$ for every $j \in \{1, ..., m\}$, the series expansion of the exponential function and $\mathbb{E}[A_n(t_j)] = 0$ for every $j \in \{1, ..., m\}$, we can write for every $\lambda \in (0, 1/e \max_{1 \leq j \leq m} \sqrt{\mathbb{E}[A_j^2]})$

$$\mathbb{E}\left[e^{\lambda A_n(t_j)}\right] = \sum_{\ell=0}^{\infty} \frac{\lambda^\ell}{\ell!} \mathbb{E}[A_n(t_j)\ell] \leq 1 + \sum_{\ell=2}^{\infty} \frac{\lambda^\ell (\ell-1)^\ell}{\ell!} \mathbb{E}[A_n(t_j)^2]^{\ell/2}.$$

Using $\sqrt{2\pi\ell}\left(\frac{\ell}{e}\right)^{\ell} < \ell!$ and $1 + a < e^a$ for $a > 0$, we have

$$\mathbb{E}\left[e^{\lambda A_n(t_j)}\right] < 1 + \sum_{\ell=2}^{\infty} \frac{\lambda^\ell e^\ell}{2\sqrt{\pi e}} \mathbb{E}[A_n(t_j)^2]^{\ell/2}$$

$$\leq 1 + \sum_{\ell=0}^{\infty} \frac{\lambda^\ell e^\ell}{2\sqrt{\pi e}} \mathbb{E}[A_n(t_j)^2]^{\ell/2}$$

$$\leq \exp\left(\sum_{\ell=0}^{\infty} \frac{\lambda^\ell e^\ell}{2\sqrt{\pi e}} \mathbb{E}[A_n(t_j)^2]^{\ell/2}\right)$$

$$= \exp\left(\frac{1}{2\sqrt{\pi e}(1 - \lambda e\sqrt{\mathbb{E}[A_n(t_j)^2]})}\right)$$

$$\leq \exp\left(\left(2\sqrt{\pi e}\left(1 - \lambda e \max_{1\leq j\leq m}\sqrt{\mathbb{E}[A_n(t_j)^2]}\right)\right)^{-1}\right).$$

The last inequality and standard convexity and monotonicity arguments (see, e.g., the bottom of page 39 in Giné and Nickl, 2015) yield

$$\mathbb{E}\left[\max_{1\leq j\leq m} A_n(t_j)\right] \leq \frac{\log\left(m \times \max_{1\leq j\leq m} \mathbb{E}\left[\exp\left(\lambda A_n(t_j)\right)\right]\right)}{\lambda}$$

$$\leq \frac{\log\left(m \exp\left[\left(2\sqrt{\pi e}\left(1 - \lambda e \max_{1\leq j\leq m}\sqrt{\mathbb{E}[A_n(t_j)^2]}\right)\right)^{-1}\right]\right)}{\lambda}$$

$$= \frac{\log m}{\lambda} + \frac{1}{\lambda}\left(2\sqrt{\pi e}\left(1 - \lambda e \max_{1\leq j\leq m}\sqrt{\mathbb{E}[A_n(t_j)^2]}\right)\right)^{-1}$$

Pick $\lambda = \gamma/[e\max_{1\leq j\leq m}\sqrt{\mathbb{E}[A_n(t_j)^2]}]$ and $\gamma \in (0,1)$. Then, minimizing over $\gamma$, we arrive at (5.7).

To bound from above $\mathbb{E}\left[\max_{1\leq j\leq m}|A_n(t_j)|\right]$, remark that it is equal to $\mathbb{E}\left[\max_{1\leq j\leq 2m}\widetilde{A}_n(t_j)\right]$, where $(\widetilde{A}_n(t_j))_{j=1}^m = (A_n(t_j))_{j=1}^m$ and $(\widetilde{A}_n(t_j))_{j=m+1}^{2m} = (-A_n(t_j))_{j=1}^m$. The previous result thus applies with $\log 2m$ instead of $\log m$. Recalling the definition of $||\cdot||_T$, we obtain

$$\mathbb{E}\left[\max_{1\leq j\leq m}|A_n(t_j)|\right] \leq e \max_{1\leq j\leq m}||t_j||_T \log 2m + \frac{\max_{1\leq j\leq m}||t_j||_T}{\sqrt{\pi}}$$
$$+ \frac{2\max_{1\leq j\leq m}||t_j||_T}{\sqrt{\pi}}\sqrt{e\log 2m + \frac{\max_{1\leq j\leq m}||t_j||_T}{\sqrt{\pi}}}. \quad (5.8)$$

Next, we bound from above the right-hand side of (5.8) by $K\int_0^{D/2}\log N(\eta, T, ||\cdot||_T)d\eta$. To do so, we simply observe that the proof of the first statement of Theorem 2.3.6 in Giné and Nickl (2015) can be replicated using (5.8) instead of the maximal inequality stemming from their Lemma 2.3.4. In our case, the stochastic process of interest is $A_n(t), t \in T$ and we choose $t_0$ equal to the null vector. Then there exists some numerical constant $K > 0$ such that for every $m \geq 1$,

$$\mathbb{E}\left[\max_{1\leq j\leq m}|A_n(t_j)|\right] \leq K\int_0^{D/2}\log N(\eta, T, ||\cdot||_T)d\eta.$$

This inequality extends to the whole set $T$ by monotone convergence, ending the proof of the lemma.

## 5.4 Covering and entropic integrals

**Lemma S11** (Properties of entropic integrals)**.**
Let $\mathcal{F}$ a class of functions with envelope $F$ such that $\int_0^\infty \zeta(\varepsilon)d\varepsilon < \infty$, with
$$\zeta(\varepsilon) = \sup_Q \sqrt{\log(N(\varepsilon||F||_{Q,2}, \mathcal{F}, ||\cdot||_{Q,2}))}.$$

1. $u \mapsto J_{\mathcal{F}}(u) = \int_0^u \zeta(\varepsilon)d\varepsilon$ is positive, non-decreasing, concave, larger than $u\zeta(u)$ for every $u > 0$ and $\sup_{u \geq 0} J_{\mathcal{F}}(u) = J_{\mathcal{F}}(2)$.

2. For every $K > 0$, $(x, y) \in [0, \infty) \times (0, \infty) \mapsto \sqrt{y} J_{\mathcal{F}}\left(K\frac{\sqrt{x}}{\sqrt{y}}\right)$ is concave.

**Lemma S12** (Covering numbers inequalities)**.**
For every $\varepsilon > 0$:

1. for every class $\mathcal{H}$, every norm $||\cdot||$ and every $\lambda > 0$: $N(\varepsilon, \mathcal{H}, \lambda||\cdot||) = N(\varepsilon/\lambda, \mathcal{H}, ||\cdot||)$.

2. for every class $\mathcal{H}$, every pair of norms $||\cdot|| \leq ||\cdot||'$: $N(\varepsilon, \mathcal{H}, ||\cdot||) \leq N(\varepsilon, \mathcal{H}, ||\cdot||')$.

3. for every $\mathcal{H} \subset \mathcal{H}'$ and every norm $||\cdot||$: $N(\varepsilon, \mathcal{H}, ||\cdot||) \leq N(\varepsilon/2, \mathcal{H}', ||\cdot||)$.

4. for every norm $||\cdot||$, every class $\mathcal{F}$ and for $\mathcal{F}_\infty = \{f : f = f_1 - f_2, (f_1, f_2) \in \mathcal{F} \times \mathcal{F}\}$:
$N(\varepsilon, \mathcal{F}_\infty, ||\cdot||) \leq N^2(\varepsilon/2, \mathcal{F}, ||\cdot||)$.

5. for every class $\mathcal{F}$ and for $\mathcal{F}_\infty^2 = \{f : f = (f_1 - f_2)^2, (f_1, f_2) \in \mathcal{F} \times \mathcal{F}\}$:
$\sup_Q N(8\varepsilon||F^2||_{Q,1}, \mathcal{F}_\infty^2, ||\cdot||_{Q,1}) \leq \sup_Q N^2(\varepsilon||F||_{Q,2}, \mathcal{F}, ||\cdot||_{Q,2})$
where the supremum is taken over the set of all finite probability measures on the domain of the functions in $\mathcal{F}$.

### 5.4.1 Proof of Lemma S11

1. $\zeta$ is nonnegative and nonincreasing. It follows that $u \mapsto J_{\mathcal{F}}(u)$ is positive, non-decreasing and concave. Furthermore, $J_{\mathcal{F}}(u) \geq \int_0^u \zeta(u)d\varepsilon = u\zeta(u)$ for every $u > 0$. For $\varepsilon \geq 2$, we have $N(\varepsilon||F||_{Q,2}, \mathcal{F}, ||\cdot||_{Q,2}) = 1$ for every probability measure $Q$. As a result, $\zeta(\varepsilon) = 0$.

2. $J$ is concave on $[0, \infty)$ which implies for $\lambda \in (0; 1)$, $(x, x') \in [0, \infty)^2$, $(y, y') \in (0, \infty)^2$

$$(\lambda y + (1-\lambda)y') J_{\mathcal{F}}\left(K\frac{\lambda x + (1-\lambda)x'}{\lambda y + (1-\lambda)y'}\right)$$
$$= (\lambda y + (1-\lambda)y') J_{\mathcal{F}}\left(\frac{\lambda y}{\lambda y + (1-\lambda)y'}\frac{Kx}{y} + \frac{(1-\lambda)y'}{\lambda y + (1-\lambda)y'}\frac{Kx'}{y'}\right)$$
$$\geq \lambda y J_{\mathcal{F}}\left(K\frac{x}{y}\right) + (1-\lambda)y' J_{\mathcal{F}}\left(K\frac{x'}{y'}\right).$$

We can therefore claim that $f(x, y) = y J_{\mathcal{F}}(K\frac{x}{y})$ is concave on $[0, \infty) \times (0, \infty)$. Moreover $f(x, y)$ is non-decreasing in $x$ as $J_{\mathcal{F}}$ is non-decreasing. We also have $f(x, y) = y \int_0^{K\frac{x}{y}} \zeta(\varepsilon)d\varepsilon = x \int_0^1 \zeta\left(K\frac{x}{y}\varepsilon\right) d\varepsilon$.

Since $\zeta$ in nonincreasing, $f$ is non-decreasing in $y$. Finally, because $u \mapsto \sqrt{u}$ is concave, we have

$$\sqrt{\lambda y + (1-\lambda)y'} J_{\mathcal{F}}\left(K\frac{\sqrt{\lambda x + (1-\lambda)x'}}{\sqrt{\lambda y + (1-\lambda)y'}}\right) = f\left(\sqrt{\lambda x + (1-\lambda)x'}, \sqrt{\lambda y + (1-\lambda)y'}\right)$$
$$\geq f\left(\lambda\sqrt{x} + (1-\lambda)\sqrt{x'}, \lambda\sqrt{y} + (1-\lambda)\sqrt{y'}\right)$$
$$\geq \lambda f\left(\sqrt{x}, \sqrt{y}\right) + (1-\lambda)f\left(\sqrt{x'}, \sqrt{y'}\right)$$
$$= \lambda\sqrt{y} J_{\mathcal{F}}\left(K\frac{\sqrt{x}}{\sqrt{y}}\right) + (1-\lambda)\sqrt{y'} J_{\mathcal{F}}\left(K\frac{\sqrt{x'}}{\sqrt{y'}}\right).$$

### 5.4.2 Proof of Lemma S12

1. A ball of radius $\varepsilon$ for the norm $\lambda||\cdot||$ is a ball of radius $\varepsilon/\lambda$ for the norm $||\cdot||$.

2. A minimal $\varepsilon$-covering for $||\cdot||'$ is also an $\varepsilon$-covering for $||\cdot||$.

3. Consider a minimal $\varepsilon/2$-covering of $\mathcal{H}'$. This is not an $\varepsilon/2$-covering of $\mathcal{H}$ in general because the centers of the covering balls need not be in $\mathcal{H}$. However, in each ball that intersects $\mathcal{H}$, we can select an element of $\mathcal{H}$ as a center of a new ball of radius $\varepsilon$. We thus obtain a new family of balls which forms an $\varepsilon$-covering of $\mathcal{H}$.

4. Let $f_1, ..., f_{N(\varepsilon/2, \mathcal{F}, ||\cdot||)}$ the centers of balls of a minimal $\varepsilon/2$-covering of $\mathcal{F}$. Consider balls of center $f_i - f_j$ and of radius $\varepsilon$ for $1 \leq i, j \leq N(\varepsilon/2, \mathcal{F}, ||\cdot||)$. The latter constitute an $\varepsilon$-covering of $\mathcal{F}_\infty$ because for $(g_1, g_2) \in \mathcal{F} \times \mathcal{F}$ we have

$$||(f_i - f_j) - (g_1 - g_2)|| \leq ||f_i - g_1|| + ||f_j - g_2||,$$

which is smaller than $\varepsilon$ for at least one pair $(i, j)$.

5. Let $f_1, ..., f_{N(\varepsilon||F||_{Q,2}, \mathcal{F}, ||\cdot||)}$ the centers of balls of a minimal $\varepsilon||F||_{Q,2}$-covering of $\mathcal{F}$ for $||\cdot||_{Q,2}$. Consider balls of center $f_i - f_j$ and radius $8\varepsilon||F^2||_{Q,1}$ for the norm $||\cdot||_{Q,1}$. For every pair $(g_1, g_2) \in \mathcal{F} \times \mathcal{F}$, the Cauchy-Schwarz inequality implies

$$||(f_i - f_j)^2 - (g_1 - g_2)^2||_{Q,1} \leq ||f_i - f_j + g_1 - g_2||_{Q,2} \times ||(f_i - f_j) - (g_1 - g_2)||_{Q,2}$$
$$\leq 4||F||_{Q,2} \times (||f_i - g_1||_{Q,2} + ||f_j - g_2||_{Q,2}),$$

which is smaller than $8\varepsilon||F||_{Q,2}^2 = 8\varepsilon||F^2||_{Q,1}$ for at least one pair $(i, j)$.



\end{document}